\documentclass[binding=0.6cm,PhD]{sapthesis}

\usepackage{microtype}
\usepackage{hyperref}
\hypersetup{pdftitle={Deformation problems in the categories of Lie groupoids and algebroids},pdfauthor={Pier Paolo La Pastina}}

\usepackage{amsthm}
\usepackage{amssymb}
\usepackage{tensor}
\usepackage{todonotes}
\usepackage[all, cmtip]{xy}
\usepackage{stmaryrd}
\usepackage{mnsymbol}
\usepackage{tikz-cd}

\title{Deformations of vector bundles in the categories of Lie algebroids and groupoids}
\author{Pier Paolo La Pastina}
\IDnumber{1751192}
\course[Mathematics]{Matematica}
\courseorganizer{Scuola Dottorale in Scienze Astronomiche, Chimiche, Fisiche, Matematiche e della Terra "Vito Volterra"}
\cycle{XXXII}
\submitdate{October 2019}
\copyyear{2019}
\advisor{Dr. Luca Vitagliano, PhD}
\authoremail{plapastina@gmail.com}

\newtheorem{theorem}{Theorem}[subsection]
\newtheorem{corollary}[theorem]{Corollary}
\newtheorem{lemma}[theorem]{Lemma}
\newtheorem{proposition}[theorem]{Proposition}
\newtheorem*{theorem*}{Theorem}

\theoremstyle{definition}
\newtheorem{definition}[theorem]{Definition}
\newtheorem{examplex}[theorem]{Example}
\newenvironment{example}
  {\pushQED{\qed}\examplex}
  {\popQED\endexamplex}
  
\theoremstyle{remark}
\newtheorem{remarkx}[theorem]{Remark}
\newenvironment{remark}
  {\pushQED{\qed}\remarkx}
  {\popQED\endremarkx}

\begin{document}

\frontmatter
\maketitle

\begin{abstract}

This thesis deals with deformations of VB-algebroids and VB-groupoids. They can be considered as vector bundles in the categories of Lie algebroids and groupoids and encompass several classical objects, including Lie algebra and Lie group representations, 2-vector spaces and the tangent and the cotangent algebroid (groupoid) to a Lie algebroid (groupoid). Moreover, they are geometric models for some kind of representations of Lie algebroids (groupoids), namely 2-term representations up to homotopy. Finally, it is well known that Lie groupoids are ``concrete'' incarnations of differentiable stacks, hence VB-groupoids can be considered as representatives of vector bundles over differentiable stacks, and VB-algebroids their infinitesimal versions.

In this work, we attach to every VB-algebroid and VB-groupoid a cochain complex controlling its deformations, their \emph{linear deformation complex}. Moreover, the deformation complex of a VB-algebroid is equipped with a DGLA structure. The basic properties of these complexes are discussed: their relationship with the deformation complexes of the total spaces and the base spaces, particular cases and generalizations. The main theoretical results are a linear van Est theorem, that gives conditions for the linear deformation cohomology of a VB-groupoid to be isomorphic to that of the corresponding VB-algebroid, and a Morita invariance theorem, that implies that the linear deformation cohomology of a VB-groupoid is really an algebraic invariant of the associated vector bundle of differentiable stacks. Finally, several examples are discussed, showing how the linear deformation cohomologies are related to other well-known cohomologies. 

\end{abstract}

\chapter*{Acknowledgments}

Research is not only about individual challenge. It is rather the result of a collective effort, where people who support us emotionally and scientifically play a fundamental role. Here I want to thank some of them, because I find their contribution to this work particularly important.

First of all, I want to thank my thesis advisor, Luca Vitagliano, for having introduced me to the world of mathematical research with patience and dedication. Moreover, I would like to thank my tutor, Andrea Sambusetti, for his guidance during the first steps of my PhD. Then, I would like to thank all my colleagues in Rome and in Salerno for having shared with me pleasures and pains of being a PhD student. A special thank goes to the geometry group in Salerno - Alfonso, Andreas, Antonio, Chiara, Luca again, Jonas and Stefano - for all our stimulating discussions, that have often gone far beyond the borders of mathematics.

I also want to mention some people I have met during my research stays. I thank Cristian Ortiz, Ivan Struchiner and all the PhD students in geometry at São Paulo University, for their warm hospitality and for some insights into research topics related to my thesis. More generally, I thank all the big community of researchers in Poisson geometry in Brazil, for having hosted me also in Rio de Janeiro and S\~{a}o Carlos and having given me the opportunity to give talks and share some ideas. In addition, I would like to thank João Nuno Mestre, for having accepted with enthusiasm to work on a common project and having hosted me at Max Planck Institute for Mathematics.

The last paragraph goes to my dear ones. I thank my family and my friends because they constantly encouraged me during difficult times, that are almost unavoidable during a PhD. Despite in these years I have been travelling quite a lot, they have always been by my side and now our bonds have become even stronger.

\tableofcontents

\mainmatter

\chapter{Introduction}

The theory of Lie groupoids and algebroids is nowadays a well-established area of research in differential geometry. While groupoids appeared in the early 20th century \cite{brandt:uber}, Lie groupoids were introduced in the 1950's by Charles Ehresmann as a tool for understanding principal bundles \cite{ehresmann:categories}. Their theory generalizes that of Lie groups: they have an infinitesimal counterpart, Lie algebroids, first defined by Jean Pradines in \cite{pradines:theorie}, and there exists a Lie differentiation functor from Lie groupoids to Lie algebroids.

It is well known that, despite its objects are very regular and exhibit a rich geometry, the category of smooth manifolds is quite ill-behaved: for example, the quotient of a smooth manifold is rarely smooth. Nevertheless, these quotients, such as orbifolds, orbit spaces of Lie group actions or leaf spaces of foliations, arise very naturally in differential geometry, so it should be possible to perform some sort of differential calculus over them. One way of doing this is by observing that they can be modeled by certain Lie groupoids, that play the role of ``desingularizations'' of such spaces. This observation can be made rigorous by the introduction of \emph{differentiable stacks} \cite{behrend:diff}: they encompass all these examples and one can show that they are equivalent to Lie groupoids up to a certain notion of equivalence, namely \emph{Morita equivalence}.

This leads to the idea that the category of Lie groupoids should be understood as a more general setting for differential geometry, and indeed several kinds of compatible geometric structures on Lie groupoids (symplectic structures, complex structures, Riemannian metrics and so on) have been considered in the last years, along with their infinitesimal counterparts. This thesis will be concerned with vector bundles in the categories of Lie groupoids and algebroids, known as \emph{VB-groupoids} and \emph{VB-algebroids}. VB-groupoids first appeared in \cite{pradines:remarque} in connection to symplectic groupoids. There is a notion of Morita equivalence for VB-groupoids \cite{delhoyo:morita} and one can interpret VB-groupoids up to Morita equivalence as vector bundles in the category of differentiable stacks. VB-algebroids appeared more recently in \cite{gracia-saz:vb}, although they are equivalent to \emph{LA-vector bundles} defined in \cite{mackenzie:double}. For example, the tangent and the cotangent bundle of a Lie groupoid are VB-groupoids. Notice that several geometric structures on manifolds can be understood as vector bundle maps between tangent and cotangent bundles, so we can think about the corresponding geometric structures on groupoids as VB-groupoid maps. Moreover, VB-groupoids provide more intrinsic models for the so-called \emph{2-term representations up to homotopy} \cite{abad:ruth2, abad:ruth}. They never appeared in the context of Lie groups and were introduced in order to find a good framework for the adjoint and the coadjoint representations of a Lie groupoid. Now we know that these two representations are only defined up to some choices: the corresponding canonical objects are precisely the tangent and the cotangent VB-groupoids. This discussion translates literally to the realm of Lie algebroids.

In particular, in this work we will be interested in the deformation theory of VB-groupoids and VB-algebroids: most of the results come from the papers \cite{lapastina:def} and \cite{lapastina:def2}, wrote in collaboration with Luca Vitagliano. We will follow the - now classical - principle that states that every reasonable deformation problem should be governed by a suitable differential graded Lie algebra (dg-Lie algebra or DGLA for short). This principle goes back to the work of Kodaira-Spencer on deformations of complex structures (see \cite{kodaira:def} for a comprehensive review) and Nijenhuis-Richardson on deformations of Lie algebras \cite{nijenhuis:coh}, and has been recently formalized in a celebrated theorem by Lurie \cite{lurie:dagx}. While the deformation theory of Lie groups and Lie algebras is classical \cite{nijenhuis:coh, nijenhuis:def, nijenhuis:def_lie}, deformations of Lie groupoids and algebroids have been studied more recently in the papers \cite{crainic:def2} by Crainic-Mestre-Struchiner and \cite{crainic:def} by Crainic-Moerdijk, respectively.

The thesis is divided in 4 chapters, the first one being this introduction. The presentation of the material is deliberately concise: we decided not to include some direct computations and to skip some details, for which we often refer to the existing literature.

In Chapter \ref{chap:prel} we recall the notions needed to develop our theory, assuming only a basic knowledge of differential geometry and, in Subsection \ref{sec:graded}, of sheaf theory. This review part is quite large, in order to keep the presentation as self-contained as possible. In Section \ref{Sec:homogeneity}, we recall the concept of homogeneity structure of a vector bundle and of linear ``objects'' (polyvector fields, differential forms) on the total space of a vector bundle: these ideas will play an essential role in the thesis. In Section \ref{sec:lie}, we define Lie groupoids and algebroids in details, giving a lot of examples and briefly discussing the integration problem. After that, we introduce representations up to homotopy of Lie groupoids and algebroids and give a partial answer to the problem of their geometric descriptions. In the last subsection, we present the formalism of graded manifolds. They are particularly interesting for us because they provide a very elegant description of Lie algebroids: indeed, the latter are equivalent to special kinds of \emph{DG-manifolds}. In Section \ref{sec:def} we deal with deformations of Lie algebroids and groupoids. Finally, in Section \ref{sec:VB} we introduce VB-algebroids and VB-groupoids and discuss their structure.

Our treatment of the deformation theory of Lie algebroids and Lie groupoids requires some explanations. Given a Lie algebroid $A \Rightarrow M$, we decided to call \emph{deformation} of $A$ any Lie algebroid structure on the vector bundle $A \to M$. It is well known that deformations of smooth vector bundles are trivial, so it is not restrictive to deform only the Lie algebroid structure. In \cite{crainic:def}, Crainic and Moerdijk introduced the DGLA $C_{\mathrm{def}}(A)$ and proved that deformations of $A$ are in one-to-one correspondence with Maurer-Cartan elements of $C_{\mathrm{def}}(A)$. We observe, moreover, that two Lie algebroid structures on $A \to M$ are equivalent, i.e.~are connected by a smooth path of Lie algebroid structures, if and only if the corresponding Maurer-Cartan elements are gauge-equivalent: thus the DGLA $C_{\mathrm{def}}(A)$ controls deformations of $A$ in the usual sense.

The situation is not so easy in the case of Lie groupoids. In \cite{crainic:def2} it is shown that, given a Lie groupoid $\mathcal G \rightrightarrows M$, its structure can be described only in terms of the source map $\mathsf s: \mathcal G \to M$ and the division map $\bar{\mathsf m}: \mathcal G \times_{\mathsf s} \mathcal G \to \mathcal G$. However, in general it is not true that deformations of the surjective submersion $\mathsf s$ are trivial, so one has to include the possibility that $\mathcal G$ and $M$ themselves are deformed. Even if $\mathcal G$ and $M$ stay fixed, a variation of $\mathsf s$ forces the domain of $\bar {\mathsf m}$ to change. It is difficult, then, to adapt the previous approach to this case. Instead, following \cite{crainic:def2}, we study \emph{families of Lie groupoids} $\tilde{\mathcal G} \rightrightarrows \tilde M \to I$ over an interval $I \subset \mathbb R$ containing $0$, and we call \emph{deformation} of $\mathcal G$ a family whose fiber over $0$ coincides with $\mathcal G$. We notice explicitly that, unlike the case of Lie algebroids, there is no known DGLA structure on the deformation complex of a Lie groupoid, so only infinitesimal deformations can be treated fully satisfactorily. Finding such a DGLA structure and unifying the two approaches would be a relevant research direction for the future.

Chapters \ref{chap:VB_alg} and \ref{chap:VB_gr} are about deformations of VB-algebroids and VB-groupoids, respectively. They are both divided in two sections: the first one discusses the general theory, the second one examples and applications. The approaches we followed here are analogous to the ones used for Lie algebroids and groupoids, for the same reason we discussed before. Given a VB-algebroid $(W \Rightarrow E; A \Rightarrow M)$, we show that its deformations are controlled by a subcomplex $C_{\mathrm{def,lin}}(W)$ of the deformation complex of $W \Rightarrow E$, defined via a certain ``linearity'' condition. In the case of a VB-groupoid $(\mathcal W \rightrightarrows E; \mathcal G \rightrightarrows M)$, with the same method one finds a subcomplex $C_{\mathrm{def,lin}}(\mathcal W) \subset C_{\mathrm{def}}(\mathcal W)$. The \emph{linear deformation complexes} of a VB-algebroid and a VB-groupoid were first introduced in \cite{etv:infinitesimal} (for different purposes from the present ones). We analyze the meaning of the first cohomology groups and the behaviour of these invariants with respect to duality. When a VB-groupoid $(\mathcal W \rightrightarrows E; \mathcal G \rightrightarrows M)$ has trivial core, the linear deformation complex has an easier description that we study in detail. We show that if, moreover, $\mathcal G$ is proper, the VB-groupoid $\mathcal W$ is cohomologically rigid. We are currently studying the case in which the core is not trivial, and we expect a similar result to be true.

Next, we show that the linear deformation cohomology of a VB-algebroid (VB-groupoid) is included in the deformation cohomology of the top algebroid (groupoid). To do this, we need to introduce the \emph{linearization map}, an important technical tool, adapting ideas from \cite{cabrera:hom}. Subsection \ref{sec:van_est} is dedicated to the proof of a van Est theorem for the linear deformation cohomology of VB-groupoids. We show that the linear deformation complex of a VB-groupoid and that of the associated VB-algebroid are intertwined by a van Est map, which is a quasi-isomorphism under certain connectedness conditions. Finally, in Subsection \ref{sec:morita} we show that the linear deformation cohomology of a VB-groupoid is Morita invariant. This is particularly important, because it means that this cohomology is really an invariant of the associated vector bundle of differentiable stacks. 

In the example sections, we show that several classical objects in algebra and differential geometry can be understood in the framework of VB-algebroids and VB-groupoids. For example, vector bundles in the category of Lie algebras (Lie groups) are equivalent to Lie algebra (Lie group) representations, and we show that indeed the linear deformation cohomologies and the classical deformation cohomologies by Nijenhuis and Richardson are related by a long exact sequence. Considering Lie algebroids (groupoids) in the category of vector spaces, we are led to \emph{LA-vector spaces} and \emph{2-vector spaces} \cite{baez:higher}. Deformations of the tangent and cotangent VB-algebroids (VB-groupoids) of a Lie algebroid (groupoid) are also discussed. More geometric examples are representations of foliations groupoids, (linear) Lie group actions on vector bundles and their infinitesimal counterparts. In the case of VB-algebroids, we have two more examples: deformations of type-1 VB-algebroids and deformations of general representations up to homotopy. In particular, in this last paragraph the description in terms of graded geometry plays a central role.

Finally, let us have a look to unsolved problems related to this thesis, that can lead to further research directions. The most interesting problem is certainly the one of finding a DGLA (or $L_\infty$) structure on the deformation complex of a Lie groupoid. Surprisingly, even in the case of Lie groups the answer is unknown. Moreover, by now there is no geometric description of general representations up to homotopy of a Lie groupoid in the literature: this issue is being addressed in the work in progress \cite{delhoyo:simp} by del Hoyo and Trentinaglia. In the case of Lie algebroids these problems can be treated, in a very natural way, by using the language of graded geometry. On the other hand, in \cite{severa:int} it is shown that a general DG-manifold integrates to a Kan simplicial manifold: one might expect, then, that developing calculus and geometry in the category of simplicial manifolds would lead to a solution in the global setting.

Another possibility would be that of studying deformations of (compatible) geometric structures on Lie groupoids. The PhD thesis \cite{cardenas:thesis} goes in this direction, studying deformations of symplectic groupoids. Moreover the author, in collaboration with João Nuno Mestre and Luca Vitagliano, is working on a project on deformations of complex groupoids. Hopefully, this will shed new light on deformations of singular spaces that often arise in complex and algebraic geometry.

\chapter{Preliminaries} \label{chap:prel}

In this chapter, we will introduce the main characters of this work: Lie groupoids, Lie algebroids and vector bundles over them. We will see (Subsections \ref{sec:lie_grpd} and \ref{sec:lie_algd}) that Lie groupoids and algebroids arise in several situations in differential geometry, including Lie groups, Lie group actions, principal bundles, Poisson manifolds, so it is worth to study them abstractly. Lie algebroids are well understood in the language of \emph{graded geometry}, which we recall in Subsection \ref{sec:graded}. Moreover, in Subsection \ref{sec:rep} we give a brief account of \emph{representations up to homotopy}, algebraic objects associated to Lie algebroids and groupoids that generalize Lie algebra and Lie group representations. Their geometric counterparts are, in the easiest cases, \emph{VB-algebroids} and \emph{VB-groupoids} (Subsections \ref{sec:VB-alg}, \ref{sec:VB-gr}). In general, representations up to homotopy of Lie algebroids correspond to special kinds of DG-vector bundles, while representations up to homotopy of Lie groupoids are being studied by del Hoyo and Trentinaglia \cite{delhoyo:simp}. In order to describe VB-algebroids and VB-groupoids, we will recall a very useful technical tool, the \emph{homogeneity structure of a vector bundle} (Section \ref{Sec:homogeneity}).

Finally, in Subsections \ref{sec:def_algd} and \ref{sec:def_grpd} we deal with deformations of Lie algebroids and groupoids. This theory generalizes the classical treatment of deformations of Lie algebras and Lie groups, as well as other deformation theories from the literature, and is the starting point for our work on VB-algebroids and VB-groupoids.

\section{The homogeneity structure of a vector bundle} \label{Sec:homogeneity}

The concept of \emph{homogeneity structure of a vector bundle} originally appeared in \cite{grab:vect}. Using this tool, it is possible to define linear polyvector fields and differential forms on the total space of a vector bundle. These ideas appeared (probably for the first time, e.g.) in \cite{iglesias:univ} and \cite{bursztyn:mult}, respectively. In the last two references, one can also find the proofs of (a version of) Propositions \ref{prop:lin_forms} and \ref{prop:lin_multi}. With respect to those references, we will offer just a slightly different point of view, in order to make the presentation consistent.

Let $p: E \to M$ be a vector bundle. The monoid $\mathbb R_{\geq 0}$ of non-negative real numbers acts on $E$ by homotheties $h_\lambda: E \to E$ (fiber-wise scalar multiplication). The action $h : \mathbb R_{\geq 0} \times E \to E$, $(\lambda, e) \mapsto h_\lambda (e)$, is called the \emph{homogeneity structure} of $E$. The homogeneity structure (together with the smooth structure) fully characterizes the vector bundle structure \cite{grab:vect}. In particular, it determines the addition. This implies that \emph{every notion that involves the linear structure of $E$ can be expressed in terms of $h$ only}: for example, a smooth map between the total spaces of two vector bundles is a bundle map if and only if it commutes with the homogeneity structures.

The homogeneity structure isolates distinguished subspaces in the algebras $\Omega (E)$ of differential forms and $\mathfrak X_{\mathrm{poly}} (E)$ of polyvector fields on the total space $E$ of the vector bundle.

\begin{definition} A differential form $\omega \in \Omega (E)$ is (\emph{homogeneous}) \emph{of weight $q$} if 
\begin{equation}\label{eq:hom_forms}
h_\lambda^* \omega = \lambda^q \omega
\end{equation}
for all $\lambda > 0$. 
A polyvector field $X \in \mathfrak X_{\mathrm{poly}} (E)$ is (\emph{homogeneous}) \emph{of weight $q$} if
\begin{equation}\label{eq:hom_multivect}
h_\lambda^* X = \lambda^q X
\end{equation}
for all $\lambda > 0$. 
\end{definition}

Clearly, for $q \geq 0$, weight $q$ functions coincide with functions on $E$ that are fiber-wise homogeneous polynomial of degree $q$, while there are no non-zero functions of weight $q < 0$. This discussion implies that, if Formula (\ref{eq:hom_forms}) holds for a form $\omega$ of a degree $0$, for any $\lambda > 0$, then it actually holds for all $\lambda \in \mathbb R$. Moreover, this is true for forms of any degree, because the coefficients in local coordinates of a weight $q$ form are fiber-wise polynomial functions of degree at most $q$. For polyvector fields, the same is true for $\lambda \neq 0$, while the pull-back $h_0^*$ is not defined on polyvector fields because $h_0$ is not a diffeomorphism.

From the functorial properties of the pull-back, it follows that the grading defined by the weight is natural with respect to all the usual operations on functions, forms and (multi)vector fields. From this remark, we easily see that \emph{there are no non-zero $k$-vector fields of weight less than $-k$} and \emph{there are no non-zero $k$-forms of weight less than $0$.}

\begin{definition} \label{def:linear}
A function (resp.~a $k$-form) on $E$ is \emph{linear} if it is of weight 1, it is \emph{core} if it is of weight 0. A $k$-vector field is \emph{linear} if it is of weight $1-k$, it is \emph{core} if it is of weight $-k$. We denote by $C^\infty_{\mathrm{lin}}(E)$, $\Omega_{\mathrm{lin}}(E)$, $\mathfrak X_{\mathrm{lin}}(E)$ and $\mathfrak X_{\mathrm{poly, lin}}(E)$ the spaces of linear functions, differential forms, vector fields and polyvector fields respectively, by $C^\infty_{\mathrm{core}}(E)$, $\Omega_{\mathrm{core}}(E)$, $\mathfrak X_{\mathrm{core}}(E)$ and $\mathfrak X_{\mathrm{poly, core}}(E)$ the spaces of core functions, differential forms, vector fields and polyvector fields respectively.  
\end{definition}

\begin{remark} Linear functions are precisely fiber-wise linear functions. While the terms ``linear function'' and ``linear polyvector field'' are already present in the literature, the terms ``core function'' and ``core polyvector field'' are being used here probably for the first time. We use them in analogy to core sections of double vector bundles, as defined in Subsection \ref{sec:DVB}. \end{remark}

The definition of linear polyvector field may sound a little strange, but it is motivated (among other things) by the following proposition.

\begin{proposition}\label{prop:lin_multi} 
Let $X \in \mathfrak X^k_{\mathrm{poly}} (E)$. The following conditions are equivalent:
\begin{enumerate} 
\item $X$ is linear;
\item $X$ takes
	\begin{enumerate}
	\item $k$ linear functions to a linear function,
	\item $k-1$ linear functions and a core function to a core function,
	\item $k-i$ linear functions and $i$ core functions to $0$, for every $i \geq 2$;
	\end{enumerate}
\item If $(x^i)$ are local coordinates on $M$ and $(u^\alpha)$ are linear fiber coordinates on $E$, $X$ is locally of the form
\begin{equation}\label{eq:mult_lin_coord}
X = X^{i \alpha_2 \dots \alpha_k} (x) \dfrac{\partial}{\partial x^i} \wedge \dfrac{\partial}{\partial u^{\alpha_2}} \wedge \dots \wedge \dfrac{\partial}{\partial u^{\alpha_k}} + X^{\alpha_1 \dots \alpha_k}_\beta (x) u^\beta \dfrac{\partial}{\partial u^{\alpha_1}} \wedge \dots \wedge \dfrac{\partial}{\partial u^{\alpha_k}}.
\end{equation}
\end{enumerate}
\end{proposition}

\proof Condition (1) clearly implies (2) because the action of $X$ on functions respects the weights. Condition (2) implies (3) because $x^i$ is a core function and $u^\beta$ is a linear. Finally, (3) implies (1): just check (\ref{eq:hom_multivect}) for the polyvector field (\ref{eq:mult_lin_coord}). \endproof

In the same way, one can prove:

\begin{proposition}\label{prop:lin_forms} 
Let $\omega \in \Omega^k (E)$, $k \geq 1$. The following conditions are equivalent:
\begin{enumerate} 
\item $\omega$ is linear;
\item $\omega$ takes
	\begin{enumerate}
	\item $k$ linear vector fields to a linear function,
	\item $k-1$ linear vector fields and a core vector field to a core function,
	\item $k-i$ linear vector fields and $i$ core vector fields to $0$, for every $i \geq 2$;
	\end{enumerate}
\item If $(x^i)$ are local coordinates on $M$ and $(u^\alpha)$ are linear fiber coordinates on $E$, $\omega$ is locally of the form
\[ 
\omega = \omega_{i_1 \dots i_k \beta} (x) u^\beta dx^{i_1} \wedge \dots \wedge dx^{i_k} + \omega_{i_1 \dots i_{k-1} \alpha} (x) dx^{i_1} \wedge \dots \wedge dx^{i_{k-1}} \wedge du^\alpha.
\]
\end{enumerate}
\end{proposition}

\begin{remark}\label{rmk:lin} In most cases, condition (2) in Propositions \ref{prop:lin_multi} and \ref{prop:lin_forms} can be simplified. Actually, if the rank of $E$ is greater than $0$, then it can be proved (by induction on $k$) that (2).(a) implies (2).(b) and (2).(c), so linear polyvector fields (resp.~differential forms) are simply those that send linear functions (resp.~vector fields) to a linear function. If $\operatorname{rank} E = 0$, then $E = 0_M$ is the zero vector bundle over $M$, $C^\infty_{\mathrm{lin}}(E) = 0$, $C^\infty_{\mathrm{core}}(E) = C^\infty (E)$ and $\mathfrak X_{\mathrm{lin}}(E) = \mathfrak X (E)$. It follows that condition (2).(a) in Proposition \ref{prop:lin_multi} is always satisfied, but every $k$-vector field takes $k$ core functions to a core function, so (2).(c) is satisfied if and only if $k = 1$. As for Proposition \ref{prop:lin_forms}, condition (2).(a) implies $\omega = 0$, so (2).(b) and (2).(c) follow trivially. 

From the coordinate expression (\ref{eq:mult_lin_coord}), it is clear that a linear $k$-vector field $X$ is completely determined by its action on $k$ linear functions and on $k-1$ linear functions and a core function. In the same way, one can prove that, if $\operatorname{rank} E > 0$, then $X$ is fully determined by its action on $k$ linear functions only. A similar property holds for linear differential forms.
\end{remark}

\begin{remark} \label{rmk:lin_vect} From the definition, it follows that a vector field $X$ on $E$ is linear if and only if it is a vector bundle morphism
\[
\xymatrix{
E \ar[d]^p \ar[r]^X & TE \ar[d]^{Tp} \\
M \ar[r] & TM}
\]
over a vector field on $M$. 

Now we want to describe core vector fields on $E$. Recall that there is a canonical isomorphism of $C^\infty (M)$-modules
\[
\Gamma(E^*) \overset{\cong}{\longrightarrow} C^\infty_{\mathrm{lin}}(E), \quad \varphi \mapsto \ell_\varphi,
\]
given by $\ell_\varphi (e) = \varphi_{p(e)} (e)$ for every $e \in E$. Then, one can verify that the $C^\infty (M)$-module $\mathfrak X_{\mathrm{core}}(E)$ is canonically isomorphic to $\Gamma(E)$. The isomorphism is given by
\begin{equation} \label{eq:vert_lift}
\Gamma(E) \overset{\cong}{\longrightarrow} \mathfrak X_{\mathrm{core}}(E), \quad \varepsilon \mapsto \varepsilon^\uparrow,
\end{equation}
where $\varepsilon^\uparrow$ is uniquely defined by $\varepsilon^\uparrow (\ell_\varphi) = p^* \langle \varphi, \varepsilon \rangle$. Here and in the following, $\langle \cdot, \cdot \rangle$ denotes the duality pairing. The vector field $\varepsilon^\uparrow$ is said to be the \emph{vertical lift} of $\varepsilon$.

The concepts of linear and core vector fields will be generalized in Subsection \ref{sec:DVB}, when we will talk about double vector bundles. \end{remark}

\begin{remark} \label{rmk:sn} Let $M$ be a smooth manifold. We recall that $\mathfrak X_{\mathrm{poly}}(M)[1]$ is a graded Lie algebra, when equipped with the \emph{Schouten-Nijenhuis bracket} $\llbracket -,- \rrbracket$. The latter is explicitly given by
\begin{equation} \label{eq:sn}
\begin{aligned}
\llbracket X_1 \wedge \dots \wedge X_k, Y_1 \wedge \dots \wedge Y_l \rrbracket = & \sum_{i,j} (-1)^{i+j} [X_i, Y_j] \wedge X_i \wedge \dots \wedge \widehat{X_i} \wedge \dots \wedge X_k \\
& \wedge Y_1 \wedge \dots \wedge \widehat{Y_j} \wedge \dots \wedge Y_l; \\
\llbracket X_1 \wedge \dots \wedge X_k, f \rrbracket = & \sum_i (-1)^{k-i} X_i (f) X_1 \wedge \dots \wedge \widehat{X_i} \wedge \dots \wedge X_k
\end{aligned}
\end{equation}
for every $X_1, \dots, X_k, Y_1, \dots, Y_l \in \mathfrak X (M)$, $f \in C^\infty(M)$.

Now, if $E \to M$ is a vector bundle, from Equation (\ref{eq:sn}) it is clear that linear polyvector fields form a graded Lie subalgebra $\mathfrak X_{\mathrm{poly, lin}}(E)$ of $\mathfrak X_{\mathrm{poly}}(E)$. \end{remark}

\section{Lie groupoids and Lie algebroids} \label{sec:lie}

\subsection{Lie groupoids} \label{sec:lie_grpd}

\begin{definition} A \emph{groupoid} $\mathcal G$ is a small category in which every arrow is invertible. \end{definition}

Every groupoid $\mathcal G$ comes with a set of arrows and a set of objects. Usually, the set of arrows is again denoted $\mathcal G$. If $X$ is the set of objects, we say that $\mathcal G$ is a \emph{groupoid over $X$} and we call $X$ the \emph{base} of $\mathcal G$. Moreover, $\mathcal G$ comes with some \emph{structure maps}:

\begin{itemize}
\item The \emph{source} and the \emph{target maps} $\mathsf{s, t}: \mathcal G \to X$ send each arrow to its source and target object. An arrow $g \in \mathcal G$ with source $x$ and target $y$ is depicted $g: x \to y$ or $x \overset{g}{\to} y$.
\item The \emph{composition} (or \emph{multiplication}) $\mathsf m: \mathcal G^{(2)} \to \mathcal G$ is defined on the set
\[
\mathcal G^{(2)} := \{ (g,h) \in \mathcal G^2: \mathsf s(g) = \mathsf t(h) \}
\]
of pairs of composable arrows and sends a pair $(g,h)$ to its composition $g \circ h = gh$.
\item The \emph{unit map} $\mathsf 1: X \to \mathcal G$ sends every object $x$ to the identity arrow $\mathsf 1_x: x \to x$.
\item The \emph{inverse map} $\mathsf i: \mathcal G \to \mathcal G$ sends an arrow $g$ to its inverse $g^{-1}$.
\end{itemize}
Finally, these maps have to satisfy the following axioms:
\begin{itemize}
\item if $x \overset{g}{\to} y \overset{h}{\to} z$, then $x \overset{hg}{\to} z$ (\emph{composition law});
\item if $x \overset{g}{\to} y \overset{h}{\to} z \overset{k}{\to} w$, then $(kh)g = k(hg)$ (\emph{associativity law});
\item for every $g: x \to y$, $g = g \mathsf 1_x = \mathsf 1_y g$ (\emph{law of units});
\item for every $g: x \to y$, $g^{-1} g = \mathsf 1_x$, $g g^{-1} = \mathsf 1_y$ (\emph{law of inverses}).
\end{itemize}

A groupoid $\mathcal G$ over $X$ is denoted $\mathcal G \rightrightarrows X$, the two arrows standing for the source and the target maps. For every $k > 0$, denote
\[
\mathcal G^{(k)} = \{ (g_1, \dots, g_k) \in \mathcal G^k: \mathsf s (g_i) = \mathsf t (g_{i+1}) \ \mathrm{for \ every} \ i < k \}
\]
the set of $k$-tuples of composable arrows of $\mathcal G$ and
\[
\mathcal G^{[k]} = \{ (g_1, \dots, g_k) \in \mathcal G^k: \mathsf s (g_i) = \mathsf s (g_{i+1}) \ \mathrm{for \ every} \ i < k \}
\]
the set of $k$-tuples of arrows of $\mathcal G$ with the same source. We will also denote $\mathsf s: \mathcal G^{[k]} \to X$ the projection to the common source. Finally, the \emph{division map} is the map
\[
\bar{\mathsf m}: \mathcal G^{[2]} \to \mathcal G, \quad (g,h) \mapsto gh^{-1}.
\]

\begin{definition} A \emph{Lie groupoid} is a groupoid $\mathcal G \rightrightarrows M$ where $\mathcal G$ and $M$ are smooth manifolds, $\mathsf s$ and $\mathsf t$ are submersions and all the structure maps are smooth. \end{definition}

The regularity condition on $\mathsf s$ and $\mathsf t$ is needed in order to ensure that $\mathcal G^{(2)}$ is also a smooth manifold, hence the definition makes sense. One can also show that the unit map $\mathsf 1: M \to \mathcal G$ is an embedding, and $M$ is usually identified with the submanifold of units of $\mathcal G$.

Now we give some standard notation and terminology. For every $x \in M$, the submanifolds
\[
\mathcal G (x,-) := \mathsf s^{-1}(x), \quad \mathcal G (-,x) := \mathsf t^{-1}(x)
\]
are called the \emph{$\mathsf s$-fiber} and the \emph{$\mathsf t$-fiber} at $x$, respectively. If $g: x \to y$, the \emph{right multiplication by $g$} is defined by
\[
R_g: \mathcal G (y,-) \to \mathcal G (x,-), \quad h \mapsto hg
\]
and the \emph{left multiplication by $g$} is defined by
\[
L_g: \mathcal G (-,x) \to \mathcal G (-,y), \quad h \mapsto gh.
\]

Since groupoids are special kinds of categories, a \emph{morphism of groupoids} $\mathcal G_1 \to \mathcal G_2$ is simply a functor. A morphism of Lie groupoids, of course, is also required to be smooth. Explicitly: 

\begin{definition} Let $\mathcal G_1 \rightrightarrows M_1$, $\mathcal G_2 \rightrightarrows M_2$ be Lie groupoids. A \emph{morphism of Lie groupoids} consists of smooth maps $\Phi: \mathcal G_1 \to \mathcal G_2$, $f: M_1 \to M_2$ such that:
\begin{itemize}
\item if $g: x \to x'$, then $\Phi(g): f(x) \to f(x')$;
\item if $(g,g') \in \mathcal G_1^{(2)}$, then $\Phi(gg') = \Phi(g) \Phi(g')$;
\item for every $x \in M_1$, $\Phi(\mathsf 1_x) = \mathsf 1_{f(x)}$;
\item for every $g \in \mathcal G_1$, $\Phi(g^{-1}) = \Phi(g)^{-1}$.
\end{itemize}
We also say that $\Phi$ is a morphism of Lie groupoids over $f$. Such a morphism is represented in the following way:
\[
\begin{array}{r}
\xymatrix{
\mathcal G_1 \ar[r]^\Phi \ar@<0.4ex>[d] \ar@<-0.4ex>[d] & \mathcal G_2 \ar@<0.4ex>[d] \ar@<-0.4ex>[d] \\
M_1 \ar[r]^f & M_2}
\end{array}.
\]
If $\Phi$ and $f$ are injective immersions, we say that $\mathcal G_1$ is a \emph{Lie subgroupoid} of $\mathcal G_2$.
\end{definition}

Besides the obvious notion of isomorphism of Lie groupoids, there is a weaker notion of equivalence that will be very important for us: \emph{Morita equivalence}.

\begin{definition}\label{def:Morita_map} A morphism of Lie groupoids $\Psi: \mathcal G_1 \to \mathcal G_2$ over a smooth map $f: M_1 \to M_2$ is a \emph{Morita map} (or a \emph{weak equivalence}) if it is
\begin{enumerate}
\item \emph{fully faithful}, i.e.~the diagram
\[
\begin{array}{c}
\xymatrix@C=45pt{\mathcal G_1 \ar[r]^{\Psi} \ar[d]_{(\mathsf s_1, \mathsf t_1)} & \mathcal G_2 \ar[d]^{(\mathsf s_2, \mathsf t_2)} \\
M_1 \times M_1 \ar[r]_{f \times f} & M_2 \times M_2}
\end{array}
\]
is a pull-back diagram, and
\item \emph{essentially surjective}, i.e.~the map $\mathcal G_2 \tensor*[_{\mathsf s_1}]{\times}{_f}  M_1 \to M_2, (g,x) \mapsto \mathsf t_2(g)$ is a surjective submersion.
\end{enumerate}

Two groupoids $\mathcal G_1$ and $\mathcal G_2$ are said to be \emph{Morita} (or \emph{weakly}) \emph{equivalent} if there exist a Lie groupoid $\mathcal H$ and Morita maps $\Psi_1: \mathcal H \to \mathcal G_1$, $\Psi_2: \mathcal H \to \mathcal G_2$.
\end{definition}

\begin{remark} Notice that Definition \ref{def:Morita_map}.(1) of fully faithful morphism is slightly different from the one in \cite{delhoyo:morita}, where an additional property is required. However, it is easy to see that for an essentially surjective morphism the two definitions are equivalent. In the literature, Morita equivalence is often expressed in terms of principal bibundles: we refer to \cite{mrcun:stability} for this notion and many more details. \end{remark}

Morita equivalent Lie groupoids share many interesting properties: we will briefly come back to this point in Subsection \ref{sec:morita}. Moreover, there is a more conceptual reason why this notion is so important. In \cite{behrend:diff}, stacks in the category of smooth manifolds, namely \emph{differentiable stacks}, are introduced, and it is proved that the category of differentiable stacks is equivalent to that of Lie groupoids up to Morita equivalence. This means that every Morita invariant is in fact an invariant of the associated stack, and one can study the geometry of stacks by means of their representative groupoids, which are usually more manageable. 

\

If $\mathcal G \rightrightarrows M$ is a Lie groupoid, there is an induced equivalence relation on $M$. If $x, y \in M$, we say that $x$ is equivalent to $y$ if there is an arrow $g: x \to y$. The equivalence class of $x$ is denoted $\mathcal O_x$ and called the \emph{orbit through $x$}, while the quotient set is denoted $M/ \mathcal G$ and called the \emph{orbit set of $\mathcal G$}. One can prove \cite{crainic:lect} that each orbit is an immersed submanifold of $M$. If the orbit set is just a point, we say that $\mathcal G$ is \emph{transitive}. If all the orbits have the same dimension, we say that $\mathcal G$ is \emph{regular}.

For every $x \in M$, the set of arrows of $\mathcal G$ starting and ending at $x$ is clearly a group. In \cite{mackenzie} it is proved that it is also a Lie group, called the \emph{isotropy group of $\mathcal G$ at $x$} and denoted $\mathcal G_x$.




\begin{example}[Lie groups] Lie groupoids $G \rightrightarrows *$ with just one object are exactly Lie groups. In fact, a Lie groupoid $\mathcal G \rightrightarrows M$ can be considered as a manifold with a ``partially defined'' smooth multiplication $\mathsf m: \mathcal G^{(2)} \to \mathcal G$, and $\mathsf m$ is defined on the whole $\mathcal G \times \mathcal G$ if and only if $\mathcal G$ is a Lie group. \end{example}

\begin{example}[Bundles of Lie groups] \label{ex:bundle} More generally, take a Lie groupoid $\mathcal G \rightrightarrows M$ where the source and the target maps coincide, and let us call this map $\pi$. Then each fiber of $\pi$ inherits a Lie group structure, and $\pi: \mathcal G \to M$ is called a \emph{bundle of Lie groups}. \end{example}

\begin{example}[Pair groupoid] Let $M$ be a manifold. The \emph{pair groupoid of $M$} is $M \times M \rightrightarrows M$, where a pair $(x,y)$ is seen as an arrow from $x$ to $y$ and the composition is given by:
\[
(y,z) \cdot (x,y) = (x,z).
\] 
Clearly, one has $\mathsf 1_x = (x,x)$ and $(x,y)^{-1} = (y,x)$. For every Lie groupoid $\mathcal G \rightrightarrows M$, there is a natural morphism of Lie groupoids
\[
\mathcal G \to M \times M, \quad g \mapsto (\mathsf s (g), \mathsf t (g)).
\]
\end{example}

\begin{example}[Action groupoid] \label{ex:action} Let $M$ be a manifold and $G$ a Lie group acting on $M$. Then we can construct the \emph{action groupoid} $G \ltimes M \rightrightarrows M$ as follows. As a set, $G \ltimes M$ is simply the product $G \times M$, and the structure maps are:
\begin{equation} \label{eq:action_grpd}
\begin{aligned}
\mathsf s(g,x) & = x, \\
\mathsf t(g,x) & = gx, \\
(h, gx) \cdot (g, x) & = (hg, x), \\
\mathsf 1_x & = (1,x), \\
(g,x)^{-1} & = (g^{-1}, gx).
\end{aligned}
\end{equation}
\end{example}

\begin{example}[Gauge groupoid] \label{ex:gauge_grpd} Let $G$ be a Lie group and $P \to M$ be a principal $G$-bundle. The \emph{gauge groupoid} of $P$ is the quotient of the pair groupoid $P \times P \rightrightarrows P$ by the diagonal action of $G$, and is denoted $P \otimes_G P \rightrightarrows M$. For every $x, y \in M$, an arrow from $x$ to $y$ in $P \otimes_G P$ can be identified with a $G$-equivariant smooth map $P_x \to P_y$. \end{example}

\begin{example}[Fundamental groupoid of a manifold] Let $M$ be a connected manifold. The \emph{fundamental groupoid of $M$} consists of all homotopy classes of paths in $M$ with fixed end points and is denoted $\Pi_1 (M)$. It can be seen as the gauge groupoid of the universal cover $\tilde M$ of $M$, that is a $\pi_1(M)$-principal bundle, hence it supports a Lie groupoid structure. The isotropy groups of $\Pi_1(M)$ are exactly the fundamental groups of $M$. \end{example}

\begin{example}[Monodromy groupoid of a foliation]\label{ex:foliation} Let $M$ be a manifold endowed with a foliation $\mathcal F$. The \emph{monodromy groupoid of $\mathcal F$}, denoted $\operatorname{Mon}(M, \mathcal F)$, consists of all leafwise homotopy classes of paths in $M$ with fixed end points. More precisely, if $x, y \in M$ belong to the same leaf $F$, then arrows $x \to y$ are smooth curves from $x$ to $y$ in $F$ modulo smooth homotopy in $F$, otherwise there is no arrow from $x$ to $y$. 

Out of a foliation, one can also construct the \emph{holonomy groupoid}. There is a deep connection between foliation theory and Lie groupoids, which is explored in detail in \cite{moerdijk:intro}. The leaf space of a foliation is rarely a manifold, can be very pathological (e.g., non-Hausdorff) and the topology alone forgets about the pre-existing smooth structure. Moreover, every leaf has its own geometry. So the quotient does not retain all the information, while the monodromy and holonomy groupoids do, at the cost of increased complexity and dimension. \end{example}

\begin{example}[General linear groupoid] Let $E \to M$ be a vector bundle. The \emph{general linear groupoid} of $E$ is a Lie groupoid $GL(E) \rightrightarrows M$, whose arrows from $x$ to $y$ are all linear isomorphisms $E_x \to E_y$. 

Now, let $\mathrm{Fr}(E) \to M$ be the frame bundle of $E$. If the rank of $E$ is $k$, then $\mathrm{Fr}(E)$ is a principal $GL(k, \mathbb R)$-bundle and its gauge groupoid is isomorphic to $GL(E)$. \end{example}

\begin{example}[Tangent prolongation Lie groupoid] Let $\mathcal G \rightrightarrows M$ be a Lie groupoid. The \emph{tangent prolongation Lie groupoid} $T \mathcal G \rightrightarrows TM$ is obtained by applying the tangent functor to all the structure maps of $\mathcal G$. \end{example}

We end this paragraph recalling a classical construction that will be needed later. Let $\mathcal G \rightrightarrows M$ be a Lie groupoid. A smooth map $\beta: M \to \mathcal G$ is a \emph{bisection} of $\mathcal G$ if it is a section of $\mathsf s$ (i.e., $\mathsf s \circ \beta = \mathrm{id}_M$) and $\mathsf t \circ \beta$ is a diffeomorphism. The value of $\beta$ at $x \in M$ is denoted $\beta_x$.

The set $\mathrm{Bis}(\mathcal G)$ of bisections of $\mathcal G$ can be given a group structure as follows. If $\beta, \beta' \in \mathrm{Bis}(\mathcal G)$, then we set:
\[
(\beta * \beta')_x = \beta_{(\mathsf t \circ \beta')(x)} \cdot \beta'_x.
\]
The group of bisections can be useful when dealing with several issues related to Lie groupoids. We refer to \cite{mackenzie} for more information about that.

\subsection{Lie algebroids} \label{sec:lie_algd}

One of the main features of differential geometry is the possibility to study geometric objects by means of their infinitesimal analogues, which are usually easier to handle. The most classical example is given by Lie groups: they differentiate to Lie algebras, that can be studied with purely algebraic methods and can give a lot of information about the corresponding groups. 

We briefly recall this procedure. Given a Lie group $G$, first one observes that right-invariant vector fields $\mathfrak X^R (G)$ on $G$ form a Lie subalgebra of $\mathfrak X (G)$ and are fully determined by their values at the identity. Then one has a natural isomorphism between $\mathfrak g := T_e G$ and $\mathfrak X^R (G)$ and can transfer the Lie algebra structure on $\mathfrak g$: this is by definition the Lie algebra of $G$.

It is natural, then, to try to extend this process to a general Lie groupoid; in this way, one obtains a \emph{Lie algebroid}. However, Lie groupoids and algebroids are much more flexible, so in general one can not expect structure or classification theorems.

\begin{definition} A \emph{Lie algebroid} over a manifold $M$ consists of a vector bundle $A \to M$ together with a bundle map $\rho: A \to TM$, called the \emph{anchor map}, and a Lie bracket $[-,-]: \Gamma(A) \times \Gamma(A) \to \Gamma(A)$ that satisfy the following Leibniz rule:
\begin{equation} \label{eq:leibniz}
[\alpha, f \beta] = \rho (\alpha) (f) \beta + f [\alpha, \beta]
\end{equation}
for every $\alpha, \beta \in \Gamma(A)$, $f \in C^\infty(M)$. A Lie algebroid $A$ over $M$ is denoted $A \Rightarrow M$. \end{definition}

Take a Lie groupoid $\mathcal G \rightrightarrows M$. In this case, right translations are defined only on $\mathsf s$-fibers, so a right-invariant vector field should take values in the kernel of $T\mathsf s$, which we denote $T^{\mathsf s} \mathcal G$. We say that a vector field $X$ on $\mathcal G$ is \emph{right-invariant}, and we write $X \in \mathfrak X^R (\mathcal G)$, if $X \in \Gamma (T^{\mathsf s} \mathcal G)$ and 
\[
TR_h (X_g) = X_{gh}
\]
for every $(g,h) \in \mathcal G^{(2)}$.

It is clear now that such a vector field is determined by its values at the units. Explicitly, if we set $A := T^{\mathsf s} \mathcal G|_M$, there is a canonical isomorphism of $C^\infty (M)$-modules
\[
\Gamma(A) \overset{\cong}{\longrightarrow} \mathfrak X^R (\mathcal G), \quad \alpha \mapsto \overrightarrow \alpha,
\] 
given by $\overrightarrow \alpha_g = TR_g (\alpha_{\mathsf t(g)})$ for every $g \in \mathcal G$. Again, one can verify that $\mathfrak X^R (\mathcal G)$ is a Lie subalgebra of $\mathfrak X (\mathcal G)$ and transfer the bracket on $\Gamma(A)$.

The last step is the construction of the anchor map: this is defined by 
\[
\rho := T\mathsf t|_A : A \to TM. 
\]
The proof of the Leibniz rule (\ref{eq:leibniz}) is a simple computation that can be found, for example, in \cite{crainic:lect}. The Lie algebroid $A \Rightarrow M$ constructed in this way is called the \emph{Lie algebroid of $\mathcal G$} and is sometimes denoted $\operatorname{Lie}(\mathcal G)$. We can also say that \emph{$\mathcal G$ integrates $A$}, or \emph{$A$ integrates to $\mathcal G$}. 

\begin{remark} For every $\alpha \in \Gamma(A)$, $\overrightarrow \alpha$ is called the \emph{right-invariant vector field determined by $\alpha$}.

Similarly, let us denote $T^{\mathsf t} \mathcal G$ the kernel of $T \mathsf t: T \mathcal G \to TM$. We say that a vector field $X$ on $\mathcal G$ is \emph{left-invariant}, and we write $X \in \mathfrak X^L (\mathcal G)$, if $X \in \Gamma (T^{\mathsf t} \mathcal G)$ and
\[
TL_g (X_h) = X_{gh}
\]
for every $(g, h) \in \mathcal G^{(2)}$. There is a canonical isomorphism of $C^\infty (M)$-modules
\[
\Gamma(A) \overset{\cong}{\longrightarrow} \mathfrak X^L (\mathcal G), \quad \alpha \mapsto \overleftarrow \alpha,
\] 
given by $\overleftarrow \alpha_g = (TL_g \circ T \mathsf i)(\alpha_{\mathsf s (g)})$. Also $\mathfrak X^L (\mathcal G)$ is a Lie subalgebra of $\mathfrak X (\mathcal G)$, and the above isomorphism induces on $\Gamma(A)$ the same Lie algebra structure as before. 

Notice that, if $G$ is a Lie group, $\mathfrak g$ its Lie algebra and $X \in \mathfrak g$, then $\overleftarrow X_g = - TL_g (X)$ for every $g \in G$, so $\overleftarrow X$ is not what is usually called the ``left invariant vector field determined by $X$''. \end{remark}

We now come to morphisms of Lie algebroids. Recall that, if $E, F$ are vector bundles over $M$, a morphism of vector bundles $\phi: E \to F$ over $M$ induces a morphism $\phi_*: \Gamma(E) \to \Gamma(F)$ between the $C^\infty (M)$-modules of sections. The notion of morphism of Lie algebroids over the same base is straightforward.

\begin{definition}\label{def:morphism} Let $A_1 \Rightarrow M$, $A_2 \Rightarrow M$ be Lie algebroids. A \emph{morphism of Lie algebroids $\phi: A_1 \to A_2$ over $M$} is a morphism of vector bundles over $M$ that preserves the bracket and the anchor:
\begin{enumerate}
\item $\phi_* [\alpha, \alpha'] = [\phi_* \alpha, \phi_* \alpha']$ for every $\alpha, \alpha' \in \Gamma(A_1)$;
\item $\rho_{A_2} \circ \phi = \rho_{A_1}$.
\end{enumerate}
\end{definition}

For example, if $A \Rightarrow M$ is a Lie algebroid, then one can verify that the anchor map $\rho: A \to TM$ is a Lie algebroid map.

\begin{remark} If two Lie algebroids have different bases, the notion of Lie algebroid morphism between them is quite involved and is described, for example, in \cite{mackenzie}. In Subsection \ref{sec:graded} we will recall another description of Lie algebroids, that allows for a much easier definition of morphism. \end{remark}

Let $A \Rightarrow M$ be a Lie algebroid, $\rho: A \to TM$ its anchor, $x \in M$. Then $\ker (\rho_x)$ is naturally a Lie algebra, denoted $\mathfrak i(A)_x$ (or simply $\mathfrak i_x$, if $A$ is clear from the context) and called the \emph{isotropy Lie algebra of $A$ at $x$}. Its bracket $[-,-]_{\mathfrak i_x}$ is defined by
\[
[\alpha_x, \beta_x]_{\mathfrak i_x} := [\alpha, \beta]_x
\]
for every $\alpha, \beta \in \Gamma(A)$. If $A$ is the Lie algebroid of a Lie groupoid $\mathcal G \rightrightarrows M$, then $\mathfrak i_x$ is the Lie algebra of $\mathcal G_x$ for every $x \in M$. 

If $A$ is regular, i.e.~$\rho$ has constant rank, then $\mathfrak i = \ker \rho$ and $\nu(A) := TM/\operatorname{im} \rho$ are vector bundles over $M$, called the \emph{isotropy bundle} and the \emph{normal bundle} of $A$ respectively. The latter is denoted simply $\nu$, if $A$ is clear from the context. If $\mathcal G \rightrightarrows M$ is a regular Lie groupoid, its Lie algebroid $A \Rightarrow M$ is regular: by definition, the \emph{isotropy bundle} of $\mathcal G$ is the isotropy bundle of $A$ and the \emph{normal bundle} of $\mathcal G$ is the normal bundle of $A$.

\

Computing the Lie algebroids corresponding to Lie groupoids in Subsection \ref{sec:lie_grpd}, we will be able to give a lot of examples.

\begin{example}[Lie algebras] Obviously, Lie algebroids over a point are Lie algebras. We recall that, according to the classical \emph{third Lie theorem}, every finite-dimensional real Lie algebra integrates to a unique simply connected Lie group. \end{example}

\begin{example}[Bundles of Lie algebras] Let $A \Rightarrow M$ be a Lie algebroid with zero anchor map. Then the bracket is $C^\infty(M)$-bilinear, so it induces a Lie algebra structure on each fiber: such an algebroid is called a \emph{bundle of Lie algebras}. Notice that, if $A \Rightarrow M$ is a regular Lie algebroid, then the isotropy bundle $\mathfrak i (A)$ is a bundle of Lie algebras in the sense of this definition. \end{example}

\begin{example}[The tangent bundle] Apart from Lie algebras, the easiest example of Lie algebroid is given by the tangent bundle of a smooth manifold $M$. Here, the bracket is the Lie bracket of vector fields and the anchor map is the identity. It turns out that $TM$ is the Lie algebroid of both the pair groupoid $M \times M$ and the fundamental groupoid $\Pi_1(M)$. \end{example}

\begin{example}[Involutive distributions] A distribution $D$ on a manifold $M$ is a subbundle of $TM$: it is called \emph{involutive} if its sections are closed under the Lie bracket. It is clear that an involutive distribution is a Lie algebroid over $M$: the bracket and the anchor are obtained by restricting the structure maps of $TM$. Notice that every Lie algebroid with injective anchor map is actually of this form.

The well-known Frobenius theorem states that every involutive distribution integrates to a foliation. As a Lie algebroid, an involutive distribution integrates to both the monodromy groupoid and the holonomy groupoid of the corresponding foliation. \end{example}

\begin{example}[Action Lie algebroid] Let $M$ be a smooth manifold and $\mathfrak g$ be a Lie algebra. An \emph{(infinitesimal) action of $\mathfrak g$ on $M$} is a Lie algebra homomorphism $\rho: \mathfrak g \to \mathfrak X(M)$. This name comes from the fact that every Lie group action induces an infinitesimal action of its Lie algebra. Indeed, if $G$ acts on $M$, then $\mathfrak g$ acts on $M$ via:
\[
\rho(X)_x := \dfrac{d}{d \epsilon} \bigg|_{\epsilon = 0} \exp (\epsilon X) \cdot x.
\]
In general, if a Lie algebra action can be obtained in this way, we say that it is \emph{integrable}.

An action of $\mathfrak g$ on $M$ gives rise to a Lie algebroid $\mathfrak g \ltimes M \Rightarrow M$, called the \emph{action Lie algebroid}. As a vector bundle, $\mathfrak g \ltimes M$ is just the trivial vector bundle over $M$ with fiber $\mathfrak g$, so $\Gamma (\mathfrak g \ltimes M) \cong C^\infty (M) \otimes \mathfrak g$. Moreover, the anchor map is obtained from $\rho$ by scalar extension, while the bracket is uniquely determined by the Leibniz rule and by
\[
[c_X, c_Y] := c_{[X,Y]},
\]
where $c_X$ is the section that is constantly equal to $X \in \mathfrak g$. If the action comes from an action of the Lie group $G$ on $M$, then $\mathfrak g \ltimes M$ is the Lie algebroid of $G \ltimes M$. \end{example}

\begin{example}[Gauge algebroid] Let $G$ be a Lie group and $\pi: P \to M$ be a principal $G$-bundle. Then $G$ acts on $TP$ and the projection $TP \to P$ is $G$-equivariant. Taking the quotient, one obtains a vector bundle $TP/G \to M$ that is canonically a Lie algebroid: the bracket is obtained observing that sections of $TP/G$ are equivalent to $G$-invariant vector fields on $P$ and the anchor map is induced by $T\pi: TP \to TM$. This algebroid is just the Lie algebroid of the gauge groupoid of Example \ref{ex:gauge_grpd}. \end{example}

\begin{example}[Derivations] \label{ex:der} Let $p: E \to M$ be a vector bundle. We will assume that the rank of $E$ is greater than $0$. A \emph{derivation of $E$ at $x \in M$} is an $\mathbb R$-linear map $\delta: \Gamma (E) \to E_x$ that satisfies
\[
\delta (f \varepsilon) = \sigma_\delta (f) \varepsilon_x + f(x) \delta \varepsilon
\]
for some $\sigma_\delta \in T_x M$ and for every $f \in C^\infty(M)$, $\varepsilon \in \Gamma(E)$. The tangent vector $\sigma_\delta$ is uniquely determined by $\delta$: it is called the \emph{symbol} of $\delta$. It is easy to see that derivations of $E$ form a vector bundle over $M$, which we denote $DE$, and the symbol map is a bundle map $\sigma: DE \to TM$. This map gives rise to a short exact sequence
\begin{equation}\label{eq:spencer}
0 \longrightarrow \operatorname{End} E \longrightarrow DE \longrightarrow TM \longrightarrow 0
\end{equation}
of vector bundles over $M$, called the \emph{Spencer sequence} of $E$.

Sections of $DE$ are called \emph{derivations} of $E$. The $C^\infty(M)$-module of derivations of $E$ is denoted $\mathfrak D (E)$ (or $\mathfrak D (E, M)$ if we want to insist on the fact that the base of the vector bundle $E$ is $M$). A derivation $\Delta$ of $E$ can be described as an $\mathbb R$-linear map $\Delta: \Gamma(E) \to \Gamma(E)$ such that there exists a (necessarily unique) vector field $\sigma_\Delta \in \mathfrak X (M)$, the \emph{symbol} of $\Delta$, satisfying
\[
\Delta(f \varepsilon) = \sigma_\Delta (f) \varepsilon+ f \Delta \varepsilon
\]
for every $f \in C^\infty(M), \varepsilon \in \Gamma(E)$.

Finally, $DE$ carries a canonical Lie algebroid structure. The bracket is given by the commutator
\[
[\Delta_1, \Delta_2] := \Delta_1 \circ \Delta_2 - \Delta_2 \circ \Delta_1,
\]
while the anchor is the symbol map. This Lie algebroid is called the \emph{gauge algebroid} of $E$. Actually, one can show that $DE$ is nothing more than the Lie algebroid of the general linear groupoid $GL(E)$: for this reason, it is also denoted $\mathfrak{gl}(E)$. This implies that $DE$ is (canonically isomorphic to) the gauge algebroid of the principal bundle $\mathrm{Fr}(E) \to M$. 

There is an alternative description of $DE \to M$ which will be useful later. First of all, notice also that there is a canonical isomorphism of Lie algebroids 
\[
DE \overset{\cong}{\longrightarrow} DE^*, \quad \delta \mapsto \delta^*,
\]
defined as follows: if $\delta \in D_x E$, then $\delta^*: \Gamma(E^*) \to E^*_x$ is uniquely determined by the condition
\[
\sigma_\delta \langle \varphi, \varepsilon \rangle = \langle \delta^* \varphi, \varepsilon_x \rangle + \langle \varphi_x, \delta \varepsilon \rangle,
\]
for all $\varepsilon \in \Gamma(E)$, $\varphi \in \Gamma(E^*)$.

Now, consider the vector bundle $Tp: TE \to TM$. We denote by $TE|_v$ the fiber $Tp^{-1} (v)$ over a tangent vector $v \in TM$. For every $x \in M$, a derivation $\delta \in D_x E$ with symbol $\sigma_\delta$ determines a linear map $\ell_\delta: E_x \to TE|_{\sigma_{\delta}}$ via
\begin{equation} \label{eq:ell_delta}
\ell_\delta (e) (\ell_\varphi) = \langle \delta^* \varphi, e \rangle,
\end{equation}
for all $e \in E_x$ and $\varphi \in \Gamma (E^\ast)$. The assignment $\delta \mapsto (\sigma_\delta, \ell_\delta)$ establishes a one-to-one correspondence between the fiber $D_x E$ of the gauge algebroid over $x$ and the space of pairs $(v, h)$ where $v \in T_x M$ and $h : E_x \to TE|_v$ is a right inverse of the projection $TE|_v \to E_x$. We will sometimes identify $\delta$ and the pair $(\sigma_\delta, \ell_\delta)$.

At the level of sections, this simply means that derivations of $E$ can be identified with linear vector fields on $E$:
\[
\mathfrak D (E) \cong \mathfrak X_{\mathrm{lin}}(E).
\]
Actually, it is very easy to describe this isomorphism globally: if $\Delta \in \mathfrak D (E)$, then the corresponding linear vector field $X_\Delta$ is implicitly defined by $[X_\Delta, \varepsilon^\uparrow] = (\Delta \varepsilon)^\uparrow$ for all $\varepsilon \in \Gamma (E)$. Derivations and linear vector fields will be fundamental in this thesis: according to the situation, we will prefer one or the other approach.

The assignment $E \mapsto DE$ is functorial in the following sense. Let
\[
\begin{array}{c}
\xymatrix{E_N \ar[r]^{\phi} \ar[d] & E \ar[d] \\
N \ar[r]^{f} & M}
\end{array}
,
\]
be a \emph{regular} vector bundle morphism, i.e.~$\phi$ is an isomorphism on each fiber, so that $E_N$ is canonically isomorphic to the pull-back bundle $f^*E$. Then there is a pull-back map $\phi^*: \Gamma(E) \to \Gamma(E_N)$ defined by
\[
(\phi^*\varepsilon)_y = \phi_y^{-1}(\varepsilon_{f (y)})
\]
for all $\varepsilon \in \Gamma (E)$ and all $y \in N$. One can use this pull-back to define a Lie algebroid morphism $D \phi: DE \to DF$. Specifically, for all $\delta \in D_y E_N$ we define $D \phi (\delta): \Gamma(E) \to E_{f(y)}$ by
\[
D \phi (\delta) = \phi \circ \delta \circ \phi^*.
\]
It is then easy to see that 
\begin{equation} \label{eq: ell_D phi}
\ell_{D \phi (\delta)} = T \phi \circ \ell_\delta \circ \phi_x^{-1} : E_x \to TE|_{Tf (\sigma_{\delta})}.
\end{equation}
Finally, the diagram
\[
\xymatrix{
DE_N \ar[r]^-{D \phi} \ar[d] & DE \ar[d]^\sigma \\
TN \ar[r]^-{Tf} & TM
}
\]
is a pull-back diagram, and this induces an isomorphism 
\begin{equation} \label{eq:DE_N}
DE_N \cong TN \mathbin{{}_{Tf} \times_{\sigma}} DE
\end{equation}
 which is sometimes useful. From now on, we will often identify $E_N$ with the pull-back $f^\ast E$. For more details about the gauge algebroid we refer to \cite{mackenzie} (see also \cite{etv:infinitesimal}). \end{example}

All the Lie algebroids that we have met so far integrate to some Lie groupoid. Although it was believed since the 1960s that a third Lie theorem for Lie algebroids was true, in 1985 the first counterexample was given in \cite{almeida:suites}. After several partial positive results were established, in 2003 Crainic and Fernandes \cite{crainic:int} characterized precisely which Lie algebroids can be integrated to a Lie groupoid.

Our last class of examples is not related to any example in Subsection \ref{sec:lie_grpd} and includes some non-integrable Lie algebroids.

\begin{example}[Poisson manifolds] \label{ex:poisson} Let $M$ be a manifold. A \emph{Poisson structure} on $M$ is a Lie bracket $\{-,-\}: C^\infty(M) \times C^\infty(M) \to C^\infty(M)$ on its algebra of functions that satisfies
\[
\{ f,gh \} = \{ f,g \} h + g \{ f,h \}
\]
for every $f, g, h \in C^\infty(M)$. A manifold endowed with a Poisson structure is called a \emph{Poisson manifold}. A Poisson bracket $\{-,-\}$ induces a bivector $\pi \in \mathfrak X^2_{\mathrm{poly}}(M)$, given by $\pi(df, dg) = \{ f,g \}$. Actually, $\mathfrak X_{\mathrm{poly}} (M)[1]$ is a graded Lie algebra with the Schouten-Nijenhuis bracket $\llbracket -,- \rrbracket$, and a Poisson bracket on $M$ is equivalent to a bivector $\pi$ on $M$ satisfying $\llbracket \pi, \pi \rrbracket = 0$. Such a structure is called a \emph{Poisson bivector} on $M$. Moreover, the differential $\partial = \llbracket \pi, - \rrbracket$ makes $\mathfrak X_{\mathrm{poly}}(M)[1]$ a differential graded Lie algebra, called the \emph{Lichnerowicz complex} of $M$ and denoted $C_\pi (M)$.

\

Poisson manifolds carry a rich and interesting structure: a detailed account of their geometry is given in the recent monograph \cite{laurent:poisson}. Their definition is motivated by mechanics and they include, for example, symplectic manifolds. 

They also have a strong connection with Lie algebroids. First of all, \emph{a Lie algebroid structure on a vector bundle $p: A \to M$ is equivalent to a linear Poisson bivector on $A^*$}: here linear is meant in the sense of Definition \ref{def:linear}. Explicitly, if $A \Rightarrow M$ is a Lie algebroid, the Poisson bracket $\{ -,- \}: C^\infty (A^*) \times C^\infty (A^*) \to C^\infty (A^*)$ is characterized by
\[
\begin{aligned}
\{ \ell_\alpha, \ell_\beta \} & = \ell_{[\alpha, \beta]} \\
\{ \ell_\alpha, p^* f \} & = p^* (\rho(\alpha)(f))
\end{aligned}
\]
for every $\alpha, \beta \in \Gamma(A)$, $f \in C^\infty(M)$.

Moreover, a Poisson structure on a manifold $M$ induces a Lie algebroid structure on its cotangent bundle $T^* M$. To see this, take a Poisson bivector $\pi$ on $M$. It induces a bundle map $\pi^\#: T^* M \to TM$. Then, the bracket $[-,-]: \Omega^1(M) \times \Omega^1(M) \to \Omega^1(M)$ is given by
\[
[\varphi, \psi] := \mathcal L_{\pi^\# (\varphi)} \psi - \mathcal L_{\pi^\# (\psi)} \varphi - d(\pi(\varphi,\psi))
\]
and the anchor map is just $\pi^\#$. With this structure, $T^* M$ is called the \emph{cotangent Lie algebroid of $M$}. If it is integrable, the cotangent Lie algebroid of a Poisson manifold is a \emph{symplectic groupoid}, i.e.~a Lie groupoid carrying a compatible symplectic form. A survey of integration of Poisson manifolds can be found in \cite{crainic:int2}. \end{example}

\subsection{Representations (up to homotopy)} \label{sec:rep}

Let $G$ be a Lie group. We know that a representation of $G$ on a vector space $V$ is a Lie group morphism $G \to GL(V)$. It is natural then to give the following definition.

\begin{definition}\label{def:grpd_rep} Let $\mathcal G \rightrightarrows M$ be a Lie groupoid, $E \to M$ be a vector bundle. A \emph{representation of $\mathcal G$ on $E$} is a morphism of Lie groupoids $\Phi: \mathcal G \to GL(E)$ over the identity of $M$. \end{definition}

To such a representation, one can associate a cochain complex $(C(\mathcal G, E), \delta_\Phi)$, the \emph{(differentiable) complex of $\mathcal G$ with coefficients in $E$}, in the following way. For every $k > 0$, $C^k(\mathcal{G},E)$ is the set of smooth maps $u: \mathcal{G}^{(k)} \rightarrow E$ such that $u(g_1, \dots, g_k) \in E_{\mathsf t(g_1)}$. In other words, if we set $\mathsf f_k: \mathcal G^{(k)} \to M, (g_1, \dots, g_k) \mapsto \mathsf t(g_1)$, then
\begin{equation} \label{eq:f_k}
C^k(\mathcal G, E) := \Gamma(\mathsf f_k^* E).
\end{equation}
The differential of $u \in C^k(\mathcal{G},E)$ is given by
\begin{equation}\label{eq:grpd_rep_0}
\begin{aligned}
\delta_\Phi u (g_1, \dots, g_{k+1}) = & \ g_1 \cdot u(g_2, \dots, g_{k+1}) \\ 
& + \sum_{i=1}^k (-1)^i u(g_1, \dots, g_i g_{i+1}, \dots, g_{k+1}) + (-1)^{k+1} u(g_1, \dots, g_k).
\end{aligned}
\end{equation}
Moreover, $C^0(\mathcal G, E) := \Gamma(E)$. For $\varepsilon \in \Gamma(E)$, $\delta_\Phi \varepsilon$ is defined by
\begin{equation}\label{eq:grpd_rep}
\delta_\Phi \varepsilon (g) = g \cdot \varepsilon_{\mathsf s(g)} - \varepsilon_{\mathsf t(g)}.
\end{equation}
The cohomology of $C (\mathcal G, E)$ is called the \emph{(differentiable) cohomology of $\mathcal{G}$ with coefficients in $E$} and denoted $H (\mathcal G, E)$. This is a direct generalization of the classical (Lie) group cohomology. In particular, elements of $H^0 (\mathcal G, E)$ are called \emph{invariant sections} of $E$.

When $E$ is the trivial line bundle with the trivial representation, the above complex is simply called the \emph{Lie groupoid complex} of $\mathcal{G}$ and denoted $(C (\mathcal{G}), \delta)$. Its cohomology is called the \emph{Lie groupoid cohomology} of $\mathcal G$ and denoted $H (\mathcal G)$. 

\

The complex $C (\mathcal G)$ carries a natural product. The product of $f_1 \in C^k (\mathcal G)$, $f_2 \in C^l(\mathcal G)$ is defined by
\begin{equation} \label{eq:prod}
\begin{aligned}
& (f_1 \star f_2)(g_1, \dots, g_{k+l}) = f_1 (g_1, \dots, g_k) f_2 (g_{k+1}, \dots, g_{k+l}) && \operatorname{if} k, l > 0, \\
& (f_1 \star f_2)(g_1, \dots, g_l) = f_1 (\mathsf t(g_1)) f_2 (g_1, \dots, g_l) && \operatorname{if} k = 0, l > 0, \\
& (f_1 \star f_2)(g_1, \dots, g_k) = f_1 (g_1, \dots, g_k) f_2(\mathsf s(g_k)) && \operatorname{if} k > 0, l = 0, \\
& f_1 \star f_2 = f_1 f_2 && \operatorname{if} k = l = 0.
\end{aligned}
\end{equation}
This product gives $C (\mathcal G)$ a structure of DG-algebra, i.e.~the following Leibniz rule holds
\[
\delta (f_1 \star f_2) = \delta f_1 \star f_2 + (-1)^{|f_1|} f_1 \star \delta f_2
\]
for every $f_1, f_2 \in C (\mathcal G)$.

The same formulas define a right $C (\mathcal G)$-DG-module structure on $(C (\mathcal G, E), \delta_\Phi)$, for every representation $E$:
\begin{equation}\label{eq:leibniz_rep}
\delta_\Phi (u \star f) = \delta_\Phi u \star f + (-1)^{|u|} u \star \delta f
\end{equation}
for every $f \in C (\mathcal G)$, $u \in C (\mathcal G, E)$. Finally, notice that from (\ref{eq:f_k}) it follows that
\[
C(\mathcal G, E) \cong C (\mathcal G) \otimes_{C^\infty(M)} \Gamma(E).
\]

Formula (\ref{eq:grpd_rep}) suggests that it may be possible to recover the representation from the DG-module structure on $C (\mathcal G, E)$. This is indeed the case, but we need one further observation.

Let $u \in C^k (\mathcal G, E)$. We say that $u$ is \emph{normalized} if $u(g_1, \dots, g_k) = 0$ whenever one of the $g_i$'s is a unit. From Equations (\ref{eq:grpd_rep_0}) and (\ref{eq:grpd_rep}), it follows that the space of normalized cochains is in fact a subcomplex of $C (\mathcal G, E)$, denoted $\widehat C (\mathcal G, E)$. The following proposition shows that this property characterizes DG-module structures on $C(\mathcal G, E)$ arising from groupoid representations.

\begin{proposition} [{\cite[Lemma 2.6]{abad:ruth2}}] The assignment $\Phi \mapsto \delta_\Phi$ establishes a bijection between representations of $\mathcal G$ on $E$ and differentials on $C (\mathcal G, E)$ satisfying the Leibniz rule (\ref{eq:leibniz_rep}) and preserving the space of normalized cochains. \end{proposition}

Every Lie groupoid carries two canonical representations, the \emph{isotropy representation} and the \emph{normal representation}. Now we will describe them in some details. 

Let $\mathcal G \rightrightarrows M$ be a Lie groupoid, $A \Rightarrow M$ be its Lie algebroid. If $\mathcal G$ is regular, we know that the isotropy bundle $\mathfrak i$ and the normal bundle $\nu$ are well-defined vector bundles over $M$. If $\mathcal G$ is not regular, we can still make sense of $\mathfrak i$ and $\nu$ as ``set-theoretic vector bundles'', i.e.~families of vector spaces over $M$. In this case, the space of sections of $\mathfrak i$ is, by definition,
\[
\Gamma(\mathfrak i) := \ker (\rho: \Gamma(A) \to \mathfrak X (M)),
\]
while the space of sections of $\nu$ is
\begin{equation}\label{eq:sections_nu}
\Gamma(\nu) := \dfrac{\mathfrak X (M)}{\operatorname{im} (\rho: \Gamma(A) \to \mathfrak X (M))}.
\end{equation}

For every $g: x \to y$ in $\mathcal{G}$, there is an obvious conjugation map $\mathrm{Ad}_g: \mathcal{G}_x \to \mathcal{G}_y$. Differentiating at units one obtains the \emph{isotropy representation} of $\mathcal{G}$ on $\mathfrak{i}$, 
\[
\mathrm{ad}_g: \mathfrak{i}_x \to \mathfrak{i}_y.
\]
Even if $\mathfrak{i}$ is not a honest vector bundle, we define invariant sections by:
\[
H^0(\mathcal{G}, \mathfrak{i}) = \Gamma(\mathfrak{i})^{\mathrm{inv}} := \{ \alpha \in \Gamma(\mathfrak{i}) : \mathrm{ad}_g (\alpha_x) = \alpha_y \ \forall \ g: x \to y \ \mathrm{in} \ \mathcal{G} \}.
\]
More generally, one can define the differentiable cohomology of $\mathcal{G}$ with coefficients in $\mathfrak{i}$ as before, where a cochain is smooth if it is smooth as an $A$-valued map. One can prove that this gives a well-defined complex $C(\mathcal G, \mathfrak i)$.

The groupoid $\mathcal G$ also acts canonically on $\nu$. Take an arrow $g: x \to y$ in $\mathcal G$ and $v \in \nu_x$, then choose a curve $g(\epsilon): x(\epsilon) \to y(\epsilon)$ such that $g(0) = g$ and $\dot x (0)$ represents $v$. We define the \emph{normal representation} by
\[
\mathrm{ad}_g (v) = \dot y (0) (\operatorname{mod} \operatorname{im} \rho).
\]
One can check that this definition does not depend on the choices involved and that the following lemma holds.

\begin{lemma} [{\cite[Lemma 4.5]{crainic:def2}}] Let $V \in \mathfrak X (M)$. Then the set-theoretic section of $\nu$
\[
x \mapsto V_x \ (\operatorname{mod} \operatorname{im} \rho)
\]
is invariant if and only if, for any $\mathsf s$-lift $X$ of $V$ and any $g: x \to y$, there exists $\eta (g) \in A_{\mathsf t(g)}$ such that
\begin{equation}\label{eq:invariance}
V_{\mathsf t(g)} = T\mathsf t(X_g) + \rho (\eta (g)).
\end{equation}
\end{lemma}
Moreover, if Equation (\ref{eq:invariance}) holds for some $\mathsf s$-lift, it holds for all of them.

Notice that every section of $\nu$, as defined in (\ref{eq:sections_nu}), induces a set-theoretic section. So it is natural to declare that a section $V (\operatorname{mod} \operatorname{im} \rho)$ of $\nu$ is invariant if, for some $\mathsf s$-lift $X$ of $V$, there exists a smooth section $\eta$ of $\mathsf t^* A \to \mathcal G$ such that Equation (\ref{eq:invariance}) holds. But the vector field $X' \in \mathfrak X (\mathcal G)$ defined by
\[
X'_g = X_g + T R_g (\eta_g)
\]
is also an $\mathsf s$-lift, so we end up with the following definition.

\begin{definition} A section $V (\operatorname{mod} \operatorname{im} \rho) \in \Gamma(\nu)$ is \emph{invariant} if there exists $X \in \mathfrak X (\mathcal G)$ that is $\mathsf s$-projectable and $\mathsf t$-projectable to $V$. We say that $X$ is an \emph{$(\mathsf s, \mathsf t)$-lift} of $V$. The space of invariant sections is denoted $H^0 (\mathcal G, \nu)$ or $\Gamma(\nu)^{\mathrm{inv}}$. \end{definition}

Our next goal will be that of enlarging the class of Lie groupoid representations. We will need the notion of \emph{representation up to homotopy}, that is essentially obtained by replacing a vector bundle with a graded vector bundle.

\begin{definition} A \emph{graded vector bundle} over a smooth manifold $M$ is a direct sum of vector bundles $E = \bigoplus_i E_i$ over $M$, where each $E_i$ is considered to be of degree $i$. We will always assume that graded vector bundles are bounded in degree. The \emph{shift by $k$} $E[k] \to M$ is the graded vector bundle defined by $E [k]_i = E_{i+k}$. The \emph{dual} of $E$ is defined by $E^* = \bigoplus_i (E^*)_i$, where $(E^*)_i := (E_{-i})^*$. If $E \to M$ and $F \to N$ are graded vector bundles, a \emph{morphism of graded vector bundles $\phi: E \to F$} is a collection of vector bundle maps $\phi_i: E_i \to F_i$ that cover the same smooth map $f: M \to N$. \end{definition}

Let $\mathcal G \rightrightarrows M$ be a Lie groupoid and $E \to M$ a graded vector bundle. We define $C (\mathcal G, E)$ by
\[
C (\mathcal G, E)^k := \bigoplus_{i+j=k} C^i(\mathcal G, E_j).
\]
With the same formulas as before, $C (\mathcal G,  E)$ is a right graded $C (\mathcal G)$-module.

\begin{definition} A \emph{representation up to homotopy} of $\mathcal G$ on $E$ is an $\mathbb R$-linear map of degree 1 $D: C (\mathcal G, E) \to C (\mathcal G, E)$ that satisfies $D^2 = 0$ and
\begin{equation} \label{eq:leibniz_ruth}
D(\varepsilon \star f) = D \varepsilon \star f + (-1)^{|\varepsilon|} \varepsilon \star \delta f
\end{equation}
for every $\varepsilon \in \Gamma(E)$, $f \in C (\mathcal G)$. Equivalently, a representation up to homotopy of $\mathcal G$ on $E$ is a structure of right $C(\mathcal G)$-DG-module on $C(\mathcal G, E)$.
\end{definition}

One clearly has $C(\mathcal G, E) \cong C(\mathcal G) \otimes_{C^\infty (M)} \Gamma (E)$. Notice that the map $D$ induces a differential on the graded $C^\infty(M)$-module $\Gamma(E) = \bigoplus_i \Gamma(E_i)$. Indeed, for every $k$, $D$ sends $\Gamma(E_k)$ to
\[
\Gamma(E_{k+1}) \oplus \Omega^1_A (E_k) \oplus \Omega^2_A (E_{k-1}) \oplus \dots.
\]
Composing with the projection onto $\Gamma(E_{k+1})$, one obtains a map $\partial: \Gamma(E_k) \to \Gamma(E_{k+1})$, and from Formula (\ref{eq:leibniz_ruth}) it follows that $\partial$ is $C^\infty(M)$-linear. It is also clear that $\partial^2 = 0$. For this reason, one sometimes assumes that a cochain complex of vector bundles $(E, \partial)$ is given, and says that $\mathcal G$ is \emph{represented on $(E, \partial)$}. If $E$ is concentrated in degrees $k, \dots, k+l-1$, one says that $E$ is an \emph{$l$-term representation up to homotopy} of $\mathcal G$.

\begin{remark} The concept of representation up to homotopy of a Lie groupoid first appeared in \cite{abad:ruth2}. There the authors show that a representation up to homotopy on the cochain complex $(E , \partial)$ induces, for every $g: x \to y$, a cochain map $\lambda_g: E_x \to E_y$. While, for plain representations, one has $\lambda_g \circ \lambda_h = \lambda_{gh}$ for every $(g,h) \in \mathcal G^{(2)}$, in general the two sides of the latter equality are only homotopic, whence the term ``up to homotopy''. \end{remark}

\

In \cite{abad:ruth2}, it is also shown that there is a good notion of adjoint representation of a Lie groupoid in this context, defined up to isomorphisms. Namely, if $\mathcal G \rightrightarrows M$ is a Lie groupoid and $A \Rightarrow M$ its Lie algebroid, $\mathcal G$ is represented on the graded vector bundle $\mathrm{Ad}(\mathcal G) = A[1] \oplus TM$ with the differential $\rho: A \to TM$. 

One can ask if representations up to homotopy can be described in a more geometric fashion. Up to now, this has been done only for 1-term and 2-term representations, via the notion of \emph{VB-groupoid}: we will recall them in Subsection \ref{sec:VB-gr}. The general case is being treated in the work in progress \cite{delhoyo:simp}.

\

We now pass to the infinitesimal picture.

\begin{definition} Let $A \Rightarrow M$ be a Lie algebroid, $E \to M$ be a vector bundle. An \emph{$A$-connection on $E$} is an $\mathbb R$-bilinear map \[
\nabla: \Gamma(A) \times \Gamma(E) \to \Gamma(E), \quad (\alpha, \varepsilon) \mapsto \nabla_\alpha \varepsilon
\]
that satisfies:
\begin{enumerate}
\item $\nabla_{f \alpha} \varepsilon = f \nabla_\alpha \varepsilon$;
\item $\nabla_\alpha (f \varepsilon) = \rho(\alpha)(f) \varepsilon + f \nabla_\alpha \varepsilon$
\end{enumerate}
for every $f \in C^\infty(M)$, $\alpha \in \Gamma(A)$, $\varepsilon \in \Gamma(E)$.

We say that $\nabla$ is \emph{flat} if $[\nabla_\alpha, \nabla_\beta] = \nabla_{[\alpha, \beta]}$ for every $\alpha, \beta \in \Gamma(A)$. A vector bundle $E \to M$ endowed with a flat $A$-connection is also called a \emph{representation} of $A$.
\end{definition}

It is easily shown that \emph{a representation $(E, \nabla)$ of $A$ is equivalent to a Lie algebroid morphism $A \to DE$}. This shows that on one hand Lie algebroid representations generalize Lie algebra representations, on the other hand they are the infinitesimal version of Lie groupoid representations. Likewise, they can be described algebraically. Set $\Omega^k_A(E) := \Gamma(\wedge^k A^* \otimes E)$ and $\Omega_A(E) = \bigoplus_k \Omega^k_A (E)$. The latter graded space is also denoted $C(A,E)$ and its elements are called \emph{$E$-valued $A$-forms} (or \emph{$A$-forms with values in $E$}). The datum of an $A$-connection $\nabla$ on $E$ induces a differential $d_\nabla$ on $\Omega_A(E)$, given by:
\begin{equation} \label{eq:d_nabla}
\begin{aligned}
d_\nabla \omega (\alpha_1, \dots, \alpha_{k+1}) := & \sum_{i=1}^{k+1} (-1)^{i+1} \nabla_{\alpha_i} (\omega(\alpha_1, \dots, \widehat{\alpha_i}, \dots, \alpha_{k+1})) \\
& + \sum_{i < j} (-1)^{i+j} \omega([\alpha_i, \alpha_j], \alpha_1, \dots, \widehat{\alpha_i}, \dots, \widehat{\alpha_j}, \dots, \alpha_{k+1}).
\end{aligned}
\end{equation}
for $\omega \in \Omega^k_A(E)$, $\alpha_1, \dots, \alpha_{k+1} \in \Gamma(A)$. We explictly notice that, for $k = 0$, one obtains
\[
d_\nabla \varepsilon (\alpha) = \nabla_\alpha \varepsilon
\]
for every $\alpha \in \Gamma(A)$, $\varepsilon \in \Gamma(E)$, hence $\nabla$ can be recovered from $d_\nabla$. Moreover, $\nabla$ is flat if and only if $d_\nabla^2 = 0$. In this case, $(\Omega_A (E), d_\nabla)$ is a cochain complex, called the \emph{Lie algebroid complex of $A$ with coefficients in the representation $E$}. Its cohomology is denoted $H(A,E)$ and called the \emph{Lie algebroid cohomology of $A$ with coefficients in $E$}. 

When $E$ is the trivial line bundle with the representation $\nabla_\alpha = \rho(\alpha)$, the above complex is simply called the \emph{Lie algebroid complex} of $A$ and denoted $(\Omega_A, d_A)$ (or $(C(A), d_A)$). Sometimes we will use $\Omega_{A, M}$ if we want to insist on $M$ being the base manifold. The complex $\Omega_A$, equipped with the wedge product, is a differential graded-commutative algebra (DGCA or CDGA, for short). Its cohomology is called the \emph{Lie algebroid cohomology} of $A$ and denoted $H(A)$. 

If $E$ is a representation of $A$, $\Omega_A(E)$ is clearly a DG-module over $\Omega_A$ via
\[
\omega \wedge (\rho \otimes \varepsilon) := (\omega \wedge \rho) \otimes \varepsilon.
\]
The above considerations show that:

\begin{proposition} [{\cite[Proposition 2.3]{abad:ruth}}] \label{prop:nabla} The assignment $\nabla \mapsto d_\nabla$ establishes a one-to-one correspondence between representations of $A$ on $E$ and DG-module structures of $\Omega_A(E)$ on $\Omega_A$. \end{proposition}

We can go one step further and give an algebraic characterization of the Lie algebroid structures on a vector bundle. Given a vector bundle $A \to M$, $\Omega_A$ is a graded-commutative algebra, and we have already observed that a Lie algebroid structure on $A$ makes $\Omega_A$ a DGCA. For $k = 0,1$, Formula (\ref{eq:d_nabla}) reads
\[
\begin{aligned}
d_A f (\alpha) & = \rho(\alpha)(f); \\
d_A \varphi (\alpha, \beta) & = \rho (\alpha) (\langle \varphi, \beta \rangle) - \rho (\beta) (\langle \varphi, \alpha \rangle) - \langle \varphi, [\alpha, \beta] \rangle.
\end{aligned}
\]
This means that the differential $d_A$ is enough to recover the whole Lie algebroid structure. Moreover, the Jacobi identity for the bracket is equivalent to $d_A^2 = 0$ and the Leibniz rule to the fact that $d_A \varphi$ is $C^\infty(M)$-linear. Hence:

\begin{proposition}\label{prop:omega_a} There is a one-to-one correspondence between Lie algebroid structures on the vector bundle $A \to M$ and DGCA structures on the graded-commutative algebra $\Omega_A$. \end{proposition}

Proposition \ref{prop:nabla} shows how to enlarge the class of representations of a Lie algebroid. Let $A \Rightarrow M$ be a Lie algebroid, $E = \bigoplus_i E_i \to M$ be a graded vector bundle.

\begin{definition} A \emph{representation up to homotopy of $A$ on $E$} is an $\Omega_A$-DG-module structure on $\Omega_A(E)$. \end{definition}

As in the Lie groupoid case, a representation up to homotopy induces a differential $\partial$ on $E$, so one also says that \emph{$A$ is represented on the cochain complex $(E, \partial)$}. If $E$ is concentrated in degrees $k, \dots, k+l-1$, one says that $E$ is an \emph{$l$-term representation up to homotopy} of $A$.

\begin{remark} The concept of representation up to homotopy of a Lie algebroid first appeared in \cite{abad:ruth}. There the authors explain the terminology ``up to homotopy'', in analogy to what happens for Lie groupoids. \end{remark}

In \cite{abad:ruth} it is also shown that every Lie algebroid $A \Rightarrow M$ is represented on the graded vector bundle $\mathrm{Ad}(A) = A[1] \oplus TM$ with the differential $\rho: A \to TM$. This representation is defined up to isomorphisms and is called the \emph{adjoint representation of $A$}: it generalizes the ordinary adjoint representation of a Lie algebra.

In the case of Lie algebroids, a geometric description for representations up to homotopy is already available. For 1-term and 2-term representations, ordinary differential geometry is sufficient: the corresponding geometric objects are \emph{VB-algebroids}. For higher representations, we will need to pass to \emph{graded geometry}, which is reviewed in the next subsection.

\subsection{Graded manifolds} \label{sec:graded}

Graded manifolds are a generalization of manifolds obtained by allowing coordinates of any degree that ``commute in the graded sense''. Coordinates on a manifold $M$ are simply locally defined smooth functions, so our idea can be formalized by replacing the sheaf of commutative algebras $C^\infty_M$ of smooth functions over $M$ with a sheaf of graded-commutative algebras on $M$, with smooth functions in degree zero. DG-manifolds are graded manifolds endowed with an additional structure, and they provide a more algebraic model for Lie algebroids.

Our standard reference for graded geometry is \cite{mehta:thesis}. In this section, vector spaces will always be real and finite-dimensional, algebras will always be associative and unital and tensor products and wedge products will be over $\mathbb R$, unless otherwise stated.

\

Let $V$ be a graded vector space. The \emph{tensor algebra on $V$}, denoted $T(V)$, is the free graded algebra on $V$. It is obtained taking the iterated tensor products on $V$
\[
T(V) := \bigoplus_i V^{\otimes i},
\]
where the tensor product is understood in the category of graded vector spaces. 

The \emph{symmetric algebra on $V$}, denoted $S(V)$, is the free graded-commutative algebra on $V$. It is the quotient of $T(V)$ by the two-sided ideal generated by elements of the form
\[
v \otimes w - (-1)^{|v| |w|} w \otimes v
\]
where $v, w \in V$ are homogeneous elements. We will denote the product of $v$ and $w$ in $S(V)$ simply by $vw$. 

Finally, the \emph{exterior} (or \emph{Grassmann) algebra} on $V$, denoted $\wedge V$, is the quotient of $T(V)$ by the two-sided ideal generated by elements of the form
\[
v \otimes w + (-1)^{|v| |w|} w \otimes v
\]
where $v, w \in V$ are homogeneous elements. We will denote the product of $v$ and $w$ in $\wedge V$ by $v \wedge w$.

Notice that $T(V)$, $S(V)$ and $\wedge V$ are all naturally bigraded vector spaces. However, we will always consider them as graded vector spaces with respect to the total degree. All the discussion also holds for graded modules over a graded ring.

\

Let $(p_i)$ be a finite sequence of non-negative integers. We denote by $\mathbb R^{(p_i)}$ the pair $(\mathbb R^{p_0}, \mathcal O^{(p_i)})$, where $\mathcal O^{(p_i)}$ is the sheaf of graded-commutative algebras $C^\infty_{\mathbb R^{p_0}} \otimes S(V^*)$ and $V$ is the graded vector space whose components are $\mathbb R^{p_i}$ in degree $i$ for every $i \neq 0$, and $0$ in other degrees.

\begin{definition} A \emph{graded manifold of dimension $(p_i)$} is a pair $\mathcal M = (M, \mathcal O_M)$, where $M$ is a smooth manifold and $\mathcal O_M$ is a sheaf of graded-commutative algebras on $M$ with the following property: for every $x \in M$ there exist an open neighbourhood $U$ of $x$, an open subset $V \subset \mathbb R^{p_0}$ and an isomorphism of sheaves $(U, \mathcal O_M|_U) \to (V, \mathbb R^{(p_i)}|_V)$ covering a diffeomorphism $U \to V$.

The manifold $M$ is the \emph{support}, the sheaf $\mathcal O_M$ is the \emph{function sheaf} of $\mathcal M$ and its sections are called \emph{(local) functions on $\mathcal M$}. The graded-commutative algebra of global functions on $\mathcal M$ is denoted $C^\infty(\mathcal M)$. If the $p_i$'s range from to $k$ to $k+l$, we say that \emph{$\mathcal M$ is concentrated in degrees $k, \dots, k+l$}. A graded manifold concentrated in non-negative degrees is called an \emph{$\mathbb N$-manifold}. \end{definition}

From the definition, it follows that $\dim M = p_0$. Moreover, there exists a canonical \emph{evaluation map} $\mathrm{ev}_M: \mathcal O_M \to C^\infty_M$ that takes functions to their zero-degree part. A \emph{morphism of graded manifolds} $f: (M, \mathcal O_M) \to (N, \mathcal O_N)$ is given by a smooth map $f_0: M \to N$ and a morphism of sheaves $f^*: \mathcal O_N \to (f_0)_* \mathcal O_M$ such that $\mathrm{ev}_M \circ f^* = f_0^* \circ \mathrm{ev}_N$. We say that $f$ is a \emph{surjection} if $f_0$ is a surjection and $f^*$ is an injection of sheaves.

\

Take $U \subset M$ a coordinate open set that also satisfies $\mathcal O_M (U) \cong C^\infty (U) \otimes S(V^*)$. If $(x^i)$ is a set of coordinates on $U$, $(u^\beta)$ is a basis of $V^*$ made of homogeneous elements, we can think of the $x^i$'s and the $u^\beta$'s as ``coordinates'' on $\mathcal M$, of zero and non-zero degree respectively. By definition, a function $f \in C^\infty (\mathcal M)$ of degree $k$ locally looks like
\begin{equation} \label{eq:f_U}
f|_U = \sum_{\beta_1 \leq \dots \leq \beta_k} f_{\beta_1 \dots \beta_k} (x) u^{\beta_1} \dots u^{\beta_k}
\end{equation}
where $f_{\beta_1 \dots \beta_k}$ is a smooth function on $U$ for every choice of indices. Moreover, the non-zero degree coordinates satisfy
\[
u^{\beta_1} u^{\beta_2} = (-1)^{|u^{\beta_1}| |u^{\beta_2}|} u^{\beta_2} u^{\beta_1}
\]
for every $\beta_1$ and $\beta_2$. Thus we have given a precise mathematical meaning to the idea of manifold with graded-commuting coordinates. From Formula (\ref{eq:f_U}), it also follows that the subsheaf of degree $0$ functions on $\mathcal M$ is naturally isomorphic to $C^\infty_M$. This implies that there is a natural projection $p: \mathcal M \to M$, defined as follows: $p_0: M \to M$ is the identity, while $p^*: C^\infty_M \to \mathcal O_M$ is the natural inclusion we have just mentioned.

\

Of course, the most trivial example of graded manifold is a manifold $M$ endowed with its sheaf of smooth functions: here we have only coordinates of degree zero. 

Graded manifolds concentrated in degree $0$ and $1$ are also easily described. Let $E \to M$ be a vector bundle. Then we can construct a graded manifold $E[1] = (M, \Omega_E)$. The notation $E[1]$ means that we think $\Gamma(E^*)$ as being concentrated in degree $1$: informally, one says that $E[1]$ is obtained from $E$ by ``shifting by 1 the degree in the fibers''.

\begin{theorem}\label{prop:e[1]} Every graded manifold concentrated in degree $0$ and $1$ with support $M$ is canonically isomorphic to $E[1]$ for some vector bundle $E \to M$. The correspondence $E \rightsquigarrow E[1]$ is an equivalence of categories between vector bundles and graded manifolds concentrated in degree $0$ and $1$. \end{theorem}

This description can be further generalized to all $\mathbb N$-manifolds. Consider a graded vector bundle $E = \bigoplus_i E_i$ over $M$, and suppose that $E_0 = 0$. Then the pair 
\begin{equation}\label{eq:graded_mfld}
(M, \Gamma (S E^*))
\end{equation}
is an $\mathbb N$-manifold, and we have:

\begin{theorem} [Batchelor, \cite{batchelor:str}] Every $\mathbb N$-manifold is isomorphic to one of the form (\ref{eq:graded_mfld}). The category of $\mathbb N$-manifolds and the category of graded vector bundles are equivalent. \end{theorem}

Notice that, unlike in Theorem \ref{prop:e[1]}, the isomorphism is in general non-canonical.

\

On graded manifolds, one can construct a differential calculus that is entirely analogous to the classical calculus on smooth manifolds. This is due to the fact that, once the algebra of smooth functions on a manifold is given, all the usual objects and operations (vector bundles, vector fields, differential forms and so on) can be described in purely algebraic terms. These algebraic definitions extend, almost literally, to the realm of graded manifolds.

\begin{definition} Let $\mathcal M$ be a graded manifold, $(q_i)$ a finite sequence of non-negative integers. A \emph{vector bundle of rank $(q_i)$ over $\mathcal M$} is a graded manifold $\mathcal E$ together with a surjection $p: \mathcal E \to \mathcal M$ and an atlas of \emph{local trivializations}, i.e.~an open cover $(U_j)$ of $M$ and a family of isomorphisms $\phi_j: \mathcal E|_{p_0^{-1}(U_j)} \to \mathcal M|_{U_j} \times \mathbb R^{(q_i)}$ such that the transition functions $\phi_j \circ \phi_k^{-1}$ are linear in the fiber coordinates. \end{definition}

It follows that, if $(x^i, u^\beta)$ are local coordinates on $\mathcal M$, of zero and non-zero degree respectively, then there are corresponding local coordinates $(x^i, u^\beta, a^\alpha, c^\gamma)$ on $\mathcal E$, where $(a^\alpha, c^\gamma)$ are \emph{linear} fiber coordinates of zero and non-zero degree respectively.


\

As in the non-graded setting, it is possible to define the \emph{homogeneity structure} of a vector bundle of graded manifolds. Let $\mathcal E \to \mathcal M$ be such a vector bundle, and suppose that $(x^i, u^\beta)$ are local coordinates on $\mathcal M$ (of zero and non-zero degree respectively), $(a^\alpha, c^\gamma)$ are linear fiber coordinates (of zero and non-zero degree respectively). For every $\lambda \geq 0$, $h_\lambda: \mathcal E \to \mathcal E$ is locally given by:
\[
(x, u, a, c) \mapsto (x, u, \lambda a, \lambda c).
\]
We say that a function $f \in C^\infty (\mathcal E)$ is \emph{(homogeneous) of weight $q$} if and only if
\[
h_\lambda^* f = \lambda^q f
\]
for all $\lambda > 0$. Functions of weight $0$ and $1$ are called \emph{core} and \emph{linear functions}, respectively. The spaces of core and linear functions on $\mathcal E$ are denoted $C^\infty_ {\mathrm{core}}(\mathcal E)$ and $C^\infty_{\mathrm{lin}}(\mathcal E)$. More generally, the direct sum of all spaces of homogeneous functions on $\mathcal E$ form a subalgebra $C^\infty_{\mathrm{poly}}(\mathcal E) \subset C^\infty (\mathcal E)$, that of \emph{(fiber-wise) polynomial functions on $\mathcal E$}. This subalgebra is naturally bigraded by the degree (inherited by $C^\infty (\mathcal E)$) and the weight. To avoid any confusion, we stress that $C^\infty_{\mathrm{poly}}(\mathcal E)$ is a graded-commutative algebra only with respect to the degree and that, if $f \in C^\infty_{\mathrm{poly}}(\mathcal E)$, then $|f|$ denotes as usual the degree. We will not need a symbol for the weight of $f$.

\

Any natural operation on graded vector spaces has an associated operation on vector bundles of graded manifolds. In particular, one can define the \emph{shift by $k$} of a vector bundle $\pi: \mathcal E \to \mathcal M$. It is denoted $\pi[k]: \mathcal E [k] \to \mathcal M$ and it is obtained by ``raising the degree of the fiber coordinates by $k$''.

\begin{definition} A \emph{homogeneous section of degree $k$} of $\mathcal E \to \mathcal M$ is a morphism of graded manifolds $\varepsilon: \mathcal M \to \mathcal E [k]$ such that $\pi[k] \circ \varepsilon = \mathrm{id}_{\mathcal M}$. The space of degree $k$ sections of $\mathcal E$ is denoted $\Gamma_k (\mathcal E)$. The \emph{space of sections} of $\mathcal E$ is $\Gamma (\mathcal E) := \bigoplus_k \Gamma_k (\mathcal E)$. \end{definition}

One can prove that $\Gamma (\mathcal E)$ is a locally free module of finite rank. Moreover, all locally free modules of finite rank over $\mathcal M$ arise in this way:

\begin{theorem} [{\cite[Theorem 2.2.21]{mehta:thesis}}] The functor $\mathcal E \rightsquigarrow \Gamma (\mathcal E)$ establishes an equivalence of categories between vector bundles over $\mathcal M$ and locally free modules of finite rank over $C^\infty (\mathcal M)$. \end{theorem}

For this reason, we will mostly think about vector bundles in the category of graded manifolds as locally free (graded) modules of finite rank. This allows us to give some definitions. If $\mathcal E \to \mathcal M$ is a vector bundle of graded manifolds, the \emph{dual bundle} of $\mathcal E$ is the vector bundle $\mathcal E^* \to \mathcal M$ whose module of sections is:
\[
\Gamma(\mathcal E^*) := \operatorname{Hom}_{C^\infty (\mathcal M)} (\Gamma(\mathcal E), C^\infty (\mathcal M)).
\]
Similarly, the \emph{endomorphism bundle} of $\mathcal E$ is the vector bundle $\operatorname{End} \mathcal E \to \mathcal M$ whose module of sections is:
\begin{equation} \label{eq:end_bundle}
\operatorname{\mathfrak{End}} \mathcal E = \Gamma(\operatorname{End} \mathcal E) := \operatorname{End}_{C^\infty (\mathcal M)} \Gamma(\mathcal E).
\end{equation}
Let $f: \mathcal M \to \mathcal N$ be a morphism of graded manifolds and $\mathcal E \to \mathcal N$ a vector bundle. The \emph{pull-back of $\mathcal E$ along $f$} is the vector bundle $f^* \mathcal E \to \mathcal M$ defined by
\begin{equation} \label{eq:pull-back}
\Gamma(f^*\mathcal E) := C^\infty(\mathcal M) \otimes_{C^\infty (\mathcal N)} \Gamma(\mathcal E).
\end{equation}
This definition also holds when the vector bundle $\mathcal E \to \mathcal N$ is replaced by a graded vector bundle $E$ over a smooth manifold $N$. 
The next theorem describes all vector bundles over $\mathbb N$-manifolds.

\begin{theorem} [{\cite[Theorem 2.1]{mehta:vaintrob}}] \label{prop:vb_graded} Let $\mathcal M$ be an $\mathbb N$-manifold, $M$ be the underlying smooth manifold and $p: \mathcal M \to M$ the canonical projection. Then every vector bundle of graded manifold $\mathcal E \to \mathcal M$ is non-canonically isomorphic to one of the form $p^* E$, for some graded vector bundle $E \to M$. \end{theorem}

We emphasize that, in the above theorem, $\mathcal E$ is allowed to have fiber coordinates in negative degrees.

\

Now we are ready to give some basic notions of calculus on graded manifolds. Let $\mathcal M$ be a graded manifold. A \emph{vector field of degree $k$ on $\mathcal M$} is a graded derivation of degree $k$ of $C^\infty (\mathcal M)$, i.e.~an $\mathbb R$-linear map $X: C^\infty (\mathcal M) \to C^\infty (\mathcal M)$ of degree $k$ that satisfies
\[
X(fg) = X(f) g + (-1)^{k|f|} f X(g)
\]
for every $f, g \in C^\infty (\mathcal M)$. The space of vector fields of degree $k$ is denoted $\mathfrak X (\mathcal M)^k$. We set $\mathfrak X (\mathcal M) = \bigoplus_k \mathfrak X (\mathcal M)^k$.

One can show that $\mathfrak X (\mathcal M)$ is a locally free graded $C^\infty (\mathcal M)$-module, so it is the module of sections of a vector bundle over $\mathcal M$, the \emph{tangent bundle} $T \mathcal M \to \mathcal M$. Moreover, $\mathfrak X (\mathcal M)$ has a natural structure of graded Lie algebra with the graded commutator:
\[
[X,Y] := X \circ Y - (-1)^{|X| |Y|} Y \circ X.
\]


The shift by 1 of $T \mathcal M$ is denoted $T[1] \mathcal M$. By definition, the \emph{algebra of differential forms of $\mathcal M$} is $\Omega (\mathcal M) := C^\infty_{\mathrm{poly}} (T[1] \mathcal M)$.

\begin{remark} \label{rmk:T[1]M} If $\mathcal M$ is an $\mathbb N$-manifold, $C^\infty_{\mathrm{poly}}(T[1] \mathcal M)$ is simply $C^\infty(T[1] \mathcal M)$. In particular one can see in coordinates that, for any smooth manifold $M$, $C^\infty (T[1] M)$ is canonically isomorphic to the usual de Rham algebra of $M$. \end{remark}

The following theorem is proved in the same way as in the non-graded setting.

\begin{theorem} Let $\mathcal M$ be a graded manifold. There is a unique $\mathbb R$-linear map of degree $1$ $d: \Omega (\mathcal M) \to \Omega (\mathcal M)$ such that
\begin{enumerate}
\item $df (X) = X(f)$ for every $f \in C^\infty(\mathcal M)$, $X \in \mathfrak X (\mathcal M)$;
\item $d(\omega \wedge \rho) = d \omega \wedge \rho + (-1)^{|\omega|} \omega \wedge d \rho$ for every $\omega, \rho \in \Omega (\mathcal M)$;
\item $d^2 = 0$.
\end{enumerate}
\end{theorem}

The map $d$ is called the \emph{de Rham differential on $\mathcal M$}.

\

If $\mathcal E \to \mathcal M$ is a vector bundle of graded manifolds, then its homogeneity structure also induces an action on $\mathfrak X (\mathcal E)$. As usual, we say that a vector field $X \in \mathfrak X (\mathcal E)$ is \emph{(homogeneous) of degree $q$} if
\[
h_\lambda^* X = \lambda^q X
\]
for every $\lambda > 0$. Vector fields of weight $0$ and $1$ are called \emph{core} and \emph{linear vector fields}, respectively. The spaces of core and linear vector fields on $\mathcal E$ are denoted $\mathfrak X_ {\mathrm{core}}(\mathcal E)$ and $\mathfrak X_{\mathrm{lin}}(\mathcal E)$.

\begin{remark} As in Subsection \ref{Sec:homogeneity}, we stress that the terms ``linear function'' and ``linear vector field'' are already present in the literature, while the terms ``core function'' and ``core vector field'' are used here for the first time. \end{remark}

\begin{remark} Let $\mathcal M$ be a graded manifold and denote by $p: T \mathcal M \to \mathcal M$ the projection. For every $X \in \mathfrak X (\mathcal M)$, its \emph{tangent lift} is the linear vector field $X_{\mathrm{tan}} \in \mathfrak X (T \mathcal M)$ uniquely determined by:
\begin{equation}\label{eq:tan_lift}
X_{\mathrm{tan}}(\ell_{df}) = \ell_{d X(f)}, \quad X_{\mathrm{tan}}(p^* f) = p^* X(f),
\end{equation}
for all $f \in C^\infty (\mathcal M)$. The tangent lift defines an inclusion:
\[
\mathrm{tan}: \mathfrak X (\mathcal M) \hookrightarrow \mathfrak X_{\mathrm{lin}}(T \mathcal M).
\]
If $M$ is simply a smooth manifold, this construction is classical. In this case, if $X \in \mathfrak X (M)$ and $(\Phi^X_\epsilon)$ is its flow, then $X_{\mathrm{tan}}$ is characterized by the fact that its flow is $(T \Phi^X_\epsilon)$.
\end{remark}

A \emph{derivation of degree $k$} of $\mathcal E$ is an $\mathbb R$-linear map $\Delta: \Gamma(\mathcal E) \to \Gamma(\mathcal E)$ of degree $k$ such that there exists a (necessarily unique) vector field $\sigma_\Delta \in \mathfrak X (\mathcal M)$ satisfying
\[
\Delta (f \varepsilon) = \sigma_\Delta (f) \varepsilon + (-1)^{k |f|} f \Delta \varepsilon
\]
for every $f \in C^\infty(\mathcal M)$, $\varepsilon \in \Gamma(\mathcal E)$. The space of derivations of degree $k$ is denoted $\mathfrak D^k (\mathcal E)$ or $\mathfrak D^k(\mathcal E, \mathcal M)$. Derivations of $\mathcal E$ form a graded $C^\infty (\mathcal M)$-module $\mathfrak D (\mathcal E) = \bigoplus_k \mathfrak D^k (\mathcal E)$. We will also denote it $\mathfrak D (\mathcal E, \mathcal M)$, if we want to insist on the base manifold $\mathcal M$. The symbol map $\sigma: \mathfrak D (\mathcal E) \to \mathfrak X (\mathcal M)$ determines a short exact sequence of graded $C^\infty (\mathcal M)$-modules, which is the graded analogue of (\ref{eq:spencer}):
\[
0 \longrightarrow \operatorname{\mathfrak{End}} \mathcal E \longrightarrow \mathfrak D (\mathcal E) \longrightarrow \mathfrak X (\mathcal M) \longrightarrow 0.
\]

\

In analogy with the non-graded case, these objects are related by several isomorphisms. As in the non-graded case, there is an isomorphism
\[
\Gamma(\mathcal E^*) \overset{\cong}{\longrightarrow} C^\infty_{\mathrm{lin}}(\mathcal E), \quad \varphi \mapsto \ell_\varphi.
\]
Moreover, a section $\varepsilon$ of $\mathcal E$ determines a core vector field $\varepsilon^\uparrow \in \mathfrak X (\mathcal E)$, its \emph{vertical lift}, uniquely defined by $\varepsilon^\uparrow (\ell_\varphi) := (-1)^{|\varepsilon||\varphi|} \langle \varphi, \varepsilon \rangle$ for every $\varphi \in \Gamma(\mathcal E^*)$. So we have an isomorphism
\[
\Gamma(\mathcal E) \overset{\cong}{\longrightarrow} \mathfrak X_{\mathrm{core}}(\mathcal E), \quad \varepsilon \mapsto \varepsilon^\uparrow.
\]
Now we study the module $\mathfrak D (\mathcal E)$ of derivations of $\mathcal E$. There is a canonical isomorphism of graded Lie algebras and graded $C^\infty (\mathcal M)$-modules 
\[
\mathfrak D(\mathcal E) \overset{\cong}{\longrightarrow} \mathfrak X_{\mathrm{lin}}(\mathcal E) , \quad \Delta \mapsto X_\Delta,
\]
implicitly defined by 
\begin{equation} \label{eq:X_Delta}
[X_\Delta, \varepsilon^\uparrow] = (\Delta \varepsilon)^\uparrow
\end{equation}
for all $\varepsilon \in \Gamma (\mathcal E)$. It follows that the symbol map determines a canonical projection $\mathfrak X_{\mathrm{lin}} (\mathcal E) \to \mathfrak X (\mathcal M)$. Finally, there is a canonical isomorphism
\begin{equation} \label{eq:graded_der}
\mathfrak D (\mathcal E) \overset{\cong}{\longrightarrow} \mathfrak D (\mathcal E^*), \quad \Delta \mapsto \Delta^*,
\end{equation}
implicitly defined by
\begin{equation} \label{eq:graded_der_2}
\sigma_\Delta \langle \varphi, \varepsilon \rangle = \langle \Delta^* \varphi, \varepsilon \rangle + (-1)^{|\varphi| |\Delta|} \langle \varphi, \Delta \varepsilon \rangle
\end{equation}
for every $\varepsilon \in \Gamma(\mathcal E)$, $\varphi \in \Gamma(\mathcal E^*)$.

\

On a smooth manifold $M$, every vector field $X$ satisfies $[X,X] = 0$. This is no longer true for a graded manifold, because the Lie bracket is graded skew-symmetric: if $X$ is a vector field of odd degree, the equation $[X,X] = 0$ is equivalent to $X^2 = 0$. This leads to the following definition.

\begin{definition} Let $\mathcal M$ be a graded manifold. A \emph{homological vector field} (or a \emph{$Q$-vector field}) on $\mathcal M$ is a degree 1 vector field $X$ on $\mathcal M$ that satisfies $[X,X] = 0$. A graded manifold equipped with a homological vector field is called a \emph{differential graded manifold} (or \emph{DG-manifold}, or \emph{$Q$-manifold}). 

A \emph{morphism of DG-manifolds} $f: (\mathcal M, X) \to (\mathcal N, Y)$ is a morphism of graded manifolds that satisfies $X \circ f^* = f^* \circ Y$. \end{definition}

Now, putting together Proposition \ref{prop:omega_a} and Theorem \ref{prop:e[1]} one obtains:

\begin{theorem} Every DG-manifold over a smooth manifold $M$ concentrated in degrees $0$ and $1$ is canonically isomorphic to $(A[1], d_A)$ for some Lie algebroid $A \Rightarrow M$. \end{theorem}

Therefore, a morphism of Lie algebroids can be defined simply as a morphism of the corresponding DG-manifolds. This agrees with Definition \ref{def:morphism} and with the general one in \cite{mackenzie}, so we have:

\begin{theorem} [\cite{vaintrob:lie}] The assignment $(A \Rightarrow M) \rightsquigarrow (A[1], d_A)$ establishes an equivalence of categories between Lie algebroids and DG-manifolds concentrated in degrees $0$ and $1$. \end{theorem}

Finally, we will give a graded-geometric description of representations up to homotopy of Lie algebroids.

\begin{definition} A \emph{DG-vector bundle} is a vector bundle of graded manifolds $\mathcal E \to \mathcal M$ such that $\mathcal E$ and $\mathcal M$ are both DG-manifolds, with homological vector fields $Q_{\mathcal E}$ and $Q_{\mathcal M}$ and, additionally, $Q_{\mathcal E}$ is linear and projects onto $Q_{\mathcal M}$. \end{definition}

For more details about DG-vector bundles see, e.g.~\cite{Vit:vv-forms}. 

\begin{definition} Let $A \Rightarrow M$ be a Lie algebroid. A \emph{Va\u \i ntrob module over $A$}, or an \emph{$A$-module}, is a DG-vector bundle $\mathcal E \to A[1]$ over $A[1]$. In other words, a Lie algebroid module over $A$ is a locally free DG-module over $\Omega_A$. \end{definition}

From Formula (\ref{eq:pull-back}), Theorem \ref{prop:vb_graded} and the definition of representation up to homotopy it follows immediately that:

\begin{theorem} [{\cite[Lemma 4.4]{mehta:vaintrob}}] Every Va\u \i ntrob module over $A$ is non-canonically isomorphic to a representation up to homotopy of $A$. \end{theorem}

\section{Deformations} \label{sec:def}

Once we have introduced Lie groupoids and algebroids, we briefly recall their deformation theory. This topic was originally devoleped in \cite{crainic:def2} and \cite{crainic:def}, respectively. We will start from the infinitesimal picture and we will add some details about equivalence of deformations which are missing in the original treatment.

We explicitly notice here that our deformation complexes are shifted by $-1$ with respect to the ones appearing in \cite{crainic:def2, crainic:def}. We did this for the deformation complex of a Lie algebroid in order to obtain a bracket of degree zero, and we modified the deformation complex of a Lie groupoid accordingly.

\subsection{Deformations of Lie algebroids} \label{sec:def_algd}

We begin with a vector bundle $E \to M$. Let $k \geq 0$.

\begin{definition} A  \emph{multiderivation} of $E$ with $k$ entries (and $C^{\infty}(M)$-multilinear symbol), also called a $k$-\emph{derivation}, is a skew-symmetric, $\mathbb{R}$-$k$-linear map
\[
c: \Gamma(E) \times \dots \times \Gamma(E) \to \Gamma(E)
\]
such that there exists a bundle map $\sigma_c: \wedge^{k-1} E \to TM$, the \emph{symbol} of $c$, satisfying the following Leibniz rule:
\[
c(\varepsilon_1, \dots, \varepsilon_{k-1}, f \varepsilon_k) = \sigma_c (\varepsilon_1, \dots, \varepsilon_{k-1})(f) \varepsilon_k + f c(\varepsilon_1, \dots, \varepsilon_k),
\]
for all $\varepsilon_1, \dots, \varepsilon_k \in \Gamma(E)$, $f \in C^{\infty}(M)$.
\end{definition}

$1$-derivations are simply derivations as defined in Example \ref{ex:der}, $2$-derivations are called \emph{biderivations}. The space of $k$-derivations is denoted $\mathfrak D^k (E)$ (or $\mathfrak D^k (E,M)$). In particular, $\mathfrak D^1 (E) = \mathfrak D (E)$. We also put $\mathfrak D^0 (E) = \Gamma (E)$ and $\mathfrak D^\bullet (E) = \bigoplus_k \mathfrak D^k (E)$. Then $\mathfrak D^{\bullet}(E)[1]$, endowed with the \emph{Gerstenhaber bracket} $\llbracket -,- \rrbracket$, is a graded Lie algebra. We recall that, for $c_1 \in \mathfrak D^k(E)$, and $c_2 \in \mathfrak D^l (E)$, the \emph{Gerstenhaber product} of $c_1$ and $c_2$ is the $\mathbb R$-$(k+l-1)$-linear map $c_1 \circ c_2$ given by
\[
(c_1 \circ  c_2) (\varepsilon_1, \dots, \varepsilon_{k+l-1}) = \sum_{\tau \in S_{l,k-1}} (-1)^{\tau} c_1(c_2 (\varepsilon_{\tau(1)}, \dots, \varepsilon_{\tau(l)}), \varepsilon_{\tau(l+1)}, \dots, \varepsilon_{\tau(l+k-1)}),
\]
for all $\varepsilon_1, \dots, \varepsilon_{k+l-1} \in \Gamma(E)$, and the Gerstenhaber bracket is defined by
\[
\llbracket c_1, c_2 \rrbracket = (-1)^{(k-1)(l-1)} c_1 \circ c_2 - c_2 \circ c_1.
\]
The graded Lie algebra $\mathfrak D^\bullet (E)[1]$ first appeared in \cite{grab:Lie}.

\begin{remark} The symbol map determines a short exact sequence of graded $C^\infty(M)$-modules:
\[
0 \longrightarrow \Omega_E (E) \longrightarrow \mathfrak D (E) \overset{\sigma}{\longrightarrow} \Omega_E (TM) [-1] \longrightarrow 0.
\]
In \cite{crainic:def} it is shown that the choice of a connection $\nabla$ on $E$ determines a splitting of this sequence. Explicitly, a left splitting is given by $c \mapsto L^\nabla_c$, with
\begin{equation} \label{eq:L_c}
L^\nabla_c (\varepsilon_1, \dots, \varepsilon_{k}) = c(\varepsilon_1, \dots, \varepsilon_{k}) + (-1)^{k+1} \sum_i (-1)^i \nabla_{\sigma_c (\varepsilon_1, \dots, \widehat \varepsilon_i, \dots, \varepsilon_k)} (\varepsilon_i).
\end{equation}
for every $c \in \mathfrak D^k (E)$.
\end{remark}

The group of vector bundle automorphisms of $E$ acts naturally on multiderivations of $E$. If $\phi: E \to E$ is an automorphism covering the diffeomorphism $\phi_M : M \to M$, then $\phi$ acts on sections of $E$ (by pull-back) via the following formula:
\[
\phi^* \varepsilon := \phi^{-1} \circ \varepsilon \circ \phi_M, \quad \varepsilon \in \Gamma (E),
\]
and it acts on higher degree multiderivations via:
\[
(\phi^* c) (\varepsilon_1, \dots, \varepsilon_{k}) := \phi^* \left(c(\phi^{-1}{}^* \varepsilon_1, \dots, \phi^{-1}{}^* \varepsilon_{k})\right)
\]
for all $\varepsilon_1, \dots, \varepsilon_{k} \in \Gamma(E)$, $c \in \mathfrak D^k (E)$. Moreover, $\phi$ acts in the obvious way on sections of the dual bundle $E^*$. It is clear that
\begin{equation}\label{eq:phi^*}
\begin{aligned}
\phi^*(f\varepsilon) & = \phi_M^* f \cdot \phi^* \varepsilon, \\
\phi^* (c(\varepsilon_1, \dots, \varepsilon_k)) & = (\phi^* c) (\phi^* \varepsilon_1, \dots, \phi^* \varepsilon_k), \\
\phi_M^* \langle \varphi, \varepsilon \rangle & = \langle \phi^* \varphi, \phi^* \varepsilon \rangle,
\end{aligned}
\end{equation}
for all $\varepsilon, \varepsilon_1, \dots, \varepsilon_k \in \Gamma (E)$, $f \in C^\infty (M)$, and $\varphi \in \Gamma (E^\ast)$. Finally, $\phi$ acts on the exterior algebras of $E$ and $E^\ast$, and it also acts on vector bundle maps $\wedge^k E \to TM$ in the obvious way.

A direct computation shows that the action of vector bundle automorphisms on multiderivations does also respect the Gerstenhaber bracket, i.e.~
\begin{equation}\label{eq:pb_Gerst}
\phi^* \llbracket c_1, c_2 \rrbracket = \llbracket \phi^* c_1, \phi^* c_2 \rrbracket
\end{equation}
for all $c_1, c_2 \in \mathfrak D^\bullet(E)$. Additionally,
\begin{equation}\label{eq:phi^*sigma}
\phi^* \sigma_c = \sigma_{\phi^* c}.
\end{equation}
for all $c \in \mathfrak D^\bullet(E)$.

If $A \Rightarrow M$ is a Lie algebroid, the Lie bracket $b_A = [-,-]$ on sections of $A$ is a biderivation and it contains the full information about $A \Rightarrow M$. Additionally, $\llbracket b_A, b_A \rrbracket = 0$ as a consequence of the Jacobi identity. We summarize this remark with the following

\begin{proposition} [{\cite[Subsection 2.2, Lemma 2]{crainic:def}}] \label{prop:def_MC} Lie algebroid structures on $A \to M$ are in one-to-one correspondence with Maurer-Cartan elements in the graded Lie algebra $\mathfrak D^\bullet (A)[1]$, i.e.~degree $1$ elements $b$ such that $\llbracket b, b \rrbracket = 0$. \end{proposition}

Now, fix a Lie algebroid structure $A \Rightarrow M$ on the vector bundle $A \to M$, and let $b_A$ be the Lie bracket on sections of $A$. Equipped with the Gerstenhaber bracket and the interior derivation $\delta := \llbracket b_A, - \rrbracket$, $\mathfrak D^\bullet (A)[1]$ is a differential graded Lie algebra (DGLA), denoted $C_{\mathrm{def}} (A)$ (or $C_{\mathrm{def}} (A, M)$ if we want to insist on the base manifold being $M$) and called the \emph{deformation complex of $A$}. The cohomology of $C_{\mathrm{def}} (A)$ is denoted $H_{\mathrm{def}} (A)$ (or $H_{\mathrm{def}} (A, M)$), and called the \emph{deformation cohomology} of $A$.


The differential $\delta: C^{k}_{\mathrm{def}}(A) \to C^{k+1}_{\mathrm{def}}(A)$ is explicitly given by
\begin{equation}\label{eq:diff_algd}
\begin{aligned}
\delta c (\alpha_0, \dots, \alpha_{k+1}) = & \sum_{i} (-1)^{i} [\alpha_i, c(\alpha_0, \dots, \widehat{\alpha_i}, \dots, \alpha_{k+1})] \\
& + \sum_{i<j} (-1)^{i+j} c([\alpha_i, \alpha_j], \alpha_0, \dots, \widehat{\alpha_i}, \dots, \widehat{\alpha_j}, \dots, \alpha_{k+1}).
\end{aligned}
\end{equation}

\begin{definition} A  \emph{deformation} of $b_A$ is a(n other) Lie algebroid structure on the vector bundle $A \to M$. \end{definition}

It is clear that $b = b_A + c$ satisfies $\llbracket b, b \rrbracket = 0$ if and only if 
\[
\delta c + \dfrac{1}{2} \llbracket c, c \rrbracket = 0,	
\]
i.e.~{$c$ is a (degree $1$) solution of the  \emph{Maurer-Cartan equation} in the DGLA $C_{\mathrm{def}}(A)$}. Hence Proposition \ref{prop:def_MC} can be rephrased saying that \emph{the assignment $c \mapsto b_A + c$ establishes a one-to-one correspondence between Maurer-Cartan elements of $C_{\mathrm{def}}(A)$ and deformations of $A$}.

Now, let $b_0, b_1$ be deformations of $b_A$. We say that $b_0$ and $b_1$ are \emph{equivalent} if there exists a \emph{fiber-wise linear} isotopy taking $b_0$ to $b_1$, i.e.~there is a smooth path of vector bundle automorphisms $\phi_\epsilon: A \to A$, $\epsilon \in [0,1]$, such that $\phi_0 = \mathrm{id}_A$ and $\phi_1^\ast b_1 = b_0$. On the other hand, two Maurer-Cartan elements $c_0, c_1$ are \emph{gauge-equivalent} if they are interpolated by a smooth path of $1$-cochains $(c_\epsilon)$, and $(c_\epsilon)$ is a solution of the following ODE: 
\begin{equation}\label{eq:c_t}
\frac{dc_\epsilon}{d\epsilon} = \delta \Delta_\epsilon + \llbracket c_\epsilon, \Delta_\epsilon \rrbracket,
\end{equation}
for some smooth path of $0$-cochains (i.e.~derivations) $(\Delta_\epsilon)$, $\epsilon \in [0,1]$. 

Notice that (\ref{eq:c_t}) is equivalent to
\begin{equation}\label{eq:b_t}
\frac{db_\epsilon}{d\epsilon} = \llbracket b_\epsilon, \Delta_\epsilon \rrbracket
\end{equation}
where $b_\epsilon = b_A + c_\epsilon$.

The following observation is original: it appeared for the first time in \cite{lapastina:def}.

\begin{proposition}\label{prop:eq_G_eq}
Let $b_0 = b_A + c_0, b_1 = b_A + c_1$ be deformations of $b_A$. If $b_0, b_1$ are equivalent, then $c_0, c_1$ are gauge-equivalent. If $M$ is compact, the converse is also true.
\end{proposition}

\proof

Suppose that $b_0$ and $b_1$ are equivalent deformations, and let $\phi_\epsilon : A \to A$, $\epsilon \in [0,1]$ be a fiber-wise linear isotopy taking $b_0$ to $b_1$. Set $b_\epsilon = \phi_\epsilon^{-1}{}^* b_0 = b_A + c_\epsilon$, and define $(\Delta_\epsilon)$ by
\begin{equation}\label{eq:xi_t}
\dfrac{d \phi_\epsilon^*}{d\epsilon} = \phi_\epsilon^* \circ \Delta_\epsilon.
\end{equation}
One can check that $(\Delta_\epsilon)$ is a time-dependent derivation of $A$. Notice that 
\[
\llbracket b_\epsilon, b_\epsilon \rrbracket = \llbracket \phi_\epsilon^{-1}{}^* b_0, \phi_\epsilon^{-1}{}^* b_0 \rrbracket = \phi_\epsilon^{-1}{}^* \llbracket b_0, b_0 \rrbracket = 0,
\]
so $b_\epsilon$ is a deformation of $b_A$ for all $\epsilon$. Moreover, we have $\phi_\epsilon^*(b_\epsilon(\alpha, \beta)) = b_0(\phi_\epsilon^* \alpha, \phi_\epsilon^* \beta)$ for all $\alpha, \beta \in \Gamma(A)$. Differentiating with respect to $\epsilon$, we obtain:
\[
\begin{aligned}
\phi_\epsilon^* \left( \Delta_\epsilon (b_\epsilon(\alpha,\beta)) + \dfrac{db_\epsilon}{d\epsilon}(\alpha,\beta) \right) & = b_0\left(\phi_\epsilon^* (\Delta_\epsilon (\alpha)), \phi_\epsilon^*\beta \right) + b_0\left(\phi_\epsilon^* \alpha, \phi_\epsilon^*(\Delta_\epsilon(\beta))\right) \\ 
& =  \phi_\epsilon^* \left(b_\epsilon(\Delta_\epsilon(\alpha),\beta) + b_\epsilon(\alpha, \Delta_\epsilon(\beta))\right),
\end{aligned}
\]
so
\[
\dfrac{db_\epsilon}{d\epsilon}(\alpha,\beta) = b_\epsilon(\Delta_\epsilon(\alpha), \beta) + b_\epsilon(\alpha, \Delta_\epsilon(\beta)) - \Delta_\epsilon(b_\epsilon(\alpha,\beta)),
\]
i.e.~(\ref{eq:b_t}), hence (\ref{eq:c_t}), holds, as desired.

Conversely, suppose that $M$ is compact and there exist a family of derivations $(\Delta_\epsilon)$ and a family of $1$-cochains $(b_\epsilon)$ such that (\ref{eq:c_t}) or, equivalently, (\ref{eq:b_t}) holds. Let $X_\epsilon$ be the symbol of $\Delta_\epsilon$, for every $\epsilon$. From compactness, $(X_\epsilon)$ is a complete time-dependent vector field on $M$, i.e.~it generates a complete one-parameter family of diffeomorphisms $(\phi_M)_\epsilon$. The time-dependent derivation $(\Delta_\epsilon)$ generates a linear isotopy $\phi_\epsilon : A \to A$, covering the complete family $(\phi_M)_\epsilon$ (and implicitly defined by the ODE (\ref{eq:xi_t})).  By linearity, $\phi_\epsilon$ is a complete family itself. We want to show that
\begin{equation}\label{eq:Xi_t}
\phi_\epsilon^*(b_\epsilon(\alpha, \beta)) = b_0(\phi_\epsilon^* \alpha, \phi_\epsilon^* \beta) \quad \alpha,\beta \in \Gamma (A).
\end{equation}
For $\epsilon = 0$ this is obviously true and the derivatives of both sides are the same because of (\ref{eq:b_t}). So we have (\ref{eq:Xi_t}), and, by taking $\epsilon = 1$, we conclude that $(\phi_\epsilon)$ is a (fiber-wise linear) isotopy taking $b_0$ to $b_1$.
\endproof

Now we will briefly describe the cohomology groups in low degrees: their interpretation extends the classical one for Lie algebras.

Let $\alpha \in C^{-1}_{\mathrm{def}}(A) = \Gamma(A)$. Then $\delta \alpha = \mathrm{ad}_\alpha$, so 
\[
H^{-1}_{\mathrm{def}}(A) = Z(\Gamma(A)).
\]
Now take $c \in C^0_{\mathrm{def}}(A) = \mathfrak D (A)$. Then $\delta c = 0$ if and only if
\[
c([\alpha, \beta]) = [c(\alpha), \beta] + [\alpha, c(\beta)]
\]
for every $\alpha, \beta \in \Gamma(A)$. This means that $c$ is also a derivation of the Lie bracket: we say that $c$ is a \emph{Lie algebroid derivation}. Moreover, 0-coboundaries are derivations of the form $\mathrm{ad}_\alpha$, $\alpha \in \Gamma(A)$: we call them \emph{inner derivations of $A$}. The quotient is, by definition, the space of \emph{outer derivations} of $A$ and is denoted $\mathrm{OutDer}(A)$, hence
\[
H^0_{\mathrm{def}}(A) = \mathrm{OutDer}(A).
\]
The first cohomology group is the most important one for our purposes. By definition, an \emph{infinitesimal deformation} of a Lie algebroid $A \Rightarrow M$ is an element $c \in C^1_{\mathrm{def}}(A)$ such that $\delta c = 0$, i.e.~a 1-cocycle in $C_{\mathrm{def}}(A)$.  {As usual in deformation theory}, this definition is motivated by the fact that, \emph{if $(c_\epsilon)$ is a smooth path of Maurer-Cartan elements starting at $0$, then $\frac{dc_\epsilon}{d\epsilon} \big|_{\epsilon=0}$ is an infinitesimal deformation of $A$}. More generally, the cocycle condition $\delta c = 0$ is just the linearization at $c = 0$ of the Maurer-Cartan equation. Hence, $1$-cocycles in $C_{\mathrm{def}} (A)$ can be seen as the (formal) tangent vectors to the variety of Maurer-Cartan elements. Similarly, $1$-coboundaries can be seen as tangent vectors to the gauge orbit through $0$. We conclude that $H^1_{\mathrm{def}}(A)$ is the \emph{formal tangent space} to the moduli space of deformations under gauge equivalence. 

Finally, $H^2_{\operatorname{def,lin}}(W)$ is the space of obstructions to extend infinitesimal deformations to formal ones; if $H^2_{\operatorname{def,lin}}(W) = 0$, every infinitesimal deformation can be extended to a formal deformation. The proof of this fact is entirely algebraic and can be performed as in the Lie algebra case.

\begin{remark} 
The deformation complex of a Lie algebroid has an efficient description in terms of graded geometry.

Let $A \Rightarrow M$ be a Lie algebroid. We recall that its structure can be encoded in the DG-manifold $(A[1], d_A)$: smooth functions on $A[1]$ coincide with the de Rham complex $C(A)$. Moreover, vector fields on $A[1]$ form a graded Lie algebra $\mathfrak X (A[1])$. It becomes a DGLA with the adjoint operator $[d_A, -]$ and there is a canonical isomorphism of DGLAs
\begin{equation} \label{eq:c_def}
C_{\mathrm{def}}(A) \overset{\cong}{\longrightarrow} \mathfrak X (A[1]), \quad c \mapsto \delta_c,
\end{equation}
that can be explicitly described as follows. Let $c \in C^k_{\mathrm{def}} (A)$ and let $\sigma_c$ be the symbol of $c$. Then $\delta_c \in \mathfrak{X}(A[1])$  {is the degree $k$ vector field that takes} $\omega \in \Omega^p_A$, to $\delta_c \omega \in \Omega^{k + p}_A$ with
\begin{equation}\label{eq:def_der}
\begin{aligned}
\delta_c \omega (\alpha_1, \dots, \alpha_{k+p}) = & \sum_{\tau \in S_{k,p}} (-1)^{\tau} \sigma_c (\alpha_{\tau(1)}, \dots, \alpha_{\tau(k)}) (\omega (\alpha_{\tau(k+1)}, \dots, \alpha_{\tau(k+p)})) \\
& - \sum_{\tau \in S_{k+1,p-1}} (-1)^{\tau} \omega (c(\alpha_{\tau(1)}, \dots, \alpha_{\tau(k+1)}), \alpha_{\tau(k+2)}, \dots, \alpha_{\tau(k+p)})
\end{aligned}
\end{equation}
where $S_{l,m}$ denotes $(l,m)$-unshuffles. Notice that $c$ can be reconstructed from $\delta_c$ by using formula (\ref{eq:def_der}) for $p = 0, 1$:
\begin{equation}\label{eq:delta_c_f}
\delta_c f (\alpha_1, \dots, \alpha_k)  = \sigma_c(\alpha_1, \dots, \alpha_{k})f,
\end{equation}
and
\[
\delta_c \varphi (\alpha_1, \dots, \alpha_{k+1})  = \sum_{i} (-1)^{k-i} \sigma_c (\alpha_1, \dots, \widehat{\alpha_i}, \dots, \alpha_{k+1})\langle \varphi, \alpha_i \rangle - \langle \varphi, c(\alpha_1, \dots, \alpha_{k+1}) \rangle,
\]
where $f \in C^{\infty}(M)$, $\varphi \in \Omega^1_{A}$, and $\alpha_1, \dots, \alpha_{k+1} \in \Gamma(A)$. From the isomorphism (\ref{eq:c_def}) it follows that $C_{\mathrm{def}}(A)$ has an induced structure of DG-module over $C(A)$.
\end{remark}
 
\begin{remark}\label{rem:def_Pois}
There is yet another description of the deformation complex of $A \Rightarrow M$. From Example \ref{ex:poisson}, we know that $A^*$ is equipped with a fiber-wise linear Poisson structure $\pi$ and $C_\pi (A^*) = \mathfrak X_{\mathrm{poly}}(A^*)[1]$ is a differential graded Lie algebra. Moreover, $\pi$ is linear, so $\mathfrak X_{\mathrm{poly,lin}}(A^*) [1]$ is a subDGLA $C_{\pi, \mathrm{lin}} (A^*)$ of $C_\pi (A^\ast)$, and there is a canonical isomorphism
\begin{equation} \label{eq:iso_poisson}
C_{\mathrm{def}} (A) \overset{\cong}{\longrightarrow} C_{\pi, \mathrm{lin}} (A^\ast).
\end{equation}
The isomorphism $c \mapsto X_c$ is defined as follows. Recall that sections of $A$ are in one-to-one correspondence with fiber-wise linear functions on $A^\ast$. Now, the polyvector field $X_c$ is uniquely determined by
\begin{equation}\label{eq:lin_multiv}
\langle X_c, d\ell_{\alpha_1}\wedge \cdots \wedge d\ell_{\alpha_k}\rangle =  \ell_{c(\alpha_1, \dots, \alpha_k)}, \quad \alpha_1, \dots, \alpha_k \in \Gamma (A).
\end{equation}
Notice that the isomorphism (\ref{eq:iso_poisson}) generalizes the one-to-one correspondence between Lie algebroid structures on the vector bundle $A \to M$ and linear Poisson structures on $A^*$: they are precisely the Maurer-Cartan elements of the two DGLAs. \end{remark}

\subsection{Deformations of Lie groupoids} \label{sec:def_grpd}

Let $\mathcal G \rightrightarrows M$ be a Lie groupoid. As in Subsection \ref{sec:lie_grpd}, we denote by $\mathsf s, \mathsf t, \mathsf 1, \mathsf m, \mathsf i$ its structure maps, by $\bar{\mathsf m}$ the division map and by $\mathcal{G}^{(k)}$ the manifold of $k$-tuples of composable arrows of $\mathcal{G}$. 

\begin{definition}\label{def:def_complex}
The \emph{deformation complex} $(C_{\mathrm{def}}(\mathcal{G}), \delta)$ of $\mathcal{G}$ is defined as follows. For $k \geq 0$, $C^k_{\mathrm{def}}(\mathcal{G})$ is the set of smooth maps $c: \mathcal{G}^{(k+1)} \to T\mathcal{G}$ such that:
\begin{enumerate}
	\item $c(g_0, \dots, g_k) \in T_{g_0} \mathcal{G}$;
	\item $(T \mathsf s \circ c) (g_0, \dots, g_k)$ does not depend on $g_0$
\end{enumerate}
for any $(g_0, \dots, g_k) \in \mathcal{G}^{(k+1)}$. Thus we define the \emph{$\mathsf s$-projection of $c$} to be 
\[
s_c: \mathcal G^{(k)} \to TM, \quad s_c(g_1, \dots, g_k) := (T\mathsf s \circ c)(g_0, \dots, g_k).
\]

The differential of $c \in C^k_{\mathrm{def}} (\mathcal{G})$ is defined by
\begin{equation}\label{eq:diff_grpd}
\begin{aligned}
\delta c (g_0, \dots, g_{k+1}) = & -T \bar{\mathsf m} (c(g_0 g_1, \dots, g_{k+1}), c(g_1, \dots, g_{k+1})) \\ 
& + \sum_{i=1}^k (-1)^{i-1} c(g_0, \dots, g_i g_{i+1}, \dots, g_{k+1}) + (-1)^k c(g_0, \dots, g_k).
\end{aligned}
\end{equation}
Moreover, $C^{-1}_{\mathrm{def}}(\mathcal{G}) := \Gamma(A)$, where $A \Rightarrow M$ is the Lie algebroid of $\mathcal{G}$, and $\delta \alpha = \overleftarrow{\alpha} + \overrightarrow{\alpha}$ for each $\alpha \in \Gamma(A)$, where $\overleftarrow{\alpha}$ and $\overrightarrow{\alpha}$ are the left-invariant and right-invariant vector fields determined by $\alpha$.
\end{definition}

Notice that Formulas (\ref{eq:prod}) now give $C_{\mathrm{def}}(\mathcal G)$ a structure of DG-module over $C(\mathcal G)$.


\begin{remark}\label{rmk:pair} We observe, for later use, that conditions (1) and (2) can be expressed in the following way. For every $k \geq 1$, define the surjective submersions
\begin{equation}\label{eq:pq}
\begin{aligned}
p_k: \mathcal G^{(k)} \to \mathcal G, \quad & (g_1, \dots, g_k) \mapsto g_1, \\
q_k: \mathcal G^{(k)} \to M, \quad & (g_1, \dots, g_k) \mapsto \mathsf s(g_1).
\end{aligned}
\end{equation}
Then an element $c \in C^k_{\mathrm{def}}(\mathcal G)$ can be seen as a section of the pull-back bundle $p^*_{k+1} T \mathcal G \to \mathcal G^{(k+1)}$ that is also $\mathsf s$-projectable, i.e.~such that there exists a section $s_c$ of $q_k^* TM \to \mathcal G^{(k)}$ fitting into the following commutative diagram
\[
\begin{array}{c}
\xymatrix{p_{k+1}^* T \mathcal G \ar[d]_{T\mathsf s} \ar[r]& \mathcal G^{(k+1)} \ar@/_1.2pc/[l]_{c} \ar[d] &  \\
q_k^* TM \ar[r] & \mathcal G^{(k)} \ar@/_1.2pc/[l]_{s_c} }
\end{array}
,
\]
where the right vertical arrow is the projection onto the last $k$ elements.
\end{remark}

If $c \in C^k_{\mathrm{def}}(\mathcal G)$, $k \geq 0$, then the $\mathsf s$-projection of $\delta c$ is given by the following formula:
\begin{equation} \label{eq: s_delta}
\begin{aligned}
s_{\delta c} (g_1, \dots, g_{k+1}) = & -T \mathsf t (c(g_1, \dots, g_{k+1})) \\ 
& + \sum_{i=1}^k (-1)^{i-1} s_c (g_1, \dots, g_i g_{i+1}, \dots, g_{k+1}) + (-1)^k s_c (g_1, \dots, g_k).
\end{aligned}
\end{equation}

Recall that there is a short exact sequence of vector bundles over $\mathcal G$:
\begin{equation}\label{eq:ses1}
0 \longrightarrow T^{\mathsf s} \mathcal G \longrightarrow T \mathcal G \overset{T\mathsf s}\longrightarrow {\mathsf s}^* TM \longrightarrow 0.
\end{equation}
Pulling back the sequence via $\mathsf 1: M \to \mathcal G$, we get a short exact sequence of vector bundles over $M$
\begin{equation}\label{eq:ses2}
0 \longrightarrow A \longrightarrow \mathsf 1^* T \mathcal G \longrightarrow TM \to 0.
\end{equation}
The latter has a canonical splitting, given by $T\mathsf 1: TM \to T \mathcal G$, so $\mathsf 1^* T \mathcal G \cong A \oplus TM$ canonically.

\begin{definition}\label{def:norm_complex}
The \emph{normalized deformation complex} $\widehat C_{\mathrm{def}}(\mathcal G)$ of $C_{\mathrm{def}}(\mathcal G)$ is defined as follows. For $k \geq 1$, $\widehat C^k_{\mathrm{def}}(\mathcal G)$ is composed by all cochains $c \in C^k_{\mathrm{def}}(\mathcal G)$ such that
\begin{equation}\label{eq:norm}
c(\mathsf 1_x, g_1, \dots, g_k) \in T_x M \subset T_{\mathsf 1_x} \mathcal G, \quad \text{and} \quad c(g_0, \dots, \mathsf 1_x, \dots, g_k) = 0,
\end{equation}
while $\widehat C^0_{\mathrm{def}}(\mathcal G)$ is composed by $0$-cochains $c$ that satisfy $c(\mathsf 1_x) = s_c(x)$. Finally, $\widehat C^{-1}_{\mathrm{def}}(\mathcal G) := \Gamma(A)$. 
\end{definition}

Notice that the first condition in (\ref{eq:norm}) implies that $c(\mathsf 1_x, g_1, \dots, g_k)$ can be identified with $(T\mathsf s \circ c)(\mathsf 1_x, g_1, \dots, g_k) = s_c (g_1, \dots, g_k)$, so we recover the definition in \cite{crainic:def2}. Moreover, we have:

\begin{proposition}[{\cite[Proposition 11.8]{crainic:def2}}] \label{prop:norm} The inclusion $\widehat C_{\mathrm{def}}(\mathcal G) \hookrightarrow C_{\mathrm{def}}(\mathcal G)$ is a quasi-isomorphism. \end{proposition}

The group of automorphisms of $\mathcal G$ acts naturally on $C_{\mathrm{def}}(\mathcal G)$ by pull-back. Explicitly, if $\Psi \in \mathrm{Aut}(\mathcal G)$, denote by $\psi$ the induced automorphism of $A$. Then the action is given by
\[
\Psi^* c := \psi^* c
\]
for $c \in C^{-1}_{\mathrm{def}}(\mathcal G)$, and by
\begin{equation} \label{eq:aut_action}
(\Psi^* c)(g_0, \dots, g_k) := T \Psi^{-1} (c(\Psi(g_0), \dots, \Psi(g_k)))
\end{equation}
for $c \in C^k_{\mathrm{def}}(\mathcal G), k \geq 0$. It is easy to check that this action preserves the differential. Indeed, if $\alpha \in \Gamma(A)$, we have
\begin{equation}\label{eq:inv_vf}
\Psi^* \overrightarrow{\alpha} = \overrightarrow{\psi^* \alpha}, \quad \Psi^* \overleftarrow{\alpha} = \overleftarrow{\psi^* \alpha}
\end{equation}
and so
\[
\Psi^* (\delta \alpha) = \Psi^* (\overleftarrow{\alpha} + \overrightarrow{\alpha}) = \overleftarrow{\psi^* \alpha} + \overrightarrow{\psi^* \alpha} = \delta (\psi^* \alpha) = \delta (\Psi^* \alpha).
\]
If $c \in C^k_{\mathrm{def}}(\mathcal G)$, $k \geq 0$, a direct computation exploiting Equation (\ref{eq:diff_grpd}) shows that $\Psi^*$ commutes with $\delta$. Finally, it is straightforward that the action preserves $\widehat C_{\mathrm{def}}(\mathcal G)$.

\

Now we recall from \cite{crainic:def2} some information about the deformation cohomology of a Lie groupoid in degrees $-1$ and $0$ and its relation with the isotropy and the normal representations, first introduced in Subsection \ref{sec:rep}.

Let $\mathcal{G} \rightrightarrows M$ be a Lie groupoid. First of all, we have:

\begin{proposition} [{\cite[Proposition 4.1]{crainic:def2}}] \label{prop:H^-1} $H^{-1}_{\mathrm{def}}(\mathcal{G}) = H^0 (\mathcal{G}, \mathfrak{i}) = \Gamma(\mathfrak{i})^{\mathrm{inv}}$. \end{proposition}

Moreover, there is an inclusion of complexes
\begin{equation} \label{eq:isotropy}
r: C (\mathcal G, \mathfrak i) \hookrightarrow C_{\mathrm{def}} (\mathcal G)[-1] , \quad u \mapsto c_u
\end{equation}
given by
\begin{equation} \label{eq:c_u}
c_u(g_1, \dots, g_k) = T R_{g_1} (u(g_1, \dots, g_k)).
\end{equation}

Remember that there is an induced tangent prolongation Lie groupoid $T \mathcal G \rightrightarrows TM$. A vector field $X \in \mathfrak X (\mathcal G)$ is called \emph{multiplicative} if it is a groupoid morphism $X: \mathcal G \to T \mathcal G$. Notice that, if $\alpha \in \Gamma (A)$, $\delta \alpha = \overrightarrow{\alpha} + \overleftarrow{\alpha}$ is a multiplicative vector field. The flow of $\delta \alpha$ consists of inner automorphisms of $\mathcal G$ and every $\delta \alpha$ is called an \emph{inner multiplicative vector field}.

\begin{proposition} [{\cite[Proposition 4.3]{crainic:def2}}] \label{prop:H^0}
\[
H^0_{\mathrm{def}}(\mathcal G) = \dfrac{\text{multiplicative vector fields on $\mathcal G$}}{\text{inner multiplicative vector fields on $\mathcal G$}}.
\]
\end{proposition}

The relation with the normal representation is given by the linear map
\[
\pi: H^0_{\mathrm{def}}(\mathcal G) \to \Gamma(\nu)^{\mathrm{inv}}
\]
that sends the class of a multiplicative vector field to the class of its projection on $M$.

\begin{lemma} [{\cite[Lemma 4.9]{crainic:def2}}] \label{prop:curvature} Let $V \ (\operatorname{mod \ im} \rho) \in \Gamma(\nu)^{\mathrm{inv}}$ and $X$ an $(\mathsf s, \mathsf t)$-lift of $V$. Then $\delta X \in C^2(\mathcal G, \mathfrak i)$ and its cohomology class does not depend on the choice of $X$. Therefore there is an induced linear map
\[
K: \Gamma(\nu)^{\mathrm{inv}} \to H^2(\mathcal G, \mathfrak i),
\]
called the \emph{curvature map}.
\end{lemma}

Finally,

\begin{proposition} [{\cite[Proposition 4.11]{crainic:def2}}] \label{prop:exact} There is an exact sequence
\[
0 \longrightarrow H^1(\mathcal G, \mathfrak i) \overset{r}{\longrightarrow} H^0_{\mathrm{def}}(\mathcal G) \overset{\pi}{\longrightarrow} \Gamma(\nu)^{\mathrm{inv}} \overset{K}{\longrightarrow} H^2(\mathcal G, \mathfrak i) \overset{r}{\longrightarrow} H^1_{\mathrm{def}}(\mathcal G).
\]
\end{proposition}

The group $H^1_{\mathrm{def}}(\mathcal G)$ is directly linked to (infinitesimal) deformations. Before discussing this relation, we need a more general definition.

\begin{definition} A \emph{family of Lie groupoids}
\[
\tilde{\mathcal G} \rightrightarrows \tilde{M} \overset{\pi}{\longrightarrow} B,
\]
consists of a Lie groupoid $\tilde{\mathcal G} \rightrightarrows \tilde M$ and a surjective submersion $\pi: \tilde M \to B$ such that $\pi \circ \tilde {\mathsf s} = \pi \circ \tilde {\mathsf t}$. In particular, for every $b \in B$, $\mathcal G_b := (\pi \circ \tilde {\mathsf s})^{-1}(b)$ is a Lie groupoid over $M_b = \pi^{-1}(b)$. \end{definition}

This definition encodes the idea of a ``smoothly varying'' Lie groupoid. If $B$ is an open interval $I$ containing $0$, we say that $\tilde{\mathcal G}$ is a \emph{deformation} of $\mathcal G_0 \rightrightarrows M_0$ and we denote the latter by $\mathcal G \rightrightarrows M$. We will often denote by $\epsilon$ the canonical coordinate on $I$. Accordingly, a deformation of $\mathcal G$ is also denoted by $(\mathcal G_\epsilon)$. The structure maps of $\mathcal G_\epsilon$ are denoted $\mathsf s_\epsilon, \mathsf t_\epsilon, \mathsf 1_\epsilon, \mathsf m_\epsilon, \mathsf i_\epsilon$, the division map is denoted $\bar {\mathsf m}_\epsilon$. A deformation $(\mathcal G_\epsilon)$ is called \emph{strict} if $\tilde{\mathcal G} \cong \mathcal G \times I$ and $\tilde M \cong M \times I$. A strict deformation is \emph{$\mathsf s$-constant} (resp.~\emph{$\mathsf t$-constant}) if $\mathsf s_\epsilon$ (resp.~$\mathsf t_\epsilon$) does not depend on $\epsilon$. A deformation which is both $\mathsf s$-constant and $\mathsf t$-constant is \emph{$(\mathsf s, \mathsf t)$-constant}. The \emph{constant deformation} is the one with $\mathcal G_\epsilon = \mathcal G$ as groupoids for all $\epsilon$.

Two deformations $(\mathcal G_\epsilon)$ and $(\mathcal G'_\epsilon)$ of $\mathcal G$ are said to be \emph{equivalent} if there exists a smooth family of groupoid isomorphisms $\Psi_\epsilon: \mathcal G_\epsilon \to \mathcal G'_\epsilon$ such that $\Psi_0 = \mathrm{id}_{\mathcal G_0}$. We say that $(\mathcal G_\epsilon)$ is \emph{trivial} if it is equivalent to the constant deformation.

Let $(\mathcal G_\epsilon)$ be a strict deformation of the Lie groupoid $\mathcal G \rightrightarrows M$. Then it is natural to look at the variation of the structure maps, in particular the multiplication. However, in general, if $(g,h) \in \mathcal G^{(2)}$, there is no guarantee that $g$ and $h$ are composable also with respect to the groupoid structure $\mathcal G_\epsilon$. We will consider this problem later: now we simply assume that $(\mathcal G_\epsilon)$ is an $(\mathsf s, \mathsf t)$-constant deformation. In this case, it makes sense to consider the tangent vector
\begin{equation}\label{eq:vector}
- \dfrac{d}{d \epsilon} \bigg|_{\epsilon = 0} \mathsf m_\epsilon(g,h) \in T_{gh} \mathcal G
\end{equation}
for any $(g,h) \in \mathcal G^{(2)}$. It is clear that (\ref{eq:vector}) is killed by both $T \mathsf s$ and $T \mathsf t$. This means that it is of the form $T R_{gh}(a)$ with $a \in A$, and moreover $a \in \mathrm{ker}(\rho) = \mathfrak i$. Hence we can define a cochain $u_0 \in C^1 (\mathcal G, \mathfrak i)$ by
\[
u_0 (g,h) := - \dfrac{d}{d \epsilon} \bigg|_{\epsilon = 0} R^{-1}_{gh} \left(\mathsf m_\epsilon(g,h)\right).
\]
Differentiating the associativity equation $\mathsf m_\epsilon (\mathsf m_\epsilon(g,h),k) = \mathsf m_\epsilon(g, \mathsf m_\epsilon(h,k))$ at $\epsilon = 0$, we find that $u_0$ is a cocycle.

Now, define $\xi_0 := r(u_0) \in C^1_{\mathrm{def}}(\mathcal G)$, where $r$ is the map (\ref{eq:isotropy}). Differentiating at 0 the identity $\mathsf m_\epsilon(\bar {\mathsf m}_\epsilon (\mathsf m_0(g,h),h),h) = \mathsf m_0(g,h)$, one obtains the following expression for $\xi_0$:
\begin{equation}\label{eq:def_coc}
\xi_0(g,h) = \dfrac{d}{d\epsilon} \bigg|_{\epsilon = 0} \bar {\mathsf m}_\epsilon (gh,h).
\end{equation}

The last computation suggests how to generalize the procedure. Let $(\mathcal G_\epsilon)$ be an $\mathsf s$-constant deformation of $\mathcal G \rightrightarrows M$. The \emph{deformation cocycle} $\xi_0 \in C^1_{\mathrm{def}}(\mathcal G)$ associated to $(\mathcal G_\epsilon)$ is defined by Formula (\ref{eq:def_coc}) (which makes sense also in the present case).

\begin{lemma} [{\cite[Lemma 5.3]{crainic:def2}}] \label{prop:s-const} $\xi_0$ is a cocycle and its cohomology class only depends on the equivalence class of the deformation. \end{lemma}

Next we interpret $\xi_0$ in terms of the groupoid $\tilde{\mathcal G}$.

\begin{proposition} [{\cite[Proposition 5.7]{crainic:def2}}] Let $\tilde{\mathcal G}$ be an $\mathsf s$-constant deformation of the Lie groupoid $\mathcal G$. Then, if we set $\xi = \delta (\frac{\partial}{\partial \epsilon}) \in C^1_{\mathrm{def}}(\tilde {\mathcal G})$, we have $\xi_0 = \xi|_{\mathcal G}$. \end{proposition}

Notice that this statement relies heavily on the fact that $(\mathcal G_\epsilon)$ is $\mathsf s$-constant: otherwise the vector field $\frac{\partial}{\partial \epsilon}$ would not be $\tilde {\mathsf s}$-projectable. In the general case, one has to find an analogue of $\frac{\partial}{\partial \epsilon}$: this leads to the concept of \emph{transverse vector field}.

\begin{definition} Let $\tilde {\mathcal G}$ be a deformation of $\mathcal G$. A \emph{transverse vector field} for $\tilde{\mathcal G}$ is a vector field $X \in \mathfrak X (\tilde{\mathcal G})$ which is $\mathsf s$-projectable to a vector field $V \in \mathfrak X (\tilde M)$ which is, in turn, $\pi$-projectable to $\frac{d}{d \epsilon}$. \end{definition}

\begin{proposition} [{\cite[Proposition 5.12]{crainic:def2}}] \label{prop:gen_def} Let $\tilde{\mathcal G}$ be a deformation of $\mathcal G$. Then:
\begin{enumerate}
\item there exist transverse vector fields for $\tilde{\mathcal G}$;
\item if $\tilde X$ is transverse, then $\delta \tilde{X}$, when restricted to $\mathcal G$, induces a cocycle $\xi_0 \in C^1_{\mathrm{def}}(\mathcal G)$;
\item the cohomology class of $\xi_0$ does not depend on the choice of $\tilde X$.
\end{enumerate}
\end{proposition}

The resulting cohomology class in $H^1_{\mathrm{def}}(\mathcal G)$ is called the \emph{deformation class} associated to the deformation $\tilde {\mathcal G}$. So, in general it is not possible to find a canonical cocycle. It was possible in the case of an $\mathsf s$-constant deformation because there was a canonical choice of a transverse vector field. Notice that from the proposition above it follows directly that the deformation class is also invariant under equivalence of deformations.

Finally, we recall a result about general families of Lie groupoids. Let $ \tilde{\mathcal G} \rightrightarrows \tilde M \overset{\pi}{\longrightarrow} B$ be a family of Lie groupoids. Then any curve $\gamma: I \to B$ induces a deformation $\gamma^* \tilde{\mathcal G}$ of $\tilde{\mathcal G}_{\gamma(0)}$. We have the following

\begin{proposition} [{\cite[Proposition 5.15]{crainic:def2}}] Let $b \in B$. For any curve $\gamma: I \to B$ with $\gamma(0) = b$, the deformation class of $\gamma^* \tilde{\mathcal G}$ at time $0$ does only depend on $\dot{\gamma}(0)$. This defines a linear map
\[
\mathrm{Var}^{\tilde{\mathcal G}}_b: T_b B \to H^1_{\mathrm{def}}(\tilde{\mathcal G}_b),
\]
called the \emph{variation map} of $\tilde{\mathcal G}$ at $b$.
\end{proposition}

\section{VB-algebroids and VB-groupoids} \label{sec:VB}

In this section, we recall the basic definitions and properties of double vector bundles, VB-algebroids and VB-groupoids that will be useful later. 

\subsection{Double vector bundles} \label{sec:DVB}

The theory of double vector bundles was developed for the first time in 1974 in the PhD thesis of Pradines \cite{pradines:fibres}, although this concept was implicit in some of the existing literature on connections. A classical treatment of double vector bundles is contained in the book \cite{mackenzie}. The recent paper \cite{meinrenken:weil} gives an excellent account of the subject, involving some new ideas.

From now on, the projection, the addition, the scalar multiplication and the zero section of any vector bundle will be called its \emph{structure maps}.

\begin{definition} A \emph{double vector bundle} (DVB for short) is a vector bundle in the category of vector bundles.  More precisely, it is a commutative square
\begin{equation}\label{eq:DVB}
\begin{array}{c}
\xymatrix{
W \ar[d]_{\tilde q} \ar[r]^{\tilde p} & E \ar[d]^{q} \\
A \ar[r]^{p} & M}
\end{array},
\end{equation}
where all four sides are vector bundles, the projection $\tilde q: W \to A$, the addition $+_A: W \times_A W \to W$, the multiplication $\lambda\, \cdot_A : W \to W$ by any scalar $\lambda \in \mathbb R$ in the fibers of $W \to A$ and the zero section $\tilde{0}^A: A \to W$ are vector bundle maps covering the projection $q: E \to M$, the addition $+ : E \times_M E \to E$, the scalar multiplication $\lambda \cdot {}: E \to E$ and the zero section $0^E: M \to E$, respectively. The manifold $W$ will be called the \emph{total space}. DVB (\ref{eq:DVB}) will be also denoted by $(W \to E; A \to M)$. \end{definition}

\begin{remark} Notice that $W$ is a vector bundle over $E$ and over $A$, so it carries two homogeneity structures. However, we will mainly use the latter and denote it simply by $h$. As the homogeneity structure (likewise the fiber-wise sum) completely determines the vector bundle structure, most of the conditions in the definition of a DVB are actually redundant and could be omitted. For instance, a DVB is, equivalently, a diagram (\ref{eq:DVB}) of vector bundles such that the two homogeneity structures on $W$ commute. In particular, the definition is symmetric in the sense that, if $(W \to E; A \to M)$ is a DVB, then so is $(W \to A; E \to M)$. \end{remark} 

Let $(W \to E; A \to M)$ be a DVB and consider the submanifold
\[
C := \operatorname{ker}(W \to E) \cap \operatorname{ker}(W \to A) \subset W.
\]
In other words, elements of $C$ are those projecting simultaneously on the (images of the) zero sections of $A$ and $E$ (which are both diffeomorphic to $M$). The fiber-wise operations of the vector bundles $W \to E$ and $W \to A$ coincide on $C$ (see \cite{mackenzie}), so they define a (unique) vector bundle structure on $C$ over $M$. The vector bundle $C \to M$ is called the \emph{core} of $(W \to E; A \to M)$. The DVB $W$ is said to be \emph{trivial-core} if $C$ is the zero vector bundle $0_M$ over $M$.

In the following, we denote by $\Gamma (W, E)$ the space of sections of $W \to E$. Sections of $C \to M$ can be naturally embedded into $\Gamma(W,E)$, via the map 
\[
\Gamma(C) \to \Gamma(W,E), \quad \chi \mapsto \overline \chi,
\]
defined by: 
\[
\overline \chi_e = \tilde{0}{}^E_e +_A \chi_{q(e)}, \quad e \in E.
\]
The image of the inclusion $\chi \mapsto \overline \chi$ is, by definition, the space $\Gamma_{\mathrm{core}}(W,E)$ of \emph{core sections} of $W \to E$. 

There is another relevant class of sections of $W \to E$: \emph{linear sections}. We say that a section of $W \to E$ is a  \emph{linear section} if it is a vector bundle map covering some section of $A \to M$. The space of linear sections of $W \to E$ is denoted $\Gamma_{\mathrm{lin}}(W,E)$.  
We will usually denote by $\tilde \alpha, \tilde \beta, \dots$ the sections in $\Gamma_{\mathrm{lin}}(W,E)$. The $C^\infty (E)$-module $\Gamma (W, E)$ is spanned by $\Gamma_{\mathrm{core}}(W,E)$ and $\Gamma_{\mathrm{lin}}(W,E)$.

Linear and core sections of $W \to E$ can be efficiently characterized using the homogeneity structure $h$. Namely, the following lemma holds.

\begin{lemma}\label{prop:char_lin}
A section $w \in \Gamma(W,E)$ is
\begin{enumerate}
\item linear if and only if $h_\lambda^* w = w$ for every $\lambda > 0$;
\item core if and only if $h_\lambda^* w = \lambda^{-1} w$ for every $\lambda > 0$.
\end{enumerate}
\end{lemma}

More generally, we say that a section $w$ of $W \to E$ is \emph{of weight $q$} if $h_\lambda^* w = \lambda^q w$ for every $\lambda > 0$. Using this terminology, linear sections are precisely sections of weight $0$ and core sections are sections of weight $-1$. It is easy to check that \emph{there are no non-zero sections of $W \to E$ of weight less than $-1$.}

\

Every DVB is non-canonically isomorphic to a DVB of the form \cite{gracia-saz:vb}
\[
\begin{array}{c}
\xymatrix{
A \times_M E \times_M C \ar[d] \ar[r] & E \ar[d] \\
A \ar[r] & M}
\end{array}.
\]
Such a DVB is called \emph{split} or \emph{decomposed} and an isomorphism $W \cong A \times_M E \times_M C$ is called a \emph{decomposition} of $W$. It follows that, if $(A_{\alpha})$ is a local frame of $A$ and $(C_\gamma)$ is a local frame of $C$, then the $A_{\alpha}$'s can be lifted to linear sections $\tilde A_{\alpha} : E \to A \times_M E \times_M C$ (setting their $C$-component to be $0$), the $C_\gamma$'s identify with some core sections $\overline C_\gamma$, and $( \tilde A_{\alpha}, \overline C_\gamma )$ is a local frame of $W \to E$. If $(x^i)$ are coordinates on $M$, $(u^\beta)$ linear coordinates on $E$ and $(a^\alpha)$, $( c^\gamma )$ linear coordinates dual to the frames $(A_\alpha)$, $(C_\gamma)$, respectively, then $(x^i, u^\beta, a^\alpha, c^\gamma)$ are coordinates on $W$ called \emph{adapted} (to the DVB structure).

Using these coordinates, we can provide a local description of linear and core sections. A section $w \in \Gamma(W,E)$ is linear if and only if it is locally of the form
\begin{equation}\label{eq:lin_coord}
w = w^{\alpha}(x) \tilde A_{\alpha} + w_\beta^\gamma(x) u^\beta \overline C_\gamma,
\end{equation}
and it is core if and only if it is locally of the form
\[
w = w^\gamma(x) \overline C_\gamma
\]
for some local functions $w^\alpha, w_\beta^\gamma$ and $w^\gamma$.

From Equation (\ref{eq:lin_coord}), one can see that $\Gamma_{\mathrm{lin}}(W,E)$ is a locally free $C^\infty(M)$-module, hence it is the module of sections of a vector bundle $\widehat A \to M$. Moreover, there is a canonical projection $\widehat \pi: \widehat A \to A$ sending each linear section to its base section, and again from Equation (\ref{eq:lin_coord}) one can see that the kernel of $\widehat \pi$ is $\operatorname{Hom}(E,C)$. So there is a short exact sequence of vector bundles over $M$:
\begin{equation}\label{eq:widehat}
0 \longrightarrow \operatorname{Hom}(E,C) \longrightarrow \widehat A \overset{\widehat \pi}{\longrightarrow} A \longrightarrow 0.
\end{equation}
The injection $\operatorname{Hom}(E,C) \cong E^* \otimes C \hookrightarrow \widehat A$ is given by $\varphi \otimes \chi \mapsto \ell_\varphi \overline \chi$. Right splittings of the short exact sequence (\ref{eq:widehat}) are called \emph{horizontal lifts} of $W$: they are in one-to-one correspondence with decompositions of $W$ \cite{gracia-saz:vb}. Explicitly, if $\phi: W \to A \times_M E \times_M C$ is a decomposition, then the induced horizontal lift is given, at the level of sections, by
\[
\Gamma(A) \to \Gamma_{\mathrm{lin}}(W,E), \quad \alpha \mapsto \widehat \alpha
\]
with $\widehat \alpha_e = \phi^{-1}(\alpha_{q(e)}, e, 0_{q(e)})$.

In particular, if $W$ is a trivial-core DVB, then $\widehat A = A$ and the sequence (\ref{eq:widehat}) has just one section, the identity. In this case, there is a canonical isomorphism 
\begin{equation} \label{eq:trivial_core_iso}
(W\to E, A \to M) \overset{\cong}{\longrightarrow} (A \times_M E \to E; A \to M), \quad w \mapsto (\tilde q (w), \tilde p (w)).
\end{equation}

\begin{remark} \label{rmk:right_inv} From the discussion above, it is also clear how to describe the vector bundle $\widehat A$ pointwise. Indeed, the map $\tilde p$ in (\ref{eq:DVB}) can be seen as a map $\tilde p: W \to p^* E$ of vector bundles over $A$, so a linear section of $W \to E$ is a vector bundle map $\tilde \alpha: p^* E \to W$ such that $\tilde p \circ \tilde \alpha = \mathrm{id}_{p^* E}$. We conclude that, for every $x \in M$, an element $\widehat a$ of $\widehat A_x$ can be seen as a pair $(a, h)$, where $a \in A_x$ and $h$ is a right inverse of the linear map $\tilde p_a : W_a \to E_x$. \end{remark}

Double vector bundles have a very interesting duality theory, which is discussed in detail in \cite{mackenzie}. Here we recall some basic concepts. Let $(W \to E;A \to M)$ be a DVB, let $C$ be its core and let $W^*_A \to A$ be the dual vector bundle of $W \to A$. Then there is a vector bundle map
\begin{equation} \label{eq:dual_DVB}
\xymatrix{
W^*_A \ar[d]_{q_A} \ar[r]^{\pi_A} & C^* \ar[d] \\
A \ar[r]^p & M}
\end{equation}
given by 
\[
\langle \pi_A(\varphi), c \rangle = \langle \varphi, \tilde 0^A_a +_A c \rangle
\]
for every $a \in A$, $\varphi \in W^*_a$, $c \in C_{p(a)}$. Moreover, the diagram (\ref{eq:dual_DVB}) is a DVB, called the \emph{dual} of $W$ over $A$, and its core is $E^*$. The operations in the fibers of $\pi_A: W^*_A \to C^*$ are given by
\[
\begin{aligned}
\langle \varphi +_{C^*} \varphi', w +_A w' \rangle & := \langle \varphi , w \rangle + \langle \varphi', w' \rangle; \\
\langle \lambda \cdot_{C^*} \varphi, \lambda \cdot_A w \rangle & := \lambda \cdot \langle \varphi, w \rangle
\end{aligned}
\]
for suitable elements.

Of course, one can also construct the dual of $W$ over $E$:
\[
\begin{array}{l}
\xymatrix{
W^*_E \ar[d]_{q_E} \ar[r]^{\pi_E} & C^* \ar[d] \\
E \ar[r] & M}
\end{array} .
\]

There is a canonical non-degenerate pairing between the vector bundles $W^*_A \to C^*$ and $W^*_E \to C^*$, given by
\[
| \varphi, \psi | := \langle \varphi, w \rangle_A - \langle \psi, w \rangle_E,
\]
where $\varphi \in W^*_A$, $\psi \in W^*_E$, $\pi_A (\varphi) = \pi_E (\psi)$ and $w$ is any element of $W$ that projects on $q_A (\varphi)$ and $q_E (\psi)$ simultaneously. Hence there is an induced isomorphism 
\begin{equation} \label{eq:beta}
\beta: W^*_E \overset{\cong}{\longrightarrow} (W^*_A)^*_{C^*}
\end{equation}
of vector bundles over $C^*$. Finally, in \cite{mackenzie} it is proved that $\beta$ is also an isomorphism of DVBs.

\begin{examplex} \label{eq:tang_double} A distinguished example of a DVB is the \emph{tangent double of a vector bundle}. If $E \to M$ is a vector bundle, then
\[
\xymatrix{
TE \ar[d] \ar[r] & E \ar[d] \\
TM \ar[r] & M}
\]
is a DVB with core canonically isomorphic to $E$. The isomorphism is explicitly given by:
\[
E \overset{\cong}{\longrightarrow} \mathrm{Core}(TE), \quad e \mapsto \frac{d}{d \lambda} \bigg|_{\lambda = 0} \lambda e.
\]
Linear and core sections of $TE \to E$ are precisely linear and core vector fields (see Remark \ref{rmk:lin_vect}). Moreover, the isomorphism $\Gamma(E) \cong \mathfrak X_{\mathrm{core}}(E)$ coincides with the one defined in (\ref{eq:vert_lift}).

The dual of $TE$ over $E$ is
\[
\begin{array}{c}
\xymatrix{
T^*E \ar[d]_\pi \ar[r] & E \ar[d] \\
E^\ast \ar[r] & M}
\end{array}.
\]
If $E = TM$, we obtain the \emph{double tangent bundle} $TTM$ of $M$. Notice that the side bundles of $TTM$ both coincide with $TM$.
\end{examplex}

We conclude this paragraph by describing another isomorphism that will be useful in the following. As before, let $p: E \to M$ be a vector bundle, and denote $q: E^* \to M$ its dual. If $\mu: E^* \times_M E \to \mathbb R$ is the canonical pairing, the \emph{tangent pairing} $T \mu: TE^* \times_{TM} TE \to \mathbb R$ is a pairing between the vector bundles $TE \to TM$ and $TE^* \to TM$, which we denote $\llangle -,- \rrangle$. There is a canonical isomorphism
\begin{equation} \label{eq:B}
B: T^* E \overset{\cong}{\longrightarrow} T^* E^*,
\end{equation}
uniquely defined by
\[
\llangle v, \xi \rrangle = \langle \varphi, \xi \rangle_E + \langle B(\varphi), v \rangle_{E^*}
\]
for every $e \in E$, $\theta \in E^*$, $\xi \in T_e E$, $v \in T_\theta E^*$ and $\varphi \in T^*_e E$ such that $Tp (\xi) = Tq (v)$ and $\pi (\varphi) = \theta$.

\subsection{VB-algebroids} \label{sec:VB-alg}

\begin{definition}\label{def:vb-alg}
A \emph{VB-algebroid} is a DVB as in (\ref{eq:DVB}), equipped with a Lie algebroid structure $W \Rightarrow E$ such that the anchor $\rho_W: W \to TE$ is a vector bundle map covering a vector bundle map $\rho_A: A \to TM$ and the Lie bracket $[-,-]_W$ on sections of $W \to E$ satisfies
\begin{equation}\label{eq:VB-alg}
\begin{aligned}
{}[\Gamma_{\mathrm{lin}}(W,E), \Gamma_{\mathrm{lin}}(W,E)]_W & \subset \Gamma_{\mathrm{lin}}(W,E), \\
[\Gamma_{\mathrm{lin}}(W,E), \Gamma_{\mathrm{core}}(W,E)]_W & \subset \Gamma_{\mathrm{core}}(W,E), \\
[\Gamma_{\mathrm{core}}(W,E), \Gamma_{\mathrm{core}}(W,E)]_W & = 0.
\end{aligned}
\end{equation}
\end{definition}

Notice that, using the grading defined above, Property (\ref{eq:VB-alg}) is equivalent to asking that the Lie bracket on $\Gamma (W, E)$ is of weight 0.

\begin{remark}
An \emph{LA-vector bundle} is a ``vector bundle in the category of Lie algebroids''. Namely, it is a DVB
\[
\begin{array}{c}
\xymatrix{
W \ar[d]_{q_W} \ar@{=>}[r]& E \ar[d]^{q} \\
A \ar@{=>}[r] & M}
\end{array},
\]
where the horizontal sides are Lie algebroids and the structure maps of $W \to A$ are Lie algebroid maps (covering the structure maps of $E \to M$). This definition first appeared in \cite{mackenzie:double}. Then, Gracia-Saz and Mehta defined VB-algebroids as in Definition \ref{def:vb-alg} and proved that the two definitions do actually coincide \cite{gracia-saz:vb}. We use the second definition because it is more useful for our purposes: in particular, it emphazises the interplay between the bracket of the top algebroid and the grading induced on $\Gamma(W,E)$ by the homogeneity structure $h$.
\end{remark}

\begin{remark}
Notice that VB-algebroids can be characterized in terms of the sole homogeneity structure. Namely, a VB-algebroid is a DVB as in (\ref{eq:DVB}) such that the rows are Lie algebroids, and, for all $\lambda > 0$, $h_\lambda$ is a Lie algebroid map. This definition is particularly efficient in proving several properties of VB-algebroids (see, e.g., \cite{grab:hom1, grab:hom2, bursztyn:vec}).
\end{remark}

\

Now we will describe some features of the VB-algebroid $(W \Rightarrow E; A \Rightarrow M)$. 

The \emph{core-anchor} of $(W \Rightarrow E; A \Rightarrow M)$ is the vector bundle map $\partial : C \to E$ defined as follows. Let $\chi$ be a section of $C$, and let $\overline \chi$ be the corresponding core section of $W \to E$. The anchor $\rho_W: W \to TE$ maps $\overline \chi$ to a core vector field $\rho_W(\overline \chi)$ on $E$; in its turn, $\rho_W (\overline \chi)$ is the vertical lift of a section $\varepsilon$ of $E$ and we define $\partial \chi := \varepsilon$.

Then, consider the vector bundle $\widehat A \to M$. From Definition \ref{def:vb-alg}, it follows that $\widehat A$ is a Lie algebroid, with the bracket and the anchor given by:
\begin{alignat*}{2}
[\tilde \alpha, \tilde \beta ]_{\widehat A} & := [\tilde \alpha, \tilde \beta]_W; \\
\rho_{\widehat A} (\tilde \alpha) & := \rho_A (\tilde \alpha_0),
\end{alignat*}
where $\tilde \alpha_0$ is the projection of $\tilde \alpha$ on $\Gamma(A)$. We refer to $\widehat A$ as the \emph{fat algebroid} of $W$. The canonical map $\pi: \widehat A \to A$ is a Lie algebroid map, so there is an induced Lie algebroid structure on $\operatorname{Hom}(E,C)$ (with zero anchor) and the sequence (\ref{eq:widehat}) becomes a short exact sequence of Lie algebroids. Notice that sections of this sequence in the category of Lie algebroids do not always exist.


There are canonical representations of $\widehat A$ on $E$ and $C$, the \emph{side representation} $\psi^s$ and the \emph{core representation} $\psi^c$, given by:
\begin{equation} \label{eq:core_rep}
\begin{aligned}
\psi^s: & \ \Gamma (\widehat A) \times \Gamma(E) \to \Gamma(E), \quad (\psi^s_{\tilde \alpha} \varepsilon)^\uparrow = [\rho_{\widehat A}(\tilde \alpha), \varepsilon^\uparrow]; \\
\psi^c: & \ \Gamma (\widehat A) \times \Gamma (C) \to \Gamma(C), \quad \overline{\psi^c_{\tilde \alpha} \chi} = [\tilde \alpha, \overline \chi]_W.
\end{aligned}
\end{equation}
Finally we notice here, for the first time, that there is a 3-term complex of vector bundles canonically associated to $W$, which we denote $E_W$. Let $DE \times_{TM} DC$ be the fiber product over the symbol maps $\sigma_E: DE \to TM$, $\sigma_C: DC \to TM$: in other words, sections of $DE \times_{TM} DC$ are pairs $(\Delta_E, \Delta_C)$ of derivations of $E$ and $C$ with the same symbol. Then $E_W$ is the complex
\begin{equation} \label{eq:E_W}
0 \longrightarrow \widehat A [1] \longrightarrow DE \times_{TM} DC \longrightarrow \operatorname{Hom}(C,E) [-1] \longrightarrow 0,
\end{equation}
where the two (non-trivial) maps are given by:
\begin{equation} \label{eq:E_W_diff}
\tilde \alpha \mapsto (\psi^s_{\tilde \alpha}, \psi^c_{\tilde \alpha}), \quad (\Delta_E, \Delta_C) \mapsto \Delta_E \circ \partial - \partial \circ \Delta_C.
\end{equation}
This complex will be useful in Subsection \ref{sec:lin_def}, where we will describe the linear deformation complex of a VB-algebroid.

\subsubsection{Alternative descriptions}\label{subsec:alt_desc}

There are several equivalent descriptions of VB-algebroids. Two of them will be particularly useful for our purposes. 

\

The first one is in terms of graded geometry. We know that Lie algebroids are equivalent to DG-manifolds concentrated in degree $0$ and $1$. For VB-algebroids there is an analogous result that we now briefly explain. Let $(W \to E; A \to M)$ be a DVB. If we shift the degree in the fibers of both $W \to E$ and $A \to M$ (and use the functoriality of the shift) we get a vector bundle of graded manifolds, denoted $W[1]_E \to A[1]$. If $(W \Rightarrow E; A \Rightarrow M)$ is a VB-algebroid, then $W[1]_E \to A[1]$ is a DG-vector bundle concentrated in degree $0$ and $1$.
 
\begin{theorem} [\cite{voronov:q}] \label{prop:grad_VB}
Correspondence $(W \Rightarrow E; A \Rightarrow M) \rightsquigarrow (W[1]_E \to A[1])$ establishes an equivalence between the category of VB-algebroids and the category of DG-vector bundles concentrated in degree 0 and 1.
\end{theorem}
 
We have already observed that, in turn, DG-vector bundles over $A[1]$ are equivalent to representations up to homotopy, so the last theorem can be formulated in the following way:

\begin{theorem} \label{thm:ruth} The category of VB-algebroids is equivalent to that of 2-term representations up to homotopy of Lie algebroids. \end{theorem}

In \cite{gracia-saz:vb} and \cite{drummond:vb} this is proved directly, without the use of graded geometry. In particular, given a VB-algebroid $(W \Rightarrow E; A \Rightarrow M)$ with core $C$, the authors construct explicitly a 2-term representation up to homotopy of $A$ via the choice of a decomposition $W \cong A \times_M E \times_M C$.

If $W$ is a trivial-core VB-algebroid, then there is a unique decomposition of $W$, so $W \cong A \times_M E$ canonically as DVBs, via the isomorphism (\ref{eq:trivial_core_iso}), and $\widehat A = A$. In this case, the VB-algebroid structure on $W$ is completely determined by the Lie algebroid structure on $A$ and the side representation $\psi^s$, so Theorem \ref{thm:ruth} specializes as follows:

\begin{corollary} The category of trivial-core VB-algebroids is equivalent to that of representations of Lie algebroids. \end{corollary}

A second description of VB-algebroids is in terms of \emph{Poisson DVBs}. Let $(W \to E; A \to M)$ be a DVB with core $C$. Then $(W^*_E \to E; C^\ast \to M)$ is a DVB. In particular the smooth function algebra $C^{\infty}(W^*_E)$ possesses a bigraded subalgebra $C^\infty _{\mathrm{poly}}(W^\ast_E) = \bigoplus_{p,q} C^\infty_{p, q} (W^\ast_E)$, where $C^{\infty}_{p,q}(W^*_E)$ is the $C^\infty (M)$-module of functions that are homogeneous of weights $p$ and $q$ with respect to the vector bundle structures $W \to E$ and $W \to C^*$, respectively. Denote by $h^{W^\ast}$ the homogeneity structure of $W^*_E \to C^*$. It is easy to see, for example in coordinates, that the isomorphism
\[
\Gamma(W,E) \to C^\infty_{1,\bullet}(W^*_E), \quad w \mapsto \ell_w
\]
satisfies
\begin{equation}\label{eq:ell}
\ell_{h_\lambda^* w} = \lambda^{-1} \left(h_\lambda^{W^\ast}\right)^\ast \ell_w
\end{equation}
for every $\lambda > 0$. In particular, it identifies $\Gamma_{\mathrm{lin}}(W,E)$ with $C^{\infty}_{1,1}(W^*_E)$ and $\Gamma_{\mathrm{core}}(W,E)$ with $C^{\infty}_{1,0}(W^*_E)$.

Now assume, additionally, that $(W \Rightarrow E; A \Rightarrow M)$ is a Lie algebroid. The Lie algebroid structure $W \Rightarrow E$ induces a Poisson bracket $\{-,-\}$ on $W^*_E$ which is fiber-wise linear with respect to $W^\ast_E \to E$, i.e.
\[
\{ C^{\infty}_{1, \bullet}(W^*_E), C^{\infty}_{1, \bullet}(W^*_E) \} \subset C^{\infty}_{1, \bullet}(W^*_E).
\]
Condition (\ref{eq:VB-alg}) is now equivalent to say that $\{-,-\}$ satisfies also
\begin{equation}\label{eq:Poisson_linear}
\{ C^{\infty}_{\bullet, 1}(W^*_E), C^{\infty}_{\bullet, 1}(W^*_E) \} \subset C^{\infty}_{\bullet, 1}(W^*_E),
\end{equation}
giving the following

\begin{theorem} [{\cite[Theorem 2.4]{mackenzie:ehresmann}}]
Consider a DVB as in (\ref{eq:DVB}). A Lie algebroid structure $W \Rightarrow E$ is part of a VB-algebroid structure $(W \Rightarrow E;A \Rightarrow M)$ if and only if the induced linear Poisson structure on $W^*_E$ satisfies (\ref{eq:Poisson_linear}), i.e.~it is also linear with respect to the vector bundle structure $W^\ast_E \to C^\ast$.
\end{theorem}

In \cite{mackenzie:ehresmann}, Mackenzie calls \emph{Poisson DVB} a DVB whose total space is equipped with a Poisson structure that is fiber-wise linear with respect to both vector bundle structures. The above discussion then shows that a VB-algebroid $(W \Rightarrow E; A \Rightarrow M)$ is equivalent to a Poisson DVB $(W^\ast_E \to E; C^\ast \to M)$. Notice that, while the definition of VB-algebroid is not symmetric with respect to the ``horizontal'' and ``vertical'' sides, the definition of Poisson DVB is so, and this has the following important consequence. A VB-algebroid structure on $(W \Rightarrow E;A \Rightarrow M)$ induces a Poisson DVB $(W^\ast_E \to E; C^\ast \to M)$ which, in turn, determines (by duality in the ``vertical'' direction), a VB-algebroid $((W^*_E)^*_{C^*} \Rightarrow C^\ast; A \Rightarrow M)$. But we know that $(W^*_E)^*_{C^*} \cong W^*_A$, so $(W^*_A \Rightarrow C^\ast; A \Rightarrow M)$ is a VB-algebroid, called the \emph{dual VB-algebroid} of $(W \Rightarrow E; A \Rightarrow M)$.

\

We close this subsection presenting two basic examples. More examples will be given in Section \ref{Sec:examples}.

\begin{examplex} \label{ex:double_tang} Let $E \to M$ be a vector bundle and consider the DVB $(TE \to E; TM \to M)$ of Example \ref{eq:tang_double}. The tangent bundles $TE \to E$ and $TM \to M$ both carry Lie algebroid structures and one can verify that they fit together into a VB-algebroid $(TE \Rightarrow E; TM \Rightarrow M)$. From the observations in Example \ref{ex:der} and Remark \ref{rmk:right_inv} it follows that the fat algebroid of $TE$ coincides with $DE$ and we obtain again that $\mathfrak D (E) \cong \mathfrak X_{\mathrm{lin}}(E)$. \end{examplex}

\begin{examplex} \label{ex:tangent_VB} Let $A \Rightarrow M$ be a Lie algebroid. Then $(TA \to TM; A \to M)$ is a DVB, whose core is canonically isomorphic to $A$ itself. In particular, any section $\alpha$ of $A$ determines a core section $\overline \alpha$ of $TA \to TM$ and a linear section: its tangent map $T \alpha$. Now we recall from \cite{mackenzie} that the Lie algebroid structure on $A$ induces a Lie algebroid structure on $TA \to TM$.

The structure maps are defined as follows. The anchor $\rho_{TA} : TA \to TTM$ is determined by
\begin{equation}\label{eq:anchor_TA}
\rho_{TA}(T\alpha) = \rho(\alpha)_{\operatorname{tan}}, \quad \rho_{TA}(\overline \alpha) = \rho(\alpha)^{\uparrow},
\end{equation}
and the bracket $[-,-]_{TA}$ on $\Gamma (TA, TM)$ is determined by
\begin{equation}\label{eq:bracket_TA}
[T\alpha, T\beta]_{TA} = T[\alpha,\beta], \quad  [T\alpha, \overline \beta]_{TA} = \overline{[\alpha, \beta]}, \quad [\overline \alpha, \overline \beta]_{TA} = 0,
\end{equation}
for all $\alpha, \beta \in \Gamma(A)$. Moreover, $(TA \Rightarrow TM; A \Rightarrow M)$ is a VB-algebroid, the \emph{tangent VB-algebroid} of $A$. The dual VB-algebroid $(T^\ast A \Rightarrow A^\ast; A \Rightarrow M)$ of the tangent VB-algebroid is called the \emph{cotangent VB-algebroid}. \end{examplex}

\subsection{VB-groupoids} \label{sec:VB-gr}

\begin{definition}\label{def:VB-gr} A \emph{VB-groupoid} is a \emph{vector bundle in the category of Lie groupoids}, i.e.~a diagram
\begin{equation}\label{eq:VB-gr}
\begin{array}{r}
\xymatrix{
\mathcal W \ar[d] \ar@<0.4ex>[r] \ar@<-0.4ex>[r] & E \ar[d] \\
\mathcal{G} \ar@<0.4ex>[r] \ar@<-0.4ex>[r] & M}
\end{array},
\end{equation}
where $\mathcal W \rightrightarrows E$ and $\mathcal{G} \rightrightarrows M$ are Lie groupoids, $\mathcal W \to \mathcal{G}$ and $E \to M$ are vector bundles and all the vector bundle structure maps (addition, multiplication, projection and zero section) are Lie groupoid maps. We denote $\tilde{\mathsf s}, \tilde{\mathsf t}, \tilde{\mathsf 1}, \tilde{\mathsf m}, \tilde{\mathsf i}$ the structure maps of $\mathcal W$, $\mathsf s, \mathsf t, \mathsf 1, \mathsf m, \mathsf i$ the structure maps of $\mathcal{G}$, $\tilde p: \mathcal W \to \mathcal{G}$, $p: E \to M$ the vector bundle projections and $\tilde 0: \mathcal{G} \to \mathcal W$, $0: M \to E$ the zero sections. The VB-groupoid (\ref{eq:VB-gr}) will be also denoted $(\mathcal W \rightrightarrows E; \mathcal{G} \rightrightarrows M)$. The groupoid $\mathcal W \rightrightarrows E$ is called the \emph{total groupoid}, $\mathcal{G} \rightrightarrows M$ is called the \emph{base groupoid}. We will sometimes say that \emph{$\mathcal W$ is a VB-groupoid over $\mathcal G$}. \end{definition}

\begin{remark} Actually, some of the conditions in the previous definition are redundant and could be omitted. Moreover, from the definition it also follows that the structure maps of $\mathcal W \rightrightarrows E$ are vector bundle maps over the respective structure maps of $\mathcal G \rightrightarrows M$. For details and other equivalent definitions of VB-groupoids, see the discussion in \cite{gracia-saz:vb2}. \end{remark} 

\begin{remark}\label{rmk:hom_str} Similarly as for DVBs and VB-algebroids, the definition of VB-groupoid can be greatly simplified using the concept of homogeneity structure. Indeed, consider a diagram of Lie groupoids and vector bundles like (\ref{eq:VB-gr}). Let us denote by $h$ the homogeneity structure of $\mathcal W \to \mathcal G$. It can be shown that such a diagram is a VB-groupoid if and only if, for every $\lambda > 0$, $h_\lambda$ is a Lie groupoid automorphism \cite{bursztyn:vec}. \end{remark}

The Lie theory of VB-groupoids is studied in \cite{bursztyn:vec}. Here the authors show that, applying the usual differentation process to the total and the base groupoids of a VB-groupoid (and the structure maps of the vector bundles), we end up with a VB-algebroid.

From Definition \ref{def:VB-gr} it follows \cite{li-bland} that the map
\[
p^R: \mathcal W \to \mathsf s^* E, \quad v \mapsto (\tilde {\mathsf s} (v), \tilde{p} (v))
\]
is a surjective submersion. Hence its kernel is a vector bundle $V^R \to \mathcal{G}$, called the \emph{right-vertical subbundle} of $\mathcal W$. Finally, the \emph{right-core} of $\mathcal W$ is $C^R := \mathsf 1^* (V^R)$. Explicitly, $C^R$ is the set of elements of $\mathcal W$ that project on the units of $\mathcal{G}$ and $\mathsf s$-project to the zero section of $M$: this is analogous to the definition of the core of a DVB.

The right-core fits in a short exact sequence of vector bundles over $\mathcal{G}$:
\begin{equation}\label{eq:coreseq}
0 \longrightarrow \mathsf t^* C^R \overset{j^R}{\longrightarrow} \mathcal W \overset{p^R}{\longrightarrow} \mathsf s^* E \longrightarrow 0,
\end{equation} 
where $j^R$ is defined by $j^R(c,g) = c \cdot \tilde {\mathsf 0}_g$. A splitting of such a sequence always exists and gives a non-canonical decomposition $\mathcal W \cong \mathsf s^* E \oplus \mathsf t^* C^R$. Additionally, over the submanifold of units of $\mathcal G$ there is a natural splitting, given by $\tilde {\mathsf 1}: E \to \mathcal W$. Hence we give the following

\begin{definition} \cite{gracia-saz:vb2} A \emph{right-horizontal lift} of the VB-groupoid (\ref{eq:VB-gr}) is a splitting $h: \mathsf s^* E \to \mathcal W$ of (\ref{eq:coreseq}) that satisfies $h(e, \mathsf 1_x) = \tilde {\mathsf 1}_e$ for all $x \in M$, and $e \in E_x$. A \emph{right-decomposition} of (\ref{eq:VB-gr}) is a direct sum decomposition $\mathcal W \cong \mathsf s^\ast E \oplus \mathsf t^\ast C^R$ that comes from a right-horizontal lift. \end{definition}

The existence of right-horizontal lifts can be proved by a partition of unity argument, hence \emph{every VB-groupoid admits a (non-canonical) right-decomposition} \cite{gracia-saz:vb2}. 


By exchanging the role of the source and the target, one can similarly define a \emph{left-core} $C^L$. The analogous short exact sequence of (\ref{eq:coreseq}) is
\begin{equation}\label{eq:coreseq_2}
0 \longrightarrow \mathsf s^* C^L \longrightarrow \mathcal W \longrightarrow \mathsf t^* E \longrightarrow 0.
\end{equation}
The splittings of (\ref{eq:coreseq_2}) that restrict to the natural splitting over the units are called \emph{left-horizontal lifts}. Moreover, the involution
\[
F: \mathcal W \to \mathcal W, \quad v \mapsto -v^{-1}
\]
induces an isomorphism of vector bundles between the right-core and the left-core.

In the following, we will always consider the right-core: it will be referred to simply as the ``core'' and denoted $C$. The \emph{core-anchor} is the vector bundle map $\partial: C \to E$ defined by $\partial c = \tilde {\mathsf t} (c)$. Following \cite{etv:infinitesimal}, we say that $\mathcal W$ is \emph{trivial-core} (or \emph{vacant}) if its core is the zero vector bundle $0_M \to M$, \emph{full-core} if its side bundle is the zero vector bundle $0_M$.

\

Most of the concepts we outlined for VB-algebroids have an analogue in the theory of VB-groupoids. For example, it is possible to contruct the global counterparts of the side and the core representations: they are called the \emph{core} and the \emph{side actions} and involve the notion of \emph{fat groupoid} of a VB-groupoid. We will not need these definitions, so we refer the interested reader to \cite{gracia-saz:vb2}. Finally, we stress that an analogue of Theorem \ref{thm:ruth} also holds for VB-groupoids:

\begin{theorem} [{\cite[Corollary 6.2]{gracia-saz:vb2}}] There is a one-to-one correspondence between isomorphism classes of VB-groupoids and isomorphism classes of 2-term representations up to homotopy. \end{theorem}

Now we review two examples of VB-groupoids that will appear later on in the paper. Other examples will be considered in Subsection \ref{Sec:examples_1}.

\begin{examplex} Let $\mathcal{G} \rightrightarrows M$ be a Lie groupoid and $T\mathcal{G} \rightrightarrows TM$ be its tangent groupoid. Then it is easy to check that
\[
\begin{array}{r}
\xymatrix{
T \mathcal{G} \ar[d] \ar@<0.4ex>[r] \ar@<-0.4ex>[r] & TM \ar[d] \\
\mathcal{G} \ar@<0.4ex>[r] \ar@<-0.4ex>[r] & M}
\end{array}
\]
is a VB-groupoid. The core of $(T\mathcal G \rightrightarrows TM; \mathcal G \rightrightarrows M)$ is, by definition, the Lie algebroid $A$ of $\mathcal{G}$ and the core-anchor is the anchor map $\rho: A \to TM$.
\end{examplex}

\begin{examplex} \label{ex:act_grpd} Let $\mathcal G \rightrightarrows M$ be a Lie groupoid and let $p: E \to M$ be a representation of $\mathcal G$. Out of these data, we can define the action groupoid $\mathcal G \ltimes E \rightrightarrows E$. As a manifold, $\mathcal G \ltimes E = \mathcal G  \tensor*[_{\mathsf s}]{\times}{_p} E$, and the structure maps are given by
\[
\begin{aligned}
\tilde{\mathsf s} (g, e) & := e, \\ 
\tilde{\mathsf t} (g, e) & := g \cdot e, \\ 
(h, ge) \cdot (g, e) & := (hg, e).
\end{aligned}
\]
It is easy to prove that
\[
\begin{array}{r}
\xymatrix{
\mathcal{G} \ltimes E \ar[d] \ar@<0.4ex>[r] \ar@<-0.4ex>[r] & E \ar[d] \\
\mathcal{G} \ar@<0.4ex>[r] \ar@<-0.4ex>[r] & M}
\end{array}
\]
is a trivial-core VB-groupoid. Actually, every trivial-core VB-groupoid arises in this way, up to a canonical isomorphism \cite{gracia-saz:vb2}.
\end{examplex}

\begin{remark} Let $(\mathcal W \rightrightarrows E; \mathcal G \rightrightarrows M)$ be a VB-groupoid, $\tilde{\mathcal G} \rightrightarrows \tilde M$ be a Lie groupoid and $\Phi: \tilde{\mathcal G} \to \mathcal G$ be a morphism of Lie groupoids over $f: \tilde M \to M$. Then $\Phi^* \mathcal W = \mathcal W \times_{\mathcal G} \tilde{\mathcal G}$ and $f^* E = E \times_M \tilde M$, so $\Phi^* \mathcal W \rightrightarrows f^* E$ is a Lie groupoid with the structure maps defined component-wise. Moreover, in \cite{bursztyn:vec} it is proved that $\Phi^* \mathcal W \rightrightarrows f^* E$ is a VB-groupoid over $\tilde{\mathcal G} \rightrightarrows \tilde M$: $(\Phi^* \mathcal W \rightrightarrows f^* E; \tilde{\mathcal G} \rightrightarrows \tilde M)$ is called the \emph{pull-back VB-groupoid along $\Phi$}. \end{remark}

Let us briefly recall how duality works for VB-groupoids. Let $\mathcal W$ be a VB-groupoid like in (\ref{eq:VB-gr}), with core $C$, and let $\mathcal W^* \to \mathcal{G}$ be the dual vector bundle of $\mathcal W \to \mathcal{G}$. Then we can define the \emph{dual VB-groupoid}
\[
\begin{array}{r}
\xymatrix{
\mathcal W^* \ar[d] \ar@<0.4ex>[r] \ar@<-0.4ex>[r] & C^* \ar[d] \\
\mathcal{G} \ar@<0.4ex>[r] \ar@<-0.4ex>[r] & M}
\end{array}
\]
as follows. The source and the target $\check{\mathsf s}, \check{\mathsf t}: \mathcal W^* \rightrightarrows C^*$ are defined by
\[
\begin{aligned}
\langle \check {\mathsf s}(\varphi), c_1 \rangle & := - \langle \varphi, 0_g \cdot c_1^{-1} \rangle, \\
\langle \check {\mathsf t}(\varphi), c_2 \rangle & := \langle \varphi, c_2 \cdot 0_g \rangle
\end{aligned}
\]
for every $g \in \mathcal{G}$, $\varphi \in \mathcal W^*_g$, $c_1 \in C_{\mathsf s(g)}$, $c_2 \in C_{\mathsf t(g)}$, while the multiplication is uniquely defined by
\[
\langle \varphi_1 \cdot \varphi_2, v_1 \cdot v_2 \rangle := \langle \varphi_1, v_1 \rangle + \langle \varphi_2, v_2 \rangle
\] 
for all $(g_1,g_2) \in \mathcal G^{(2)}$, $(\varphi_1, \varphi_2) \in (\mathcal W^*)^{(2)}$, $\varphi_i \in \mathcal W^*_{g_i}$, $(v_1, v_2) \in \mathcal W^{(2)}$, $v_i \in \mathcal W_{g_i}$. By construction, the dual of a trivial-core VB-groupoid is a full-core VB-groupoid. One can also verify that the core of $\mathcal W^*$ is (canonically isomorphic to) $E^*$. For details and proofs see \cite{mackenzie}. 

\begin{examplex} Let $\mathcal G \rightrightarrows M$ be a Lie groupoid, and let $A \Rightarrow M$ be its Lie algebroid. The dual of the tangent VB-groupoid is the \emph{cotangent VB-groupoid}
\[
\begin{array}{r}
\xymatrix{
T^* \mathcal{G} \ar[d] \ar@<0.4ex>[r] \ar@<-0.4ex>[r] & A^* \ar[d] \\
\mathcal{G} \ar@<0.4ex>[r] \ar@<-0.4ex>[r] & M}
\end{array}.
\]
\end{examplex}

Finally, we describe the linear complex and the VB-complex of a VB-groupoid $(\mathcal W \rightrightarrows E, \mathcal{G} \rightrightarrows M)$ \cite{gracia-saz:vb2}. 

We know that the total groupoid comes with its Lie groupoid complex $(C (\mathcal W), \delta)$. It is easy to check that there is an induced vector bundle structure on $\mathcal W^{(k)} \to \mathcal{G}^{(k)}$, so there is a natural subcomplex $C_{\mathrm{lin}}(\mathcal W)$ of $C(\mathcal W)$, whose $k$-cochains are functions on $\mathcal W^{(k)}$ that are linear over $\mathcal{G}^{(k)}$. Inside $C_{\mathrm{lin}}(\mathcal W)$, there is a distinguished subcomplex $C_{\mathrm{proj}}(\mathcal W)$ consisting of \emph{left-projectable linear cochains} \cite{gracia-saz:vb2}. By definition, a linear cochain $f \in C^k_{\mathrm{lin}}(\mathcal W)$ is left-projectable if
\begin{enumerate}
\item $f(\tilde 0_g, v_2, \dots, v_k) = 0$ for every $(\tilde 0_g, v_2, \dots, v_k) \in \mathcal W^{(k)}$;
\item $f(\tilde 0_g \cdot v_1, \dots, v_k) = f(v_1, \dots, v_k)$ for $(v_1, \dots, v_k) \in \mathcal W^{(k)}$ and $g \in \mathcal G$ such that $\tilde {\mathsf t}(v_1) = 0_{\mathsf s(g)}$.
\end{enumerate}

The complexes $C_{\mathrm{lin}}(\mathcal W)$ and $C_{\mathrm{proj}}(\mathcal W)$ are called the \emph{linear complex} and the \emph{VB-complex} of $\mathcal W$, respectively. Their cohomologies are denoted $H_{\mathrm{lin}}(\mathcal W)$ and $H_{\mathrm{proj}}(\mathcal W)$ and called the \emph{linear cohomology} and the \emph{VB-cohomology} of $\mathcal W$.

\begin{remark} We are adopting the terminology of \cite{delhoyo:morita}. In \cite{gracia-saz:vb2} and \cite{crainic:def2}, instead, the VB-complex and the VB-cohomology of $\mathcal W$ are defined to be $C_{\mathrm{proj}}(\mathcal W^*)$ and $H_{\mathrm{proj}}(\mathcal W^*)$, respectively. \end{remark}

Notice that $C (\mathcal G)$ can be identified with the subcomplex of $C (\mathcal W)$ of fiber-wise constant cochains. Then, one can check that the product (\ref{eq:prod}) gives $C_{\mathrm{lin}}(\mathcal W)$ and $C_{\mathrm{proj}}(\mathcal W)$ a $C (\mathcal G)$-DG-module structure.

It turns out that the linear cohomology and the VB-cohomology are isomorphic:

\begin{lemma} [{\cite[Lemma 3.1]{cabrera:hom}}] \label{prop:proj;lin} The inclusion $C_{\mathrm{proj}}(\mathcal W) \hookrightarrow C_{\mathrm{lin}}(\mathcal W)$ induces an isomorphism of $H (\mathcal G)$-modules in cohomology. \end{lemma}

For our purposes, the VB-complex is particularly important because it gives another description of the deformation complex of a Lie groupoid $\mathcal G$. Indeed, we have:

\begin{proposition} [{\cite[Proposition 3.9]{crainic:def2}}] \label{prop:cotangent}
There is an isomorphism of $C (\mathcal G)$-DG-modules
\[
\phi: C_{\mathrm{def}}(\mathcal G) \overset{\cong}{\longrightarrow} C_{\mathrm{proj}} (T^* \mathcal G)[1]
\]
given by
\[
\phi(c)(\theta_0, \dots, \theta_k) = \langle \theta_0, c(g_0, \dots, g_k) \rangle
\]
for all $c \in C^k_{\mathrm{def}}(\mathcal G)$ and $(\theta_0, \dots, \theta_k) \in (T^* \mathcal G)^{(k+1)}$ such that $\theta_i \in T^*_{g_i} \mathcal G$.
\end{proposition}

\chapter{Deformations of VB-algebroids} \label{chap:VB_alg}

In this chapter we develop a deformation theory for VB-algebroids. Given a VB-algebroid $(W \Rightarrow E; A \Rightarrow M)$, in Subsection \ref{sec:lin_def} we show that the deformation complex of $W \Rightarrow E$ contains a distinguished subcomplex, the \emph{linear deformation complex} of $W$, that controls deformations of the VB-algebroid structure. We compute its low-degree cohomology groups and we describe it in terms of graded geometry and Poisson double vector bundles. Moreover, we recall from \cite{etv:infinitesimal} the ``classical description'' of the linear deformation complex and we deduce that the linear deformation complex is non-canonically isomorphic to a 3-term representation up to homotopy. The relationship between the linear deformation complex and the deformation complex of the base algebroid $A$ is also discussed. In Subsection \ref{sec:linearization} we introduce the \emph{linearization map}, adapting a technique from \cite{cabrera:hom}. This technical tool is used to show that the linear deformation cohomology of $W$ is embedded in the deformation cohomology of $W$ (as a Lie algebroid).

In Section \ref{Sec:examples} we present several applications of these results. In Subsection \ref{sec:VB_alg} we show that \emph{VB-algebras}, i.e.~vector bundles in the category of Lie algebras, are nothing more than Lie algebra representations, and we relate their linear deformation cohomologies with the classical deformation cohomology of a Lie algebra representation. In Subsection \ref{sec:LA-vect} we compute the linear deformation cohomology of \emph{LA-vector spaces}, i.e.~Lie algebroids in the category of vector spaces. Subsection \ref{sec:tangent-VB} is about the tangent and the cotangent VB-algebroids of a given Lie algebroid. In Subsections \ref{sec:partial_conn} and \ref{sec:Lie_vect} we describe other VB-algebroids arising from geometric situations and we show that their linear deformation cohomologies are strongly related to other well-known cohomologies. Subsection \ref{sec:type_1} is about a particular kind of VB-algebroid arising from a classification by Gracia-Saz and Mehta \cite{gracia-saz:vb}, while finally Subsection \ref{sec:RUTH} extends our discussion to a general representation up to homotopy of a Lie algebroid.

\section{The general theory}

\subsection{The linear deformation complex of a VB-algebroid} \label{sec:lin_def}

In this subsection we introduce the \emph{linear deformation complex of a VB-algebroid}. Actually, the whole discussion in Section \ref{sec:def_algd} extends to VB-algebroids. We skip most of the proofs: they can be carried out in a very similar way as for plain Lie algebroids.

We begin with a DVB $(W \to E; A \to M)$. Denote by $\mathfrak D^\bullet (W, E)$ the space of multiderivations of the vector bundle $W \to E$. As in Subsection \ref{sec:DVB}, denote by $h$ the homogeneity structure of $W \to A$. The action of $h$ induces a grading on the space of multiderivations.

\begin{definition} A multiderivation $c \in \mathfrak D^\bullet (W, E)$ is \emph{homogeneous of weight $q$} (or, simply, of \emph{weight} $q$) if  $h_\lambda^* c = \lambda^q c$ for every $\lambda > 0$. A multiderivation is \emph{linear} if it is of weight $0$, and it is \emph{core} if it is of weight $-1$.

We denote by $\mathfrak D^\bullet_q (W,E)$ the space of multiderivations of weight $q$, and by $\mathfrak D^\bullet_{\mathrm{lin}}(W,E)$ and $\mathfrak D^\bullet_{\mathrm{core}}(W,E)$, respectively, the spaces of linear and core multiderivations. \end{definition}

As $\Gamma_\mathrm{core}(W,E)$ and $\Gamma_\mathrm{lin}(W,E)$ generate $\Gamma(W,E)$ as a $C^\infty(E)$-module, a multiderivation is completely characterized by its action, and the action of its symbol, on linear and core sections. From Equation (\ref{eq:phi^*}) and the fact that there are no non-zero sections of weight less than $-1$, it then follows \emph{there are no non-zero multiderivations of weight less than $-1$.} Moreover:

\begin{proposition}\label{prop:lin_der}
Let $c$ be a $k$-derivation of $W \to E$. Then $c$ is linear if and only if all the following conditions are satisfied:
\begin{enumerate}
\item $c(\tilde \alpha_1, \dots, \tilde \alpha_k)$  is a linear section,
\item $c (\tilde \alpha_1, \dots, \tilde \alpha_{k-1}, \overline \chi_1)$ is a core section,
\item $c (\tilde \alpha_1, \dots, \tilde \alpha_{k-i}, \overline \chi_1, \dots, \overline \chi_i) = 0$,
\item $\sigma_{c} (\tilde \alpha_1, \dots, \tilde \alpha_{k-1})$ is a linear vector field,
\item $\sigma_{c} (\tilde \alpha_1, \dots, \tilde \alpha_{k-2}, \overline \chi_1)$ is a core vector field,
\item $\sigma_{c}(\tilde \alpha_1, \dots, \tilde \alpha_{k-i-1}, \overline \chi_1, \dots, \overline \chi_i) = 0$
\end{enumerate}
for all linear sections $\tilde \alpha_1, \dots, \tilde \alpha_k$, all core sections $\overline \chi_1, \dots, \overline \chi_i$ of $W \to E$, and all $i \geq 2$.
\end{proposition}

\proof This follows immediately from Lemma \ref{prop:char_lin} (see also \cite{etv:infinitesimal}). \endproof

In particular, a linear $k$-derivation is uniquely determined by its action on $k$ linear sections and on $k-1$ linear sections and a core section, and by the action of its symbol on $k-1$ linear sections and on $k-2$ linear sections and a core section. This observation will be made more precise in Theorem \ref{thm:classical}.

\begin{remark} If the rank of $E$ is greater than $0$, condition (1) in Proposition \ref{prop:lin_der} alone characterizes linear multiderivations. Indeed, suppose that $c \in \mathfrak D^k (W,E)$ takes linear sections to linear sections. Then the derivation $c(\tilde \alpha_1, \dots, \tilde \alpha_{k-1}, -)$ does the same for every $\tilde \alpha_1, \dots, \tilde \alpha_{k-1} \in \Gamma_{\mathrm{lin}}(W,E)$. It is easy to see in coordinates that its symbol takes linear functions to linear functions and, from Remark \ref{rmk:lin}, this means that $\sigma_c (\tilde \alpha_1, \dots, \tilde \alpha_{k-1})$ is a linear vector field. Now, from Formula (\ref{eq:phi^*sigma}), we get that $h_\lambda^* c = c$ and $\sigma_{h_\lambda^* c} = \sigma_c$ on linear sections, for every $\lambda > 0$. Finally, we observe that, away from the zero section of $E \to M$, $\Gamma(W,E)$ is generated over $C^\infty(E)$ by linear sections. When $\operatorname{rank} E > 0$, the zero section is a nowhere dense subset in $E$ and we conclude that $h_\lambda^* c = c$. 

If the rank of $E$ is zero, condition (1) is not enough. One can check that, in this case, every derivation is linear, but a biderivation that takes linear sections to linear sections does not necessarily take two core sections to 0. This is similar to what we have already observed for polyvector fields, in Section \ref{Sec:homogeneity}. \end{remark}

It immediately follows from (\ref{eq:pb_Gerst}) that $\mathfrak D^\bullet_{\mathrm{lin}} (W, E)[1]$ is a graded Lie subalgebra of $\mathfrak D^\bullet (W, E)[1]$. The following proposition is then straightforward.

\begin{proposition}\label{prop:def_VB}
VB-algebroid structures on the DVB $(W \to E; A \to M)$ are in one-to-one correspondence with Maurer-Cartan elements in $\mathfrak D^\bullet_{\mathrm{lin}} (W, E)[1]$.
\end{proposition}

Now, fix a VB-algebroid structure $(W \Rightarrow E; A \Rightarrow M)$ on the DVB $(W \to E; A \to M)$, and denote by $b_W$ the Lie bracket on sections of $W \to E$. We also denote by $C_{\mathrm{def}} (W, E)$ the deformation complex of the top algebroid $W \Rightarrow E$. It is clear that $b_W$ is a linear biderivation of $W \to E$, i.e.~$b_W \in \mathfrak D^2_{\mathrm{lin}} (W, E)$. Hence $\mathfrak D_{\mathrm{lin}} (W, E)[1]$ is a subDGLA of $C_{\mathrm{def}} (W, E)$, denoted $C_{\mathrm{def}, \mathrm{lin}}(W)$, and called the \emph{linear deformation complex} of $W \Rightarrow E$. Its cohomology is denoted $H_{\mathrm{def}, \mathrm{lin}} (W)$ and called the \emph{linear deformation cohomology} of $(W \Rightarrow E; A \Rightarrow M)$.

\begin{definition} A \emph{linear deformation} of $b_W$ (or simply a \emph{deformation}, if this does not lead to confusion) is a(n other) VB-algebroid structure on the DVB $(W \to E; A \to M)$. \end{definition}

Exactly as for Lie algebroids, Proposition \ref{prop:def_VB} is equivalent to say that \emph{the assignment $c \mapsto b_W + c$ establishes a one-to-one correspondence between Maurer-Cartan elements of $C_{\mathrm{def, lin}}(W)$ and deformations of $W$}.

Let $b_0, b_1$ be linear deformations of $b_W$. We say that $b_0$ and $b_1$ are  \emph{equivalent} if there exists a DVB isotopy taking $b_0$ to $b_1$, i.e.~a smooth path of DVB automorphisms $\phi_\epsilon: W \to W$, $\epsilon \in [0,1]$, such that $\phi_0 = \mathrm{id}_W$ and $\phi_1^* b_1 = b_0$. On the other hand, two Maurer-Cartan elements $c_0, c_1$ in $C_{\mathrm{def}, \mathrm{lin}}(W)$ are \emph{gauge-equivalent} if they are interpolated by a smooth path of $1$-cochains $c_\epsilon \in C_{\mathrm{def}, \mathrm{lin}}(W)$, and $(c_\epsilon)$ is a solution of the following ODE
\[
\dfrac{dc_\epsilon}{d\epsilon} = \delta \Delta_\epsilon + \llbracket c_\epsilon, \Delta_\epsilon \rrbracket,
\]  
for some smooth path of $0$-cochains $\Delta_\epsilon \in C_{\mathrm{def}, \mathrm{lin}}(W)$, $\epsilon \in [0,1]$. Equivalently,
\[
\frac{db_\epsilon}{d\epsilon} = \llbracket b_\epsilon, \Delta_\epsilon \rrbracket,
\]
where $b_\epsilon = b_W + c_\epsilon$.

\begin{proposition}
The DGLA $C_{\mathrm{def}, \mathrm{lin}}(W)$ controls deformations of the VB-algebroid $(W \Rightarrow E; A \Rightarrow M)$ in the following sense. Let $b_0 = b_W + c_0$, $b_1 = b_W + c_1$ be linear deformations of $b_W$. If $b_0$ and $b_1$ are equivalent, then $c_0$ and $c_1$ are gauge-equivalent. If $M$ is compact, the converse is also true. \end{proposition}

\proof
The proof is similar to that of Proposition \ref{prop:eq_G_eq}, with linear derivations replacing derivations and DVB automorphisms replacing vector bundle automorphisms. We only need to be careful when using the compactness hypothesis. Recall from \cite{etv:infinitesimal} that a linear derivation generates a flow by DVB automorphisms. In particular, if $(\Delta_\epsilon)$ is a time-dependent linear derivation of $W \to E$, then its symbol $X_\epsilon = \sigma (\Delta_\epsilon) \in \mathfrak X (E)$ is a linear vector field, hence it generates a flow by vector bundle automorphisms of $E$. From the compactness of $M$, it follows that $(X_\epsilon)$, hence the flow of $(\Delta_\epsilon)$, is complete.
\endproof 
 
\begin{remark}
An \emph{infinitesimal deformation} of $(W \Rightarrow E; A \Rightarrow M)$ is an element $c \in C^1_{\mathrm{def}, \mathrm{lin}}(W)$ such that $\delta c = 0$, i.e.~$c$ is a 1-cocycle in $C_{\mathrm{def}, \mathrm{lin}}(W)$. If $(c_\epsilon)$ is a smooth path of Maurer-Cartan elements starting at $0$, then $\frac{dc_\epsilon}{d\epsilon} \big|_{\epsilon=0}$ is an infinitesimal deformation of $(W \Rightarrow E; A \Rightarrow M)$. Similarly as for Lie algebroids, $H^1_{\mathrm{def}, \mathrm{lin}}(W)$ is the formal tangent space to the moduli space of linear deformations under gauge equivalence. It also follows from standard deformation theory arguments that $H^2_{\mathrm{def}, \mathrm{lin}}(W)$ contains obstructions to the extension of an infinitesimal linear deformation to a formal one. Finally, we interpret $0$-degree deformation cohomologies. It easily follows from the definition that $0$-cocycles in $C_{\mathrm{def}, \mathrm{lin}} (A)$ are \emph{infinitesimal multiplicative} (IM) derivations of $(W \Rightarrow E; A \Rightarrow M)$ i.e.~derivations of $W \to E$ generating a flow by VB-algebroid automorphisms \cite{etv:infinitesimal}. Among those, $1$-cocycles are \emph{inner IM derivations}, i.e.~IM derivations of the form $[\tilde \alpha, -]$ for some linear section $\tilde \alpha$ of $W \to E$. So $H^0_{\mathrm{def}, \mathrm{lin}}(W)$ consists of \emph{outer IM derivations}. See \cite{etv:infinitesimal} for more details.
\end{remark}

\subsubsection{Alternative descriptions}

Let $(W\Rightarrow E; A \Rightarrow M)$ be a VB-algebroid. Then $W \Rightarrow E$ is a Lie algebroid and from (\ref{eq:c_def}) we get $C_{\mathrm{def}}(W) \cong \mathfrak{X}(W[1]_E)$. Moreover, it is easy to see from Formula (\ref{eq:def_der}) and Proposition \ref{prop:lin_der} that a deformation cochain $c \in C_{\mathrm{def}}(W)$ is linear if and only if the corresponding vector field $\delta_c \in \mathfrak{X}(W[1]_E)$ is a linear vector field with respect to the vector bundle structure $W[1]_E \to A[1]$. So there is a canonical isomorphism of DGLAs
\[
C_{\mathrm{def}, \mathrm{lin}}(W) \overset{\cong}{\longrightarrow} \mathfrak{X}_{\mathrm{lin}}(W[1]_E), \quad c \mapsto \delta_c.
\]
As linear vector fields are equivalent to derivations, we also get a canonical isomorphism
\begin{equation}\label{eq:der_grad}
C_{\mathrm{def}, \mathrm{lin}}(W)  \overset{\cong}{\longrightarrow} \mathfrak D(W[1]_E, A[1]), \quad c \mapsto \Delta_c.
\end{equation}
This implies that $C_{\mathrm{def,lin}}(W)$ has an induced structure of DG-module over $C(A)$. From Equation (\ref{eq:X_Delta}), it follows that the latter isomorphism is characterized by
\[
(\Delta_c w)^\uparrow = [\delta_c, w^\uparrow]
\]
for every $w \in \Gamma(W[1]_E, A[1])$.

\

The next description of the linear deformation complex is in terms of Poisson DVBs. Recall that the total space of the dual vector bundle $W^*_E \to E$ is equipped with a fiber-wise linear Poisson structure $\pi$, and let $C_{\pi}(W^*_E)$ be the associated Lichnerowicz complex. As discussed in Remark \ref{rem:def_Pois}, the deformation complex of $W \Rightarrow E$ can be seen as the subcomplex $C_{\pi, \mathrm{lin}}(W^*_E) \subset C_{\pi}(W^*_E)$ consisting of polyvector fields that are linear with respect to the vector bundle structure $W \to E$. Since $(W \to E; A \to M)$ is a DVB, $W^\ast_E$ is also a vector bundle over $C^\ast$, and we can consider the subspace $C_{\pi, \mathrm{lin, lin}} (W^*_E)\subset C_{\pi, \mathrm{lin}}(W^*_E)$ of cochains that are linear with respect to both vector bundle structures $W^\ast_E \to E$ and $W^\ast_E \to C^\ast$. As $(W^\ast_E \to E, C^\ast \to M)$ is a Poisson DVB, $C_{\pi, \mathrm{lin, lin}} (W^*_E)$ is actually a subDGLA.

Now, let $c \in C^{k-1}_{\mathrm{def}}(W)$ and let $X_c$ be the corresponding polyvector field on $W^*_E$. From Equation (\ref{eq:ell}), for every $\lambda > 0$,
\begin{equation} \label{eq:lambda^{k-1}}
X_{h_\lambda^* c} = \lambda^{k-1} \left(h^{W^\ast}_\lambda\right)^* X_c,
\end{equation}
so $c$ is linear if and only if $X_c \in C_{\pi, \mathrm{lin, lin}} (W^*_E)$. We have proved that:
\begin{proposition}\label{theor:PoissonDVB}
Correspondence $c \mapsto X_c$ establishes an isomorphism of DGLAs $C_{\mathrm{def}, \mathrm{lin}}(W) \cong C_{\pi,\mathrm{lin,lin}}(W^*_E)$.
\end{proposition}

As in Subsection \ref{sec:VB-alg}, let $\widehat A$ denote the fat algebroid of $(W \Rightarrow E; A \Rightarrow M)$. In \cite{etv:infinitesimal}, it is shown that linear deformation cochains can be described in terms of more classical data.

\begin{theorem} [{\cite[Theorem 2.33]{etv:infinitesimal}}] \label{thm:classical}
  \label{theor:linear_cochains_alg} \quad
  \begin{enumerate}
  \item For $k>0$, there is a one-to-one correspondence $\tilde c \mapsto (c_{\widehat A}, c_E, c_C, \mathrm{d}_{\widehat A}) $ between cochains in $C^k_{\mathrm{def}, \mathrm{lin}}(W)$ and $4$-tuples consisting of
    \begin{itemize}
    \item a $k$-cochain $c_{\widehat A}$ in the deformation complex of
      the fat algebroid $\widehat A$,
    \item a vector bundle morphism $c_E : \wedge^k \widehat A \to
      D E$,
      \item a vector bundle morphism $c_C : \wedge^k \widehat A \to D C$ and
      \item a vector bundle morphism $\mathrm{d}_{\widehat A} : \wedge^{k-1} \widehat A \to \operatorname{Hom} (C,E)$, such that
    \end{itemize}
    \begin{equation*}
      \label{eq:derivation_valued}
      \sigma \circ c_E = \sigma \circ c_C = \sigma_{c_{\widehat A}} 
    \end{equation*}
    where $\sigma : D E \to TM$ (respectively, $\sigma : D C \to TM$) is the symbol map, and additionally
\[
\begin{aligned}
      c_{\widehat A} (\tilde \alpha_0, \dots, \tilde \alpha_{k-1}, \phi) 
      & = 
      c_C (\tilde \alpha_0, \dots, \tilde \alpha_{k-1}) \circ \phi - \phi \circ c_E (\tilde \alpha_0, \dots, \tilde \alpha_{k-1}), 
      \\
      c_E (\tilde \alpha_0, \dots, \tilde \alpha_{k-2}, \phi) 
      & = 
      - \mathrm{d}_{\widehat A} (\tilde \alpha_0, \dots, \tilde \alpha_{k-2}) \circ \phi, 
      \\
      c_C (\tilde \alpha_0, \dots, \tilde \alpha_{k-2}, \phi) 
      & = 
      - \phi \circ \mathrm{d}_{\widehat A} (\tilde \alpha_0, \dots, \tilde \alpha_{k-2}), 
      \\
      \mathrm{d}_{\widehat A} (\tilde \alpha_0, \dots, \tilde \alpha_{k-3}, \phi) 
      & = 0,
    \end{aligned}
\] 
    for all $\tilde \alpha_0, \dots, \tilde \alpha_{k-1} \in \Gamma (\widehat A)$, $\phi \in \operatorname{\mathfrak{Hom}}(E,C)$. The 4-tuple $(c_{\widehat A}, c_E, c_C, \mathrm d_{\widehat A})$ is given by:
\begin{equation} \label{eq:4-tuple}
\begin{aligned}
c_{\widehat A} (\tilde \alpha_0, \dots, \tilde \alpha_k) & = c(\tilde \alpha_0, \dots, \tilde \alpha_k); \\
X_{c_E (\tilde \alpha_0, \dots, \tilde \alpha_{k-1})} & = \sigma_c (\tilde \alpha_0, \dots, \tilde \alpha_{k-1}); \\
\overline{c_C (\tilde \alpha_0, \dots, \tilde \alpha_{k-1})(\chi)} & = c(\tilde \alpha_0, \dots, \tilde \alpha_{k-1}, \overline \chi) ; \\
\mathrm d_{\widehat A} (\tilde \alpha_0, \dots, \tilde \alpha_{k-2})(\chi)^\uparrow & = \sigma_c (\tilde \alpha_0, \dots, \tilde \alpha_{k-2}, \overline \chi)
\end{aligned}
\end{equation}
for every $\tilde \alpha_0, \dots, \tilde \alpha_k \in \Gamma(\widehat A)$, $\varepsilon \in \Gamma(E)$, $\chi \in \Gamma(C)$.

   If $\tilde c$ corresponds to $(c_{\widehat A}, c_E, c_C,
   \mathrm{d}_{\widehat A})$, then $\delta \tilde c$ corresponds to
   $(\delta c_{\widehat A}, c_E', c_C',\mathrm{d}'_{\widehat A})$, where
\begin{equation} \label{eq:c'_E}
\begin{aligned}
        c'_E (\tilde \alpha_0, \dots, \tilde \alpha_k) = & \sum_{i = 0}^k (-1)^i [\psi^s_{\tilde \alpha_i}, c_E (\tilde \alpha_0, \dots, \widehat{\tilde \alpha}_i, \dots, \tilde \alpha_k)] \\
        & + \sum_{i<j}(-1)^{i+j} c_E ([\tilde \alpha_i, \tilde \alpha_j], \tilde \alpha_0, \dots, \widehat{\tilde \alpha}_i, \dots, \widehat{\tilde \alpha}_j, \dots, \tilde \alpha_k) \\
        & - (-1)^{k+1} \psi^s_{c_{\widehat A}(\tilde \alpha_0, \dots, \tilde \alpha_k)},
\end{aligned}
\end{equation}
$c'_C$ is defined by an analogous formula and
\[
\begin{aligned}
        & \mathrm{d}'_{\widehat A} (\tilde \alpha_0, \dots, \tilde \alpha_{k-1}) \\
        & = \sum_{i = 0}^{k-1} (-1)^{i} \left( \psi^s_{\tilde{\alpha}_i} \circ \mathrm{d}_{\widehat A} (\tilde \alpha_0, \dots, \widehat{\tilde \alpha}_i, \dots, \tilde \alpha_{k-1}) - \mathrm{d}_{\widehat A} (\tilde \alpha_0, \dots, \widehat{\tilde \alpha}_i, \dots, \tilde \alpha_{k-1}) \circ \psi^c_{\tilde{\alpha}_i} \right) \\
        & \ + \sum_{i<j} (-1)^{i+j} \mathrm{d}_{\widehat A}([\tilde \alpha_i, \tilde \alpha_j], \tilde \alpha_0, \dots, \widehat{\tilde \alpha}_i, \dots,
        \widehat{\tilde \alpha}_j, \dots, \tilde \alpha_{k-1}) \\
        & \ + (-1)^{k+1} \left( c_E(\tilde \alpha_0, \dots, \tilde \alpha_{k-1}) \circ \partial - \partial \circ c_C (\tilde \alpha_0, \dots, \tilde \alpha_{k-1})\right),
\end{aligned}
\]
for all $\tilde \alpha_0, \dots, \tilde \alpha_k \in \Gamma (\widehat A)$.
  \item For $k=0$, there is a one-to-one correspondence $\tilde c \mapsto (c_{\widehat A}, c_E, c_C)$ between $1$-cochains in $C_{\mathrm{def,lin}}(W)$ and triples $(c_{\widehat A},c_E,c_{C})$, where $c_{\widehat A}$, $c_E$, $c_C$ are derivations of the vector bundles $\widehat A$, $E$, $C$ respectively, such that
\[
\sigma(c_E)=\sigma(c_C)=\sigma(c_{\widehat A})
\]
and 
\[
c_{\widehat A}(\phi)=c_C \circ \phi - \phi \circ c_E
\]
for all $\phi \in \operatorname{\mathfrak{Hom}}(E,C)$. The triple $(c_{\widehat A}, c_E, c_C)$ is given by
\begin{equation} \label{eq:triple}
c_{\widehat A} (\tilde \alpha) = c(\tilde \alpha), \quad X_{c_E} = \sigma_c, \quad \overline{c_C (\chi)} = c(\overline \chi)
\end{equation}
for every $\tilde \alpha \in \Gamma(\widehat A)$, $\chi \in \Gamma(C)$. If $\tilde c$ corresponds to $(c_{\widehat A},c_E,c_C)$, then $\delta \tilde c$ corresponds to the $4$-tuple $(\delta c_{\widehat A},c'_E, c'_C,\mathrm{d}'_{\widehat A})$, where
\[
\mathrm{d}'_{\widehat A}=\partial \circ c_C-c_E \circ \partial
\]
and, for all $\tilde \alpha \in \Gamma(\widehat{A})$,
\[
c'_E(\tilde \alpha) = [\psi^s_{\tilde \alpha},c_E]+\psi^s_{c_{\widehat A}(\tilde \alpha)},\qquad c'_C(\tilde \alpha)=[\psi^c_{\tilde \alpha},c_C]+\psi^c_{c_{\widehat A}(\tilde \alpha)}.
\]
\item Cochains in $C_{\mathrm{def,lin}}(W)$ of degree $-1$ are the same as sections of $\widehat A$.  If $\tilde \alpha \in \Gamma(\widehat A) = C^{-1}_{\mathrm{def, lin}}(W)$, then $\delta c$ corresponds to the triple $(\delta \tilde \alpha, \psi^s_{\tilde \alpha}, \psi^c_{\tilde \alpha})$.
  \end{enumerate}
\end{theorem}

This theorem implies that $(C_{\mathrm{def,lin}}(W), \delta)$ admits another description as a representation up to homotopy. More precisely, if we choose a horizontal lift $h: A \to \widehat A$ and a connection $\nabla$ on $\widehat A$, we get an isomorphism of DG-modules over $C(A)$
\[
C_{\mathrm{def,lin}}(W) \cong C(A,E_W),
\]
where $E_W$ is the complex (\ref{eq:E_W}), in the following way. Take $\tilde c \in C^k_{\mathrm{def,lin}}(W)$ and consider the associated 4-tuple $(c_{\widehat A}, c_E, c_C, \mathrm{d}_{\widehat A})$ of Theorem \ref{theor:linear_cochains_alg}. Then we can define an element $(\check c_{\widehat A}, \check c_E, \check c_C, \check{\mathrm{d}}_{\widehat A}) \in C^k(A, E_W)$ by
\[
\begin{aligned}
\check c_{\widehat A} & := L^\nabla_{c_{\widehat A}} \circ \wedge^{k+1} h; \\
\check c_E & := c_E \circ \wedge^k h; \\
\check c_C & := c_C \circ \wedge^k h; \\
\check{\mathrm{d}}_{\widehat A} & := \mathrm{d}_{\widehat A} \circ \wedge^{k-1} h,
\end{aligned}
\]
where $L^\nabla_{c_{\widehat A}}$ is defined by Formula (\ref{eq:L_c}). One can verify that $\tilde c \mapsto (\check c_{\widehat A}, \check c_E, \check c_C, \check{\mathrm{d}}_{\widehat A})$ is indeed an isomorphism of $C(A)$-modules and induces a differential on $C(A, E_W)$ that restricts to the differential (\ref{eq:E_W_diff}) in low degrees.

\

Suppose now that $W$ is a trivial-core VB-algebroid. Then Theorem \ref{thm:classical} specializes as follows:
\begin{corollary}
 \begin{enumerate}
  \item For $k>0$, there is a one-to-one correspondence $\tilde c \mapsto (c_A, c_E) $ between cochains in $C^k_{\mathrm{def}, \mathrm{lin}}(W)$ and pairs $(c_A, c_E)$, where $c_A \in C^k_{\mathrm{def}}(A)$, $c_E: \wedge^k A \to DE$ is a vector bundle map and
    \begin{equation*}
      \sigma \circ c_E = \sigma_{c_A}.
    \end{equation*}
The pair $(c_A, c_E)$ is given by Formulas (\ref{eq:4-tuple}).
If $\tilde c$ corresponds to $(c_A, c_E)$, then $\delta \tilde c$ corresponds to $(\delta c_{\widehat A}, c_E')$, where $c'_E$ is given by Formula (\ref{eq:c'_E}).
  \item There is a one-to-one correspondence $\tilde c \mapsto (c_A, c_E)$ between $C^0_{\mathrm{def,lin}}(W)$ and $\Gamma (DA \times_{TM} DE)$. If $\tilde c$ corresponds to $(c_A, c_E)$, then $\delta \tilde c$ corresponds to $(\delta c_A, c'_E)$, where
\[
c'_E(\alpha) = [\psi^s_{\alpha},c_E]+\psi^s_{c_A (\alpha)}
\]
for all $\alpha \in \Gamma(A)$. The pair $(c_A, c_E)$ is given by Formulas (\ref{eq:triple}).
\item Cochains in $C_{\mathrm{def,lin}}(W)$ of degree $-1$ are the same as sections of $A$.  If $\alpha \in \Gamma(A) = C^{-1}_{\mathrm{def, lin}}(W)$, then $\delta c$ corresponds to the pair $(\delta \alpha, \psi^s_{\alpha})$.
  \end{enumerate}
\end{corollary}

\subsubsection{Deformations of $A$ from linear deformations of $W$}

There is a natural surjection $C_{\mathrm{def}, \mathrm{lin}}(W) \to C_{\mathrm{def}}(A)$ which is easily described in the graded geometric picture: it is just the projection 
\[
\mathfrak{X}_{\mathrm{lin}}(W[1]_E) \to \mathfrak{X}(A[1]),
\]
of linear vector fields on the base. Equivalently, it is the symbol map
\[
\sigma : \mathfrak{D}(W[1]_E, A[1]) \to \mathfrak{X}(A[1]).
\]
In particular, we get a short exact sequence of DGLAs
\begin{equation}\label{eq:SESDGLAs}
0 \longrightarrow \operatorname{\mathfrak{End}}( W[1]_E) \longrightarrow \mathfrak{X}_{\mathrm{lin}}(W[1]_E) \longrightarrow \mathfrak{X}(A[1]) \longrightarrow 0.
\end{equation}
Equivalently, there is a short exact sequence
\[
0 \longrightarrow \operatorname{\mathfrak{End}}(W[1]_E) \longrightarrow C_{\mathrm{def},\mathrm{lin}}(W) \longrightarrow C_{\mathrm{def}}(A) \longrightarrow 0.
\]
In degree $-1$, this is just the sequence (\ref{eq:widehat}). Notice that the subDGLA $\operatorname{\mathfrak{End}}( W[1]_E)$ \emph{controls deformations of $(W \Rightarrow E; A \Rightarrow M)$ that fix $A \Rightarrow M$}, i.e.~deformations of $W$ that fix $b_A$ (the Lie algebroid structure on $A$) identify with Maurer-Cartan elements in $\operatorname{\mathfrak{End}}( W[1]_E)$. Finally, we obtain a long exact sequence
\begin{equation}\label{eq:LES_VB}
\longrightarrow H^k(\operatorname{\mathfrak{End}}(W[1]_E)) \longrightarrow H^k_{\mathrm{def}, \mathrm{lin}}(W) \longrightarrow H^k_{\mathrm{def}}(A) \longrightarrow H^{k+1}(\operatorname{\mathfrak{End}}(W[1]_E)) \longrightarrow \cdots
\end{equation}
connecting the linear deformation cohomology of $W$ with the deformation cohomology of $A$.

\begin{remark} We will not need a description of the subcomplex $\operatorname{\mathfrak{End}}( W[1]_E ) \subset \mathfrak{X}_{\mathrm{lin}}(W[1]_E)$ in terms of more classical data. However, we stress that this description exists in analogy with Theorem \ref{thm:classical}. \end{remark}

\subsubsection{Deformations of the dual VB-algebroid}

We conclude this section noticing that the linear deformation complex of a VB-algebroid is canonically isomorphic to that of its dual. Let $(W \Rightarrow E; A \Rightarrow M)$ be a VB-algebroid with core $C$, and let $(W^\ast_A \Rightarrow C^\ast ; A \Rightarrow M)$ be the dual VB-algebroid.

\begin{theorem} \label{thm:dual_vb_algd}
There is a canonical isomorphism of DGLAs 
\begin{equation} \label{eq:dual_VB_alg}
C_{\mathrm{def}, \mathrm{lin}}(W) \cong C_{\mathrm{def}, \mathrm{lin}}(W^*_A).
\end{equation}
\end{theorem}

\proof Recall from \cite{mackenzie} that the isomorphism (\ref{eq:beta})
\[
(W^\ast_E \to C^\ast; E \to M) \cong ((W^\ast_A)^\ast_{C^\ast} \to C^\ast; E \to M).
\]
is an isomorphism of Poisson DVBs. Using Proposition \ref{theor:PoissonDVB} we now see that there are isomorphisms of DGLAs
\[
C_{\mathrm{def}, \mathrm{lin}}(W) \cong C_{\pi,\mathrm{lin, lin}}(W^*_E) \cong C_{\pi,\mathrm{lin, lin}}((W^*_A)^*_{C^*}) \cong C_{\mathrm{def}, \mathrm{lin}}(W^*_A).
\] 
\endproof

\begin{remark} One can see in coordinates that the dual of the DG-vector bundle $W[1]_E \to A[1]$ is the DG-vector bundle $W^*_A [1]_{C^*} \to A[1]$. In this language, isomorphism (\ref{eq:dual_VB_alg}) translates into $\mathfrak D(W[1]_E, A[1]) \cong \mathfrak D (W^*_A [1]_{C^*}, A[1])$, so it is just an instance of (\ref{eq:graded_der}). In particular, it is implicitly defined by (\ref{eq:graded_der_2}). \end{remark}

\subsection{The linearization map} \label{sec:linearization}

Let $(W \Rightarrow E; A \Rightarrow M)$ be a VB-algebroid. We have shown that deformations of the VB-algebroid structure are controlled by a subDGLA $C_{\mathrm{def}, \mathrm{lin}} (W)$ of the deformation complex $C_{\mathrm{def}}(W)$ of the top Lie algebroid $W \Rightarrow E$. In the next section, we show that the inclusion $C_{\mathrm{def}, \mathrm{lin}} (W) \hookrightarrow C_{\mathrm{def}}(W)$ induces an inclusion $H_{\mathrm{def}, \mathrm{lin}} (W) \hookrightarrow H_{\mathrm{def}}(W)$ in cohomology. In particular, given an infinitesimal linear deformation that is trivial as infinitesimal deformation of the Lie algebroid $W \Rightarrow A$, i.e.~it is connected to the zero deformation by an infinitesimal isotopy of vector bundle maps, then it is also trivial as infinitesimal linear deformation, i.e.~it is also connected to the zero deformation by an infinitesimal isotopy of DVB maps.

The key idea is adapting to the present setting the ``homogenization trick'' of \cite{cabrera:hom}. Let $E \to M$ be a vector bundle. In their paper, Cabrera and Drummond consider the following natural projections from $C^\infty (E)$ to its $C^\infty (M)$-submodules $C^\infty_q(E)$ (of weight $q$ homogeneous functions):

\begin{equation}\label{eq:proj_fun}
\mathrm{pr}_q: C^\infty (E) \to C^\infty_q (E), \quad f \mapsto \dfrac{1}{q!} \dfrac{d^q}{d \lambda^q} \bigg|_{\lambda = 0} h_\lambda^* f.
\end{equation}

Notice that $\mathrm{pr}_q(f)$ is just the degree $q$ part of the (fiber-wise) Taylor series of $f$. In the following, we adopt the notations from Section \ref{Sec:homogeneity} and denote
\begin{equation}\label{eq:core,lin}
\begin{aligned}
\mathrm{core} := & \ \mathrm{pr}_0: C^\infty (E) \to C^\infty_{\mathrm{core}} (E), \quad f \mapsto f_{\mathrm{core}} = h_0^* f, \\
\mathrm{lin} := & \ \mathrm{pr}_1: C^\infty (E) \to C^\infty_\mathrm{lin} (E), \quad f \mapsto f_{\mathrm{lin}} =  \dfrac{d}{d \lambda} \bigg|_{\lambda = 0} h_\lambda^* f.
\end{aligned}
\end{equation}

Formula (\ref{eq:proj_fun}) does not apply directly to multiderivations. To see why, let $(W \to E, A \to M)$ be a DVB, let $h$ be the homogeneity structure of $W \to A$, and let $c \in \mathfrak D^\bullet (W,E)$. Then the curve $\lambda \mapsto h_\lambda^* c$ is not defined in $0$. Actually, $\lambda = 0$ is a ``pole of order $1$'' for $h_\lambda^\ast c$. 

In order to formalize this idea, let $(c_\lambda)$ be a curve in $\mathfrak D^k (W,E)$ defined around $0$. Then the limit $\lim_{\lambda \to 0} c_\lambda$ is the element in $\mathfrak D^k (W,E)$ given by:
\[
\big( \lim_{\lambda \to 0} c_\lambda \big) (w_1, \dots, w_k)_e := \lim_{\lambda \to 0} (c_\lambda (w_1, \dots, w_k)_e).
\]
Notice that the limit on the right is well-defined because $c_\lambda (w_1, \dots, w_k)_e \in W_e$ for every $\lambda$ and $W_e$ is a real finite-dimensional vector space, so it is endowed with a natural topology.

\begin{proposition} \label{prop:c_core}
The limit
\[
\lim_{\lambda \to 0} \lambda \cdot h_\lambda^* c
\]
exists and defines a core multiderivation $c_{\mathrm{core}}$.
\end{proposition}

\proof The existence of the limit can be shown in coordinates: for the sake of simplicity, we will sketch the proof for $w \in \mathfrak D^0 (W,E) = \Gamma (W,E)$. If $(x^i, a^\alpha, u^\beta, c^\gamma)$ are adapted coordinates on $W$ as in Subsection \ref{sec:DVB}, and $(A_\alpha)$, $(C_\gamma)$ are frames dual to the linear fiber coordinates $(a^\alpha)$, $(c^\gamma)$ (on $A$ and $C$, respectively), then $w$ locally looks like
\begin{equation} \label{eq:w}
w = w^{\alpha}(x,u) \tilde A_{\alpha} + w^\gamma(x,u) \overline C_\gamma
\end{equation}
and $h_\lambda$ reads
\[
h_\lambda (x,u,a,c) = (x, \lambda u, a, \lambda c)
\]
for every $\lambda > 0$. Now, a direct computation shows that
\[
\lambda \cdot h_\lambda^* w = \lambda \cdot w^\alpha (x, \lambda u) \tilde A_\alpha + w^\gamma (x, \lambda u) \overline C_\gamma
\]
and it is clear that the limit exists. Finally, for every $\mu \neq 0$,
\[
h_\mu^* c_{\mathrm{core}} = h_\mu^* \left( \lim_{\lambda \to 0} \lambda \cdot h_\lambda^* c \right) = \lim_{\lambda \to 0} \lambda \cdot h^*_\mu h^*_\lambda c = \lim_{\lambda \to 0} \mu^{-1} \left( \lambda \mu \cdot h^*_{\lambda \mu} c \right) = \mu^{-1} c_{\mathrm{core}}. 
\]
\endproof

Next proposition can be proved in the same way.

\begin{proposition}\label{prop:c_0}
The limit
\[
\lim_{\lambda \to 0} \big(h_\lambda^* c - \lambda^{-1} \cdot c_{\mathrm{core}}\big)
\]
exists and defines a linear multiderivation $c_{\mathrm{lin}}$. \end{proposition}

So far we have defined maps
\begin{equation}\label{eq:projs}
\begin{aligned}
\mathrm{core} & {} :  \mathfrak D^\bullet (W,E) \to \mathfrak D_{\mathrm{core}}^\bullet (W,E), \quad  c \mapsto \lim_{\lambda \to 0} \lambda \cdot  h_\lambda^* c \\
\mathrm{lin} & {} : \mathfrak D^\bullet (W,E) \to \mathfrak D_\mathrm{lin}^\bullet (W,E), \quad  c \mapsto \lim_{\lambda \to 0} \big(h_\lambda^* c - \lambda^{-1} \cdot c_{\mathrm{core}}\big)
\end{aligned}
\end{equation}
that split the inclusions in $\mathfrak D^\bullet (W,E)$. We call the latter the \emph{linearization map}.

\begin{remark} Once we have removed the singularity at 0, we can proceed as in (\ref{eq:proj_fun}) and define the projections on homogeneous multiderivation of positive weights $q > 0$:
\[
\mathrm{pr}_q: \mathfrak D^\bullet (W,E) \to \mathfrak D^\bullet_q (W,E), \quad c \mapsto \dfrac{1}{q!} \dfrac{d^q}{d \lambda^q} \bigg|_{\lambda = 0} (h_\lambda^* c - \lambda^{-1} \cdot c_{\mathrm{core}}).
\]
\end{remark}

Now, let $(W \Rightarrow E, A \Rightarrow M)$ be a VB-algebroid. Then we have a linearization map
\[
\mathrm{lin}: C_{\mathrm{def}}(W) \to C_{\mathrm{def,lin}}(W).
\]

\begin{theorem}[Linearization of deformation cochains] \label{theor:lin}
The linearization map is a cochain map splitting the inclusion $C_{\mathrm{def}, \mathrm{lin}}(W) \hookrightarrow C_{\mathrm{def}}(W)$. In particular there is a direct sum decomposition
\[
C_{\mathrm{def}}(W) \cong C_{\mathrm{def}, \mathrm{lin}}(W) \oplus \ker(\mathrm{lin}).
\]
of cochain complexes. Hence, the inclusion of linear deformation cochains into deformation cochains induces an injection
\begin{equation}\label{eq:incl}
H_{\mathrm{def}, \mathrm{lin}}(W) \hookrightarrow H_{\mathrm{def}}(W).
\end{equation}
\end{theorem}

\proof We only have to prove that the linearization preserves the differential $\delta = \llbracket b_W, - \rrbracket$ (here, as usual, $b_W$ is the Lie bracket on sections of $W \Rightarrow E$). Using the fact that $b_W$ is linear, we have that $\delta$ commutes with $h_\lambda^*$. From (\ref{eq:diff_algd}) it is obvious that $\delta$ preserves limits. So
\[
(\delta c)_{\mathrm{core}}  = \lim_{\lambda \to 0} \lambda \cdot h_\lambda^* (\delta c) = \lim_{\lambda \to 0} \lambda \cdot \delta(h_\lambda^* c) = \delta \left(\lim_{\lambda \to 0} \lambda \cdot h_\lambda^* c \right) = \delta c_{\mathrm{core}},
\]
and
\[
(\delta c)_\mathrm{lin}  = \lim_{\lambda \to 0} \left(h_\lambda^* (\delta c) - \lambda^{-1} \delta(c_{\mathrm{core}})\right) = \delta \left( \lim_{\lambda \to 0} (h_\lambda^* c - \lambda^{-1} \cdot c_{\mathrm{core}}) \right) = \delta c_{\mathrm{lin}},
\]
as desired. \endproof

The inclusion (\ref{eq:incl}) can be used to transfer vanishing results from deformation cohomology of the Lie algebroid $W \Rightarrow E$ to the linear deformation cohomology of the VB-algebroid $(W \Rightarrow E; A \Rightarrow M)$.  For example, if $H^0_{\mathrm{def}}(W) = 0$, every Lie algebroid derivation of $W \Rightarrow E$ is inner, and hence every IM derivation of the VB-algebroid $W$ is inner. Similarly, if $W \Rightarrow E$ has no non-trivial infinitesimal deformations, so does $(W \Rightarrow E; A \Rightarrow M)$, and so on.

For example, consider a vector bundle $E \to M$. Then (see Example \ref{ex:double_tang}) $(TE \Rightarrow E; TM \Rightarrow M)$ is a VB-algebroid, and we have:

\begin{proposition}\label{prop:TE}
The linear deformation cohomology of $(TE \Rightarrow E; TM \Rightarrow M)$ is trivial.
\end{proposition}

\proof
From Theorem \ref{theor:lin}, $H_{\mathrm{def}, \mathrm{lin}} (TE)$ embeds into the deformation cohomology $H_{\mathrm{def}} (TE)$ of the tangent algebroid $TE \Rightarrow E$, which is trivial (see, e.g., \cite{crainic:def}).
\endproof


\begin{remark}\label{rem:lin_sec}
Let $(W \to E; A \to M)$ be a DVB. Then $\Gamma(W,E) = \mathfrak D^0 (W,E)$. For $0$-derivations, maps (\ref{eq:projs}) read
\begin{equation}\label{eq:lin_core}
\begin{aligned}
\mathrm{core} & {} : \Gamma(W,E)  \to \Gamma_{\mathrm{core}}(W,E), \\
\mathrm{lin} & {} : \Gamma(W,E)  \to \Gamma_{\mathrm{lin}}(W,E).
\end{aligned}
\end{equation}
Now, let $w \in \Gamma(W,E)$. If $(x^i, a^\alpha, u^\beta, c^\gamma)$ are adapted coordinates on $W$ as before and $w$ is given by Equation (\ref{eq:w}), the
\begin{equation}\label{eq:eq:w_core,lin}
w_{\mathrm{core}}  = w^\gamma (x,0) \overline C_\gamma, \quad 
w_{\mathrm{lin}} = w^\alpha (x,0) \tilde A_\alpha + \dfrac{\partial w^\gamma}{\partial u^\beta} (x,0) u^\beta \overline C_\gamma.
\end{equation}
It is now easy to see, using e.g.~(\ref{eq:eq:w_core,lin}), that the maps (\ref{eq:lin_core}) satisfy the following identities:
\[
\begin{array}{ll}
(f \overline \chi)_\mathrm{core} = f_\mathrm{core} \overline \chi, &  (f \overline \chi)_\mathrm{lin} = f_\mathrm{lin} \overline \chi, \\
(f \tilde \alpha)_\mathrm{core} = 0, &  (f \tilde \alpha)_\mathrm{lin} = f_\mathrm{core} \tilde \alpha,\\
(\kappa w)_\mathrm{core} = \kappa w_\mathrm{core}, &  (\kappa w)_\mathrm{lin} = \kappa w_\mathrm{lin},\\
(\ell w)_\mathrm{core} = 0, &  (\ell w)_\mathrm{lin} = \ell w_\mathrm{core},
\end{array}
\]
for all  $f \in C^\infty(E)$, $w \in \Gamma(W,E)$, $\overline \chi \in \Gamma_{\mathrm{core}} (W, E)$, $\tilde \alpha \in \Gamma_{\mathrm{lin}} (W, E)$, $\kappa \in C^\infty_{\mathrm{core}} (E)$, and $\ell \in C^\infty_{\mathrm{lin}}(E)$.


When $W =TE$ is the tangent double of $E$, the maps (\ref{eq:lin_core}) specialize further and give maps
\begin{equation} \label{eq:lin_vect_fields}
\begin{aligned}
\mathrm{core} & {} :  \mathfrak X (E) \to \mathfrak X_{\mathrm{core}}(E) \\
\mathrm{lin} & {} : \mathfrak X (E) \to \mathfrak X_{\mathrm{lin}}(E).
\end{aligned}
\end{equation}
that interact with maps (\ref{eq:core,lin}) according to the following identities
\begin{equation}\label{eq:core,lin,rem2}
\begin{aligned}
X_\mathrm{core}(f_\mathrm{core}) &= 0,\\
 X_\mathrm{core}(f_\mathrm{lin}) & = X(f_\mathrm{lin})_\mathrm{core}, \\
X_\mathrm{lin}(f_\mathrm{core}) & = X(f_\mathrm{core})_\mathrm{core}, \\
X_\mathrm{lin}(f_\mathrm{lin}) & = X(f_\mathrm{lin})_\mathrm{lin},
\end{aligned}
\end{equation}
for all $f \in C^\infty (E)$, and $X \in \mathfrak X (E)$.

The first one of (\ref{eq:core,lin,rem2}) holds because $X_\mathrm{core}(f_\mathrm{core})$ is necessarily a function of weight $-1$. For the second one, compute:
\[
\begin{aligned}
X_\mathrm{core}(f_\mathrm{lin}) & = \lim_{\lambda \to 0} \lambda \cdot h_\lambda^* X (f_\mathrm{lin}) \\
& = \lim_{\lambda \to 0} \lambda (h_\lambda^* \circ X \circ h_{\lambda^{-1}}^*)(f_\mathrm{lin}) \\
& = \lim_{\lambda \to 0} h_\lambda^* (X(f_\mathrm{lin})) = X(f_\mathrm{lin})_\mathrm{core}.
\end{aligned}
\]
The third one can be proved in a similar way. Finally, using the previous identities, we have
\[
\begin{aligned}
X_\mathrm{lin}(f_\mathrm{lin}) & = \lim_{\lambda \to 0} \big(h_\lambda^* X - \lambda^{-1} X_\mathrm{core} \big) \big(f_\mathrm{lin} \big) \\
& = \lim_{\lambda \to 0} \left((h_\lambda^* \circ X \circ h_{\lambda^{-1}}^*) (f_\mathrm{lin}) - \lambda^{-1} X(f_\mathrm{lin})_\mathrm{core} \right) \\
& = \lim_{\lambda \to 0} \lambda^{-1} \left(h_\lambda^* (X(f_\mathrm{lin})) - h_0^*(X (f_\mathrm{lin})) \right) \\
& = \dfrac{d}{d\lambda} \bigg|_{\lambda = 0} h_\lambda^* (X (f_\mathrm{lin})) \\
& = X(f_\mathrm{lin})_\mathrm{lin}
\end{aligned}
\]
as announced.
\end{remark}

Remember from Subsection \ref{sec:lin_def} that a linear deformation cochain $c \in C^k_{\mathrm{def}, \mathrm{lin}}(W)$ is completely determined by its action on $k$ linear sections and on $k-1$ linear sections and a core section, and the action of its symbol on $k-1$ linear sections and on $k-2$ linear sections and a core section. We conclude this subsection providing a slightly more explicit description of the linearization map (\ref{eq:projs}) in terms of these restricted actions. This allows you to reconstruct $c_{\mathrm{lin}}$ in terms of the ``classical data'' of Theorem \ref{thm:classical}.

\begin{proposition} Let $c \in \mathfrak D^k (W,E)$. Then $c_{\mathrm{lin}}$ is completely determined by the following identities:
	\begin{enumerate}
		\item $c_{\mathrm{lin}}(\tilde \alpha_1, \dots, \tilde \alpha_k) = c(\tilde \alpha_1, \dots, \tilde \alpha_k)_{\mathrm{lin}}$,
		\item $c_{\mathrm{lin}}(\tilde \alpha_1, \dots, \tilde \alpha_{k-1}, \overline \chi) = c(\tilde \alpha_1, \dots, \tilde \alpha_{k-1}, \overline \chi)_{\mathrm{core}}$,
		\item $\sigma_{c_{\mathrm{lin}}}(\tilde \alpha_1, \dots, \tilde \alpha_{k-1}) = \sigma_c (\tilde \alpha_1, \dots, \tilde \alpha_{k-1})_{\mathrm{lin}}$,
		\item $\sigma_{c_{\mathrm{lin}}}(\tilde \alpha_1, \dots, \tilde \alpha_{k-2}, \overline \chi) = \sigma_c (\tilde \alpha_1, \dots, \tilde \alpha_{k-2}, \overline \chi)_{\mathrm{core}}$,
	\end{enumerate}
for all $\tilde \alpha_1, \dots, \tilde \alpha_k \in \Gamma_{\mathrm{lin}}(W,E)$, $\overline \chi \in \Gamma_{\mathrm{core}}(W,E)$.
\end{proposition}

\proof We first compute
\[
\begin{aligned}
c_\mathrm{core}(\tilde \alpha_1, \dots, \tilde \alpha_k) & = \lim_{\lambda \to 0} \lambda \cdot (h_\lambda^* c) (\tilde \alpha_1, \dots, \tilde \alpha_k) \\
& = \lim_{\lambda \to 0} \lambda \ h_\lambda^*(c(\tilde \alpha_1, \dots, \tilde \alpha_k)) \\
& = c(\tilde \alpha_1, \dots, \tilde \alpha_k)_{\mathrm{core}}.
\end{aligned}
\]
Then
\[
\begin{aligned}
c_\mathrm{lin}(\tilde \alpha_1, \dots, \tilde \alpha_k) & = \lim_{\lambda \to 0} \big( h_\lambda^* c - \lambda^{-1} \cdot c_\mathrm{core} \big) (\tilde \alpha_1, \dots, \tilde \alpha_k) \\ 
& = \lim_{\lambda \to 0} \left( h_\lambda^*(c(\tilde \alpha_1, \dots, \tilde \alpha_k)) - \lambda^{-1} c(\tilde \alpha_1, \dots, \tilde \alpha_k)_\mathrm{core} \right) \\
& = c(\tilde \alpha_1, \dots, \tilde \alpha_k)_\mathrm{lin}.
\end{aligned}
\]
Identity (2) in the statement can be proved in a similar way. To prove (3) first notice that
\[
\sigma_{c_{\mathrm{core}}} = \sigma_{\lim_{\lambda \to 0} \lambda \cdot h_\lambda^\ast c} = \lim_{\lambda \to 0} \lambda \cdot \sigma_{h_\lambda^\ast c} = \lim_{\lambda \to 0} \lambda \cdot h_\lambda^\ast \sigma_{ c},
\]
where we used (\ref{eq:phi^*sigma}). Hence
\[
\begin{aligned}
\sigma_{c_{\mathrm{core}}} (\tilde \alpha_1, \dots, \tilde \alpha_{k-1}) & =  \lim_{\lambda \to 0} (\lambda \cdot h_\lambda^\ast \sigma_{ c} )(\tilde \alpha_1, \dots, \tilde \alpha_{k-1}) \\
& = \lim_{\lambda \to 0} (\lambda \cdot h_\lambda^\ast (\sigma_{ c} (\tilde \alpha_1, \dots, \tilde \alpha_{k-1}))) \\
& = \sigma_{ c} (\tilde \alpha_1, \dots, \tilde \alpha_{k-1})_{\mathrm{core}}.
\end{aligned}
\]
Similarly,
\[
\begin{aligned}
\sigma_{c_{\mathrm{lin}}} & = \sigma_{\lim_{\lambda \to 0} \left(h_\lambda^\ast c - \lambda^{-1} \cdot c_{\mathrm{core}} \right)} \\
& = \lim_{\lambda \to 0} \sigma_{h_\lambda^\ast c - \lambda^{-1} \cdot c_{\mathrm{core}}} \\
&= \lim_{\lambda \to 0} \left(h_\lambda^\ast\sigma_{ c } - \lambda^{-1} \sigma_{c_{\mathrm{core}}}\right),
\end{aligned}
\]
hence
\[
\begin{aligned}
\sigma_{c_{\mathrm{lin}}} (\tilde \alpha_1, \dots, \tilde \alpha_{k-1}) & = \lim_{\lambda \to 0} \left((h_\lambda^\ast\sigma_{ c })(\tilde \alpha_1, \dots, \tilde \alpha_{k-1}) - \lambda^{-1} \sigma_{c_{\mathrm{core}}}(\tilde \alpha_1, \dots, \tilde \alpha_{k-1})\right) \\
& = \lim_{\lambda \to 0} \left(h_\lambda^\ast (\sigma_{ c }(\tilde \alpha_1, \dots, \tilde \alpha_{k-1})) - \lambda^{-1} \sigma_{c}(\tilde \alpha_1, \dots, \tilde \alpha_{k-1})_{\mathrm{core}}\right) \\
& = \sigma_{c} (\tilde \alpha_1, \dots, \tilde \alpha_{k-1})_{\mathrm{lin}}
\end{aligned}.
\]
Identity (4) can be proved in a similar way. \endproof

Finally, we will give another description of the linearization map. To do this, we need to define the linearization map for polyvector fields on the total space of a vector bundle, in analogy to what we have just done for multiderivations. So, let $E \to M$ be a vector bundle and denote by $h$ its homogeneity structure. If $X \in \mathfrak X^k_{\mathrm{poly}} (E)$, then the curve $\lambda \mapsto h_\lambda^* X$ has a ``pole of order $k$ in $0$'', in the sense of the following

\begin{proposition} 
The limit
\begin{equation} \label{eq:X_core}
\lim_{\lambda \to 0} \lambda^k \cdot h_\lambda^* X
\end{equation}
exists and defines a core polyvector field $X_{\mathrm{core}}$.
\end{proposition}

\proof Let $(x^i)$ be local coordinates on $M$ and $(u^\alpha)$ linear fiber coordinates on $E$. Then $X$ locally looks like
\[
X = \sum_{l = 0}^k \sum_{\substack{i_1 < \dots < i_l \\ \alpha_{l+1} < \dots < \alpha_k}} X^{i_1 \dots i_l \alpha_{l+1} \dots \alpha_k} (x,u) \dfrac{\partial}{\partial x^{i_1}} \wedge \dots \wedge \dfrac{\partial}{\partial x^{i_l}} \wedge \dfrac{\partial}{\partial u^{\alpha_{l+1}}} \wedge \dots \wedge \dfrac{\partial}{\partial u^{\alpha_k}},
\]
so
\[
\lambda^k \cdot h_\lambda^* X = \sum_{l = 0}^k \lambda^l \sum_{\substack{i_1 < \dots < i_l \\ \alpha_{l+1} < \dots < \alpha_k}} X^{i_1 \dots i_l \alpha_{l+1} \dots \alpha_k} (x, \lambda u) \dfrac{\partial}{\partial x^{i_1}} \wedge \dots \wedge \dfrac{\partial}{\partial x^{i_l}} \wedge \dfrac{\partial}{\partial u^{\alpha_{l+1}}} \wedge \dots \wedge \dfrac{\partial}{\partial u^{\alpha_k}}.
\]
From this formula it is clear that the limit exists. Finally, for every $\mu \neq 0$,
\[
h_\mu^* X_{\mathrm{core}} = h_\mu^* \left( \lim_{\lambda \to 0} \lambda^k \cdot h_\lambda^* X \right) = \lim_{\lambda \to 0} \lambda^k \cdot h^*_\mu h^*_\lambda X = \lim_{\lambda \to 0} \mu^{-k} \left( (\lambda \mu)^k \cdot h^*_{\lambda \mu} X \right) = \mu^{-k} X_{\mathrm{core}}. 
\]
\endproof

In the same way, one can prove the next proposition.

\begin{proposition} 
The limit
\begin{equation} \label{eq:X_lin}
\lim_{\lambda \to 0} \lambda^{k-1} \cdot (h_\lambda^* X - \lambda^{-k} X_{\mathrm{core}})
\end{equation}
exists and defines a linear polyvector field $X_{\mathrm{lin}}$.
\end{proposition}

It is clear that, in the case $k = 1$, Equations (\ref{eq:X_core}) and (\ref{eq:X_lin}) agree with the maps defined in (\ref{eq:lin_vect_fields}).

\

Now, let $(W \Rightarrow E, A \Rightarrow M)$ be a VB-algebroid. Recall that $C_{\mathrm{def}} (W)$ is isomorphic to the space $C_{\pi, \mathrm{lin}} (W^\ast_E)$ of linear polyvector fields on $W^\ast_E$, and, under this isomorphism, $C_{\mathrm{def}, \mathrm{lin}} (W)$ identifies with the subspace $C_{\pi, \mathrm{lin,lin}} (W^*_E)$ of polyvector fields that are linear with respect to both vector bundle structures: $W^\ast_E \to E$ and $W^\ast_E \to C^\ast$. Hence we have a diagram of DGLAs
\[
\begin{array}{c}
\xymatrix{C_{\mathrm{def}}(W) \ar[r]^{\mathrm{lin}} \ar[d]_\cong & C_{\mathrm{def,lin}}(W) \ar[d]^\cong \\
C_{\pi, \mathrm{lin}} (W^*_E) \ar[r]^{\mathrm{lin}_{C^*}} & C_{\pi, \mathrm{lin,lin}} (W^*_E)}
\end{array}
\]
where $\mathrm{lin}_{C^*}$ is the linearization map for polyvector fields on $W^*_E$ with respect to the vector bundle structure over $C^*$. Exploiting Equations (\ref{eq:lin_multiv}) and (\ref{eq:lambda^{k-1}}), one easily proves that the above diagram is commutative.

\section{Examples and applications}\label{Sec:examples}
 
In this section we provide several examples. Examples in Subsections \ref{sec:VB_alg}, \ref{sec:partial_conn} and \ref{sec:Lie_vect} parallel the analogous examples in \cite{crainic:def}, connecting our linear deformation cohomology to known cohomologies. Examples in Subsections \ref{sec:LA-vect}, \ref{sec:tangent-VB} and \ref{sec:type_1} are specific to VB-algebroids. Finally, in Subsection \ref{sec:RUTH} we discuss deformations of representations up to homotopy, with a special emphasis on 2-term representations up to homotopy.

\subsection{VB-algebras}\label{sec:VB_alg}

A \emph{VB-algebra} is a \emph{vector bundle object in the category of Lie algebras}. In other words, it is a VB-algebroid of the form
\[
\begin{array}{c}
\xymatrix{\mathfrak k \ar@{=>}[r] \ar[d]_p & 0 \ar[d] \\
 \mathfrak g \ar@{=>}[r]  & \ast}
 \end{array}.
\]

In particular, $\mathfrak k$ and $\mathfrak g$ are Lie algebras. Let $C := \ker p$ be the core of $(\mathfrak k \Rightarrow 0; \mathfrak g \Rightarrow \ast)$. The map $p$ has a canonical section in the category of Lie algebras, the zero section $0: \mathfrak g \to \mathfrak k$. By standard Lie algebra theory, $\mathfrak g$ acts on the core $C$ and there is a canonical isomorphism of Lie algebras
\[
\mathfrak k \overset{\cong}{\longrightarrow} \mathfrak g \ltimes C, \quad \xi \mapsto (p(\xi), \xi - 0_{p(\xi)}),
\]
where $\mathfrak g \ltimes C$ is the semidirect product Lie algebra. Thus VB-algebras are equivalent to Lie algebra representations. Notice that the $\mathfrak g$-action is just the core representation defined in (\ref{eq:core_rep}).


Let $\operatorname{End} C$ denote endomorphisms of the vector space $C$. In the present case, the short exact sequence (\ref{eq:SESDGLAs}) reads
\begin{equation}\label{eq:seq_vb-alg}
0 \longrightarrow C(\mathfrak g, \operatorname{End} C) \longrightarrow C_{\mathrm{def,lin}}(\mathfrak k) \longrightarrow C_{\mathrm{def}}(\mathfrak g) \longrightarrow 0,
\end{equation}
where $C(\mathfrak g, \operatorname{End} C) = \wedge \mathfrak g^* \otimes \operatorname{End} C$ is the Chevalley-Eilenberg complex of $\mathfrak g$ with coefficients in the induced representation $\operatorname{End} C$, and $C_{\mathrm{def}}(\mathfrak g) = (\wedge \mathfrak g^* \otimes \mathfrak g)[1]$ is the Chevalley-Eilenberg complex with coefficients in the adjoint representation. From the classical theory of Nijenhuis and Richardson \cite{nijenhuis:coh, nijenhuis:def, nijenhuis:def_lie}, the latter controls deformations of $\mathfrak g$, while the former controls deformations of the representation of $\mathfrak g$ on $C$.

The sequence (\ref{eq:seq_vb-alg}) has a natural splitting in the category of graded Lie algebras. Namely, there is an obvious graded Lie algebra map
\[
C_{\mathrm{def}}(\mathfrak g) \to C_{\mathrm{def,lin}}(\mathfrak k), \quad c \mapsto \tilde c
\]
given by
\[
\tilde c (v_0 + \chi_0, \dots, v_k + \chi_k) := c (v_0, \dots, v_k)
\]
for all $c \in C^k_{\mathrm{def}}(\mathfrak g) = \wedge^{k+1} \mathfrak g^\ast \otimes \mathfrak g$, and all $v_i + \chi_i \in \mathfrak k = \mathfrak g \oplus C$, $i = 0, \dots, k$. It is clear that the inclusion $C_{\mathrm{def}}(\mathfrak g) \to C_{\mathrm{def,lin}}(\mathfrak k)$ splits the projection $C_{\mathrm{def,lin}}(\mathfrak k) \to C_{\mathrm{def}}(\mathfrak g)$. Hence 
\begin{equation}\label{eq:split_VB_alg}
C_{\mathrm{def,lin}}(\mathfrak k) \cong C(\mathfrak g, \operatorname{End} C) \oplus C_{\mathrm{def}}(\mathfrak g).
\end{equation}
as graded Lie algebras. This isomorphism, which we denote $c \mapsto (c_1, c_2)$, is just given by the restrictions:
\[
\begin{aligned}
c_1 (v_0, \dots, v_{k-1})(\chi) & := c(v_0, \dots, v_{k-1}, \chi); \\
c_2 (v_0, \dots, v_k) & := c(v_0, \dots, v_k).
\end{aligned}
\]
However, (\ref{eq:split_VB_alg}) is not a DGLA isomorphism. We now describe the differential $\delta$ in $C_{\mathrm{def,lin}}(\mathfrak k)$ in terms of the splitting
(\ref{eq:split_VB_alg}). First of all, denote by $\theta : \mathfrak g \to \operatorname{End} C$ the action of $\mathfrak g$ on $C$, and let
\[
\Theta : \wedge \mathfrak g^\ast \otimes \mathfrak g \to \wedge \mathfrak g^\ast \otimes \operatorname{End} C
\]
be the map obtained from $\theta$ by extension of scalars. Explicitly, 
\[
\Theta(c)(v_0, \dots, v_k)(\chi) = (-1)^k c(v_0, \dots, v_k) \cdot \chi
\]
for every $c \in C^k_{\mathrm{def}}(\mathfrak g)$. From the properties of the action, $\Theta$ is actually a cochain map
\[
\Theta : C_{\mathrm{def}}(\mathfrak g)[-1] \to C(\mathfrak g, \operatorname{End} C).
\]
The isomorphism (\ref{eq:split_VB_alg}) identifies the differential in $C_{\mathrm{def,lin}}(\mathfrak k)$ with that of the mapping cone $\operatorname{Cone} (\Theta)$:
\[
(C_{\mathrm{def,lin}}(\mathfrak k), \delta) \cong \operatorname{Cone} (\Theta)
\]
as cochain complexes. To see this, take $c \in C^k_{\mathrm{def,lin}}(\mathfrak h)$, $(c_1, c_2)$ its image under the isomorphism (\ref{eq:split_VB_alg}), $((\delta c)_1, (\delta c)_2)$ the image of $\delta c$ under (\ref{eq:split_VB_alg}). Then it is clear that $(\delta c)_1 = \delta c_1$. Moreover, if $v_0, \dots, v_k \in \mathfrak g, \chi \in C$, we compute:

\begin{alignat*}{3}
(\delta c)_2 (v_0, \dots, v_k)(\chi) & = && \delta c (v_0, \dots, v_k, \chi) \\
& = \ && \sum_{i=0}^k (-1)^i [v_i, c(v_0, \dots, \widehat v_i, \dots, v_k, \chi)] \\
& && + (-1)^{k+1} [\chi, c(v_0, \dots, v_k)] \\
& && + \sum_{i<j} (-1)^{i+j} c([v_i, v_j], v_0, \dots, \widehat v_i, \dots, \widehat v_j, \dots, v_k, \chi) \\
& && + \sum_{i=0}^k (-1)^{i+k+1} c([v_i, \chi], v_0, \dots, \widehat{v_i}, \dots, v_k) \\
& = && \sum_{i=0}^k (-1)^i v_i \cdot c_2(v_0, \dots, \widehat{v_i}, \dots, v_k)(\chi) \\
& && + (-1)^k c_1(v_0, \dots, v_k) \cdot \chi \\
& && + \sum_{i<j} (-1)^{i+j} c_2([v_i, v_j], v_0, \dots, \widehat{v_i}, \dots \widehat{v_j}, \dots, v_k) \cdot \chi \\
& && + \sum_{i=0}^k (-1)^{i+k+1} c(v_i \cdot \chi, v_0, \dots, \widehat{v_i}, \dots, v_k) \\
& = && \sum_{i=0}^k (-1)^i v_i \cdot c_2(v_0, \dots, \widehat{v_i}, \dots, v_k)(\chi) \\
& && + (-1)^k c_1(v_0, \dots, v_k) \cdot \chi \\
& && + \sum_{i<j} (-1)^{i+j} c_2([v_i, v_j], v_0, \dots, \widehat{v_i}, \dots, \widehat{v_j}, \dots, v_k)(\chi) \\
& && + \sum_{i=0}^k (-1)^{i+1} c_2(v_0, \dots, \widehat{v_i}, \dots, v_k)(v_i \cdot \chi) \\
\end{alignat*}
\begin{alignat*}{3}
& = \ && \Theta(c)(v_0, \dots, v_k)(\chi) \\
& && + \sum_{i=0}^k (-1)^i (v_i \cdot c_2(v_0, \dots, \widehat{v_i}, \dots, v_k)(\chi) \\
& && - c_2(v_0, \dots, \widehat{v_i}, \dots, v_k)(v_i \cdot \chi)) \\
& && + \sum_{i<j} (-1)^{i+j} c_2([v_i, v_j], v_0, \dots, \widehat{v_i}, \dots, \widehat{v_j}, \dots, v_k)(\chi) \\
& = \ && \Theta(c)(v_0, \dots, v_k)(\chi) + \delta c_2 (v_0, \dots, v_k)(\chi).
\end{alignat*}

In this computation we used Equation (\ref{eq:core_rep}) and the fact that the induced action of $\mathfrak g$ on $\operatorname{End} C$ is given by:
\[
(v \cdot \phi)(c) = v \cdot \phi(c) - \phi(v \cdot c).
\]
Finally, notice that the long exact cohomology sequence of the mapping cone is just (\ref{eq:LES_VB}).

\subsection{LA-vector spaces}\label{sec:LA-vect}
An \emph{LA vector space} is a \emph{Lie algebroid object in the category of vector spaces}. In other words, it is a VB-algebroid of the form
\[
\begin{array}{c}
\xymatrix{V_1 \ar@{=>}[r] \ar[d] & V_0 \ar[d] \\
0 \ar@{=>}[r]  & \ast}
\end{array}.
\]
In particular, $V_1$ and $V_0$ are vector spaces. Now, let $C := \ker (V_1 \to V_0)$ be the core of $(V_1 \Rightarrow V_0; 0 \Rightarrow *)$. It easily follows from the definition of VB-algebroid that $V_1$ identifies canonically with the direct sum $C \oplus V_0$ and all the structure maps are completely determined by a linear map $\partial : C \to V_0$. Specifically, sections of $V_1 \to V_0$ are the same as smooth maps $V_0 \to C$, and given a basis $(C_\gamma)$ of $C$, the Lie bracket on maps $V_0 \to C$ is given by
\begin{equation}\label{eq:LA-vect1}
[f, g] = f^\gamma (\partial C_\gamma)^\uparrow g - g^\gamma (\partial C_\gamma)^\uparrow f, 
\end{equation}
where $f = f^\gamma C_\gamma$ and $g = g^\gamma C_\gamma$. It follows that the anchor $\rho : V_1 \to TV_0$ is given on sections by
\begin{equation}\label{eq:LA-vect2}
\rho (f) = f^\gamma (\partial C_\gamma)^\uparrow.
\end{equation}

Linear deformations of $(V_1 \Rightarrow V_0; 0 \Rightarrow *)$ are the same as deformations of $\partial$ as a linear map. Let us describe the linear deformation complex explicitly. 
As the bottom Lie algebroid is trivial, $C_{\mathrm{def}, \mathrm{lin}} (V_1)$ consists of graded endomorphisms $\operatorname{End} (C[1] \oplus V_0)$ of the graded vector space $C[1] \oplus V_0$. From (\ref{eq:LA-vect1}) and (\ref{eq:LA-vect2}) the differential $\delta$ in $\operatorname{End} (C[1] \oplus V_0)$ is just the commutator with $\partial$, meaning that the deformation cohomology consists of homotopy classes of graded cochain maps $(C[1] \oplus V_0, \partial) \to (C[1] \oplus V_0, \partial)$. More explicitly, $(\operatorname{End} (C[1] \oplus V_0), \delta)$ is concentrated in degrees $-1, 0, 1$. Namely, it is
\[
0 \longrightarrow \operatorname{Hom} (V_0, C)[1] \overset{\delta_0}{\longrightarrow}  \operatorname{End} (C) \oplus \operatorname{End} (V_0)  \overset{\delta_1}{\longrightarrow}  \operatorname{Hom} (C, V_0)[-1] \longrightarrow 0,
\]
where $\delta_0$ and $\delta_1$ are given by:
\[
\begin{aligned}
\delta_0 \gamma & = (\gamma \circ \partial, \partial \circ \gamma), \\
\delta_1 (\gamma_1, \gamma_2) & = \partial \circ \gamma_1 - \gamma_2 \circ \partial,
\end{aligned}
\]
where $\gamma \in \operatorname{Hom} (V_0, C)$, $\gamma_1 \in \operatorname{End}(C)$ and $\gamma_2 \in \operatorname{End}(V_0)$. We conclude immediately that 
\[
H_{\mathrm{def}, \mathrm{lin}} (V_1) = \operatorname{End} (\operatorname{coker} \partial \oplus \ker \partial [1])
\]
i.e.
\[
\begin{aligned}
H^{-1}_{\mathrm{def}, \mathrm{lin}} (V_1) & = \operatorname{Hom} (\operatorname{coker} \partial, \ker \partial),\\
H^0_{\mathrm{def}, \mathrm{lin}} (V_1) & = \operatorname{End} (\operatorname{coker} \partial) \oplus \operatorname{End} (\ker \partial), \\
H^1_{\mathrm{def}, \mathrm{lin}} (V_1) & = \operatorname{Hom}(\ker \partial, \operatorname{coker} \partial). 
\end{aligned}
\]
This shows, for instance, that infinitesimal deformations of a linear map $\partial: C \to V_0$ are all trivial if and only if $\partial$ is injective or surjective, as expected.

\subsection{Tangent and cotangent VB-algebroids}\label{sec:tangent-VB}

Let $A \Rightarrow M$ be a Lie algebroid.   Then $(TA \Rightarrow TM; A \Rightarrow M)$ is a VB-algebroid, as described in Example \ref{ex:tangent_VB}, and here we want to discuss its linear deformation cohomology. We use the graded geometric description: linear deformation cochains are linear vector fields on the total space of the vector bundle $TA[1]_{TM} \to A[1]$.
\begin{lemma}\label{lem:iota}
	Let $T A[1]$ be the tangent bundle of $A[1]$ and let $\tau : TA[1] \to A[1]$ be the projection. There is a canonical isomorphism of vector bundles of graded manifolds
	\[
	\xymatrix{ TA [1]_{TM} \ar[rr]^\iota \ar[dr]& & TA[1] \ar[dl] \\
		& A[1] &}
	\]
	uniquely determined by the following condition: 
	\begin{equation}\label{eq:iota}
	\langle \iota^\ast \ell_{d \omega}, T\alpha_1 \wedge \cdots \wedge T\alpha_k \rangle = \ell_{d \langle \omega , \alpha_1 \wedge \cdots \wedge \alpha_k \rangle}
	\end{equation}
	for all $\omega \in C^k (A)$, all sections $\alpha_1, \dots, \alpha_k \in \Gamma (A)$, and all $k$. Additionally
	\begin{equation}\label{eq:iota_2}
	\langle \iota^\ast \ell_{d \omega}, T\alpha_1 \wedge \cdots \wedge T\alpha_{k-1} \wedge \overline \alpha_k \rangle = \tau^\ast \langle \omega , \alpha_1 \wedge \cdots \wedge \alpha_k \rangle.
	\end{equation}
\end{lemma}
Formulas (\ref{eq:iota}) and (\ref{eq:iota_2}) in the statement require some explanations. The expression $d \omega$ in the left hand side should be interpreted as a $1$-form on $A[1]$, the de Rham differential of the function $\omega$, and $\ell_{d \omega}$ is the associated fiber-wise linear function on $T A[1]$. The pull-back of $\ell_{d\omega}$ along $\iota^\ast$ is a function on $TA[1]_{TM}$, i.e.~a $C^\infty (TM)$-valued, skew-symmetric multilinear map on sections of $TA \to TM$. The right hand side of (\ref{eq:iota}) is the fiber-wise linear function on $TM$ corresponding to the $1$-form $d \langle \omega , \alpha_1 \wedge \cdots \wedge \alpha_k \rangle$ on $M$. Here we interpret $\omega$ as a skew-symmetric multilinear map on sections of $A$.

\proof[Proof of Lemma \ref{lem:iota}]
	Let $(x^i)$ be coordinates on $M$, let $(u_\alpha)$ be a local basis of $\Gamma (A)$, and let $(u^\alpha)$ be the associated fiber-wise linear coordinates on $A$. These data determine coordinates $(x^i, \tilde u{}^\alpha)$ on $A[1]$ in the obvious way. In particular the $x^i$ have degree $0$ and the $\tilde u{}^\alpha$ have degree $1$. We also consider standard coordinates $(x^i, u^\alpha, \dot{x}{}^i, \dot{u}{}^\alpha)$ induced by $(x^i, u^\alpha)$ on $TA$. Notice that $(u^\alpha, \dot{u}{}^\alpha)$ are fiber-wise linear coordinates with respect to the vector bundle structure $TA \to TM$. More precisely, they are the fiber-wise linear coordinates associated to the local basis $(T u_\alpha, \overline u_\alpha)$ of $\Gamma (TA, TM)$. Next we denote by $(x^i, \dot{x}{}^i, \tilde u{}^\alpha, \tilde{\dot{u}}{}^\alpha)$ the induced coordinates on $TA[1]_{TM}$. They have degree $0,0,1,1$ respectively. Finally, we denote $(x^i, \tilde u{}^\alpha, X^i, \tilde U{}^\alpha)$ the standard coordinates on $TA[1]$ induced by $(x^i, \tilde u{}^\alpha)$. Define $\iota$ by putting
	\[
	\iota^\ast X^i = \dot{x}{}^i, \quad \text{and} \quad \iota^\ast \tilde U{}^\alpha = \tilde{\dot{u}}{}^\alpha.
	\]
	A direct computation exploiting the appropriate transition maps reveals that $\iota$ is globally well-defined. Now we prove (\ref{eq:iota}). We work in coordinates. Take a degree $k$ function $\omega = f_{\alpha_1 \cdots \alpha_k} (x) \tilde u{}^{\alpha_1} \cdots \tilde u{}^{\alpha_k}$ on $A[1]$. A direct computation shows that
	\begin{equation}\label{eq:iota_l}
	\iota^\ast \ell_{d \omega} = \frac{\partial f_{\alpha_1 \cdots \alpha_k} }{\partial x^i}  \tilde u{}^{\alpha_1} \cdots \tilde u{}^{\alpha_k} \dot{x}{}^i + k f_{\alpha_1 \cdots \alpha_k} \tilde u{}^{\alpha_1} \cdots  \tilde u{}^{\alpha_{k-1}} \tilde{\dot{u}}{}^{\alpha_k}.
	\end{equation}
	Now, let $\alpha_1, \dots, \alpha_k \in \Gamma (A)$, and $a = 1, \dots, k$. If $\alpha_a$ is locally given by $\alpha_a = g_a^\alpha (x) u_\alpha$, then
	\[
	T\alpha_a = \frac{\partial g^\alpha_a}{\partial x^i} \dot{x}{}^i \overline u{}_\alpha + g^\alpha_a Tu_\alpha,
	\]
	and, from (\ref{eq:iota_l}),
	\[
	\begin{aligned}
	\langle \iota^\ast \ell_{d \omega}, T \alpha_1 \wedge \cdots \wedge T \alpha_k \rangle 
	&=  k! \left( \frac{\partial f_{\alpha_1 \cdots \alpha_k} }{\partial x^i}  g_1^{\alpha_1} \cdots g_k^{\alpha_k}+  f_{\alpha_1 \cdots \alpha_k} g_1^{\alpha_1} \cdots  g_{k-1}^{\alpha_{k-1}} \frac{\partial g^{\alpha_k}_k}{\partial x^i} \right)\dot{x}{}^i \\
	& = k! \frac{\partial}{\partial x^i} \left( f_{\alpha_1 \cdots \alpha_k} g_1^{\alpha_1} \cdots g_k^{\alpha_k} \right) \dot{x}{}^i  = \ell_{d \langle \omega , \alpha_1 \wedge \cdots \wedge \alpha_k \rangle} .
	\end{aligned}
	\]
	Identity (\ref{eq:iota_2}) is proved in a similar way. To see that there is no other vector bundle isomorphism $\iota : TA[1]_{TM} \to TA[1]$ with the same property
	(\ref{eq:iota}) notice that $X^i = \ell_{dx^i}$ and $\tilde U{}^\alpha = \ell_{d\tilde u^\alpha}$. Now use (\ref{eq:iota}) to show that $\iota^\ast X^i = \dot{x}{}^i$ and $\iota^\ast \tilde U{}^\alpha = \tilde{\dot{u}}{}^\alpha$, necessarily.
\endproof

In the following we will understand the isomorphism $\iota$ of Lemma \ref{lem:iota}, and identify $TA[1]_{TM}$ with $TA[1]$. Now, recall that $C_{\mathrm{def}, \mathrm{lin}} (TA) = \mathfrak X_{\mathrm{lin}} (TA[1])$ fits in the short exact sequence of DGLAs:
\begin{equation}\label{eq:SES_TA}
0 \longrightarrow \operatorname{\mathfrak{End}}( TA[1]) \longrightarrow \mathfrak X_{\mathrm{lin}} (TA[1]) \longrightarrow \mathfrak X (A[1]) \longrightarrow 0.
\end{equation}
The tangent lift 
\begin{equation}\label{eq:tan}
\operatorname{tan} : \mathfrak X (A[1]) \hookrightarrow \mathfrak X (TA[1]), \quad X \mapsto X_{\mathrm{tan}}
\end{equation}
splits the sequence (\ref{eq:SES_TA}) in the category of DGLAs. As $\mathfrak X (A[1]) = C_{\mathrm{def}}(A)$, we immediately have the following
\begin{proposition}
	For every Lie algebroid $A \Rightarrow M$ there is a direct sum decomposition
	\[
	H_{\mathrm{def}, \mathrm{lin}} (TA) = H_{\mathrm{def}, \mathrm{lin}} (T^\ast A) = H (\operatorname{\mathfrak{End}} (TA[1])) \oplus H_{\mathrm{def}} (A).
	\]
\end{proposition}
In the last part of the subsection we describe the inclusion (\ref{eq:tan}) in terms of deformation cochains. This generalizes (\ref{eq:anchor_TA}) and (\ref{eq:bracket_TA}) to possibly higher cochains. In particular, in the next proposition $c_{\mathrm{tan}}$ is described in terms of classical data (see Theorem \ref{thm:classical}). Using the canonical isomorphisms $C_{\mathrm{def}, \mathrm{lin}} (TA) = \mathfrak X_{\mathrm{lin}} (TA[1])$, and $C_{\mathrm{def}} (A) = \mathfrak X (A[1])$ we get an inclusion
\[
\operatorname{tan} : C_{\mathrm{def}} (A) \hookrightarrow C_{\mathrm{def}, \mathrm{lin}} (TA), \quad c \mapsto c_{\mathrm{tan}}.
\]

\begin{proposition}\label{prop:tan_VB_class}
 Let $c \in C^{k-1}_{\mathrm{def}}(A)$. Then $c_{\operatorname{tan}} \in C^{k-1}_{\mathrm{def}, \mathrm{lin}}(TA)$ satisfies:
	\begin{enumerate}
		\item $c_{\operatorname{tan}}(T\alpha_1, \dots, T\alpha_k) = T c(\alpha_1, \dots, \alpha_k)$,
		\item $c_{\operatorname{tan}}(T\alpha_1, \dots, T\alpha_{k-1}, \overline \alpha_k) = \overline{c(\alpha_1, \dots, \alpha_k)}$,
		\item $\sigma_{c_{\operatorname{tan}}} (T\alpha_1, \dots, T\alpha_{k-1}) = \sigma_c (\alpha_1, \dots, \alpha_{k-1})_{\mathrm{tan}}$,
		\item $\sigma_{c_{\operatorname{tan}}} (T\alpha_1, \dots, T\alpha_{k-2}, \overline \alpha_{k-1})  = \sigma_c (\alpha_1, \dots, \alpha_{k-1})^\uparrow$,
	\end{enumerate}
	for all $\alpha_1, \dots, \alpha_k \in \Gamma(A)$. Identities (1)-(4) (together with the fact that $c_{\mathrm{tan}}$ is a linear cochain) determine $c_{\mathrm{tan}}$ completely.
\end{proposition}

\proof 
	We begin with (3). Recall that the tangent lift $X_{\mathrm{tan}}$ of a vector field $X$ is completely determined by (\ref{eq:tan_lift}). So, let $X \in \mathfrak{X}(A[1])$ be the graded vector field corresponding to $c$ (hence $X_{\mathrm{tan}} \in \mathfrak X (TA[1]_{TM})$ is the graded vector field corresponding to $c_{\mathrm{tan}}$), let $f \in C^\infty (M)$ and let $\alpha_1, \dots, \alpha_{k-1} \in \Gamma (A)$. Using (\ref{eq:delta_c_f}), compute
\[
\begin{aligned}
\sigma_{c_{\operatorname{tan}}}(T\alpha_1, \dots, T\alpha_{k-1}) \ell_{df} & = \langle X_{\operatorname{tan}}(\ell_{df}), T\alpha_1 \wedge \cdots \wedge T\alpha_{k-1} \rangle \\
& = \langle \ell_{d (X(f))}, T\alpha_1 \wedge \cdots \wedge T\alpha_{k-1}\rangle.
\end{aligned}
\]
From (\ref{eq:iota}),
\[
\begin{aligned}
\langle \ell_{d (X(f))}, T\alpha_1 \wedge \cdots \wedge T\alpha_{k-1}\rangle & = \ell_{d\langle X(f), \alpha_1\wedge \cdots \wedge \alpha_{k-1}\rangle} \\
& = \ell_{d (\sigma_c (\alpha_1, \dots, \alpha_{k-1}) f)} \\ 
& = \sigma_c (\alpha_1, \dots, \alpha_{k-1})_{\mathrm{tan}} \ell_{df}.
\end{aligned}
\]
As $\sigma_{c_{\operatorname{tan}}}(T\alpha_1, \dots, T\alpha_{k-1}) $ and $\sigma_c (\alpha_1, \dots, \alpha_{k-1})_{\mathrm{tan}} $ are both linear and they both project onto $\sigma_c (\alpha_1, \dots, \alpha_{k-1})$, this is enough to conclude that $\sigma_{c_{\operatorname{tan}}}(T\alpha_1, \dots, T\alpha_{k-1}) = \sigma_c (\alpha_1, \dots,\alpha_{k-1})_{\mathrm{tan}} $. Identity (4) can be proved in a similar way from (\ref{eq:iota_2}) using that
	\[
	\sigma_c (\alpha_1, \dots, \alpha_{k-1})^\uparrow \ell_{df} = \tau^\ast \langle df, \sigma_c (\alpha_1, \dots, \alpha_{k-1}) \rangle = \tau^\ast (\sigma_c (\alpha_1, \dots, \alpha_{k-1}) f). 
	\]

	We now prove (1). Both sides of the identity are linear sections of $TA \to TM$ and one can easily check in local coordinates that a linear section $\tilde \alpha$ is completely determined by pairings of the form $\langle \ell_{d\varphi}, \tilde \alpha \rangle$. Here, $\varphi$ is a section of $A^\ast \to M$ seen as a degree 1 function on $A[1]$, $d\varphi$ is its de Rham differential, and $\ell_{d\varphi}$ is the associated degree 1 fiber-wise linear function on $TA[1]$, which, in turn, can be interpreted as a $1$-form on the algebroid $TA \Rightarrow TM$, as in Lemma \ref{lem:iota}.
	
	So, take $\varphi \in \Gamma(A^\ast)$, $c \in C^{k-1}_{\operatorname{def}}(A)$, $\alpha_1, \dots, \alpha_k \in \Gamma(A)$, and compute
	\[
	\begin{aligned}
	& \langle \ell_{d\varphi}, c_{\operatorname{tan}}(T\alpha_1, \dots, T\alpha_k) \rangle \\
	& = \langle X_{\operatorname{tot}} (\ell_{d\varphi}), T \alpha_1 \wedge \dots \wedge T \alpha_k \rangle  - \sum_i (-1)^{k-i} \sigma_{c_{\operatorname{tan}}} (T \alpha_1, \dots, \widehat{T\alpha_i}, \dots, T\alpha_k) \langle \ell_{d \varphi}, T \alpha_i \rangle \\
	& = \langle \ell_{d(X(\varphi))}, T\alpha_1 \wedge \dots \wedge T\alpha_k \rangle - \sum_i (-1)^{k-i} \sigma_c (\alpha_1, \dots, \widehat{\alpha_i}, \dots, \alpha_k)_{\operatorname{tan}} \ell_{d \langle \omega, \alpha_i \rangle} \\
	& = \ell_{d \langle X(\varphi), \alpha_1 \wedge \dots \wedge \alpha_k \rangle} - \sum_i (-1)^{k-i} \ell_{d( \sigma_c (\alpha_1, \dots, \widehat{\alpha_i}, \dots, \alpha_k) \langle \omega, \alpha_i \rangle)} \\
	& = \ell_{d \langle \varphi, c(\alpha_1, \dots, \alpha_k) \rangle} = \langle \ell_{d \varphi}, Tc(\alpha_1, \dots, \alpha_k) \rangle,
	\end{aligned}
	\]
	where we used, in particular, (\ref{eq:def_der}), the first one in (\ref{eq:tan_lift}), Identity (3) in the statement and (\ref{eq:iota}). So (1) holds.
	
	Identity (2) can be proved in a similar way using (\ref{eq:def_der}), both identities (\ref{eq:tan_lift}), Identity (4), and (\ref{eq:iota_2}). 
\endproof

\begin{remark}
	Proposition \ref{prop:tan_VB_class} shows, in particular, that the Lie bracket $b_A$ on $\Gamma (A)$ and the Lie bracket $b_{TA}$ in $\Gamma (TA, TM)$ are related by $b_{TA} = (b_A)_{\mathrm{tan}}$.
\end{remark}

\subsection{Partial connections} \label{sec:partial_conn}

Let $M$ be a manifold, $D \subset TM$ an involutive distribution, and let $\mathcal F$ be the integral foliation of $D$. In particular $D \Rightarrow M$ is a Lie algebroid with injective anchor. A flat (partial) $D$-connection $\nabla$ in a vector bundle $E \to M$ defines a VB-algebroid 
\[
\begin{array}{c}
\xymatrix{\mathsf H \ar[d] \ar@{=>}[r] & E \ar[d] \\
	D \ar@{=>}[r] & M}
\end{array},
\]
where $\mathsf H \subset TE$ is the \emph{horizontal distribution} determined by $D$. Notice that the core of $(\mathsf H \Rightarrow E; D \Rightarrow M)$ is trivial, and every VB-algebroid with injective (base) anchor and trivial core arises in this way. Hence, (small) deformations of $(\mathsf H \Rightarrow E; D \Rightarrow M)$ are the same as simultaneous deformations of the foliation $\mathcal F$ and the flat partial connection $\nabla$. 

We now discuss the linear deformation cohomology. Denote by $q : E \to M$ the projection. First of all, the de Rham complex of $D \Rightarrow M$ is the same as leaf-wise differential forms $\Omega (\mathcal F)$ with the leaf-wise de Rham differential $d_{\mathcal F}$. Hence, the deformation complex of $D$ consists of derivations of $\Omega (\mathcal F)$ (the differential being the graded commutator with $d_{\mathcal F}$). As the core of $(\mathsf H \Rightarrow E; D \Rightarrow M)$ is trivial, $\mathsf H$ is canonically isomorphic to $q^* E$, via the isomorphism (\ref{eq:trivial_core_iso}). It easily follows that
 the linear deformation complex $(C_{\mathrm{def}, \mathrm{lin}}(\mathsf H), \delta)$ consists of derivations of the graded module $\Omega (\mathcal F, E)$ of $E$-valued, leaf-wise differential forms, and the differential $\delta$ is the commutator with the (leaf-wise partial) connection differential $d^\nabla_{\mathcal F}$. The kernel of $C_{\mathrm{def}, \mathrm{lin}}(\mathsf H) \to C_{\mathrm{def}} (D)$ consists of graded $\Omega(\mathcal F)$-linear endomorphisms of $\Omega (\mathcal F, E)$. The latter are the same as $\operatorname{End} E$-valued leaf-wise differential forms $\Omega (\mathcal F, \operatorname{End} E)$, and the restricted differential is the connection differential (corresponding to the induced connection in $\operatorname{End} E$).

Denote by $\nu = TM / D$ the normal bundle to $\mathcal F$ and by $TM \to \nu, X \mapsto \overline X$ the projection. The normal bundle is canonically equipped with the Bott (partial) connection $\nabla^{\mathrm{Bott}}$, given by
\[
\nabla^{\mathrm{Bott}}: \Gamma(D) \times \Gamma(\nu) \to \Gamma(\nu), \quad \nabla^{\mathrm{Bott}}_X \overline Y = \overline{[X,Y]}.
\]
Moreover, there is a deformation retraction, hence a quasi-isomorphism, $\pi : C_{\mathrm{def}}(D) \to \Omega (\mathcal F, \nu)$ that maps a deformation cochain $c \in C^k_{\mathrm{def}}(D)$ to the composition $\pi (c)$ of the symbol $\sigma_c : \wedge^k D \to TM$ followed by the projection $TM \to \nu$ \cite{crainic:def}. A similar construction can be applied to linear deformation cochains. To see this, first notice that derivations of $E$ modulo covariant derivatives along $\nabla$, $\mathfrak D (E) / \operatorname{im} \nabla$, are sections of a vector bundle $\tilde \nu \to M$. Additionally, $\tilde \nu$ is canonically equipped with a flat partial connection, also called the \emph{Bott connection} and denoted $\nabla^{\mathrm{Bott}}$, defined by
\[
\nabla^{\mathrm{Bott}}_{X} (\Delta (\operatorname{mod} \operatorname{im} \nabla)) = [\nabla_X, \Delta] (\operatorname{mod} \operatorname{im} \nabla)
\]
for all $\Delta \in \mathfrak D (E)$, and $X \in \Gamma (D)$. The symbol map $\sigma : \mathfrak D (E) \to \mathfrak X (M)$ descends to a surjective vector bundle map $\tilde \nu \to \nu$, intertwining the Bott connections. As $\operatorname{End} E \cap \operatorname{im} \nabla = 0$, we have $\ker (\tilde \nu \to \nu) = \operatorname{End} E $. In other words, there is a short exact sequence of vector bundles with partial connections:
\[
0 \longrightarrow \operatorname{End} E  \longrightarrow \tilde \nu  \longrightarrow \nu  \longrightarrow 0.
\]

Now, we define a surjective cochain map $\tilde \pi : C_{\mathrm{def}, \mathrm{lin}}(\mathsf H) \to \Omega (\mathcal F, \tilde \nu)$. Let $\tilde c \in C^k_{\mathrm{def,lin}}(\mathsf H)$ be a linear deformation cochain. Its symbol $\sigma_{\tilde c}$ maps linear sections of $\mathsf H \to E$ to linear vector field on $E$. As $\mathsf H \cong q^\ast D$, linear sections identify with plain sections of $D$. Accordingly $\sigma_{\tilde c}$ can be seen as a $\mathfrak D (E)$-valued $D$-form. Take this point of view and denote by $\tilde \pi (\tilde c) : \wedge^k D \to \tilde \nu$ the composition of $\sigma_{\tilde c}$ followed by the projection $\mathfrak D (E) \to \Gamma (\tilde \nu)$.  

Summarizing, we have the following commutative diagram
\[
\begin{array}{c}
\xymatrix{
	0 \ar[r] & \Omega (\mathcal F, \operatorname{End} E) \ar[r] \ar@{=}[d] & C_{\mathrm{def}, \mathrm{lin}} (\mathsf H) \ar[r] \ar[d]^{\tilde \pi} & C_{\mathrm{def}} (D) \ar[r]  \ar[d]^{\pi}  & 0 \\
	0 \ar[r] & \Omega (\mathcal F,  \operatorname{End} E) \ar[r] & \Omega (\mathcal F,  \tilde \nu ) \ar[r] & \Omega (\mathcal F,  \nu ) \ar[r] & 0
}
\end{array}.
\]
The rows are short exact sequences of DG-modules, and the vertical arrows are DG-module surjections. Additionally, $\pi$ is a quasi-isomorphism. Hence, it immediately follows from the Snake Lemma and the Five Lemma that $\tilde \pi$ is a quasi-isomorphism as well. We have thus proved the following

\begin{proposition} The map $\tilde \pi$ induces a canonical isomorphism of graded vector spaces between the linear deformation cohomology of the VB-algebroid $(\mathsf H \Rightarrow E; D \Rightarrow M)$, and the leaf-wise cohomology with coefficients in $\tilde \nu$:
	\[
	H_{\mathrm{def}, \mathrm{lin}} (\mathsf H) \cong H (\mathcal F, \tilde \nu).
	\]
\end{proposition}

\subsection{Lie algebra actions on vector bundles}\label{sec:Lie_vect}

Let $\mathfrak g$ be a (finite dimensional, real) Lie algebra acting on a vector bundle $E \to M$ by infinitesimal vector bundle automorphisms. In particular $\mathfrak g$ acts on $M$ and there is an associated action Lie algebroid $\mathfrak g \ltimes M \Rightarrow M$. Additionally, $\mathfrak g$ acts on the total space $E$ by linear vector fields; equivalently, there is a Lie algebra homomorphism $\mathfrak g \to \mathfrak D (E)$ covering the (infinitesimal) action $\mathfrak g \to \mathfrak X (M)$. It follows that $(\mathfrak g \ltimes E \Rightarrow E; \mathfrak g \ltimes M \Rightarrow M)$ is a VB-algebroid. We want to discuss linear deformation cohomologies of $\mathfrak g \ltimes E \Rightarrow E$. We begin reviewing Crainic and Moerdijk' remarks on the deformation cohomology of $\mathfrak g \ltimes M \Rightarrow M$ \cite{crainic:def} providing a graded geometric interpretation. 

By definition, $C^\infty(M)$ is a $\mathfrak g$-module. Following \cite{crainic:def}, we denote by $\mathfrak g_M := C^\infty(M) \otimes \mathfrak g$ the tensor product of $C^\infty(M)$ and the adjoint representation. Now, denote the action of $\mathfrak g$ on $C^\infty(M)$ by
\[
\mathfrak g \to \mathfrak X (M), \quad X \mapsto \alpha_X,
\]
and notice that $\mathfrak g$ also acts on $\mathfrak X (M)$ via
\[
\mathfrak g \to \operatorname{Der} \mathfrak X (M), \quad X \mapsto [\alpha_X, -].
\]
The deformation complex $C_{\mathrm{def}}(\mathfrak g \ltimes M)$ consists of vector fields on $(\mathfrak g \times M)[1]$. Notice that 
\[
C^\infty((\mathfrak g \times M)[1]) = \Omega_{\mathfrak g \times M} \cong \Gamma (\wedge \mathfrak g^* \times M, M) \cong C^\infty (M) \otimes \wedge \mathfrak g^* = C^\infty (M) \otimes C^\infty (\mathfrak g[1]).
\]
We will understand this isomorphism. It follows that there is a canonical projection $\pi_{\mathfrak g}: (\mathfrak g \times M)[1] \to \mathfrak g [1]$: it is the obvious morphism $M \to *$ on the underlying smooth manifolds and acts on functions by
\[
\pi_{\mathfrak g}^*: C^\infty (\mathfrak g[1]) \to C^\infty (M) \otimes C^\infty (\mathfrak g[1]), \quad \omega \mapsto 1 \otimes \omega.
\]
Now, let $\mathfrak X_{\mathrm{rel}}(\pi_{\mathfrak g})$ be the graded $C^\infty((\mathfrak g \times M)[1])$-module of $\pi_{\mathfrak g}$-relative vector fields, i.e. the module of $\mathbb R$-linear maps $X: C^\infty (\mathfrak g[1]) \to C^\infty((\mathfrak g \times M)[1])$ such that
\begin{equation} \label{eq:X_rel}
X(\omega \wedge \rho) = X(\omega) \cdot \pi_{\mathfrak g}^* \rho + (-1)^{|X| |\omega|} \pi_{\mathfrak g}^* \omega \cdot X(\rho)
\end{equation}
for every $\omega, \rho \in C^\infty (\mathfrak g[1])$. Then composition on the right with $\pi_{\mathfrak g}^*$ establishes a projection from vector fields on $(\mathfrak g \times M)[1]$ to $\mathfrak X_{\mathrm{rel}} (\pi_{\mathfrak g})$:
\begin{equation}\label{eq:proj}
\mathfrak X ( (\mathfrak g \times M)[1]) \to \mathfrak X_{\mathrm{rel}} (\pi_{\mathfrak g}), \quad X \mapsto X \circ \pi_{\mathfrak g}^\ast.
\end{equation}
The kernel of projection (\ref{eq:proj}) consists of $\pi_\mathfrak g$-vertical vector fields $\mathfrak X^{\pi_\mathfrak g} ((\mathfrak g \times M)[1])$. Denote by $d_{\mathfrak g} \in \mathfrak X ((\mathfrak g \times M)[1])$ the homological vector field on $(\mathfrak g \times M)[1]$. The graded commutator $\delta := [d_{\mathfrak g}, -]$ preserves $\pi_\mathfrak g$-vertical vector fields. Hence there is a short exact sequence of cochain complexes:
\begin{equation}\label{SES:gM}
0 \longrightarrow \mathfrak X^{\pi_\mathfrak g} ((\mathfrak g \times M)[1]) \longrightarrow \mathfrak X ((\mathfrak g \times M)[1]) \longrightarrow \mathfrak X_{\mathrm{rel}} (\pi_{\mathfrak g}) \longrightarrow 0.
\end{equation}
Now, $\mathfrak X ((\mathfrak g \times M)[1])$ is exactly the deformation complex of $\mathfrak g \ltimes M$. Moreover, from the Leibniz rule (\ref{eq:X_rel}) it follows that a $\pi_{\mathfrak g}$-relative vector field $X$ of degree $k$ is completely determined by its action on $\mathfrak g^*$, i.e. by the restriction $X: \mathfrak g^* \to C^\infty(M) \otimes \wedge^{k+1} \mathfrak g^*$ Therefore, the restriction gives a canonical isomorphism
\[
\mathfrak X_{\mathrm{rel}}(\pi_{\mathfrak g}) \overset{\cong}{\longrightarrow} C(\mathfrak g, \mathfrak g_M)[1].
\]
Finally, there is a canonical isomorphism
\[
C(\mathfrak g, \mathfrak X (M)) \overset{\cong}{\longrightarrow} \mathfrak X^{\pi_{\mathfrak g}} ((\mathfrak g \times M)[1]), \quad \omega \otimes X \mapsto \widehat{\omega \otimes X}
\]
given by $\widehat{\omega \otimes X}(\rho \otimes f) = (\omega \wedge \rho) \otimes X(f)$. So there is a short exact sequence of cochain complexes
\[
0 \longrightarrow C (\mathfrak g, \mathfrak X (M)) \longrightarrow C_{\mathrm{def}}(\mathfrak g \ltimes M) \longrightarrow C(\mathfrak g, \mathfrak g_M)[1] \longrightarrow 0,
\]
and a long exact cohomology sequence
\[
\cdots \longrightarrow H^k (\mathfrak g, \mathfrak X (M)) \longrightarrow H^k_{\mathrm{def}}(\mathfrak g \ltimes M) \longrightarrow H^{k+1} (\mathfrak g, \mathfrak g_M) \longrightarrow  H^{k+1}(\mathfrak g, \mathfrak X (M)) \longrightarrow \cdots.
\]
Observe that this is exactly the long exact sequence first described in \cite{crainic:def}.

We now pass to $\mathfrak g \ltimes E$. The linear deformation complex $C_{\mathrm{def}, \mathrm{lin}} (\mathfrak g \ltimes E)$ consists of linear vector fields on $(\mathfrak g \times E)[1]_E$. Similarly as above, we consider the projection $\tilde \pi_{\mathfrak g}: (\mathfrak g \times E)[1]_E \to \mathfrak g [1]$. Composition on the right with the pull-back $\tilde \pi_{\mathfrak g}^\ast$ establishes a projection:
\[
\mathfrak X_{\mathrm{lin}}( (\mathfrak g \times E)[1]_E) \to \mathfrak X_{\mathrm{rel}} (\pi_{\mathfrak g}), \quad X \mapsto X \circ \tilde \pi{}^\ast_{\mathfrak g}
\]
(beware, the range consists of $\pi_{\mathfrak g}$-relative, not $\tilde \pi_{\mathfrak g}$-relative vector fields) whose kernel consists of $\tilde \pi_\mathfrak g$-vertical linear vector fields $\mathfrak X^{\tilde \pi_\mathfrak g}_{\mathrm{lin}} ((\mathfrak g \times E)[1]_E)$. Hence there is a short exact sequence of cochain complexes:
\begin{equation}\label{SES:gE}
0 \longrightarrow \mathfrak X^{\tilde \pi_\mathfrak g}_{\mathrm{lin}} ((\mathfrak g \times E)[1]_E) \longrightarrow \mathfrak X_{\mathrm{lin}} ((\mathfrak g \times E)[1]_E) \longrightarrow \mathfrak X_{\mathrm{rel}} (\pi_{\mathfrak g}) \longrightarrow 0.
\end{equation}
Using the projection $\mathfrak X_{\mathrm{lin}} ((\mathfrak g \times E)[1]_E) \to \mathfrak X ((\mathfrak g \times M)[1])$, we can combine sequences (\ref{SES:gM}) and (\ref{SES:gE}) in an exact diagram
\[
\begin{array}{c}
\xymatrix@C=15pt@R=15pt{
& 0 \ar[d] & 0 \ar[d] & & \\
0 \ar[r] & \operatorname{\mathfrak{End}}((\mathfrak g \times E)[1]_E) \ar[d] \ar@{=}[r] & \operatorname{\mathfrak{End}} ((\mathfrak g \times E)[1]_E) \ar[r] \ar[d] & 0 \ar[d] \\
0 \ar[r] & \mathfrak X^{\tilde \pi_\mathfrak g}_{\mathrm{lin}} ((\mathfrak g \times E)[1]_E) \ar[r] \ar[d] & \mathfrak X_{\mathrm{lin}} ((\mathfrak g \times E)[1]_E) \ar[r] \ar[d] & \mathfrak X_{\mathrm{rel}} (\pi_{\mathfrak g}) \ar[r] \ar@{=}[d] & 0 \\
0 \ar[r] & \mathfrak X^{\pi_\mathfrak g} ((\mathfrak g \times M)[1]) \ar[r] \ar[d] & \mathfrak X((\mathfrak g \times M)[1]) \ar[r] \ar[d] & \mathfrak X_{\mathrm{rel}} (\pi_{\mathfrak g}) \ar[r] \ar[d] & 0 \\
 & 0 & 0 & 0 &
}
\end{array}.
\]
Now, $\mathfrak X_{\mathrm{lin}}((\mathfrak g \times E)[1]_E)$ is the linear deformation complex of $(\mathfrak g \ltimes E \Rightarrow E ; \mathfrak g \ltimes M \Rightarrow M)$. Moreover, $\mathfrak X^{\tilde \pi_\mathfrak g}_{\mathrm{lin}} ((\mathfrak g \times E)[1]_E)$ is canonically isomorphic to the Chevalley-Eilenberg cochain complex of $\mathfrak g$ with coefficients in $\mathfrak D (E)$. The isomorphism 
\[
C (\mathfrak g , \mathfrak D (E)) \overset{\cong}{\longrightarrow} \mathfrak X^{\tilde \pi_\mathfrak g}_{\mathrm{lin}} ((\mathfrak g \times E)[1]_E)
\]
maps a cochain $\omega \otimes \Delta$ to the vector field $\tilde \pi{}^\ast_{\mathfrak g}(\omega) X_{\Delta}$, where $X_\Delta$ is the unique $\tilde \pi_{\mathfrak g}$-vertical vector field on $(\mathfrak g \times E)[1]_E$ projecting on the linear vector field on $E$ corresponding to derivation $\Delta$. Finally, one can verify in coordinates that this isomorphism sends $C(\mathfrak g, \operatorname{\mathfrak{End}} E)$ to $\operatorname{\mathfrak{End}}((\mathfrak g \times E)[1]_E)$.

We conclude that there is an exact diagram of cochain complexes
\[
\begin{array}{c}
\xymatrix@C=15pt@R=15pt{
& 0 \ar[d] & 0 \ar[d] & & \\
0 \ar[r] & C (\mathfrak g, \operatorname{\mathfrak{End}} E) \ar[d] \ar@{=}[r] & C (\mathfrak g, \operatorname{\mathfrak{End}} E) \ar[r] \ar[d] & 0 \ar[d] \\
0 \ar[r] & C (\mathfrak g, \mathfrak{D} (E)) \ar[r] \ar[d] & C_{\mathrm{def}, \mathrm{lin}} (\mathfrak g\ltimes E) \ar[r] \ar[d] & C (\mathfrak g, \mathfrak g_M) [1] \ar[r] \ar@{=}[d] & 0 \\
0 \ar[r] & C (\mathfrak g, \mathfrak X(M)) \ar[r] \ar[d] & C_{\mathrm{def}}(\mathfrak g \ltimes M) \ar[r] \ar[d] & C (\mathfrak g, \mathfrak g_M) [1] \ar[r] \ar[d] & 0 \\
 & 0 & 0 & 0 &
}
\end{array}.
\]
This proves the following
\begin{proposition}
Let $\mathfrak g$ be a Lie algebra acting on a vector bundle $E \to M$ by infinitesimal vector bundle automorphisms. The linear deformation cohomology of the VB-algebroid $(\mathfrak g \ltimes E \Rightarrow E, \mathfrak g \ltimes M \Rightarrow M)$ fits in the exact diagram:
\[
\begin{array}{c}
\xymatrix@C=15pt@R=15pt{
& \vdots \ar[d] &  \vdots \ar[d] &  \vdots \ar[d] &  \vdots \ar[d] & \\
\ar[r] & H^k (\mathfrak g, \operatorname{\mathfrak{End}} E) \ar@{=}[r] \ar[d]& H^k(\mathfrak g, \operatorname{\mathfrak{End}} E) \ar[r] \ar[d]& 0 \ar[r] \ar[d] &  H^{k+1}(\mathfrak g, \operatorname{\mathfrak{End}} E) \ar[r] \ar[d]& \\
\ar[r] & H^k (\mathfrak g, \mathfrak D (E)) \ar[r] \ar[d]& H^k_{\mathrm{def}, \mathrm{lin}}(\mathfrak g \ltimes E) \ar[r] \ar[d]& H^{k +1} (\mathfrak g, \mathfrak g_M) \ar[r] \ar@{=}[d]&  H^{k+1}(\mathfrak g, \mathfrak D (E)) \ar[r] \ar[d]& \\
\ar[r] & H^k (\mathfrak g, \mathfrak X (M)) \ar[r] \ar[d] & H^k_{\mathrm{def}}(\mathfrak g \ltimes M) \ar[r] \ar[d]& H^{k +1}(\mathfrak g, \mathfrak g_M) \ar[r] \ar[d]&  H^{k+1}(\mathfrak g, \mathfrak X (M)) \ar[r] \ar[d] & \\
\ar[r] & H^{k+1} (\mathfrak g, \operatorname{\mathfrak{End}} E) \ar@{=}[r] \ar[d]& H^{k+1}(\mathfrak g, \operatorname{\mathfrak{End}} E) \ar[r] \ar[d]& 0 \ar[r] \ar[d] &  H^{k+2}(\mathfrak g, \operatorname{\mathfrak{End}} E) \ar[r] \ar[d]& \\
& \vdots  &  \vdots  &  \vdots  &  \vdots &
}
\end{array}.
\]
\end{proposition}

\subsection{VB-algebroids of type 1}\label{sec:type_1}

Let $(W \Rightarrow E; A \Rightarrow M)$ be a VB-algebroid with core $C$. According to a definition by Gracia-Saz and Mehta \cite{gracia-saz:vb}, $(W \Rightarrow E; A \Rightarrow M)$ is \emph{of type $1$} (resp.~\emph{of type $0$}) if the core-anchor $\partial: C \to E$ is an isomorphism (resp.~is the zero map). More generally, $(W \Rightarrow E; A \Rightarrow M)$ is \emph{regular} if the core-anchor has constant rank. In this case $(W \Rightarrow E; A \Rightarrow M)$ is the direct sum of a VB-algebroid of type 0 and a VB-algebroid of type 1, up to isomorphisms. So VB-algebroids of type 0 and 1 are the building blocks of regular VB-algebroids. In this subsection we discuss linear deformation cohomologies of VB-algebroids of type 1.

Let $(W \Rightarrow E; A \Rightarrow M)$ be a VB-algebroid of type 1, and denote by $p : E \to M$ the projection. Gracia-Saz and Mehta show that $(W \Rightarrow E; A \Rightarrow M)$ is canonically isomorphic to the VB-algebroid $(p^! A \Rightarrow E; A \Rightarrow M)$ \cite{gracia-saz:vb}, where $p^! A \Rightarrow E$ is the \emph{pull-back Lie algebroid} defined in the following way. Its total space $p^! A$ is the fibered product $p^! A := TE \tensor*[_{Tp}]{\times}{_\rho} A$. Hence, sections of $p^! A \to E$ are pairs $(X, \alpha)$, where $X$ is a vector field on $E$ and $\alpha$ is a section of the pull-back bundle $p^{\ast }A \to E$, with the additional property that $Tp (X_{e}) = \rho (\alpha_{p(e)})$ for all $e \in E$. Moreover, the anchor $p^! A \to TE$ is the projection $(X, \alpha) \mapsto X$, and the Lie bracket is uniquely defined by
\[
\left [  (X, p^{\ast }\alpha ), (Y, p^{\ast }\beta) \right ]  =
\left (  [X,Y], p^{\ast }[\alpha, \beta]\right ),
\]
on sections of the special form $(X, p^{\ast }\alpha)$, $(Y, p^{\ast} \beta )$, with $\alpha , \beta \in \Gamma (A)$. Finally, $(p^! A \Rightarrow E; A \Rightarrow M)$ is a VB-algebroid, and every VB-algebroid of type $1$ arises in this way (up to isomorphisms).

As $E \to M$ is a vector bundle, it has contractible fibers. So, according to \cite{spar:def}, $p^! A \Rightarrow E$ and $A \Rightarrow M$ share the same deformation cohomology. As an immediate consequence, we get that the canonical map $C_{\mathrm{def}, \mathrm{lin}} (p^! A) \to C_{\mathrm{def}} (A)$ induces an injection in cohomology. We want to show that, even more, it is a quasi-isomorphism. To do this it is enough to prove that the kernel $\operatorname{\mathfrak{End}} (p^! A [1]_E)$ of $C_{\mathrm{def}, \mathrm{lin}} (p^! A) \to C_{\mathrm{def}} (A)$ is acyclic. We use graded geometry again. So, consider the pull-back diagram
\[
\begin{array}{c}
\xymatrix{p^! A \ar[r] \ar[d]&  TE \ar[d]^{Tp} \\
A \ar[r]^-{\rho} & TM
} 
\end{array}.
\]
All vertices are vector bundles and, shifting by one the degree in their fibers, we get a pull-back diagram of DG-manifolds:
\[
\begin{array}{c}
\xymatrix{p^! A[1]_E \ar[r] \ar[d]_-{\tilde p}&  T[1]E \ar[d]^-{Tp} \\
A[1] \ar[r]^-{\rho} & T[1]M
} 
\end{array}.
\]
From Equations (\ref{eq:end_bundle}), (\ref{eq:pull-back}) and Remark \ref{rmk:T[1]M}, it follows that
\[
\begin{aligned}
\operatorname{\mathfrak{End}} (p^! A[1]_E) & = \operatorname{End}_{C(A)} \Gamma(p^! A[1]_E) \\
& = \operatorname{End}_{C(A)} (C(A) \otimes_{\Omega(M)} \Gamma (T[1] E)) \\
& \cong C(A) \otimes_{\Omega(M)} \operatorname{\mathfrak{End}} (T[1]E).
\end{aligned}
\]
Now, there is a canonical contracting homotopy $H: \mathfrak X (T[1]E) \to \mathfrak X (T[1]E)[-1]$: this implies, in particular, that the deformation cohomology of $TE \Rightarrow E$ is trivial. In order to define $H$, we recall some basics of \emph{Frölicher-Nijenhuis calculus} \cite{frolicher:theory}. Every $U \in \Omega^k(E,TE)$ defines an \emph{insertion operator} $i_U \in \operatorname{Der}^{k-1} \Omega(E) = \mathfrak X (T[1]E)^{k-1}$ and a \emph{Lie derivative} $\mathcal L_U = [d, i_U] \in \mathfrak X (T[1]E)^k$. Moreover, every $X \in \mathfrak X(T[1]E)$ can be written as $X = i_{U_X} + \mathcal L_{V_X}$ for uniquely defined $U_X, V_X \in \Omega (E,TE)$. Hence, the contracting homotopy $H$ is defined by
\[
H(X) = (-1)^{|X|} i_{V_X}.
\]
One can check that $H$ preserves both $\mathfrak X_{\mathrm{lin}}(T[1]E)$ and $\mathfrak{End} (T[1]E)$.


Finally, we define a contracting homotopy $h: \operatorname{\mathfrak{End}} (p^! A[1]_E) \to \operatorname{\mathfrak{End}}(p^! A[1]_E)[-1]$ by putting $h (\omega \otimes \Phi) := (-1)^{|\omega|} \omega \otimes H(\Phi)$, for all $\omega \in C(A)$, $\Phi \in \operatorname{\mathfrak{End}} (T[1] E)$. Summarizing, we proved the following

\begin{proposition} Let $(W \Rightarrow E; A \Rightarrow M)$ be a VB-algebroid of type 1. Then the canonical surjection $C_{\mathrm{def}, \mathrm{lin}}(W) \to C_{\mathrm{def}}(A)$ is a quasi-isomorphism. In particular, $H_{\mathrm{def}, \mathrm{lin}}(W) \cong H_{\mathrm{def}}(A)$. \end{proposition}

Morally, deforming a VB-algebroid of type 1 is the same as deforming its base Lie algebroid.

\subsection{Deformations of representations up to homotopy} \label{sec:RUTH}

We have already observed (Theorem \ref{thm:ruth}) that VB-algebroids are equivalent to 2-term representations up to homotopy of Lie algebroids. In this subsection, we want to have a look to deformations of general representations up to homotopy and relate them to our linear deformation complex.

\

Let $(E, D)$ be a representation up to homotopy of $A$. Notice preliminarily that $D$ actually encodes both the Lie algebroid structure on $A$ and the representation up to homotopy structure on $E$. Graded derivations of $\Omega_A (E)$, endowed with the graded commutator and the differential $\delta = [D, -]$, form a DGLA that we denote $C_{\mathrm{def}}(A,E, D)$. 

\begin{definition}
A \emph{simultaneous deformation} of the Lie algebroid structure on $A$ and the representation up to homotopy structure on $E$ is a pair consisting of a deformation of $A$ and a representation up to homotopy structure on $E$ with respect to the new Lie algebroid.
\end{definition}

Let $\gamma$ be a degree $1$ derivation of $\Omega_A (E)$ and put $\Delta = D + \gamma$. Then $\Delta$ is a simultaneous deformation of $(b_A, D)$ if and only if 
\[
\delta \gamma + \frac{1}{2}[\gamma, \gamma] = 0,
\]
i.e., simultaneous deformations of $(b_A, D)$ are in one-to-one correspondence with Maurer-Cartan elements of $C_{\mathrm{def}}(A, E, D)$.

As usual, there is a projection $\sigma: C_{\mathrm{def}}(A, E, D) \to C_{\mathrm{def}}(A)$: the symbol. Its kernel is the DGLA $\operatorname{\mathfrak{End}} (\Omega_A (E))$ of $\Omega_A$-linear endomorphisms of $\Omega_A (E)$, and we get a short exact sequence of DGLAs:
\[
0 \longrightarrow \operatorname{\mathfrak{End}} (\Omega_A (E)) \longrightarrow C_{\mathrm{def}}(A, E, D) \overset{\sigma}{\longrightarrow} C_{\mathrm{def}}(A) \longrightarrow 0.
\]
Deformations of $(b_A, D)$ that fix $b_A$ identify with Maurer-Cartan elements in the subDGLA $\operatorname{\mathfrak{End}} (\Omega_A (E))$.

Finally, let $(W \Rightarrow E; A \Rightarrow M)$ be a VB-algebroid with core $C$ and let $\tilde{E} =  C[1] \oplus E$. Choose an isomorphism of DVBs $W \cong A \times_M E \times_M C$, and let $D$ be the corresponding representation up to homotopy of $A$ on $\tilde{E}$. Then the module of sections of $W[1]_E \to A[1]$ is canonically isomorphic to $\Omega_A (\tilde E)$, hence from (\ref{eq:der_grad}) we get
\[
C_{\mathrm{def}, \mathrm{lin}}(W) \cong C_{\mathrm{def}}(A, \tilde E, D),
\]
as DGLAs, showing that our theory fits into this general framework.

\chapter{Deformations of VB-groupoids} \label{chap:VB_gr}

In the last chapter we pass to the global picture and study deformations of VB-groupoids. We follow a similar procedure to the one in Chapter \ref{chap:VB_alg}, and indeed most of the results there have a global counterpart. The connection between the infinitesimal and the global levels is explored in Subsection \ref{sec:van_est}, where we define a van Est map for the linear deformation cohomology. Moreover, in Subsection \ref{sec:morita} we deal with Morita invariance, that is specific to VB-groupoids. However, there are some differences that here we want to stress. In the case of VB-groupoids there is no known DGLA structure on the linear deformation complex and no description that parallels the graded-geometric one. Also a description of the linear deformation complex as a 3-term representation up to homotopy is not available at the moment. However, these structures should exist and further research should be pursued in this direction.

Now we describe our results in details. In Subsection \ref{sec:lin_def_1} we introduce the \emph{linear deformation complex} of a VB-groupoid $(\mathcal W \rightrightarrows E; \mathcal G \rightrightarrows M)$ as a subcomplex of the deformation complex of $\mathcal W \rightrightarrows E$, and we show that it controls infinitesimal deformations of the VB-groupoid structure. We describe it in terms of the cotangent groupoid $T^* \mathcal W$ and we study its relationship with the deformation complex of the base groupoid. When the core of the VB-groupoid is trivial, the linear deformation complex admits a much simpler description. Finally, we compute low-degree cohomology groups and we show that dual VB-groupoids have isomorphic linear deformation cohomologies. Subsection \ref{sec:linearization_1} is devoted to the construction of a \emph{linearization map}, essentially using results from Subsection \ref{sec:linearization}: this map is used to show that the linear deformation cohomology of a VB-groupoid is embedded into the deformation cohomology of the top groupoid. In Subsection \ref{sec:van_est} we define a van Est map connecting the linear deformation complex of a VB-groupoid and that of its VB-algebroid and we study its properties, while in Subsection \ref{sec:morita} we show that Morita equivalent VB-groupoids share the same linear deformation cohomology. The latter result shows that the linear deformation cohomology is really an invariant of the associated vector bundle of differentiable stacks.

In Section \ref{Sec:examples_1} we deal with examples. In Subsection \ref{sec:VB_grp} we show that \emph{VB-groups}, i.e.~vector bundles in the category of Lie groups, are equivalent to Lie group representations and we study the relationship between the deformation complexes of the two objects. In Subsection \ref{sec:2-vect} we compute the linear deformation cohomology of \emph{2-vector spaces}, i.e.~Lie groupoids in the category of vector spaces, and we show that it is isomorphic to its infinitesimal counterpart via the van Est map. Subsection \ref{sec:tangent_VB_grpd} is about the tangent and the cotangent VB-groupoids of a given Lie groupoid. Finally, in Subsections \ref{sec:fol_grpd} and \ref{sec:Lie_grp_vect} we describe VB-groupoids arising from representations of foliation groupoids and linear Lie group actions, and we relate their linear deformation cohomologies to other well-known cohomologies.

\section{The general theory}

\subsection{The linear deformation complex of a VB-groupoid} \label{sec:lin_def_1}

In this subsection we introduce the \emph{linear deformation complex of a VB-groupoid}. Let $(\mathcal W \rightrightarrows E; \mathcal G \rightrightarrows M)$ be a VB-groupoid, and let $(W \Rightarrow E; A \Rightarrow M)$ be its VB-algebroid. In this chapter, unless otherwise stated, $h$ will be the homogeneity structure of the vector bundle $\mathcal W \to \mathcal G$. By Remark \ref{rmk:hom_str}, $h_\lambda$ is a Lie groupoid automorphism for every $\lambda > 0$, so it acts on the deformation complex $C_{\mathrm{def}}(\mathcal W)$ of $\mathcal W \rightrightarrows E$. We say that a deformation cochain $\tilde c$ is \emph{linear} if
\begin{equation}\label{eq:lin_coch}
h_\lambda^* \tilde c = \tilde c
\end{equation}
for every $\lambda > 0$. Hence, linear cochains are those which are invariant under the action induced by the homogeneity structure.

We know from Subsection \ref{sec:def_grpd} that $h_\lambda^*$ commutes with $\delta$, for all $\lambda > 0$, so \emph{linear deformation cochains form a subcomplex of $C_{\mathrm{def}}(\mathcal W)$}. We denote the latter by $C_{\mathrm{def,lin}}(\mathcal W)$ and we call it the \emph{linear deformation complex} of $\mathcal W$. Its cohomology is called the \emph{linear deformation cohomology} of $\mathcal W$ and denoted $H_{\mathrm{def,lin}}(\mathcal W)$. Formula (\ref{eq:prod}) also shows that \emph{$C_{\mathrm{def}}(\mathcal W)$ is a $C (\mathcal G)$-module and $C_{\mathrm{def,lin}}(\mathcal W)$ is a $C(\mathcal G)$-submodule of $C_{\mathrm{def}}(\mathcal W)$}.

The action induced by the homogeneity structure $h$ on $C^{-1}_{\mathrm{def}}(\mathcal W)$ coincides, by definition, with the action induced by the homogeneity structure of $W \to A$ on $\Gamma(W,E)$, so $C^{-1}_{\mathrm{def,lin}}(\mathcal W)$ is simply the space $\Gamma_{\mathrm{lin}}(W,E)$ of linear sections of $W \to E$. For $k \geq 0$, Equation (\ref{eq:lin_coch}) is equivalent to saying that $\tilde c: \mathcal W^{(k+1)} \to T \mathcal W$ intertwines the homogeneity structures of $\mathcal W^{(k+1)} \to \mathcal G^{(k+1)}$ and $T \mathcal W \to T \mathcal G$. Again by Remark \ref{rmk:hom_str}, this means that $\tilde c$ is a vector bundle map over some map $c: \mathcal G^{(k+1)} \to T \mathcal G$. In this way we recover the definition in \cite{etv:infinitesimal}. 

For $k \geq 0$, a linear $k$-cochain $\tilde c$ can also be seen as an $\tilde {\mathsf s}$-projectable, linear section of the DVB $(p^*_{k+1} T \mathcal W \to \mathcal W^{(k+1)}; p^*_{k+1} T \mathcal G  \to \mathcal G^{(k+1)})$:
\[
\begin{array}{c}
\xymatrix{\tilde p^*_{k+1} T \mathcal W \ar[r] \ar[d] & \mathcal W^{(k+1)} \ar[d] \ar@/_1.25pc/[l]_-{\tilde c}\\
p^*_{k+1} T \mathcal G \ar[r] & \mathcal G^{(k+1)} \ar@/_1.25pc/[l]_-{ c}}
\end{array}
.
\]
Here we denote $\tilde p_k: \mathcal W^{(k)} \to \mathcal W$, and $\tilde q_k: \mathcal W^{(k)} \to E$ the maps (\ref{eq:pq}) for $\mathcal W$:
\[
\begin{aligned}
\tilde p_k: \mathcal W^{(k)} \to \mathcal W, \quad & (w_1, \dots, w_k) \mapsto w_1, \\
\tilde q_k: \mathcal W^{(k)} \to E, \quad & (w_1, \dots, w_k) \mapsto \mathsf s(w_1).
\end{aligned}
\]

There is another way to describe the linear deformation complex of a VB-groupoid. Let $(\mathcal W \rightrightarrows E; \mathcal G \rightrightarrows M)$ be a VB-groupoid, $C$ be its core, and let $(W \Rightarrow E; A \Rightarrow M)$ be its VB-algebroid. Consider the cotangent VB-groupoid of $\mathcal W \rightrightarrows E$, $(T^* \mathcal W \rightrightarrows W^*_E; \mathcal W \rightrightarrows E)$ (here, as usual, we denote by $W^\ast_E \to E$ the dual of the vector bundle $W \to E$). 

Actually, $T^* \mathcal W \rightrightarrows W^*_E$ is the top groupoid of another VB-groupoid. To see this, first take the dual of $\mathcal W$, that is $(\mathcal W^* \rightrightarrows C^*; \mathcal G \rightrightarrows M)$. Then, the cotangent VB-groupoid of $\mathcal W^* \rightrightarrows C^*$ is $(T^* \mathcal W^* \rightrightarrows (W^*_A)^*_{C^*}; \mathcal W^* \rightrightarrows C^*)$. Now, the isomorphism of DVBs (\ref{eq:B}) for $\mathcal W \to \mathcal G$ reads
\[ 
B: T^*\mathcal W \overset{\cong}{\longrightarrow} T^* \mathcal W^*
\]
and in \cite{mackenzie} it is proved that $B$ is also an isomorphism of Lie groupoids, covering the isomorphism of DVBs (\ref{eq:beta}):
\[
\beta : W^\ast_E \overset{\cong}{\longrightarrow} (W^\ast_A)^\ast_{C^\ast}.
\]

Combining all these maps, we obtain a diagram:

\begin{equation}\label{eq:dvb-grpd}
\begin{array}{c}
\xymatrix@C=10pt@R=15pt{
 & T^\ast \mathcal W \ar[dd]|!{[dl];[dr]}{\hole} \ar@<0.4ex>[rr] \ar@<-0.4ex>[rr] \ar[dl]^-{B}_-{\cong}& & W^\ast_E \ar[dl]^-{\beta}_-{\cong} \ar[dd] \\
 T^\ast \mathcal W^\ast \ar[dd] \ar@<0.4ex>[rr] \ar@<-0.4ex>[rr] & & (W^*_A)^*_{C^*} \ar[dd]& \\
 & \mathcal W^\ast \ar@<0.4ex>[rr]|!{[ur];[dr]}{\hole} \ar@<-0.4ex>[rr]|!{[ur];[dr]}{\hole} \ar@{=}[dl]& & C^\ast \ar@{=}[dl] \\
 \mathcal W^\ast \ar@<0.4ex>[rr] \ar@<-0.4ex>[rr] & & C^* & 
}
\end{array}.
\end{equation}

The maps $B$ and $\beta$ are isomorphisms of vector bundles over the identity, so the back face in the diagram (\ref{eq:dvb-grpd}) is also a VB-groupoid.

It follows that inside $C (T^* \mathcal W)$ there are two distinguished subcomplexes, those of cochains that are linear over $\mathcal W$ and over $\mathcal W^*$: we denote them by $C_{\mathrm{lin},\bullet}(T^* \mathcal W)$ and $C_{\bullet, \mathrm{lin}}(T^* \mathcal W)$ respectively. Moreover, denote $C_{\mathrm{proj},\bullet}(T^* \mathcal W)$ the subcomplex of left-projectable linear cochains over $\mathcal W$ and define $C_{\mathrm{lin,lin}}(T^* \mathcal W) := C_{\mathrm{lin}, \bullet}(T^* \mathcal W) \cap C_{\bullet, \mathrm{lin}}(T^* \mathcal W)$, $C_{\mathrm{proj,lin}}(T^* \mathcal W) := C_{\mathrm{proj}, \bullet}(T^* \mathcal W) \cap C_{\bullet, \mathrm{lin}}(T^* \mathcal W)$. We denote their cohomologies $H_{\mathrm{lin,lin}}(\mathcal W)$ and $H_{\mathrm{proj,lin}}(\mathcal W)$ respectively.

From Proposition \ref{prop:cotangent}, there is an isomorphism of $C (\mathcal W)$-modules
\begin{equation}\label{eq:iso}
C_{\mathrm{def}}(\mathcal W) \cong C_{\mathrm{proj},\bullet} (T^* \mathcal W)[1].
\end{equation}

It is easy to check that this isomorphism takes linear deformation cochains to cochains on $T^* \mathcal W$ that are linear over $\mathcal W^*$. So we get the following

\begin{proposition}\label{prop:cotangent_VB} There is an isomorphism of $C (\mathcal G)$-modules
\begin{equation}\label{eq:iso_lin}
C_{\mathrm{def,lin}}(\mathcal W) \cong C_{\mathrm{proj,lin}}(T^* \mathcal W)[1].
\end{equation}
\end{proposition}

For later use, we notice that a ``linear version'' of Lemma \ref{prop:proj;lin} holds. Namely, we have

\begin{lemma}\label{prop:proj,lin;lin,lin} The inclusion $C_{\mathrm{proj,lin}}(T^* \mathcal W) \hookrightarrow C_{\mathrm{lin,lin}}(T^* \mathcal W)$ induces an isomorphism in cohomology. \end{lemma}

\proof The proof of \cite[Lemma 3.1]{cabrera:hom} works identically in our setting without significant modifications. \endproof

\subsubsection{Deformations of $\mathcal G$ from linear deformations of $\mathcal W$} 
We have the following

\begin{proposition} [{\cite[Lemma 2.27]{etv:infinitesimal}}] \label{prop:proj} If $\tilde c: \mathcal W^{(k+1)} \to T\mathcal W$ belongs to $C^k_{\mathrm{def,lin}}(\mathcal W)$, then its projection $c: \mathcal{G}^{(k+1)} \to T\mathcal{G}$ belongs to $C^k_{\mathrm{def}}(\mathcal{G})$ and $\delta \tilde c$ projects to $\delta c$. \end{proposition}

It follows that there exists a natural cochain map:
\begin{equation}\label{eq:VB_pr}
C_{\mathrm{def,lin}}(\mathcal W) \to C_{\mathrm{def}}(\mathcal{G}).
\end{equation}
In degree $k = -1$, this is simply the projection $\Gamma_{\mathrm{lin}}(W,E) \to \Gamma(A)$ and we obtain the short exact sequence (\ref{eq:widehat}):
\[
0 \longrightarrow \mathfrak{Hom} (E,C) \longrightarrow \Gamma_{\mathrm{lin}}(W,E) \longrightarrow \Gamma(A) \longrightarrow 0,
\]
where, as usual, $\mathfrak{Hom} (E,C)$ is the $C^\infty(M)$-module of vector bundle morphisms $E \to C$. 

We now show that the map (\ref{eq:VB_pr}) is surjective for all $k \geq 0$. Let $c \in C^k_{\mathrm{def}}(\mathcal G)$. By Remark \ref{rmk:pair}, we have a diagram:
\[
\begin{array}{c}
\xymatrix{\mathcal G^{(k+1)} \ar[r]^{c} \ar[d] & p_{k+1}^* T \mathcal G \ar[d]^{T\mathsf s} \\
\mathcal G^{(k)} \ar[r]^{s_c} & q_k^* TM}
\end{array}
.
\]
We can lift $s_c$ to a linear section $s_{\tilde c}$ of $\tilde q_k^* TE \to \mathcal W^{(k)}$, so we obtain the following diagram
\begin{equation} \label{eq:fibration_dvb}
\begin{array}{c}
\xymatrix@C=9pt@R=12pt{
 & \tilde p^*_{k+1} T \mathcal W \ar[dd]|!{[dl];[dr]}{\hole} \ar[rr] \ar[dl] & & \tilde q_k^* TE \ar[dl] \ar[dd] \\
 \mathcal W^{(k+1)} \ar[dd] \ar[rr]  & & \mathcal W^{(k)} \ar@/_0.8pc/[ur]_-{s_{\tilde c}} \ar[dd]& \\
 & p^*_{k+1} T \mathcal G \ar[rr]|!{[ur];[dr]}{\hole} \ar[dl]& & q_k^* TM \ar[dl] \\
 \mathcal G^{(k+1)} \ar[rr]  \ar@/_0.8pc/[ur]_-{c}& & \mathcal G^{(k)} \ar@/_0.8pc/[ur]_-{s_{ c}}& 
}
\end{array},
\end{equation}
where the left and the right faces are DVBs and the horizontal arrows form a surjective DVB morphism. We would like to show that there exists a linear section $\tilde c: \mathcal W^{(k+1)} \to \tilde p^*_{k+1} T \mathcal W$ that projects on $c$. We will solve the problem in an abstract setting, in the following way.

Let $(W_1 \to E_1; A_1 \to M_1)$, $(W_2 \to E_2; A_2 \to M_2)$ be DVBs. A \emph{fibration of DVBs} is a DVB morphism 
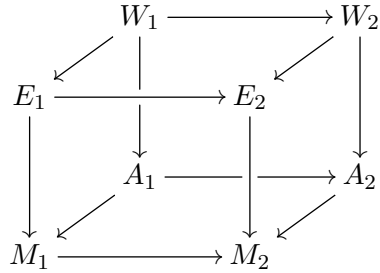
\[
\xymatrix@C=20pt@R=15pt{
& W_1 \ar[dl] \ar[rr] \ar[dd]|!{[dl];[dr]}{\hole} & & W_2 \ar[dl] \ar[dd] \\
E_1 \ar[rr] \ar[dd] & & E_2 \ar[dd] \\
& A_1 \ar[dl] \ar[rr]|!{[ur];[dr]}{\hole} & & A_2 \ar[dl] \\
M_1 \ar[rr] & & M_2 }
\]
where all the horizontal maps are surjective submersions.

\begin{lemma}\label{lemma:section} Consider a fibration of DVBs:
\[
\begin{array}{c}
\xymatrix@C=20pt@R=15pt{
& W_1 \ar[dl] \ar[rr]^(.6){\phi} \ar[dd]|!{[dl];[dr]}{\hole} & & W_2 \ar[dl] \ar[dd] \\
E_1 \ar[rr]^(.65){\phi_E} \ar[dd] & & E_2 \ar@/_0.7pc/[ur]|(.5){\tilde \alpha_2} \\
& A_1 \ar[dl] \ar[rr]|!{[ur];[dr]}{\hole} & & A_2 \ar[dl] \\
M_1 \ar[rr] \ar@/_0.7pc/[ur]|(.6){\alpha_1} & & M_2 \ar[uu] \ar@/_0.7pc/[ur]|(.6){\alpha_2} }
\end{array}
.
\]
Let $\tilde \alpha_2$ be a linear section of $W_2 \to E_2$ and let $\alpha_1$ be a section of $A_1$ that projects on the same section $\alpha_2$ of $A_2$. Then there exists a linear section of $W_1 \to E_1$ that projects simultaneously on $\tilde \alpha_2$ and $\alpha_1$. \end{lemma}

\proof We can decompose the morphism $\phi$ as follows:
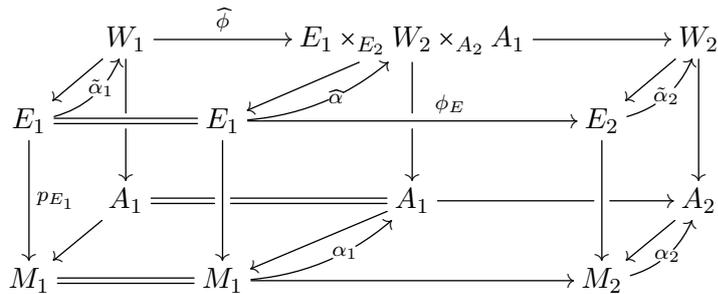
\[
\xymatrix@C=15pt@R=15pt{
& W_1 \ar[dl] \ar[rr]^-{\widehat \phi} \ar[dd]|!{[dl];[dr]}{\hole} & & E_1 \times_{E_2} W_2 \times_{A_2} A_1 \ar[dl] \ar[dd]|!{[dl];[dr]}{\hole} \ar[rr] & & W_2 \ar[dd] \ar[dl] \\
E_1 \ar@{=}[rr] \ar@/_0.7pc/[ur]|(.6){\tilde \alpha_1} \ar[dd]^-{p_{E_1}} & & E_1 \ar@/_0.7pc/[ur]|(.55){\widehat{\alpha}} \ar[rr]^(.6){\phi_E} \ar[dd] & & E_2 \ar[dd] \ar@/_0.7pc/[ur]|-{\tilde \alpha_2} \\
& A_1 \ar[dl] \ar@{=}[rr]|!{[ur];[dr]}{\hole} & & A_1 \ar[dl] \ar[rr]|!{[ur];[dr]}{\hole} & & A_2 \ar[dl] \\
M_1 \ar@{=}[rr] & & M_1 \ar@/_0.7pc/[ur]|(.6){\alpha_1} \ar[rr] & & M_2 \ar@/_0.7pc/[ur]|-{\alpha_2}
}
\]
Define
\[
\widehat \alpha: E_1 \to E_1 \times_{E_2} W_2 \times_{A_2} A_1, \quad \widehat \alpha_{e} = (e, \tilde \alpha_2 |_{\phi_E(e)}, \alpha_1 |_{p_{E_1}(e)}).
\]
It is clear that $\widehat \alpha$ is a well-defined linear section of $E_1 \times_{E_2} W_2 \times_{A_2} A_1 \to E_1$, that projects simultaneously on $\tilde \alpha_2$ and $\alpha_1$, so the problem is reduced to find a splitting of $\widehat \phi$. This can be done first locally, and then globally via the choice of a partition of unity on $M_1$. \endproof

Now, the diagram (\ref{eq:fibration_dvb}) is a fibration of DVBs, so we get the existence of the linear section $\tilde c: \mathcal W^{(k+1)} \to \tilde p^*_{k+1} T \mathcal W$, as desired.

Summarizing, there is a canonical short exact sequence of cochain maps:
\begin{equation}\label{eq:ses_VB}
0 \longrightarrow \ker \Pi \longrightarrow C_{\mathrm{def,lin}}(\mathcal W) \overset{\Pi}{\longrightarrow} C_{\mathrm{def}} (\mathcal{G})\longrightarrow 0.
\end{equation}

Next we compute $\ker \Pi$. By definition, a $k$-cochain $\tilde c: \mathcal W^{(k+1)} \to T\mathcal W$ is killed by $\Pi$ if and only if it takes values in the vertical bundle $T^{\tilde p} \mathcal W$ of $\tilde p: \mathcal W \to \mathcal G$. It is easy to check that $T^{\tilde p} \mathcal W \rightrightarrows T^p E$ is a subgroupoid of $T\mathcal W \rightrightarrows TE$. Moreover, $T^{\tilde p} \mathcal W \cong \mathcal W \times_{\mathcal G} \mathcal W$ and $T^p E \cong E \times_M E$ canonically as vector bundles. Under these isomorphisms, $T^{\tilde p} \mathcal W$ is identified with the groupoid $\mathcal W \times_{\mathcal G} \mathcal W \rightrightarrows E \times_M E$ with the component-wise structure maps. We will understand this identification.

It is clear that $\tilde c (w_0, \dots, w_k)$ has $w_0$ as first component, so one can think of $\tilde c$ as a vector bundle map $\mathcal W^{(k+1)} \to \mathcal W$ covering $p_{k+1}: \mathcal G^{(k+1)} \to \mathcal G$. Then $\ker \Pi$ is given by elements $\tilde c \in \mathfrak{Hom} (\mathcal W^{(k+1)}, p_{k+1}^* \mathcal W)$ such that $\tilde {\mathsf s} (c(w_0, \dots, w_k))$ does not depend on $w_0$ for any $(w_0, \dots, w_k) \in \mathcal W^{(k+1)}$. Finally, we observe that, on $\ker \Pi$, the differential is simply given by
\[
\begin{aligned}
\delta \tilde c (w_0, \dots, w_{k+1})  = &-\tilde{\bar{\mathsf m}} (\tilde c(w_0 w_1, \dots, w_{k+1}), \tilde c(w_1, \dots, w_{k+1})) \\ 
& + \sum_{i=1}^k (-1)^{i-1} \tilde c(w_0, \dots, w_i w_{i+1}, \dots, w_{k+1}) + (-1)^{k} \tilde c(w_0, \dots, w_k).
\end{aligned}
\]

\subsubsection{Trivial-core VB-groupoids}

Let $(\mathcal W \rightrightarrows E; \mathcal G \rightrightarrows M)$ be a trivial-core VB-groupoid. We want to show that its linear deformation complex has a particularly simple shape. First recall that, in Example \ref{ex:act_grpd}, we have observed that the total groupoid of a trivial-core VB-groupoid is canonically isomorphic to the action groupoid $\mathcal G \ltimes E \rightrightarrows E$  associated to a representation of the base groupoid $\mathcal G \rightrightarrows M$ on the side bundle $p: E \to M$. Therefore, without loss of generality, we will assume that $\mathcal W = \mathcal G \ltimes E$. As a vector bundle over $\mathcal G$, it is the pull-back $\mathsf s^\ast E$, and we denote it also by $E_{\mathcal G}$. We denote the action of an arrow $g: x \to y$ on an element $e \in E_x$ by $g \cdot e$, and the induced action of the groupoid $T \mathcal G \rightrightarrows TM$ on the vector bundle $TE \to TM$ by $\cdot_T$.

Consider a linear cochain $\tilde c \in C^k_{\mathrm{def,lin}}(E_{\mathcal G})$. By definition, it gives a commutative diagram:
\begin{equation}\label{diag:c_tilde}
\begin{array}{c}
\xymatrix@C=7pt@R=7pt{
& E_{\mathcal G}^{(k+1)} \ar[dl] \ar[rr]^{\tilde c} \ar[dd]|!{[dl];[dr]}{\hole} & & TE_{\mathcal G} \ar[dl] \ar[dd] \\
E_{\mathcal G}^{(k)} \ar[rr]^(.6){s_{\tilde c}} \ar[dd] & & TE \ar[dd] \\
& \mathcal G^{(k+1)} \ar[dl] \ar[rr]^(.4)c |!{[ur];[dr]}{\hole} & & T \mathcal G  \ar[dl] \\
\mathcal G^{(k)} \ar[rr]^(.6){s_c} & & TM  }
\end{array}.
\end{equation}
But, for every $k$, there is a canonical isomorphism 
\begin{equation} \label{eq:nerve_act}
E_{\mathcal G}^{(k)} \overset{\cong}{\longrightarrow} \mathcal G^{(k)} \mathbin{{}_{\mathsf s} \times_p} E, \quad ((g_1, e_1), \dots, (g_k, e_k)) \longmapsto ((g_1, \dots, g_k); e_k)
\end{equation}
where $\mathsf s$ is the projection $\mathsf s: \mathcal G^{(k)} \to M$, $(g_1, \dots, g_k) \mapsto \mathsf s(g_k)$. We will often understand this isomorphism. We also have that $TE_{\mathcal G} \cong T \mathcal G \mathbin{{}_{T \mathsf s} \times_{Tp}} TE$. So, we get the following alternative description of (\ref{diag:c_tilde}):
\begin{equation} \label{diag:c_tilde_2}
\begin{array}{c}
\xymatrix@C=0pt@R=7pt{
& \mathcal G^{(k+1)} \mathbin{{}_{\mathsf s} \times_p} E \ar[dl] \ar[rr]^{\tilde c} \ar[dd]|!{[dl];[dr]}{\hole} & & T \mathcal G \mathbin{{}_{T \mathsf s} \times_{Tp}} TE \ar[dl] \ar[dd] \\
\mathcal G^{(k)} \mathbin{{}_{\mathsf s} \times_p} E \ar[rr]^(.65){s_{\tilde c}} \ar[dd] & & TE \ar[dd] \\
& \mathcal G^{(k+1)} \ar[dl] \ar[rr]^(.4)c |!{[ur];[dr]}{\hole} & & T \mathcal G  \ar[dl] \\
\mathcal G^{(k)} \ar[rr]^(.65){s_c} & & TM  }
\end{array},
\end{equation}
where the vertical arrows, except for the front right one, are projections onto the first factor. In particular, $\tilde c$ is fully determined by $c$ and $s_{\tilde c}$. Set $\tilde c_1 := c$ and observe that, for every $(g_1, \dots, g_k) \in \mathcal G^{(k)}$, $s_{\tilde c} ((g_1, \dots, g_k); -)$ is a linear map $E_{\mathsf s (g_k)} \to TE|_{s_c (g_1, \dots, g_k)}$. Consider the linear map
\begin{equation} \label{eq:c_2}
\tilde c_2 (g_1, \dots, g_k): E_{\mathsf t (g_1)} \to TE|_{s_c (g_1, \dots, g_k)}, \quad e \mapsto s_{\tilde c} ((g_1, \dots, g_k); (g_1 \dots g_k)^{-1} \cdot e).
\end{equation}
It is easy to see that $\tilde c_2 (g_1, \dots, g_k)$ splits the projection $ TE|_{s_c (g_1, \dots, g_k)} \to E_{\mathsf t (g_1)}$. Hence, the pair $(s_c (g_1, \dots, g_k), \tilde c_2 (g_1, \dots, g_k))$ corresponds to a derivation in $D_{\mathsf t (g_1)} E$ (via the inverse to the correspondence $\delta \mapsto (\sigma_\delta, \ell_\delta)$). We denote by $\tilde c_2 (g_1, \dots, g_k)$ again the latter derivation. In this way, we have defined a map $\tilde c_2: \mathcal G^{(k)} \to DE$ such that
\begin{enumerate}
		\item[(TC1)] $\tilde c_2 (g_1, \dots, g_k) \in D_{\mathsf t (g_1)} E$ for every $(g_1, \dots, g_k) \in \mathcal G^{(k)}$; 
		\item[(TC2)] $\sigma \circ \tilde c_2 = s_{\tilde c_1}$.
\end{enumerate}
Conversely, from a pair $(\tilde c_1, \tilde c_2)$ with $\tilde c_1 \in C^k_{\mathrm{def}}(\mathcal G)$ and $\tilde c_2 : \mathcal G^{(k)} \to DE$ satisfying (TC1) and (TC2) above, we can reconstruct a linear deformation cochain $\tilde c \in C^k_{\mathrm{def}, \mathrm{lin}} (E_{\mathcal G})$ in the obvious way. From Proposition \ref{prop:proj}, it follows that
\begin{equation}\label{eq:delta_c_1}
(\delta \tilde c)_1 = \delta (\tilde c_1)
\end{equation}
Now we will prove that, for every $\tilde c \in C^k_{\mathrm{def}, \mathrm{lin}} (E_{\mathcal G})$, $(\delta \tilde c)_2$ is given by the following formula
\begin{equation} \label{eq:delta_c_2}
\begin{aligned}
& (\delta \tilde c)_2 (g_1, \dots, g_{k+1}) (e) \\
& =  - \tilde c_1(g_1, \dots, g_{k+1}) \cdot_T \left(\tilde c_2 (g_2, \dots, g_{k+1}) (g_1^{-1} e)\right) \\
& \quad + \sum_{i=1}^k (-1)^{i-1} \tilde c_2 (g_1, \dots, g_i g_{i+1}, \dots, g_{k+1})(e) + (-1)^k \tilde c_2 (g_1, \dots, g_k)(e).
\end{aligned}
\end{equation}
Indeed, let $(g_1, \dots, g_{k+1}) \in \mathcal G^{(k+1)}$, $e \in E_{\mathsf t (g_1)}$. Then, applying Formula (\ref{eq: s_delta}), one obtains:
\begin{alignat*}{3}
(\delta \tilde c)_2 (g_1, \dots, g_{k+1}) (e) & = && \ s_{\delta \tilde c} ((g_1, \dots, g_{k+1}); (g_1 \dots g_{k+1})^{-1} e) \\
& = && -T \tilde{\mathsf t} (\tilde c((g_1, \dots, g_{k+1}); (g_1 \dots g_{k+1})^{-1} e)) \\
& && + \sum_{i=1}^k (-1)^{i-1} s_{\tilde c} ((g_1, \dots, g_i g_{i+1}, \dots, g_{k+1}); (g_1 \dots g_{k+1})^{-1} e) \\
& && + (-1)^k s_{\tilde c} ((g_1, \dots, g_k); (g_1 \dots g_k)^{-1} e) \\
& = && -T \tilde{\mathsf t} (\tilde c_1(g_1, \dots, g_{k+1}), \tilde c_2 (g_2, \dots, g_{k+1}) (g_1^{-1} e)) \\
& && + \sum_{i=1}^k (-1)^{i-1} \tilde c_2 (g_1, \dots, g_i g_{i+1}, \dots, g_{k+1})(e) \\
& && + (-1)^k \tilde c_2 (g_1, \dots, g_k)(e) \\
& = && - \tilde c_1(g_1, \dots, g_{k+1}) \cdot_T \left(\tilde c_2 (g_2, \dots, g_{k+1}) (g_1^{-1} e)\right) \\
& && + \sum_{i=1}^k (-1)^{i-1} \tilde c_2 (g_1, \dots, g_i g_{i+1}, \dots, g_{k+1})(e) \\
& && + (-1)^k \tilde c_2 (g_1, \dots, g_k)(e).
\end{alignat*}

The above discussion proves the following

\begin{lemma} \label{lem:trivial_core} Let $(\mathcal G \ltimes E \rightrightarrows E; \mathcal G \rightrightarrows M)$ be a trivial-core VB-groupoid and let $k \geq 0$. The assignment $\tilde c \mapsto (\tilde c_1, \tilde c_2)$ establishes an isomorphism between $C^k_{\mathrm{def,lin}}(E_{\mathcal G})$ and the space of pairs $(\tilde c_1, \tilde c_2)$ with $\tilde c_1 \in C_{\mathrm{def}}^k(\mathcal G)$, and $\tilde c_2 : \mathcal G^{(k)} \to DE$ satisfying \emph{(TC1)} and \emph{(TC2)} above. Under this isomorphism the differential $\delta \tilde c$ corresponds to the pair $((\delta \tilde c)_1, (\delta \tilde c)_2)$ given by Formulas (\ref{eq:delta_c_1}) and (\ref{eq:delta_c_2}).
\end{lemma}

Now, take $\tilde c \in C^k_{\mathrm{def,lin}}(E_{\mathcal G})$. From Formula (\ref{eq:DE_N}), we get
\[
DE_{\mathcal G} \cong T \mathcal G \mathbin{{}_{T \mathsf s} \times_{\sigma}} DE,
\]
so we can define a map $\widehat{\tilde c} : \mathcal G^{(k+1)} \to DE_{\mathcal G}$ by putting
\[
\widehat{\tilde c} (g_0, \dots, g_k) := \left(\tilde c_1 (g_0, \dots, g_k), \tilde c_2 (g_1, \dots, g_k) \right).
\]
Observe that 
\begin{enumerate}
\item[(TC3)] $\widehat{\tilde c} (g_0, \dots, g_k) \in D_{g_0} E_{\mathcal G}$, for all $(g_0, \dots, g_k) \in \mathcal G^{(k+1)}$;
\item[(TC4)] there exists a (necessarily unique) smooth map $\mathcal G^{(k)} \to DE$ making the following diagram commutative:
\[
\begin{array}{c}
\xymatrix{\mathcal G^{(k+1)} \ar[r]^{\widehat{\tilde c}} \ar[d] & DE_{\mathcal G} \ar[d]^{D \tilde {\mathsf s}} \\
\mathcal G^{(k)} \ar[r] & DE}
\end{array}
,
\]
where the map on the left is the projection onto the last $k$ arrows.
\end{enumerate}
It is clear that $\tilde c_1 (g_0, \dots, g_k)$ is the symbol of $\widehat{\tilde c}(g_0, \dots, g_k)$. Conversely, given a map $\widehat{\tilde c} : \mathcal G^{(k+1)} \to DE_{\mathcal G}$ satisfying (TC3) and (TC4) we can reconstruct $\tilde c_1, \tilde c_2$ and hence $\tilde c$. 

Now we will describe the map $\widehat{\delta \tilde c} : \mathcal G^{(k+2)} \to DE_{\mathcal G}$. First of all, notice that there is a canonical isomorphism
\[
E_{\mathcal G}^{[2]} \cong \mathcal G^{[2]} \mathbin{{}_{\mathsf s} \times_p} E, \quad ((g,e),(h,e)) \mapsto (g,h,e),
\]
where, as in Subsection \ref{sec:lie_grpd}, we denote $\mathsf s: \mathcal G^{[2]} \to M$ the projection to the common source. Under this isomorphism, the map
\[
\begin{array}{c}
\xymatrix{E_{\mathcal G}^{[2]} \ar[r]^{\tilde{\bar{\mathsf m}}} \ar[d]_\cong & E_{\mathcal G} \ar[d]^\cong \\
\mathcal G^{[2]} \mathbin{{}_{\mathsf s} \times_p} E \ar[r] & \mathcal G \mathbin{{}_{\mathsf s} \times_p} E}
\end{array}
\]
translates into 
\[
(g, h, e) \mapsto (gh^{-1}, he).
\]
Therefore, it is a regular vector bundle map, so there is an induced map $D \tilde{\bar{\mathsf m}}: D E_{\mathcal G}^{[2]} \to D E_{\mathcal G}$. Notice that
\begin{equation} \label{eq:DE_G_2}
DE_{\mathcal G}^{[2]} \cong DE_{\mathcal G} \mathbin{{}_{D\tilde{\mathsf s}} \times_{D\tilde{\mathsf s}}}DE_{\mathcal G}
\end{equation}
and $DE_{\mathcal G} \mathbin{{}_{D\tilde{\mathsf s}} \times_{D\tilde{\mathsf s}}}DE_{\mathcal G} \cong (T \mathcal G \mathbin{{}_{T \mathsf s} \times_{\sigma}} DE) \mathbin{{}_{D\tilde{\mathsf s}} \times_{D\tilde{\mathsf s}}} (T \mathcal G \mathbin{{}_{T \mathsf s} \times_{\sigma}} DE)$. On the other hand, $\tilde {\mathsf s}: E_{\mathcal G} \cong \mathcal G \mathbin{{}_{\mathsf s} \times_p} E \to E$ is just the second projection, so $D \tilde{\mathsf s}: DE_{\mathcal G} \cong T \mathcal G \mathbin{{}_{T \mathsf s} \times_{\sigma}} DE \to DE$ is also the second projection. We conclude that
\begin{equation} \label{eq:DE_G_2_bis}
DE_{\mathcal G}^{[2]} \cong T \mathcal G^{[2]} \mathbin{{}_{T \mathsf s} \times_{\sigma}} DE.
\end{equation}

The next step will be finding an explicit formula for the map $D \tilde{\bar{\mathsf m}}$. To do this, take a derivation $\mathsf d \in D_{(g_1,g_2)} E_{\mathcal G}^{[2]}$. In particular $(g_1,g_2) \in \mathcal G^{[2]}$, so $\mathsf s (g_1) = \mathsf s (g_2)$. According to Equation (\ref{eq:DE_G_2_bis}), $\mathsf d$ can be seen as a pair $(\sigma_{\mathsf d}, \check {\mathsf d}) \in T_{(g_1,g_2)} \mathcal G^{[2]} \times D_{\mathsf s (g_1)} E$. In its turn, $\sigma_{\mathsf d}$ can be seen as a pair $((\sigma_{\mathsf d})_1, (\sigma_{\mathsf d})_2) \in T_{g_1} \mathcal G \times_{T_{\mathsf s (g_1)} M} T_{g_2} \mathcal G$.

Moreover, we know that $D \tilde{\bar{\mathsf m}}(\mathsf d)$ can be identified with the pair $(\sigma_{D \tilde{\bar{\mathsf m}}(\mathsf d)}, \ell_{D \tilde{\bar{\mathsf m}}(\mathsf d)})$, where $\sigma_{D \tilde{\bar{\mathsf m}}(\mathsf d)} = T \bar{\mathsf m} (\sigma_\mathsf d)$ and $\ell_{D \tilde{\bar{\mathsf m}}(\mathsf d)} : E_{\mathcal G}|_{g_1 g_2^{-1}} \to TE_{\mathcal G}|_{T \bar{\mathsf m}(\sigma_\mathsf d)}$ is the linear map defined in (\ref{eq:ell_delta}).
From Formula (\ref{eq: ell_D phi}) we get:
\[
\ell_{D \tilde{\bar{\mathsf m}}(\mathsf d)} = T \tilde{\bar{\mathsf m}} \circ \ell_\mathsf d \circ \tilde{\bar{\mathsf m}}_{(g_1,g_2)}^{-1}.
\]
Unravelling the definitions, one can see that there is a commutative diagram
\[
\begin{array}{c}
\xymatrix{E_{\mathcal G}|_{g_1 g_2^{-1}} \ar[r]^{\tilde{\bar{\mathsf m}}_{(g_1,g_2)}^{-1} \hspace{0.3 cm} } \ar[d]^\cong & E_{\mathcal G}^{[2]}|_{(g_1,g_2)} \ar[d]^{\cong} \ar[r]^{\ell_{\mathsf d}} & TE_{\mathcal G}^{[2]}|_{\sigma_{\mathsf d}} \ar[r]^{T \tilde{\bar{\mathsf m}}} \ar[d]^\cong & TE_{\mathcal G}|_{T \bar{\mathsf m}(\sigma_{\mathsf d})} \ar[d]^\cong \\
E_{\mathsf t (g_2)} \ar@{=}[d] \ar[r]^{h^{-1} \cdot} & E_{\mathsf s (g_2)} \ar@{=}[d] \ar[r]^{\ell_{\check {\mathsf d}}} & TE|_{T \mathsf s (\sigma_{\mathsf d})} \ar@{=}[d] \ar[r] & TE|_{T(\mathsf s \circ \bar{\mathsf m})(\sigma_{\mathsf d})} \ar@{=}[d] \\
E_{\mathsf t (g_2)} \ar[r]^{g_2^{-1} \cdot} & E_{\mathsf s (g_2)} \ar[r]^{\ell_{\check {\mathsf d}}} & TE|_{T \mathsf s ((\sigma_{\mathsf d})_2)} \ar[r]^{(\sigma_{\mathsf d})_2 *} & TE|_{T \mathsf t ((\sigma_{\mathsf d})_2)}}
\end{array}
\]
Hence
\begin{equation} \label{eq:Dm}
\ell_{D \tilde{\bar{\mathsf m}}(\mathsf d)}(e) = (\sigma_{\mathsf d})_2 * \ell_{\check {\mathsf d}} (g_2^{-1} e).
\end{equation}

We are finally ready to prove that $\widehat{\delta \tilde c} : \mathcal G^{(k+2)} \to DE_{\mathcal G}$ is given by
\begin{equation}\label{eq:delta_c_hat}
\begin{aligned}
\widehat{\delta \tilde c}  (g_0, \dots, g_{k+1}) = & -D \tilde{\bar{\mathsf m}} \left(\widehat{\tilde c}(g_0 g_1, \dots, g_{k+1}), \widehat{\tilde c}(g_1, \dots, g_{k+1})\right) \\ 
& + \sum_{i=1}^k (-1)^{i-1} \widehat{\tilde c}(g_0, \dots, g_i g_{i+1}, \dots, g_{k+1}) + (-1)^k \widehat{\tilde c}(g_0, \dots, g_k),
\end{aligned}
\end{equation}
where we understand the isomorphism (\ref{eq:DE_G_2}). Indeed, it is clear that
\begin{equation} \label{eq:sigma_D}
\begin{aligned}
\sigma_{D \bar{\mathsf m} (\widehat{\tilde c} (g_0 g_1, \dots, g_{k+1}), \widehat{\tilde c} (g_1, \dots, g_{k+1}))} & = T \bar{\mathsf m} (\sigma_{(\widehat{\tilde c} (g_0 g_1, \dots, g_{k+1}), \widehat{\tilde c} (g_1, \dots, g_{k+1}))}) \\
& = T \bar{\mathsf m} (\tilde c_1 (g_0 g_1, \dots, g_{k+1}), \tilde c_1 (g_1, \dots, g_{k+1})).
\end{aligned}
\end{equation}
Moreover, from Formula (\ref{eq:Dm}) it follows that, for every $(g_0, \dots, g_{k+1}) \in \mathcal G^{(k+2)}$, 
\begin{equation} \label{eq:ell_Dm}
\ell_{D \tilde{\bar{\mathsf m}} \left(\widehat{\tilde c}(g_0 g_1, \dots, g_{k+1}), \widehat{\tilde c}(g_1, \dots, g_{k+1})\right)}(e) = \tilde c_1 (g_1, \dots, g_{k+1}) \cdot_T \tilde c_2 (g_2, \dots, g_{k+1})(g_1^{-1} e).
\end{equation}
From (\ref{eq:delta_c_1}) and (\ref{eq:sigma_D}) it follows that, taking the symbol of the right hand side in (\ref{eq:delta_c_hat}), we get $(\delta \tilde c)_1$. Analogously, from (\ref{eq:delta_c_2}) and (\ref{eq:ell_Dm}) it follows that (the linear map induced by) the right hand side of (\ref{eq:delta_c_hat}) is exactly $(\delta \tilde c)_2$. Summarizing, we have the following

\begin{corollary} \label{cor:trivial_core}
Let $(\mathcal G \ltimes E \rightrightarrows E; \mathcal G \rightrightarrows M)$ be a trivial-core VB-groupoid and let $k \geq 0$. The assignment $\tilde c \mapsto \widehat{\tilde c}$ establishes an isomorphism between $C^k_{\mathrm{def,lin}}(E_{\mathcal G})$ and the space of maps $\widehat{\tilde c} : \mathcal G^{(k+1)} \to DE_{\mathcal G}$ satisfying \emph{(TC3)} and \emph{(TC4)}. Under this isomorphism, the differential $\delta \tilde c$ corresponds to the map $\widehat{\delta \tilde c} : \mathcal G^{(k+2)} \to DE_{\mathcal G}$ given by Formula (\ref{eq:delta_c_hat}).
\end{corollary}
Notice the analogy between (\ref{eq:delta_c_hat}) and (\ref{eq:diff_grpd}).

Finally, in this case, the projection (\ref{eq:VB_pr}) is the map $\tilde c \mapsto \tilde c_1$. From the condition (TC2), a cochain in the kernel of this projection is equivalent to a map $\tilde c_2: \mathcal G^{(k)} \to DE$ such that $\sigma \circ \tilde c_2 = 0$, so $\tilde c_2$ takes values in $\operatorname{End} E$. It follows that the sequence (\ref{eq:ses_VB}) takes the form:
\begin{equation} \label{eq:coch_2}
0 \longrightarrow C(\mathcal G, \operatorname{End} E) \longrightarrow C_{\mathrm{def,lin}}(\mathcal G \ltimes E) \longrightarrow C_{\mathrm{def}}(\mathcal G) \to 0,
\end{equation}
where $C(\mathcal G, \operatorname{End} E)$ is the Lie groupoid complex of $\mathcal G$ with coefficients in the representation of $\mathcal G$ on $\operatorname{End} E$ induced by that on $E$. Therefore, there is a long exact sequence in cohomology:
\begin{equation} \label{eq:les_trivial_core}
\dots \longrightarrow H^k (\mathcal G, \operatorname{End} E) \longrightarrow H^k_{\mathrm{def,lin}}(\mathcal G \ltimes E) \longrightarrow H^k_{\mathrm{def}}(\mathcal G) \longrightarrow H^{k+1} (\mathcal G, \operatorname{End} E) \longrightarrow \dots .
\end{equation}

We conclude this paragraph proving a rigidity result. Recall that a Lie groupoid $\mathcal G \rightrightarrows M$ is \emph{proper} if the map $(\mathsf s, \mathsf t): \mathcal G \to M \times M$ is proper.

\begin{proposition} \label{prop:rigidity} Let $(\mathcal G \ltimes E \rightrightarrows E; \mathcal G \rightrightarrows M)$ be a trivial-core VB-groupoid, and suppose that $\mathcal G$ is proper. Then $H^k_{\mathrm{def,lin}}(\mathcal G \ltimes E) = 0$ for $k \geq 1$. \end{proposition}

\proof Consider the long exact sequence (\ref{eq:les_trivial_core}). In degrees $k \geq 1$, both the Lie groupoid cohomology with coefficients \cite{crainic:coh} and the deformation cohomology \cite{crainic:def2} of a proper Lie groupoid vanish. It follows that $H^k_{\mathrm{def,lin}}(\mathcal G \ltimes E) = 0$ for $k \geq 1$. \endproof

At the end of Subsection \ref{sec:lin_def_1}, we will prove an analogous result for full-core VB-groupoids. (Corollary \ref{cor:full-core_rigidity}).

\subsubsection{Low-degree cohomology groups}

Next we describe low-degree cohomology groups. The entire discussion of Subsection \ref{sec:def_algd} goes through, with minor changes.

Let $(\mathcal W \rightrightarrows E; \mathcal G \rightrightarrows M)$ be a VB-groupoid, and let $(W \Rightarrow E; A \Rightarrow M)$ be its VB-algebroid. Consider the isotropy bundle $\mathfrak i$ of $\mathcal W$ and its sections. It is natural to define \emph{linear sections} of $\mathfrak i$ by
\[
\Gamma_{\mathrm{lin}}(\mathfrak i) := \Gamma(\mathfrak i) \cap \Gamma_{\mathrm{lin}}(W,E)
\]
and \emph{linear invariant sections} by
\[
H^0_{\mathrm{lin}}(\mathcal W, \mathfrak i) = \Gamma_{\mathrm{lin}}(\mathfrak i)^{\mathrm{inv}} := \Gamma(\mathfrak i)^{\mathrm{inv}} \cap \Gamma_{\mathrm{lin}}(W,E).
\]
The following proposition is easily proved as in \cite[Proposition 4.1]{crainic:def2}.

\begin{proposition}\label{prop:H^-1_lin} Let $(\mathcal W \rightrightarrows E; \mathcal G \rightrightarrows M)$ be a VB-groupoid. Then $H^{-1}_{\mathrm{def,lin}}(\mathcal W) \cong H^0_{\mathrm{lin}}(\mathcal W, \mathfrak i) = \Gamma_{\mathrm{lin}}(\mathfrak i)^{\mathrm{inv}}$. \end{proposition}

Now we consider the complex $C (\mathcal W, \mathfrak i)$. It is clear that there is a distinguished subcomplex $C_{\mathrm{lin}}(\mathcal W, \mathfrak i)$, the one of cochains that are vector bundle maps $\mathcal W^{(k)} \to W$ covering some map $\mathcal G^{(k)} \to A$. From Formula (\ref{eq:c_u}) it follows that the map (\ref{eq:isotropy}) takes $C_{\mathrm{lin}}(\mathcal W, \mathfrak i)$ to $C_{\mathrm{def,lin}}(\mathcal W)$, so we have a map
\[
r: C_{\mathrm{lin}}(\mathcal W, \mathfrak i) \hookrightarrow C_{\mathrm{def,lin}}(\mathcal W).
\]

\begin{proposition}[see also {\cite[Section 2.3]{etv:infinitesimal}}]
\[
H^0_{\mathrm{def,lin}}(\mathcal W) = \dfrac{\text{linear multiplicative vector fields on $\mathcal W$}}{\text{inner linear multiplicative vector fields on $\mathcal W$}}.
\]
\end{proposition}

\proof Let $X \in C^0_{\mathrm{def,lin}}(\mathcal W)$. Then $X$ is a linear vector field, and from Proposition \ref{prop:H^0} it follows that $X$ is closed if and only if it is multiplicative, and it is exact if and only if is inner multiplicative, as desired. \endproof

Recall that also the normal bundle $\nu$ of $\mathcal W$ is defined. Observe that the anchor $\rho: W \to TE$ of the Lie algebroid $W$ is a morphism of DVBs, hence it takes linear sections to linear vector fields. We set $\rho_{\mathrm{lin}}: \Gamma_{\mathrm{lin}}(W,E) \to \mathfrak X_{\mathrm{lin}}(E)$ and
\[
\Gamma_{\mathrm{lin}}(\nu) := \mathfrak X_{\mathrm{lin}}(E)/ \operatorname{im} \rho_{\mathrm{lin}}.
\]
Following the discussion in Subsection \ref{sec:def_algd}, we declare that a section $V (\operatorname{mod} \operatorname{im} \rho_{\mathrm{lin}}) \in \Gamma_{\mathrm{lin}}(\mathcal W)$ is \emph{invariant} if it possesses an $(\mathsf s, \mathsf t)$-lift $X \in \mathfrak X_{\mathrm{lin}}(\mathcal W)$. The space of invariant linear sections is denoted $H^0_{\mathrm{lin}} (\mathcal W, \nu)$ or $\Gamma_{\mathrm{lin}}(\nu)^{\mathrm{inv}}$. Observing that the projection on $E$ of a linear multiplicative vector field is linear, we obtain a linear map
\[
\pi: H^0_{\mathrm{def,lin}}(\mathcal W) \to \Gamma_{\mathrm{lin}}(\nu)^{\mathrm{inv}}.
\]

From Lemma \ref{prop:curvature}, Proposition \ref{prop:exact} and their proofs in \cite{crainic:def2}, the ``linear versions'' follow immediately.

\begin{lemma} Let $V (\operatorname{mod} \operatorname{im} \rho_{\mathrm{lin}}) \in \Gamma_{\mathrm{lin}}(\nu)^{\mathrm{inv}}$ and $X \in \mathfrak X_{\mathrm{lin}}(\mathcal W)$ an $(\mathsf s, \mathsf t)$-lift of $\mathcal W$. Then $\delta X \in C^2_{\mathrm{lin}}(\mathcal W, \mathfrak i)$ and its cohomology class does not depend on the choice of $X$, hence there is an induced linear map
\[
K: \Gamma_{\mathrm{lin}}(\nu)^{\mathrm{inv}} \to H^2_{\mathrm{lin}}(\mathcal W, \mathfrak i).
\]
\end{lemma}

\begin{proposition} There is an exact sequence
\[
0 \longrightarrow H^1_{\mathrm{lin}}(\mathcal W, \mathfrak i) \overset{r}{\longrightarrow} H^0_{\mathrm{def,lin}}(\mathcal W) \overset{\pi}{\longrightarrow} \Gamma_{\mathrm{lin}}(\nu)^{\mathrm{inv}} \overset{K}{\longrightarrow} H^2_{\mathrm{lin}}(\mathcal W, \mathfrak i) \overset{r}{\longrightarrow} H^1_{\mathrm{def,lin}}(\mathcal W).
\]
\end{proposition}

\subsubsection{Deformations}

Here we discuss deformations of a VB-groupoid $(\mathcal W \rightrightarrows E; \mathcal G \rightrightarrows M)$ and their relation with the $1$-cohomology group $H^1_{\mathrm{def,lin}}(\mathcal W)$. Let $B$ be a smooth manifold.

\begin{definition} A \emph{family of VB-groupoids} over $B$ is a diagram
\[
\begin{array}{r}
\xymatrix{
\tilde{\mathcal W} \ar[d] \ar@<0.4ex>[r] \ar@<-0.4ex>[r] & \tilde E \ar[d] \ar[r] & B \ar@{=}[d] \\
\tilde{\mathcal{G}} \ar@<0.4ex>[r] \ar@<-0.4ex>[r] & \tilde M \ar[r] & B }
\end{array}
\]
such that the first square is a VB-groupoid and the rows are families of Lie groupoids. In particular, for every $b \in B$, 
\[
\begin{array}{r}
\xymatrix{
\mathcal W_b \ar[d] \ar@<0.4ex>[r] \ar@<-0.4ex>[r] & E_{b} \ar[d] \\
\mathcal{G}_{b} \ar@<0.4ex>[r] \ar@<-0.4ex>[r] & M_{b}}
\end{array}
\]
is a VB-groupoid.

If $B$ is an open interval $I$ containing 0, the family is said to be a \emph{deformation} of $\mathcal W_0$ and the latter is denoted simply by $(\mathcal W \rightrightarrows E; \mathcal G \rightrightarrows M)$. A deformation of $\mathcal W$ is also denoted $(\mathcal W_\epsilon)$, where $\epsilon$ denotes the canonical coordinate on $I$.
\end{definition}

The structure maps of $\mathcal W_\epsilon$ are denoted $\tilde {\mathsf s}_\epsilon, \tilde {\mathsf t}_\epsilon, \tilde {\mathsf 1}_\epsilon, \tilde {\mathsf m}_\epsilon, \tilde {\mathsf i}_\epsilon$, the division map is denoted $\tilde{\bar {\mathsf m}}_\epsilon$. Strict, $\mathsf s$-constant, $\mathsf t$-constant, $(\mathsf s ,\mathsf t)$-constant and constant deformations are defined as in the groupoid case. Two deformations $(\mathcal W_\epsilon)$ and $(\mathcal W'_\epsilon)$ of $\mathcal W$ are called \emph{equivalent} if there exists a smooth family of VB-groupoid isomorphisms $\Psi_\epsilon: \mathcal W_\epsilon \to \mathcal W'_\epsilon$ such that $\Psi_0 = \mathrm{id}$. We say that $(\mathcal W_\epsilon)$ is \emph{trivial} if it is equivalent to the constant deformation.

We proceed as in Subsection \ref{sec:def_grpd}: first we study $(\mathsf s, \mathsf t)$-constant deformations, then $\mathsf s$-constant deformations and finally general deformations.

\begin{proposition} [$(\mathsf s, \mathsf t)$-constant deformations] Let $(\mathcal W_\epsilon)$ be an $(\mathsf s, \mathsf t)$-constant deformation of the VB-groupoid $\mathcal W$. Then formula
\[
u_0 (w, w') := - \dfrac{d}{d \epsilon} \bigg|_{\epsilon = 0} R^{-1}_{w w'}(\tilde {\mathsf m}_\epsilon(w,w'))
\]
defines a cocycle $u_0 \in C^1_{\mathrm{lin}}(\mathcal W, \mathfrak i)$ (here $R$ denotes right translation in $\mathcal W$). Its image $\xi_0$ in $C^1_{\mathrm{def,lin}}(\mathcal W)$ is
\[
\xi_0(w,w') = \dfrac{d}{d \epsilon} \bigg|_{\epsilon = 0} \tilde{\bar {\mathsf m}}_\epsilon (w w',w').
\]
\end{proposition}

\proof We only need to show that $u_0$ is linear. Remembering that $\tilde{\mathsf m}_\epsilon$ is linear for every $\epsilon$, we compute
\[
\begin{aligned}
u_0( \lambda (w,w')) & = - \dfrac{d}{d \epsilon} \bigg|_{\epsilon = 0} R_{\lambda(ww')}^{-1}(\tilde {\mathsf m}_\epsilon (\lambda (w,w')) \\
& = - \dfrac{d}{d\epsilon} \bigg|_{\epsilon = 0} (\lambda \tilde{\mathsf m}_\epsilon (w,w')) \cdot (\lambda(ww'))^{-1} \\
& = - \dfrac{d}{d\epsilon} \bigg|_{\epsilon = 0} (\lambda \tilde{\mathsf m}_\epsilon (w,w')) \cdot (\lambda (ww')^{-1}) \\
& = - \dfrac{d}{d\epsilon} \bigg|_{\epsilon = 0} \lambda (\tilde{\mathsf m}_\epsilon (w,w') \cdot (ww')^{-1}) \\
& = \lambda \cdot_W \bigg( -\dfrac{d}{d\epsilon} \bigg|_{\epsilon = 0} R^{-1}_{ww'}(\tilde{\mathsf m}_\epsilon (w,w')) \bigg) \\
& = \lambda \cdot_W u_0 (w,w').
\end{aligned}
\]
\endproof

We now pass to $\mathsf s$-constant deformations.

\begin{proposition} [$\mathsf s$-constant deformations] Let $(\mathcal W_\epsilon)$ be an $\mathsf s$-constant deformation of the VB-groupoid $\mathcal W$. Then the formula
\[
\xi_0(w,w') = \dfrac{d}{d \epsilon} \bigg|_{\epsilon = 0} \tilde{\bar {\mathsf m}}_\epsilon (ww',w').
\]
defines a cocycle $\xi_0$ in $C^1_{\mathrm{def,lin}}(\mathcal W)$ and its cohomology class does only depend on the equivalence class of the deformation.
\end{proposition}

\proof We know that $\xi_0$ is a linear cocycle. Moreover, if $\Psi_\epsilon: \mathcal W_\epsilon \to \mathcal W'_\epsilon$ is an equivalence of deformations, $X = \frac{d \Psi_\epsilon}{d \epsilon}|_{\epsilon = 0} \in C^0_{\mathrm{def,lin}}(\mathcal W)$ and we can proceed as in the proof of \cite[Lemma 5.3]{crainic:def2}. Indeed, given $(w,w') \in \mathcal W^{(2)}$, one has
\[
\bar{\mathsf m}'_\epsilon (\Psi_\epsilon (ww'), \Psi_\epsilon (w')) = \Psi_\epsilon (\bar{\mathsf m}_\epsilon (ww', w'))
\]
for every $\epsilon$. Differentiating at $\epsilon = 0$, we get
\[
\dfrac{d}{d \epsilon} \bigg|_{\epsilon = 0} \tilde{\bar{\mathsf m}}'_\epsilon(ww', w') + T \tilde{\bar{\mathsf m}} (X(ww'), X(w')) = X(w) + \dfrac{d}{d \epsilon} \bigg|_{\epsilon = 0} \tilde{\bar{\mathsf m}}_\epsilon (ww',w'),
\]
so $(\xi_0 - \xi'_0)(w,w') = \delta X (w,w')$ as desired. \endproof

In order to discuss the general case, we need to prove a linear version of Proposition \ref{prop:gen_def}. Recall from \cite{abad:ruth2} that an \emph{Ehresmann connection} on a Lie groupoid $\mathcal G \rightrightarrows M$ is a splitting of the short exact sequence (\ref{eq:ses1}) that restricts to the canonical splitting (\ref{eq:ses2}) over $M$. Notice that an Ehresmann connection on $\mathcal G$ is exactly the same as a right-horizontal lift of the VB-groupoid $T\mathcal G$. In particular, such a connection always exists.

Now let $(\mathcal W \rightrightarrows E, \mathcal G \rightrightarrows M)$ be a VB-groupoid. Then there is a diagram of morphisms of DVBs
\[
\begin{array}{c}
\xymatrix@C=9pt@R=12pt{& 0 \ar[rr] & & T^{\tilde {\mathsf s}} \mathcal W \ar[dl] \ar[rr] \ar[dd]|!{[dl];[dr]}{\hole} & & T \mathcal W \ar[dl] \ar[dd]|!{[dl];[dr]}{\hole} \ar[rr]^(.6){T \tilde {\mathsf s}} & & \tilde {\mathsf s}^* TE \ar[dl] \ar[dd]|!{[dl];[dr]}{\hole} \ar[rr] & & 0 \\
0 \ar[rr] & & T^{\mathsf s} \mathcal G \ar[rr] \ar[dd] & & T \mathcal G \ar[rr]^(.6){T \mathsf s} \ar[dd] & & \mathsf s^* TM \ar[rr] \ar[dd] & & 0 \\
& & & \mathcal W \ar[dl] \ar@{=}[rr]|!{[ur];[dr]}{\hole} & & \mathcal W \ar[dl] \ar@{=}[rr]|!{[ur];[dr]}{\hole} & & \mathcal W \ar[dl] & & \\
& & \mathcal G \ar@{=}[rr] & & \mathcal G \ar@{=}[rr] & & \mathcal G & &}
\end{array}
\]
where the top rows are short exact sequences of vector bundles.

\begin{definition} A \emph{linear Ehresmann connection} on $\mathcal W$ is a morphism of DVBs 
\[
(\tilde{\mathsf s}^* TE \to \mathsf s^* TM; \mathcal W \to \mathcal G) \longrightarrow (T \mathcal W \to T \mathcal G; \mathcal W \to \mathcal G)
\]
such that the maps $\tilde{\mathsf s}^* TE \to T \mathcal W$ and $\mathsf s^* TM \to T \mathcal G$ are Ehresmann connections on $\mathcal W$ and $\mathcal G$, respectively. \end{definition}

\begin{lemma} On every VB-groupoid $(\mathcal W \rightrightarrows E; \mathcal G \rightrightarrows M)$ there exists a linear Ehresmann connection. \end{lemma}

\proof First of all, choose local coordinates on $\mathcal G$ adapted to the submersion $\mathsf s: \mathcal G \to M$. Up to a translation, one can also assume that they are adapted to the immersion $\mathsf 1: M \to \mathcal G$. Using a right-decomposition $\mathcal W \cong \mathsf s^* E \oplus \mathsf t^* C$, we find fiber coordinates on $\mathcal W$ with analogous properties. Now it is easy to see that linear Ehresmann connections exist locally, and one can conclude with a partition of unity argument. \endproof

\begin{proposition} Let 
\[
\begin{array}{r}
\xymatrix{
\tilde{\mathcal W} \ar[d] \ar@<0.4ex>[r] \ar@<-0.4ex>[r] & \tilde E \ar[d] \ar[r]^{\tilde \pi} & I \ar@{=}[d] \\
\tilde{\mathcal{G}} \ar@<0.4ex>[r] \ar@<-0.4ex>[r] & \tilde M \ar[r] & I }
\end{array}
\]

be a deformation of $\mathcal W$. Then:
\begin{enumerate}
\item there exist transverse linear vector fields for $\tilde{\mathcal W}$;
\item if $\tilde X$ is a transverse linear vector field, then $\delta \tilde{X}$, when restricted to $\mathcal W$, induces a cocycle $\xi_0 \in C^1_{\mathrm{def,lin}}(\mathcal W)$;
\item the cohomology class of $\xi_0$ does not depend on the choice of $\tilde X$.
\end{enumerate}
\end{proposition}

\proof 
\
\begin{enumerate}
\item Take a vector field $Y$ on $\tilde M$ that projects on $\frac{d}{d\epsilon}$. Choosing a linear connection on $\tilde E \to \tilde M$, one can lift it to a linear vector field $\tilde Y$ on $\tilde E$, that obviously projects on $\frac{d}{d\epsilon}$. Now the choice of a linear Ehresmann connection on $\mathcal W$ gives a transverse linear vector field $\tilde X$ on $\mathcal W$ that projects on $\tilde Y$, as desired.
\item Let $\tilde X$ be a transverse linear vector field. Then $\delta \tilde X \in C^1_{\mathrm{def,lin}}(\tilde{\mathcal W})$ and, by Proposition \ref{prop:gen_def}, it restricts to $\mathcal W$, so it belongs to $C^1_{\mathrm{def,lin}}(\mathcal W)$. 
\item This can be proved as in \cite[Proposition 5.12]{crainic:def2}. Let $\tilde X'$ be another transverse linear vector field for $\tilde{\mathcal W}$. Then $\tilde Y = \tilde X - \tilde X'$ is killed by $T(\tilde \pi \circ \tilde{\mathsf s})$, so $\tilde Y|_{\mathcal W}$ is tangent to $\mathcal W$. It follows that $\tilde Y |_{\mathcal W} \in C^0_{\mathrm{def,lin}}(\mathcal W)$ and $\xi_0 - \xi'_0 = \delta \tilde Y |_{\mathcal W}$, as desired. \endproof
\end{enumerate}

The cohomology class $[\xi_0] \in H^1_{\mathrm{def,lin}}(\mathcal W)$ is called the \emph{linear deformation class} associated to the deformation $\tilde {\mathcal W}$. From the last proposition, it follows directly that this class is also independent of the equivalence class of the deformation. 

In particular, an infinitesimal deformation is trivial if and only if its cohomology class is zero. From Proposition \ref{prop:rigidity} it follows that \emph{every infinitesimal deformation of a trivial-core VB-groupoid with proper base is trivial}.

\begin{remark}\label{rem:ind_def_base}
Clearly, a deformation $\tilde{\mathcal W}$ of the VB-groupoid $\mathcal W$ induces a deformation $\tilde{\mathcal G}$ of the base groupoid $\mathcal G$. Now, it is easy to check that the projection (\ref{eq:VB_pr}) sends the linear deformation class of $\tilde{\mathcal W}$ to the deformation class of $\tilde{\mathcal G}$. Indeed, let $\tilde X$ be a transverse linear vector field on $\tilde{\mathcal W}$ and $X$ its projection on $\tilde{\mathcal G}$. Then $X$ is transverse again. But the linear deformation class of $\tilde {\mathcal W}$ is $[\delta \tilde X|_{\mathcal W}]$ and the projection (\ref{eq:VB_pr}) sends it to $[\delta X|_{\mathcal G}]$, which is the deformation class of $\tilde{\mathcal G}$.
\end{remark}

Finally, we discuss the variation map associated to deformations of VB-groupoids. Let $(\tilde{\mathcal W} \rightrightarrows \tilde E; \tilde{\mathcal G} \rightrightarrows \tilde M)$ be a family of VB-groupoids over a smooth manifold $B$. Then any curve $\gamma: I \to B$ induces a deformation $\gamma^* \tilde{\mathcal W}$ of $\tilde{\mathcal W}_{\gamma(0)}$, and we have:

\begin{proposition}[The linear variation map] Let $b \in B$. For any curve $\gamma: I \to B$ with $\gamma(0) = b$, the deformation class of $\gamma^* \tilde{\mathcal W}$ at time 0 does only depend on $\dot{\gamma}(0)$. This defines a linear map
\[
\mathrm{Var}^{\tilde{\mathcal W}}_{\mathrm{lin}, b}: T_b B \to H^1_{\mathrm{def,lin}}(\tilde{\mathcal W}_b),
\]
called the \emph{linear variation map} of $\tilde{\mathcal W}$ at $b$, that makes the following diagram commutative:
\[
\xymatrix@C=45pt{
T_b B \ar[r]^-{\mathrm{Var}^{\tilde{\mathcal W}}_{\mathrm{lin}, b}} \ar[dr]_-{\mathrm{Var}^{\tilde{\mathcal G}}_b} & H^1_{\mathrm{def,lin}}(\tilde{\mathcal W}_b) \ar[d] \\
& H^1_{\mathrm{def}}(\tilde{\mathcal G}_b)
}
\]
\end{proposition}

\proof The first statement is proved as in \cite[Proposition 5.15]{crainic:def2}, while the second statement trivially follows from Remark \ref{rem:ind_def_base}. \endproof

\subsubsection{Deformations of the dual VB-groupoid}

We conclude this section noticing that the linear deformation cohomology of a VB-groupoid is canonically isomorphic to that of its dual. Notice that this result is weaker than the analogous one obtained for VB-algebroids (Theorem \ref{thm:dual_vb_algd}), that gives an isomorphism even at the level of complexes.

\begin{theorem} \label{thm:dual_vb_grpd} Let $(\mathcal W \rightrightarrows E; \mathcal G \rightrightarrows M)$ be a VB-groupoid. Then there is a canonical isomorphism
\begin{equation}\label{eq:dual_def}
H_{\mathrm{def,lin}}(\mathcal W) \cong H_{\mathrm{def,lin}}(\mathcal W^*)
\end{equation}
of $H(\mathcal G)$-modules.
\end{theorem}

\proof Using Proposition \ref{prop:cotangent_VB} and Lemma \ref{prop:proj,lin;lin,lin}, we get:
\[
H_{\mathrm{def,lin}}(\mathcal W)\cong H_{\mathrm{proj,lin}}(T^* \mathcal W)[1]  \cong H_{\mathrm{lin,lin}}(T^* \mathcal W)[1].
\]
For the same reason, $H_{\mathrm{def,lin}}(\mathcal W^*) \cong H_{\mathrm{lin,lin}}(T^* \mathcal W^*)[1]$. But we have already noticed that $T^* \mathcal W \cong T^* \mathcal W^*$ as double vector bundles and Lie groupoids, so we obtain (\ref{eq:dual_def}). \endproof

As the dual of a trivial-core VB-groupoid is a full-core VB-groupoid, from Proposition \ref{prop:rigidity} we immediately get:

\begin{corollary} \label{cor:full-core_rigidity} Let $(\mathcal W \rightrightarrows 0_M; \mathcal G \rightrightarrows M)$ be a full-core VB-groupoid with proper base. Then $H^k_{\mathrm{def,lin}}(\mathcal W) = 0$ for $k \geq 1$. In particular, every infinitesimal deformation of $\mathcal W$ is trivial. \end{corollary}

\subsection{The linearization map} \label{sec:linearization_1}

Let $(\mathcal W \rightrightarrows E; \mathcal G \rightrightarrows M)$ be a VB-groupoid. We have shown that deformations of the VB-groupoid structure are controlled by a subcomplex $C_{\mathrm{def}, \mathrm{lin}} (\mathcal W)$ of the deformation complex $C_{\mathrm{def}}(\mathcal W)$ of the top Lie groupoid $\mathcal W \rightrightarrows E$. In this section, we prove that there is a canonical splitting of the inclusion $C_{\mathrm{def,lin}}(\mathcal W) \hookrightarrow C_{\mathrm{def}}(\mathcal W)$ in the category of cochain complexes, the \emph{linearization map}. This will imply, in particular, that the inclusion induces an injection in cohomology $H_{\mathrm{def,lin}}(\mathcal W) \hookrightarrow H_{\mathrm{def}}(\mathcal W)$.

The procedure we are going to describe is analogous to the one we used in Subsection \ref{sec:linearization} to define the linearization of sections of a DVB. As usual, denote by $h$ the homogeneity structure of $\mathcal W \to \mathcal G$. For every $\lambda > 0$, $h_\lambda$ is a groupoid automorphism of $\mathcal W \rightrightarrows E$. If $k = -1$, by definition the action induced by $h$ on $C^{-1}_{\mathrm{def}}(\mathcal W) = \Gamma(W,E)$ coincides with that induced by the homogeneity structure of $W \to A$. 

Let now $k \geq 0$ and $\tilde c \in C^k_{\mathrm{def}}(\mathcal W)$. Then, by Equation (\ref{eq:aut_action}), we have:
\[
(h_\lambda^* \tilde c) (w_0, \dots, w_k) = Th_\lambda^{-1}(\tilde c(\lambda (w_0, \dots, w_k))).
\]
On the other hand, $C^k_{\mathrm{def}}(\mathcal W) \subset \Gamma(\tilde p_{k+1}^* T \mathcal W, \mathcal W^{(k+1)})$. For every $\lambda > 0$, let $\tilde h_\lambda$ be the multiplication by $\lambda$ in the fibers of the vector bundle $\tilde p_{k+1}^* T \mathcal W \to p_{k+1}^* T \mathcal G$. Then $\tilde h_\lambda$ acts on sections of $\tilde p_{k+1}^* T \mathcal W \to \mathcal W^{(k+1)}$ by pull-back. We have $\tilde p^*_{k+1} T \mathcal W \cong \mathcal W^{(k+1)} \times_{\mathcal W} T \mathcal W$, $p^*_{k+1} T \mathcal G \cong \mathcal G^{(k+1)} \times_{\mathcal G} T \mathcal G$ and
\[
\tilde h_\lambda = h_\lambda^{(k)} \times T h_\lambda.
\]
Therefore, if we interpret $\tilde c$ as a section of $\tilde p^*_{k+1} T \mathcal W \to \mathcal W^{(k+1)}$, we get:
\[
\begin{aligned}
(\tilde h_\lambda^* \tilde c) (w_0, \dots, w_k) & = (h_\lambda^{(k)} \times T h_\lambda)^{-1} (\lambda (w_0, \dots, w_k), \tilde c ( \lambda (w_0, \dots, w_k))) \\
& = ((w_0, \dots, w_k), Th_\lambda^{-1} (\tilde c (\lambda (w_0, \dots, w_k)))).
\end{aligned}
\]
This shows that the two actions agree on $C^k_{\mathrm{def}}(\mathcal W)$. As a consequence, by Propositions \ref{prop:c_core}, \ref{prop:c_0} the limits
\[
\begin{aligned}
\tilde c_{\mathrm{core}} := & \lim_{\lambda \to 0} \big( \lambda \cdot h_\lambda^* \tilde c \big) \\
\tilde c_{\mathrm{lin}} := & \lim_{\lambda \to 0} \big( h_\lambda^* \tilde c - \lambda^{-1} \cdot \tilde c_{\mathrm{core}} \big) \\
\end{aligned}
\]
are well-defined. The latter equation defines a linear map
\[
\mathrm{lin}: C_{\mathrm{def}}(\mathcal W) \to C_{\mathrm{def,lin}}(\mathcal W), \quad \tilde c \mapsto \tilde c_{\mathrm{lin}},
\]
which we call the \emph{linearization map}.

\begin{theorem} [Linearization of deformation cochains] \label{prop:linearization}
The linearization map is a cochain map splitting the inclusion $C_{\mathrm{def}, \mathrm{lin}}(\mathcal W) \hookrightarrow C_{\mathrm{def}}(\mathcal W)$. In particular there is a direct sum decomposition
\[
C_{\mathrm{def}}(\mathcal W) \cong C_{\mathrm{def}, \mathrm{lin}}(\mathcal W) \oplus \ker(\mathrm{lin}).
\]
of cochain complexes. Hence, the inclusion of linear deformation cochains into deformation cochains induces an injection
\[
H_{\mathrm{def}, \mathrm{lin}}(\mathcal W) \hookrightarrow H_{\mathrm{def}}(\mathcal W).
\]
\end{theorem}

\proof We only need to prove that the linearization map respects the differential. As $h_\lambda$ is an automorphism of $\mathcal W \rightrightarrows E$ for every $\lambda > 0$, we have that $h_\lambda^*$ commutes with $\delta$, as explained in Subsection \ref{sec:def_grpd}. It is also clear that $\delta$ preserves limits, so we compute
\[
\begin{aligned}
(\delta \tilde c)_{\mathrm{core}} = & \lim_{\lambda \to 0} \big( \lambda \cdot h_\lambda^* (\delta \tilde c) \big) = \lim_{\lambda \to 0} \big( \lambda \cdot \delta (h_\lambda^* \tilde c) \big) = \delta \big( \lim_{\lambda \to 0} \lambda \cdot \delta (h_\lambda^* \tilde c) \big) = \delta \tilde c_{\mathrm{core}}, \\
(\delta \tilde c)_{\mathrm{lin}} = & \lim_{\lambda \to 0} \big( h_\lambda^* (\delta \tilde c) - \lambda^{-1} (\delta \tilde c)_{\mathrm{core}} \big) = \delta \bigg( \lim_{\lambda \to 0} (h_\lambda^* \tilde c - \lambda^{-1} \delta \tilde c_{\mathrm{core}}) \bigg) = \delta \tilde c_{\mathrm{lin}}
\end{aligned}
\]
and we are done. \endproof 

\begin{remark}\label{rmk:cotangent} Applying the isomorphisms (\ref{eq:iso}) and (\ref{eq:iso_lin}), we obtain that the submodule $H_{\mathrm{proj,lin}}(T^* \mathcal W)$ is a direct summand of $H_{\mathrm{proj},\bullet}(T^* \mathcal W)$: it identifies with classes in $H_{\mathrm{proj},\bullet}(T^* \mathcal W)$ which can be represented by cochains that are linear over $\mathcal W^*$. \end{remark}

Finally, we discuss a first consequence of Theorem \ref{prop:linearization}. We call $\widehat C_{\mathrm{def,lin}}(\mathcal W) := C_{\mathrm{def,lin}}(\mathcal W) \cap \widehat C_{\mathrm{def}}(\mathcal W)$ the \emph{linear normalized deformation complex} of $\mathcal W$.

\begin{proposition} The inclusion $\widehat C_{\mathrm{def,lin}}(\mathcal W) \hookrightarrow C_{\mathrm{def,lin}}(\mathcal W)$ is a quasi-isomorphism. \end{proposition}

\proof Take $c \in C^k_{\mathrm{def,lin}}(\mathcal W)$. By Proposition \ref{prop:norm}, there exist $\widehat c \in \widehat C^k_{\mathrm{def}}(\mathcal W)$ and $c' \in C^{k-1}_{\mathrm{def}}(\mathcal W)$ such that $c - \widehat c = \delta c'$. Applying the linearization map, we get $c - \widehat c_{\mathrm{lin}} = \delta c'_{\mathrm{lin}}$ and $\widehat c_{\mathrm{lin}} \in \widehat C^k_{\mathrm{def,lin}}(\mathcal W)$, as desired. \endproof

Other applications of the linearization map will be considered in the next sections.

\subsection{The van Est map} \label{sec:van_est}

The van Est theorem is a classical result relating the differentiable cohomology of a Lie group and the Chevalley-Eilenberg cohomology of its Lie algebra \cite{est:group_coh, est:alg_coh}. It was later extended to differentiable cohomology \cite{weinstein:extensions, crainic:coh} and deformation cohomology \cite{crainic:def2} of a Lie groupoid, and to the VB-cohomology of a VB-groupoid \cite{cabrera:hom}. In this subsection, we want to prove an analogous theorem for the linear deformation cohomology of a VB-groupoid.

Let $\mathcal G \rightrightarrows M$ be a Lie groupoid, and let $A \Rightarrow M$ be its Lie algebroid. The normalized deformation complex of $\mathcal G$ and the deformation complex of $A$ are intertwined by the \emph{van Est map}, defined as follows. Given a section $\alpha \in \Gamma(A)$, we define a map $R_\alpha: \widehat C^k_{\mathrm{def}}(\mathcal G) \to \widehat C^{k-1}_{\mathrm{def}}(\mathcal G)$ by
\begin{equation}\label{eq:R_alpha_1}
R_\alpha(c) = [c, \overrightarrow{\alpha}]|_M
\end{equation}
if $k = 0$, and
\begin{equation}\label{eq:R_alpha_2}
R_\alpha(c) (g_1, \dots, g_k) = (-1)^{k} \dfrac{d}{d \epsilon} \bigg|_{\epsilon = 0} c(g_1, \dots, g_k, \Phi^\alpha_\epsilon(\mathsf s (g_k))^{-1})
\end{equation}
if $k > 0$, where $\Phi^\alpha_\epsilon (x) = \Phi^{\overrightarrow{\alpha}}_\epsilon (\mathsf 1_x)$ for every $x \in M$: the image of $\mathsf 1_x$ under the flow $(\Phi^{\overrightarrow{\alpha}}_\epsilon)$ of the right invariant vector field $\overrightarrow{\alpha}$ associated to $\alpha$. Then the van Est map 
\begin{equation}\label{eq:VE}
\mathrm{VE}: \widehat C_{\mathrm{def}}(\mathcal G) \to C_{\mathrm{def}}(A)
\end{equation}
is given by:
\begin{equation}\label{eq:VE_formula}
\mathrm{VE}(c)(\alpha_0, \dots, \alpha_k) = \sum_{\tau \in S_{k+1}} (-1)^\tau (R_{\alpha_{\tau(k)}} \circ \dots \circ R_{\alpha_{\tau(0)}}) (c).
\end{equation}


\begin{theorem}[{\cite[Theorem 10.1]{crainic:def2}}] \label{prop:van_est}  The van Est map (\ref{eq:VE}) is a cochain map. Moreover, if $\mathcal G$ has $k$-connected $\mathsf s$-fibers, it induces an isomorphism in cohomology in all degrees $p < k$. \end{theorem}

We are going to prove an analogous theorem for the linear deformation complex of a VB-groupoid. To do this, we need a simple preliminary lemma.
>F
Let $\mathcal G \rightrightarrows M$ be a Lie groupoid, and let $A \Rightarrow M$ be its Lie algebroid.  We know that the group $\mathrm{Aut}(\mathcal G)$ of automorphisms of $\mathcal G$ acts on $\widehat C_{\mathrm{def}}(\mathcal G)$. Clearly, it also acts on $C_{\mathrm{def}}(A)$, by
\[
\Psi^* c := \psi^* c.
\]
As expected, we have the following

\begin{lemma}\label{lem:VE_equiv} The van Est map (\ref{eq:VE}) is equivariant with respect to the action of $\mathrm{Aut}(\mathcal G)$. \end{lemma}

\proof We will prove that
\begin{equation}\label{eq:R_equiv}
\Psi^*(R_\alpha (c)) = R_{\psi^* \alpha}(\Psi^* c)
\end{equation}
for all $c \in \widehat C^k_{\mathrm{def}}(\mathcal G)$, $\alpha \in \Gamma(A)$, $\Psi \in \mathrm{Aut}(\mathcal G)$. If $k = 0$, we have
\[
\begin{aligned}
\Psi^*(R_\alpha(c)) & = \Psi^*([c, \overrightarrow \alpha]|_M) = (\Psi^*[c, \overrightarrow \alpha])|_M \\
& = [\Psi^* c, \Psi^* \overrightarrow \alpha]|_M = [\Psi^*c, \overrightarrow{\psi^* \alpha}]|_M = R_{\psi^* \alpha}(\Psi^* c).
\end{aligned}
\]
Now, let $k > 0$. Then
\begin{equation}\label{eq:Psi^*R}
\begin{aligned}
\Psi^*(R_\alpha(c)) (g_1, \dots, g_k) & = T \Psi^{-1} (R_\alpha (c) (\Psi(g_1), \dots, \Psi(g_k))) \\ 
& = T\Psi^{-1} \bigg( (-1)^{k} \dfrac{d}{d \epsilon} \bigg|_{\epsilon = 0} c(\Psi(g_1), \dots, \Psi(g_k), \Phi^\alpha_\epsilon(\mathsf s (\Psi(g_k)))^{-1}) \bigg)
\end{aligned}
\end{equation} 
Using (\ref{eq:inv_vf}), we compute:
\[
\begin{aligned}
\Phi^\alpha_\epsilon(\mathsf s(\Psi(g_k))) & = \Phi^{\overrightarrow{\alpha}}_\epsilon(\mathsf 1_{\mathsf s(\Psi(g_k))}) = \Phi^{\overrightarrow{\alpha}}_\epsilon(\Psi(\mathsf 1_{\mathsf s(g_k)})) \\ & = \Psi(\Phi^{\Psi^*\overrightarrow{\alpha}}_\epsilon(\mathsf 1_{\mathsf s(g_k)})) = \Psi(\Phi_\epsilon^{\overrightarrow{\psi^* \alpha}}(\mathsf 1_{\mathsf s(g_k)})) = \Psi(\Phi^{\psi^* \alpha}_\epsilon(\mathsf s(g_k))).
\end{aligned}
\]
So, from (\ref{eq:Psi^*R}) we have:
\[
\begin{aligned}
\Psi^*(R_\alpha(c)) (g_1, \dots, g_k) & = (-1)^{k}  \dfrac{d}{d \epsilon} \bigg|_{\epsilon = 0} \Psi^{-1}(c(\Psi(g_1), \dots, \Psi(g_k), \Psi(\Phi_\epsilon^{\psi^* \alpha} (\mathsf s(g_k)))^{-1})) \\
& = (-1)^{k} \dfrac{d}{d\epsilon} \bigg|_{\epsilon = 0} (\Psi^* c)(g_1, \dots, g_k, \Phi_\epsilon^{\psi^* \alpha}(\mathsf s(g_k))^{-1}) \\
& = R_{\psi^* \alpha}(\Psi^* c)(g_1, \dots, g_k).
\end{aligned}
\]
Finally, by applying repeatedly formula (\ref{eq:R_equiv}) in (\ref{eq:VE_formula}), we obtain
\[
\Psi^*(\mathrm{VE}(c)) = \mathrm{VE}(\Psi^* c)
\]
for every $c \in \widehat C_{\mathrm{def}}(\mathcal G)$, as desired. \endproof

Now we are ready for the main theorem of this subsection. Notice that the first part of the following statement has been already proved in \cite{etv:infinitesimal} in the special case where the rank of $E$ is greater than zero. Here we provide an alternative proof exploiting Lemma \ref{lem:VE_equiv} which is valid in all cases.

\begin{theorem}[Linear van Est map] \label{prop:lin_van_est} Let $(\mathcal W \rightrightarrows E; \mathcal G \rightrightarrows M)$ be a VB-groupoid, and let $(W \Rightarrow E; A \Rightarrow M)$ be its VB-algebroid. Then the van Est map for the Lie groupoid $\mathcal W \rightrightarrows E$ restricts to a cochain map
\[
\mathrm{VE}: \widehat C_{\mathrm{def,lin}}(\mathcal W) \to C_{\mathrm{def,lin}}(W),
\]
which we call the \emph{linear van Est map}. If $\mathcal G$ has $k$-connected $\mathsf s$-fibers, this map induces an isomorphism in cohomology in all degrees $p < k$. \end{theorem}

\proof As before, we denote by $h$ the homogeneity structure of $\mathcal W \to \mathcal G$. Then the last proposition shows that
\[
h_\lambda^* (\mathrm{VE}(\tilde c)) = \mathrm{VE}(h_\lambda^* \tilde c)
\]
for every $\tilde c \in \widehat C_{\mathrm{def}}(\mathcal W)$, $\lambda > 0$. But the VB-algebroid automorphism corresponding to $h_\lambda$ is exactly the one induced by the homogeneity structure of $W \to A$, so the last equation implies that the van Est map preserves linear cochains.

Now we prove the second part of the theorem. First, we have to observe that the van Est map commutes with linearization, i.e.
\begin{equation}\label{eq:VE_lin}
\mathrm{VE}(\tilde c)_{\mathrm{lin}} = \mathrm{VE}(\tilde c_{\mathrm{lin}}).
\end{equation}
To see this, take $\tilde c \in C^k_{\mathrm{def,lin}}(\mathcal W)$ and $w_0, \dots, w_k \in \Gamma(W,E)$ and compute
\[
\begin{aligned}
\mathrm{VE}(\tilde c)_{\mathrm{lin}} (w_0, \dots, w_k) & = \lim_{\lambda \to 0} \big( h_\lambda^* (\mathrm{VE}(\tilde c)) - \lambda^{-1} \mathrm{VE}(\tilde c_{\mathrm{core}}) \big) (w_0, \dots, w_k) \\
& = \lim_{\lambda \to 0} \mathrm{VE}(h_\lambda^* \tilde c - \lambda^{-1} \tilde c_{\mathrm{core}}) (w_0, \dots, w_k) \\
& = \lim_{\lambda \to 0} \sum_{\tau \in S_{k+1}} (-1)^\tau (R_{w_{\tau(k)}} \circ \dots \circ R_{w_{\tau(0)}}) (h_\lambda^* \tilde c - \lambda^{-1} \tilde c_{\mathrm{core}}).
\end{aligned}
\]
We would like to swap the limit and the derivatives that appear in definitions (\ref{eq:R_alpha_1}) and (\ref{eq:R_alpha_2}). This is ultimately possible because of smoothness,
and we get
\[
\begin{aligned}
\mathrm{VE}(\tilde c)_{\mathrm{lin}}(w_0, \dots, w_k) & = \sum_{\tau \in S_{k+1}} (-1)^\tau (R_{w_{\tau(k)}} \circ \dots \circ R_{w_{\tau(0)}}) \big( \lim_{\lambda \to 0} (h_\lambda^* \tilde c - \lambda^{-1} \tilde c_{\mathrm{core}}) \big) \\
& = \sum_{\tau \in S_{k+1}} (-1)^\tau (R_{w_{\tau(k)}} \circ \dots \circ R_{w_{\tau(0)}}) (\tilde c_{\mathrm{lin}}) \\
& = \mathrm{VE}(\tilde c_{\mathrm{lin}}) (w_0, \dots, w_k)
\end{aligned}
\]
as desired.

Finally, suppose that $\mathcal G$ has $k$-connected $\mathsf s$-fibers. Then $\mathcal W$ has $k$-connected $\tilde{\mathsf s}$-fibers (they are vector bundles over the $\mathsf s$-fibers of $\mathcal G$). Take $p < k$. We want to prove that the induced map
\[
\mathrm{VE}: H^p_{\mathrm{def,lin}}(\mathcal W) \to H^p_{\mathrm{def,lin}}(W)
\]
is an isomorphism. If $\tilde c \in C^p_{\mathrm{def,lin}}(\mathcal W)$ is closed, we denote $[\tilde c]$ its class in $H^p_{\mathrm{def}}(\mathcal W)$ and $[\tilde c]_{\mathrm{lin}}$ its class in $H^p_{\mathrm{def,lin}}(\mathcal W)$. We use an analogous notation for $W$.

First, suppose that $[\mathrm{VE}(\tilde c)]_{\mathrm{lin}} = 0$. Then $[\mathrm{VE}(\tilde c)] = 0$ and, by Theorem \ref{prop:van_est}, $[\tilde c] = 0$, i.e.~$\tilde c = \delta \tilde \gamma$ for some $\tilde \gamma \in C^{p-1}_{\mathrm{def}}(\mathcal W)$. Applying the linearization map, we get $\tilde c = \delta \tilde \gamma_{\mathrm{lin}}$ and $\tilde \gamma_{\mathrm{lin}} \in \widehat C^{p-1}_{\mathrm{def,lin}}(\mathcal W)$, i.e.~$\mathrm{VE}$ is injective in degree $p$ cohomology. To conclude, take $[c]_{\mathrm{lin}} \in H^p_{\mathrm{def,lin}}(W)$. Then $[c] \in H^p_{\mathrm{def}}(W)$, so $[c] = [\mathrm{VE}(\tilde c)]$, i.e.~$c - \mathrm{VE}(\tilde c) = \delta \gamma$ for some $\tilde c \in \widehat C^p_{\mathrm{def}}(\mathcal W), \gamma \in C^{p-1}_{\mathrm{def}}(W)$. Applying again the linearization map and using (\ref{eq:VE_lin}), we get $c - \mathrm{VE}(\tilde c_{\mathrm{lin}}) = \delta \gamma_{\mathrm{lin}}$, i.e.~$\mathrm{VE}$ is also surjective in degree $p$ cohomology, as desired. \endproof

\subsection{Morita invariance} \label{sec:morita}

The notion of Morita equivalence of VB-groupoids first appears in \cite{delhoyo:morita}. In that reference, the authors prove that the VB-cohomologies of Morita equivalent VB-groupoids are isomorphic. As a corollary, they give a conceptual and very simple proof of the fact, first appeared in \cite{crainic:def2}, that Morita equivalent Lie groupoids have isomorphic deformation cohomologies. This second result means that the deformation cohomology of a Lie groupoid is in fact an invariant of the associated differentiable stack.

In this paragraph, we want to prove an analogous result for the linear deformation cohomology of a VB-groupoid. Let $(\mathcal W_1 \rightrightarrows E_1; \mathcal G_1 \rightrightarrows M_1)$ and $(\mathcal W_2 \rightrightarrows E_2; \mathcal G_2 \rightrightarrows M_2)$ be VB-groupoids. A VB-groupoid morphism $\Psi: \mathcal W_1 \to \mathcal W_2$ is a \emph{VB-Morita map} \cite{delhoyo:morita} if the Lie groupoid morphism $\Psi$ is a Morita map. The VB-groupoids $\mathcal W_1$ and $\mathcal W_2$ are \emph{Morita equivalent} if there exist a VB-groupoid $\mathcal V$ and VB-Morita maps $\mathcal V \to \mathcal W_1$, $\mathcal V \to \mathcal W_2$.

Here are some basic properties of VB-Morita maps.

\begin{proposition}[{\cite[Corollary 3.7]{delhoyo:morita}}]\label{prop:VBM_1} Let $\Psi: \mathcal G_1 \to \mathcal G_2$ be a Morita map and let $\mathcal W$ be a VB-groupoid over $\mathcal G_2$. Then the canonical map $\Psi^* \mathcal W \to \mathcal W$ is VB-Morita. \end{proposition}

\begin{proposition}[{\cite[Corollary 3.8]{delhoyo:morita}}]\label{prop:VBM_2} Let $\Psi: \mathcal G_1 \to \mathcal G_2$ be a Morita map. Then its tangent map $T \Psi: T \mathcal G_1 \to T \mathcal G_2$ is a VB-Morita map. \end{proposition} 

\begin{proposition}[{\cite[Corollary 3.9]{delhoyo:morita}}]\label{prop:VBM_3} A map $\Psi: \mathcal W_1 \to \mathcal W_2$ over the identity is VB-Morita if and only if its dual is so. \end{proposition}

Morita invariance of the VB-cohomology is expressed by the following theorem.

\begin{theorem} [{\cite[Theorem 4.2]{delhoyo:morita}}] Let $\Psi: \mathcal W_1 \to \mathcal W_2$ be a VB-Morita map. Then $\Psi^*: H_{\mathrm{proj}}(\mathcal W_2) \to H_{\mathrm{proj}}(\mathcal W_1)$ is an isomorphism. \end{theorem}

Now we are ready to prove Morita invariance of the linear deformation cohomology.

\begin{theorem} Let $\Psi: \mathcal W_1 \to \mathcal W_2$ be a VB-Morita map. Then $H_{\mathrm{def,lin}}(\mathcal W_1) \cong H_{\mathrm{def,lin}}(\mathcal W_2)$. \end{theorem}

\proof It is enough to show that $H_{\mathrm{proj,lin}}(T^* \mathcal W_1) \cong H_{\mathrm{proj,lin}}(T^* \mathcal W_2)$. To do this we will use Propositions \ref{prop:VBM_1}--\ref{prop:VBM_3} and linearization.

Recall that both $T^* \mathcal W_1$ and $T^* \mathcal W_2$ have two VB-groupoid structures, as discussed in Subsection \ref{sec:lin_def_1}, and observe that $\Psi^* (T^* \mathcal W_2)$ possesses also two VB-groupoid structures, that fit in the following commuting diagram:
\[
\begin{array}{c}
\xymatrix@C=9pt@R=12pt{
& \Psi^* (T^* \mathcal W_2) \ar[dl] \ar@<0.4ex>[rr] \ar@<-0.4ex>[rr] \ar[dd]|!{[dl];[dr]}{\hole} & & \Psi^* (W_2^* |_{E_2}) \ar[dl] \ar[dd] \\
\mathcal W_1 \ar@<0.4ex>[rr] \ar@<-0.4ex>[rr] \ar[dd] & & E_1 \ar[dd]\\
& \Psi^* \mathcal W_2^* \ar[dl] \ar@<0.4ex>[rr]|!{[ur];[dr]}{\hole} \ar@<-0.4ex>[rr]|!{[ur];[dr]}{\hole} & & \Psi^* C_2^* \ar[dl] \\
\mathcal G_1 \ar@<0.4ex>[rr] \ar@<-0.4ex>[rr] & & M_1 
}
\end{array}.
\]

We denote by $C_{\mathrm{proj},\bullet}(\Psi^* (T^* \mathcal W_2))$ the VB-complex of the VB-groupoid upstairs and by $C_{\mathrm{proj,lin}}(\Psi^* (T^* \mathcal W_2))$ its subcomplex of cochains that are linear with respect to the vertical vector bundle structures. We denote their cohomologies by $H_{\mathrm{proj},\bullet}(\Psi^* (T^* \mathcal W_2))$ and $H_{\mathrm{proj,lin}}(\Psi^* (T^* \mathcal W_2))$ respectively. Now, from \cite{cabrera:hom} we know that there is a linearization map
\[
\mathrm{lin}: C(\Psi^*(T^* \mathcal W_2)) \to C_{\mathrm{lin}}(\Psi^* (T^* \mathcal W_2)), \quad c \mapsto c_{\mathrm{lin}},
\]
that is a cochain map splitting the inclusion $C_{\mathrm{lin}}(\Psi^*(T^* \mathcal W_2)) \hookrightarrow C(\Psi^*(T^* \mathcal W_2))$. If $\cdot_{\Psi^* \mathcal W^*_2}$ denotes the scalar multiplication in the fibers of $\Psi^* (T^* \mathcal W_2) \to \Psi^* \mathcal W_2^*$, $c \in C^k (\Psi^* (T^* \mathcal W_2))$, then $c_{\mathrm{lin}}$ is simply given by
\[
c_{\mathrm{lin}}(\theta_1, \dots, \theta_k) = \dfrac{d}{d \lambda} \bigg|_{\lambda=0} c( \lambda \cdot_{\Psi^* \mathcal W^*_2} (\theta_1, \dots, \theta_k)).
\]
From the definition, it is clear that lin restricts to a cochain map
\[
\mathrm{lin}: C_{\mathrm{proj},\bullet}(\Psi^*(T^* \mathcal W_2)) \to C_{\mathrm{proj,lin}}(\Psi^* (T^* \mathcal W_2)),
\]
so $H_{\mathrm{proj,lin}}(\Psi^* (T^* \mathcal W_2))$ embeds in $H_{\mathrm{proj},\bullet}(\Psi^* (T^* \mathcal W_2))$ as a direct summand.

The map $\Psi$ is VB-Morita, so, from Proposition \ref{prop:VBM_2}, $T \Psi: T \mathcal W_1 \to T \mathcal W_2$ is VB-Morita again. It follows from Propositions \ref{prop:VBM_1} and \ref{prop:VBM_3} that the dual map $(T \Psi)^*: \Psi^* (T^* \mathcal W_2) \to T^* \mathcal W_1$ is VB-Morita as well.  Remember that also the canonical map $\Psi^* (T^* \mathcal W_2) \to T^* \mathcal W_2$ is VB-Morita (Proposition \ref{prop:VBM_1} again). As a result, we get isomorphisms in VB-cohomology:
\begin{equation}\label{eq:proj_coh}
H_{\mathrm{proj},\bullet}(T^* \mathcal W_1) \overset{\cong}{\longrightarrow} H_{\mathrm{proj},\bullet}(\Psi^*(T^* \mathcal W_2)) \overset{\cong}{\longleftarrow} H_{\mathrm{proj},\bullet}(T^* \mathcal W_2).
\end{equation}
Now, the maps 
\[
T^* \mathcal W_1 \overset{(T\Psi)^*}{\longleftarrow} \Psi^* (T^* \mathcal W_2) \longrightarrow T^* \mathcal W_2
\]
are also DVB morphisms. This implies, on one hand, that the maps (\ref{eq:proj_coh}) preserve linear cohomologies, on the other hand that they commute with the respective linearization maps. From this last property we deduce, as in Theorem \ref{prop:lin_van_est}, that the maps (\ref{eq:proj_coh}) induce isomorphisms on linear cohomologies, so
\[
H_{\mathrm{proj,lin}}(T^* \mathcal W_1) \cong H_{\mathrm{proj,lin}}(\Psi^* (T^* \mathcal W_2)) \cong H_{\mathrm{proj,lin}}(T^* \mathcal W_2)
\]
as desired. \endproof

\section{Examples and applications}\label{Sec:examples_1}

In this section we provide several examples. Examples in Subsections \ref{sec:VB_grp}, \ref{sec:fol_grpd} and \ref{sec:Lie_grp_vect} parallel the analogous examples in \cite{crainic:def2}, connecting our linear deformation cohomology to known cohomologies, while examples in Subsections \ref{sec:2-vect} and \ref{sec:tangent_VB_grpd} are specific to VB-groupoids.

\subsection{VB-groups and their duals}\label{sec:VB_grp}

A \emph{VB-group} is a \emph{vector bundle object in the category of Lie groups}. In other words, it is a VB-groupoid of the form
\[
\begin{array}{r}
\xymatrix{
K \ar[d]_p \ar@<0.4ex>[r] \ar@<-0.4ex>[r] & 0 \ar[d] \\
G \ar@<0.4ex>[r] \ar@<-0.4ex>[r] & \ast}
\end{array}.
\]

In particular, $K$ and $G$ are Lie groups. Let $C := \ker p$ be the core of $(K \rightrightarrows 0; G \rightrightarrows \ast)$. The map $p$ has a canonical section in the category of Lie groups, the zero section $0: G \to K$. By standard Lie group theory, $G$ acts on the core $C$ via
\[
g \star c = 0_g \cdot c \cdot 0_g^{-1}
\]
and there is a canonical isomorphism of Lie groups
\[
K \overset{\cong}{\longrightarrow} G \ltimes C, \quad k \mapsto (p(k), k \cdot 0_{p(k)}^{-1}),
\]
where $G \ltimes C$ is the semidirect product Lie group. Thus VB-groups are equivalent to Lie group representations.


It is then natural to study the relationship between the linear deformation complex of $K$ and the classical complex $C(G, \operatorname{End} C)$ (of the Lie group $G$ with coefficients in the representation $\operatorname{End} C$)  that controls deformations of the $G$-module $C$ \cite{nijenhuis:def}. To do this, we notice that the dual of $K$ is the VB-groupoid
\[
\begin{array}{r}
\xymatrix{
G \ltimes C^* \ar[d] \ar@<0.4ex>[r] \ar@<-0.4ex>[r] & C^* \ar[d] \\
G \ar@<0.4ex>[r] \ar@<-0.4ex>[r] & \ast}
\end{array},
\]
i.e.~it is the action VB-groupoid associated to the dual representation of $G$ on $C^*$. In particular, it is a trivial-core VB-groupoid, and the short exact sequence (\ref{eq:coch_2}) reads:
\[
0 \longrightarrow C(G, \operatorname{End} C^*) \longrightarrow C_{\mathrm{def,lin}}(G \ltimes C^*) \longrightarrow C_{\mathrm{def}}(G) \longrightarrow 0.
\]
But $C(G, \operatorname{End} C^*) \cong C(G, \operatorname{End} C)$ canonically, so the latter is recovered as the subcomplex of $C_{\mathrm{def,lin}}(G \ltimes C^*)$ controlling deformations of the representation $C^*$ that fix the Lie group structure on $G$. Moreover, by Theorem \ref{thm:dual_vb_grpd} $H_{\mathrm{def,lin}}(G \ltimes C^*) \cong H_{\mathrm{def,lin}}(K)$, so the long exact sequence (\ref{eq:les_trivial_core}) reads:
\[
\cdots \longrightarrow H^k(G, \operatorname{End} C) \longrightarrow H^k_{\mathrm{def}, \mathrm{lin}}(K) \longrightarrow H^k_{\mathrm{def}}(G) \longrightarrow H^{k+1}(G, \operatorname{End} C) \longrightarrow \cdots .
\]

Finally, suppose that $G$ is compact. Then Proposition \ref{prop:rigidity} applies, and we get the following

\begin{corollary} Let $(K \rightrightarrows 0; G \rightrightarrows *)$ be a VB-group with compact base. Then $H^k_{\mathrm{def,lin}}(K) = 0$ for $k \geq 1$. \end{corollary}

\subsection{2-vector spaces}\label{sec:2-vect}

A \emph{2-vector space} is a \emph{(Lie) groupoid object in the category of vector spaces}. In other words, it is a VB-groupoid of the form
\[
\begin{array}{r}
\xymatrix{
V_1 \ar[d] \ar@<0.4ex>[r] \ar@<-0.4ex>[r] & V_0 \ar[d] \\
\ast \ar@<0.4ex>[r] \ar@<-0.4ex>[r] & \ast}
\end{array}.
\]

If $\mathsf s$ and $\mathsf t$ are the source and the target maps of $V_1 \rightrightarrows V_0$ and $C = \ker \mathsf s$, then $\mathsf 1: V_0 \to V_1$ is a right inverse of $\mathsf s$ and so there is a canonical isomorphism of vector spaces:
\begin{equation} \label{eq:2-vect_0}
V_1 \overset{\cong}{\longrightarrow} C \oplus V_0, \quad v \mapsto (v - \mathsf 1_{\mathsf s (v)}, \mathsf s (v)).
\end{equation}
Denote $\partial = \mathsf t|_C: C \to V_0$. Then $C$ acts on $V_0$ by
\[
c \cdot v = \partial c + v
\]
and one can consider the associated action groupoid $C \ltimes V_0 \rightrightarrows V_0$. Notice that $C$ does not act by linear isomorphisms, but by translations. In \cite{baez:higher} it is proved that (\ref{eq:2-vect_0}) is also an isomorphism of groupoids $V_1 \cong C \ltimes V_0$: from now on, we will identify $V_1$ with $C \ltimes V_0$.

Now we describe the linear deformation complex of $(V_1 \rightrightarrows V_0; * \rightrightarrows *)$. First of all remember that, for every finite-dimensional real vector space $V$, $TV \cong V \oplus V$ canonically. The structure maps of $V_1$ are given by (\ref{eq:action_grpd}). It follows that, for every $(c,v) \in V_1$,
\begin{equation} \label{eq:2-vect_maps}
\begin{aligned}
L_{(c,v)} & : V_1 (-,v) \to V_1(-, c \cdot v), & (c',v') & \mapsto (c+c',v'); \\
R_{(c,v)} & : V_1 (c \cdot v, -) \to V_1 (v,-), & (c',c \cdot v) & \mapsto (c'+c,v); \\
TL_{(c,v)} & : T V_1 (-,v) \to T V_1(-, c \cdot v), & (c',v', \dot c, \dot v) & \mapsto (c+c',v', \dot c, \dot v); \\
TR_{(c,v)} & : T V_1 (c \cdot v, -) \to T V_1 (v,-), & (c',cv, \dot c, 0) & \mapsto (c'+c, v, \dot c, 0).
\end{aligned}
\end{equation}
Moreover, we compute:
\begin{equation} \label{eq:2-vect_maps_2}
\begin{aligned}
\mathsf i & : V_1 \to V_1 & (c,v) & \mapsto (-c, c \cdot v); \\
\bar{\mathsf m} & : V_1^{[2]} \to V_1 & ((c_1,v),(c_2,v)) & \mapsto (c_1-c_2, c_2 \cdot v); \\
T \mathsf i & : TV_1 \to TV_1, & (c,v, \dot c, \dot v) & \mapsto (-c, c \cdot v, - \dot c, \dot c \cdot \dot v); \\
T \bar{\mathsf m} & : TV_1^{[2]} \to TV_1 & ((c_1,v,\dot c_1, \dot v),(c_2,v,\dot c_2, \dot v)) & \mapsto (c_1-c_2, c_2 \cdot v, \dot c_1 - \dot c_2, \dot c_2 \cdot \dot v).
\end{aligned}
\end{equation}

As $\mathsf s: C \oplus V_0 \to V_0$ is the first projection, $\ker T \mathsf s$ can be identified with the subset $C \oplus V_0 \oplus C \oplus 0$ of $C \oplus V_0 \oplus C \oplus V_0$, so the Lie algebroid of $V_1 \rightrightarrows V_0$ can be identified with $0 \oplus V_0 \oplus C \oplus 0 \to V_0$.

Suppose that $\gamma \in C^{-1}_{\mathrm{def,lin}}(V_1) = \Gamma_{\mathrm{lin}}(V_1 \Rightarrow V_0)$. Then $\gamma$ has the form
\begin{equation} \label{eq:gamma_1}
\gamma: V_0 \to 0 \oplus V_0 \oplus C \oplus 0, \quad v \mapsto (0, v, \gamma_1(v), 0)
\end{equation}
for some linear map $\gamma_1: V_0 \to C$. Hence we have a canonical isomorphism
\begin{equation} \label{eq:2-vect_iso_1}
C^{-1}_{\mathrm{def,lin}}(V_1) \overset{\cong}{\longrightarrow} \operatorname{Hom}(V_0,C), \quad \gamma \mapsto \gamma_1.
\end{equation}
Let now $\gamma \in C^k_{\mathrm{def,lin}}(V_1)$, $k \geq 0$. As $V_1 \rightrightarrows V_0$ is an action groupoid, Equation (\ref{eq:nerve_act}) yields a canonical isomorphism $(C \ltimes V_0)^{(k)} \cong C^k \oplus V_0$. Its inverse is given by:
\begin{equation} \label{eq:2-vect_iso}
\begin{aligned}
C^k \oplus V_0 & \overset{\cong}{\longrightarrow} (C \ltimes V_0)^{(k)} \\ 
(c_1, \dots, c_k, v) & \mapsto ((c_1, (c_2 + \dots + c_k)v), (c_2, (c_3 + \dots + c_k)v), \dots, (c_k, v)).
\end{aligned}
\end{equation}
In what follows, we will understand this isomorphism. This allows us to identify $\gamma$ with a linear map $\gamma: C^{k+1} \oplus V_0 \to C \oplus V_0 \oplus C \oplus V_0$, and conditions 1 and 2 of Definition \ref{def:def_complex} mean that $\gamma$ takes the form
\begin{equation} \label{eq:gamma_1,2}
\gamma(c_0, \dots, c_k, v) = (c_0, (c_1 + \dots + c_k) \cdot v, \gamma_1 (c_0, \dots, c_k, v), \gamma_2 (c_1, \dots, c_k, v))
\end{equation}
for some linear maps $\gamma_1: C^{k+1} \oplus V \to C$, $\gamma_2: C^k \oplus V \to C$. So there is a canonical isomorphism
\begin{equation} \label{eq:2-vect_iso_2}
C^k_{\mathrm{def,lin}}(V_1) \overset{\cong}{\longrightarrow} \operatorname{Hom}(C^{k+1} \oplus V_0, C) \oplus \operatorname{Hom}(C^k \oplus V_0, C), \quad \gamma \mapsto (\gamma_1, \gamma_2).
\end{equation}
We are left with describing the differential in terms of isomorphisms (\ref{eq:2-vect_iso_1}) and (\ref{eq:2-vect_iso_2}). Let us start with $\gamma \in C^{-1}_{\mathrm{def,lin}}(V_1)$. For every $(c,v) \in V_1$,
\[
\begin{aligned}
\delta \gamma (c,v) & = \overleftarrow{\gamma}_{(c,v)} + \overrightarrow{\gamma}_{(c,v)} \\
& = TR_{(c,v)} (\gamma(\mathsf t (c,v))) + (TL_{(c,v)} \circ T \mathsf i)(\gamma(\mathsf s (c,v))) \\
& = TR_{(c,v)}(0,c \cdot v, \gamma_1 (c \cdot v), 0) + (TL_{(c,v)} \circ T \mathsf i)(0,v, \gamma_1(v),0) \\
& = (c,v,\gamma_1(c \cdot v),0) + (c,v,-\gamma_1(v), (\partial \circ \gamma_1)(v)) \\
& = (c,v,\gamma_1(c \cdot v) - \gamma_1(v), (\partial \circ \gamma_1)(v)) \\
& = (c,v,(\gamma_1 \circ \partial)(c), (\partial \circ \gamma_1)(v)).
\end{aligned}
\]
Here we used (\ref{eq:2-vect_maps}), (\ref{eq:2-vect_maps_2}) and (\ref{eq:gamma_1}). We deduce that
\[
\delta \gamma_1 = (\gamma_1 \circ \partial, \partial \circ \gamma_1).
\]
Now let $\gamma \in C^k_{\mathrm{def,lin}}(V_1)$, $k \geq 0$. Using isomorphism (\ref{eq:2-vect_iso}), Formula (\ref{eq:diff_grpd}) translates into:
\begin{alignat*}
\delta \gamma (c_0, \dots, c_{k+1}, v) & = && - T \bar{\mathsf m} (\gamma(c_0 + c_1, c_2, \dots, c_{k+1}, v), \gamma(c_1, \dots, c_{k+1}, v)) \\
& && + \sum_{i=1}^k (-1)^{i-1} \gamma(c_0, \dots, c_i + c_{i+1}, \dots, c_{k+1}, v) + (-1)^k \gamma(c_0, \dots, c_k, v).
\end{alignat*}
Applying (\ref{eq:2-vect_maps_2}) and (\ref{eq:gamma_1,2}), we get:
\[
\begin{aligned}
& T \bar{\mathsf m} (\gamma(c_0 + c_1, c_2, \dots, c_{k+1}, v), \gamma(c_1, \dots, c_{k+1}, v)) \\
& = T \bar{\mathsf m}((c_0+c_1, (c_2 + \dots + c_{k+1}) \cdot v, \gamma_1 (c_0+c_1, c_2, \dots, c_{k+1}, v), \gamma_2 (c_2, \dots, c_{k+1}, v)), \\
& \quad (c_1, (c_2 + \dots + c_{k+1}) \cdot v, \gamma_1 (c_1, \dots, c_{k+1}, v), \gamma_2 (c_2, \dots, c_{k+1}, v))) \\
& = (c_0, (c_1 + \dots + c_{k+1}) \cdot v, \gamma_1 (c_0, 0, \dots, 0), \gamma_1 (c_1, \dots, c_{k+1}, v) \cdot \gamma_2 (c_2, \dots, c_{k+1}, v)).
\end{aligned}
\]
It is clear then that $\delta (\gamma_1, \gamma_2) = (\Gamma_1, \Gamma_2)$, with
\[
\begin{aligned}
\Gamma_1 (c_0, \dots, c_{k+1}, v) = & -\gamma_1 (c_0, 0, \dots, 0) \\
& + \sum_{i=1}^k (-1)^{i-1} \gamma_1 (c_0, \dots, c_i + c_{i+1}, \dots, c_{k+1}, v)  \\
& + (-1)^k \gamma_1 (c_0, \dots, c_k, c_{k+1} \cdot v),
\end{aligned}
\]
and 
\[
\begin{aligned}
\Gamma_2 (c_1, \dots, c_{k+1}, v) = & - \gamma_1 (c_1, \dots, c_{k+1}, v) \cdot \gamma_2 (c_2, \dots, c_{k+1}, v) \\
& + \sum_{i=1}^k (-1)^{i-1} \gamma_2 (c_1, \dots, c_i + c_{i+1}, \dots, c_{k+1}, v) \\
& + (-1)^k \gamma_2 (c_1, \dots, c_k, c_{k+1} \cdot v).
\end{aligned}
\]

We further notice that, if $\gamma \in C^k_{\mathrm{def,lin}}(V_1)$, $k \geq 0$, then $\gamma = (\gamma_1, \gamma_2)$ belongs to the normalized deformations subcomplex $\widehat C_{\mathrm{def,lin}}(V_1)$ if and only if
\begin{center}
\begin{tabular}{cc}
$\gamma_1 (c_0, \dots, \underset{i}{0}, \dots, c_k, v) = 0$ & for every $i \geq 0$, \\
$\gamma_2 (c_1, \dots, \underset{i}{0}, \dots, c_k, v) = 0$ & for every $i \geq 1$. \\ 
\end{tabular}
\end{center}
It follows that $\widehat C_{\mathrm{def,lin}}(V_1)$ reduces to
\begin{equation} \label{eq:2-vect_comp}
0 \longrightarrow \mathrm{Hom}(V_0, C)[1] \overset{\delta_0}{\longrightarrow} \operatorname{End} C \oplus \operatorname{End} V_0 \overset{\delta_1}{\longrightarrow} \mathrm{Hom}(C, V_0)[-1] \longrightarrow 0
\end{equation}
with
\[
\begin{aligned}
\delta_0 \gamma & = (\gamma \circ \partial, \partial \circ \gamma), \\
\delta_1 (\gamma_1, \gamma_2) & = \gamma_2 \circ \partial - \partial \circ \gamma_1. \\
\end{aligned}
\]
So the linear deformation cohomology of $V_1$ is:
\[
\begin{aligned}
H^{-1}_{\mathrm{def}, \mathrm{lin}} (V_1) & = \operatorname{Hom} (\operatorname{coker} \partial, \ker \partial),\\
H^0_{\mathrm{def}, \mathrm{lin}} (V_1) & = \operatorname{End} (\operatorname{coker} \partial) \oplus \operatorname{End} (\ker \partial), \\
H^1_{\mathrm{def}, \mathrm{lin}} (V_1) & = \operatorname{Hom}(\ker \partial, \operatorname{coker} \partial).
\end{aligned}
\]

Finally, notice that the VB-algebroid of $V_1$ is an LA-vector space of the form $(V_1 \Rightarrow V_0; 0 \Rightarrow *)$. In Subsection \ref{sec:LA-vect} we showed that the linear deformation complex of $(V_1 \Rightarrow V_0; 0 \Rightarrow *)$ is again (\ref{eq:2-vect_comp}), and one can see, for example in coordinates, that the van Est map
\begin{equation*}
\mathrm{VE}: \widehat C_{\mathrm{def,lin}}(V_1 \rightrightarrows V_0) \to C_{\mathrm{def,lin}}(V_1 \Rightarrow V_0)
\end{equation*}
is simply the identity.

\subsection{Tangent and cotangent VB-groupoids} \label{sec:tangent_VB_grpd}

Let $\mathcal G \rightrightarrows M$ be a Lie groupoid. We want to relate the linear deformation cohomology of $T \mathcal G$ with the deformation cohomology of $\mathcal G$. First recall that there is a projection
\begin{equation}\label{eq:proj_TG}
\Pi: C_{\mathrm{def,lin}}(T \mathcal G) \to C_{\mathrm{def}}(\mathcal G).
\end{equation}

We want to show that there is a canonical right inverse of this projection. To do this, we recall from \cite{mackenzie} the definition of the \emph{canonical involution} of the double tangent bundle of a smooth manifold. Let $N$ be a smooth manifold and let $p: TN \to N$ be the canonical projection. Given a smooth map $u: \mathbb R^2 \to N$, we define its second mixed derivative 
\[
\dfrac{\partial^2 u}{\partial \epsilon \partial \eta} \bigg|_{(0,0)} \in TTN
\]
in the following way: we consider the first derivative $\frac{\partial u}{\partial \eta}(-,0): \mathbb R \to TN$ as a curve in $TN$ and then we take its derivative at $\epsilon = 0$. One can show in coordinates that every tangent vector in $TTN$ can be represented in this way.

The canonical involution is the DVB morphism $J: TTN \to TTN$ defined by:
\[
J \bigg( \dfrac{\partial^2 u}{\partial \epsilon \partial \eta} \bigg|_{(0,0)} \bigg) := \dfrac{\partial^2 u}{\partial \eta \partial \epsilon} \bigg|_{(0,0)}.
\]
It is not difficult to see that $J$ is a DVB automorphism of $TTN$ that exchanges the two projections $Tp: TTN \to T \mathcal G$ and $p_{TN}: TTN \to TN$:
\begin{equation}\label{eq:J}
J \circ Tp = p_{TN} \circ J.
\end{equation}
Moreover, suppose that $f: N_1 \to N_2$ is a smooth map, and denote $J_1$, $J_2$ the canonical involutions of $TTN_1$ and $TTN_2$, respectively. Then one can verify in coordinates that 
\[
TTf \bigg( \dfrac{\partial^2 u}{\partial \epsilon \partial \eta} \bigg|_{(0,0)} \bigg) = \dfrac{\partial^2}{\partial \epsilon \partial \eta} \bigg|_{(0,0)} (f \circ u)
\]
and
\begin{equation} \label{eq:TTf}
J_2 \circ TTf = TTf \circ J_1.
\end{equation}

Now define an inclusion
\[
\iota: C_{\mathrm{def}}(\mathcal G) \to C_{\mathrm{def,lin}}(T \mathcal G), \quad  c \mapsto \iota_c := J \circ Tc.
\]
The inclusion $\iota$ is well-defined, i.e., if $c \in C^k_{\mathrm{def}}(\mathcal G)$, then $\iota_c \in C^k_{\mathrm{def,lin}}(T \mathcal G)$. To see this we show that properties 1 and 2 of Definition \ref{def:def_complex} hold. So, let $\epsilon \mapsto (g_0 (\epsilon), \dots, g_k (\epsilon))$ be a curve in $\mathcal G^{(k+1)}$ defined around 0. Then $(p \circ c)(g_0 (\epsilon), \dots, g_k (\epsilon)) = g_1 (\epsilon)$. Differentiating at $\epsilon = 0$ we get
\[
(Tp \circ Tc)(\dot g_0(0), \dots, \dot g_k(0)) = \dot g_0(0).
\]
Applying $J$ and remembering Equation (\ref{eq:J}) we get
\[
p_{T\mathcal G} (\iota_c (\dot g_0(0), \dots, \dot g_k(0))) = \dot g_0(0),
\]
i.e.~$\iota_c (\dot g_0(0), \dots, \dot g_k(0)) \in T_{\dot g_0(0)} T \mathcal G$ as desired.

Now notice that property 2 of Definition \ref{def:def_complex} can be expressed in the following way: if $q: \mathcal G^{(k+1)} \to \mathcal G^{(k)}$ is the map that forgets the first arrow, then $T\mathsf s \circ c$ descends to a map $\mathcal G^{(k)} \to TM$:
\[
\begin{array}{c}
\xymatrix@C=45pt{
& \mathcal{G}^{(k+1)} \ar[r]^-{T \mathsf s \circ c} \ar[d]_-q & TM \\
& \mathcal G ^{(k)} \ar[ur]}
\end{array}.
\]
Differentiating, we obtain that also $TT \mathsf s \circ Tc: T \mathcal G^{(k+1)} \to TTM$ descends to a map $T \mathcal G^{(k)} \to TTM$. Applying $J$ and remembering that it commutes with $TT \mathsf s$, we obtain the same statement for $TT \mathsf s \circ \iota_c$, as desired. Moreover, $\iota_c$ is obviously linear.

Next, we show that $\iota$ is a cochain map. Let $k \geq 0$, $c \in C^k_{\mathrm{def}}(\mathcal G)$ and $\epsilon \mapsto (g_0(\epsilon), \dots, g_{k+1}(\epsilon))$ be a curve in $\mathcal G^{(k+2)}$ defined around $0$. Then, using Equation (\ref{eq:TTf}), we get:
\[
\begin{aligned}
& \iota_{\delta c} (\dot g_0 (0), \dots, \dot g_{k+1} (0)) \\
& = (J \circ T(\delta c))(\dot g_0 (0), \dots, \dot g_{k+1} (0)) \\
& = J \bigg( \frac{d}{d\epsilon} \bigg|_{\epsilon = 0} \delta c (g_0 (\epsilon), \dots, g_{k+1} (\epsilon)) \bigg) \\
& = J \bigg( \frac{d}{d \epsilon} \bigg|_{\epsilon = 0} ( - T \bar{\mathsf m} (c(g_0 (\epsilon) g_1 (\epsilon), \dots, g_{k+1} (\epsilon)), c(g_1 (\epsilon), \dots, g_{k+1}(\epsilon))) \\
& \quad + \sum_{i=1}^{k+1} (-1)^{i+1} c(g_0(\epsilon), \dots, g_i (\epsilon) g_{i+1} (\epsilon), \dots, g_{k+1}(\epsilon)) + (-1)^{k+1} c(g_0(\epsilon), \dots, g_k(\epsilon))) \bigg) \\
& = J ( - TT \bar{\mathsf m} (Tc(\dot g_0 (0) \dot g_1 (0), \dots, \dot g_{k+1} (0)), Tc(\dot g_1 (0), \dots, \dot g_{k+1}(0))) \\
& \quad + \sum_{i=1}^{k+1} (-1)^{i+1} Tc (\dot g_0(0), \dots, \dot g_i (0) \dot g_{i+1} (0), \dots, \dot g_{k+1}(0)) + (-1)^{k+1} Tc(\dot g_0(0), \dots, \dot g_k(0))) \\
& = - TT \bar{\mathsf m} ((J \circ Tc)(\dot g_0 (0) \dot g_1 (0), \dots, \dot g_{k+1} (0)), (J \circ Tc)(\dot g_1 (0), \dots, \dot g_{k+1}(0))) \\
& \quad + \sum_{i=1}^{k+1} (-1)^{i+1} (J \circ Tc) (\dot g_0(0), \dots, \dot g_i (0) \dot g_{i+1} (0), \dots, \dot g_{k+1}(0)) \\
& \quad + (-1)^{k+1} (J \circ Tc) (\dot g_0(0), \dots, \dot g_k(0)) \\
& = \delta \iota_c (\dot g_0 (0), \dots, \dot g_{k+1} (0)).
\end{aligned}
\]
A similar computation holds for cochains of degree $-1$.

Finally, we observe that the following diagram commutes:
\[
\begin{array}{c}
\xymatrix{
 T \mathcal{G}^{(k+1)} \ar[r]^-{Tc} \ar[d] & TT \mathcal G \ar[d]^-{p_{T\mathcal G}} \ar[r]^-{J} & TT \mathcal G \ar[d]^-{Tp} \\
 \mathcal G ^{(k+1)} \ar[r]^-{c} & T \mathcal G \ar@{=}[r] & T \mathcal G 
}
\end{array}.
\]
This shows that $\iota$ inverts the projection (\ref{eq:proj_TG}). It follows that
\[
C_{\mathrm{def,lin}}(T \mathcal G) \cong C_{\mathrm{def}}(\mathcal G) \oplus \ker \Pi
\]
as cochain complexes, hence
\[
H_{\mathrm{def,lin}}(T\mathcal G) = H_{\mathrm{def,lin}}(T^* \mathcal G) = H_{\mathrm{def}}(\mathcal G) \oplus H (\ker \Pi).
\]

\subsection{Representations of foliation groupoids}\label{sec:fol_grpd}

A \emph{foliation groupoid} is a Lie groupoid whose anchor map is injective. This condition ensures that the connected components of the orbits of the groupoid are the leaves of a regular foliation of the base manifold, whence the name. On the other hand, foliation groupoids encompass several classical groupoids associated to a foliated manifold, such as the holonomy and the monodromy groupoids.

By Example \ref{ex:act_grpd}, representations of a foliation groupoid $\mathcal G$ are equivalent to trivial-core VB-groupoids over $\mathcal G$. Here we want to study the linear deformation cohomology of such VB-groupoids. First of all, let $\mathcal G \rightrightarrows M$ be a foliation groupoid, let $A \Rightarrow M$ be its Lie algebroid, $\rho: A \to TM$ the (injective) anchor map, $\nu = TM/ \operatorname{im} \rho$ the normal bundle and let $\pi: TM \to \nu$ be the projection. In this case, $\nu$ has constant rank and the normal representation is a plain representation of $\mathcal G$ on $\nu$. Moreover, we recall from \cite{crainic:def2} that the map
\[
p: C_{\mathrm{def}}(\mathcal G) \to C (\mathcal G, \nu), \quad c \mapsto \pi \circ s_c
\]
is a surjective quasi-isomorphism.

Consider a representation $E \to M$ of $\mathcal G$ and construct the associated trivial-core VB-groupoid $(\mathcal G \ltimes E \rightrightarrows E; \mathcal G \rightrightarrows M)$. At the infinitesimal level, there is an induced representation of $A$ on $E$, i.e.~an $A$-flat connection $\nabla: A \to DE$, and $\nabla$ is injective because $\rho$ is so. Therefore, the cokernel $\tilde \nu = DE/ \operatorname{im} \nabla$ is a vector bundle over $M$. Denote by
\[
\tilde \pi: DE \to \tilde \nu, \quad \delta \mapsto \bar \delta
\]
the projection.

We want to show that, in this situation, $\mathcal G$ acts on $\tilde \nu$. To see this, notice that the group $\mathrm{Bis}(\mathcal G)$ of bisections of $\mathcal G$ acts on $\Gamma(E)$ via
\begin{equation} \label{eq:bis_act_1}
(\beta \star \varepsilon)_x = \beta_{(\mathsf t \circ \beta)^{-1} (x)} \cdot \varepsilon_{(\mathsf t \circ \beta)^{-1} (x)},
\end{equation}
so it also acts on derivations of $E$ by
\begin{equation} \label{eq:bis_act_2}
(\beta \bullet \Delta)(\varepsilon) = \beta \star (\Delta( \beta^{-1} \star \varepsilon)).
\end{equation}
If $\beta$ is a local bisection around $x$ and $\beta(x) = g: x \to y$, these formulas still make sense: Equation (\ref{eq:bis_act_1}) shows that the action of $\beta$ takes local sections around $y$ to local sections around $x$, and Equation (\ref{eq:bis_act_2}) shows that $\beta$ acts on derivations locally defined around $x$ (to give a derivation locally defined around $y$).

Now, let $g: x \to y$ be an arrow in $\mathcal G$ and $\delta \in D_x E$. Choose a local bisection $\beta$ of $\mathcal G$ passing through $g$ and a derivation $\Delta \in \mathfrak D (E)$ such that $\Delta_x = \delta$. Our action is then defined by
\[
g \cdot \bar \delta = \overline{\beta \bullet \Delta}|_y.
\]
Let us show that the definition does not depend on the choice of $\Delta$. If $\Delta \in \mathfrak D (E)$ and $\Delta_x = 0$, then, for every $\varepsilon \in \Gamma(E)$, we have:
\[
\begin{aligned}
(\beta \bullet \Delta)_y (\varepsilon) & = (\beta \star \Delta(\beta^{-1} \star \varepsilon))_y \\
& = \beta_{(\mathsf t \circ \beta)^{-1}(y)} \cdot \Delta_{(\mathsf t \circ \beta)^{-1}(y)} (\beta^{-1} \star \varepsilon) \\
& = \beta_x \cdot \Delta_x (\beta^{-1} \star \varepsilon) = 0.
\end{aligned}
\]
The definition is also independent of the choice of $\beta$, though the proof is a bit more involved. This is equivalent to say that, if $\beta \in \mathrm{Bis}(\mathcal G)$ and $\beta_z = \mathsf 1_z$ for some $z \in M$, then there exists $\alpha \in \Gamma(A)$ such that 
\[
(\Delta - \beta \bullet \Delta - \nabla_\alpha)_z = 0.
\]
Consider the vector field $\sigma_{\Delta - \beta \bullet \Delta}$. We have
\begin{equation}\label{eq:xi}
\xi := (\sigma_{\Delta - \beta \bullet \Delta})_z = (\sigma_\Delta - (\mathsf t \circ \beta)_* (\sigma_\Delta))_z = \sigma_{\Delta_z} - T(\mathsf t \circ \beta)(\sigma_{\Delta_z}).
\end{equation}
But $\mathsf t \circ \beta$ preserves the orbits of $\mathcal G$, hence it maps a sufficiently small neighborhood of $z$ in the leaf $\mathcal L_z$ of $\operatorname{im} \rho$ through $z$ to itself. It then follows from (\ref{eq:xi}) that $\xi$ kills all the functions that are constant along $\mathcal L_z$, hence it belongs to $\operatorname{im} \rho$.

Now, let $\alpha \in \Gamma(A)$ be any section such that $\rho (\alpha_z) = \xi$ and put $D := \Delta - \beta \bullet \Delta - \nabla_\alpha$. By construction, $\sigma_{D_z} = 0$, so it suffices to show that $D_z$ vanishes on $\nabla$-flat sections. If $\varepsilon$ is such a section, then
\begin{equation}\label{eq:D_z}
D_z \varepsilon = (\Delta - \beta \bullet \Delta)_z \varepsilon = \Delta_z (\varepsilon - \beta^{-1} \star \varepsilon).
\end{equation}
But the hypothesis $\nabla \varepsilon = 0$ implies that $\varepsilon$ is invariant under the action of $\mathrm{Bis} (\mathcal G)$, at least locally around $z$, and the claim follows from (\ref{eq:D_z}).

Notice that the symbol map $\sigma: DE \to TM$ descends to a $\mathcal G$-equivariant map $\tilde \nu \to \nu$, and the fact that $\operatorname  {End} E \cap \operatorname{im} \nabla = 0$ implies that its kernel is again $\operatorname{End} E$. Hence we have a short exact sequence of vector bundles with $\mathcal G$-action:
\[
0 \longrightarrow \operatorname{End} E \longrightarrow \tilde \nu \longrightarrow \nu \longrightarrow 0.
\]
In turn, this induces a short exact sequence of DG-modules:
\[
0 \longrightarrow C(\mathcal G, \operatorname{End} E) \longrightarrow C(\mathcal G, \tilde \nu) \longrightarrow C(\mathcal G, \nu) \longrightarrow 0.
\]
But $\mathcal G \ltimes E$ is a trivial-core VB-groupoid, so we also have the sequence (\ref{eq:coch_2}):
\[
0 \longrightarrow C(\mathcal G, \operatorname{End} E) \longrightarrow C_{\mathrm{def,lin}}(\mathcal G \ltimes E) \longrightarrow C_{\mathrm{def}}(\mathcal G) \longrightarrow 0.
\] 
We are looking for a map relating the two sequences. If $\tilde c \in C^k_{\mathrm{def,lin}}(\mathcal G \ltimes E)$ and $\tilde c_2: \mathcal G^{(k)} \to DE$ is the map (\ref{eq:c_2}), we can simply define:
\begin{equation} \label{eq:tilde_p}
\tilde p: C_{\mathrm{def,lin}}(\mathcal G \ltimes E) \to C(\mathcal G, \tilde \nu), \quad \tilde c \mapsto \tilde \pi \circ \tilde c_2.
\end{equation}
This is a cochain map and we obtain the following commutative diagram:
\[
\begin{array}{c}
\xymatrix@C=15pt@R=20pt{
	0 \ar[r] & C(\mathcal G, \operatorname{End} E) \ar[r] \ar@{=}[d] & C_{\mathrm{def,lin}}(\mathcal G \ltimes E) \ar[r] \ar[d]^{\tilde p} & C_{\mathrm{def}}(\mathcal G) \ar[r]  \ar[d]^{p}  & 0 \\
	0 \ar[r] & C(\mathcal G, \operatorname{End} E) \ar[r] & C(\mathcal G, \tilde \nu) \ar[r] & C(\mathcal G, \nu) \ar[r] & 0
}
\end{array}.
\]
The rows are short exact sequences of DG-modules and the vertical arrows are DG-module surjections. Additionally, $p$ is a quasi-isomorphism. Hence, it immediately follows from the Snake Lemma and the Five Lemma that $\tilde p$ is a quasi-isomorphism as well. We have thus proved the following
\begin{proposition} The map (\ref{eq:tilde_p}) induces a canonical isomorphism between the linear deformation cohomology of the VB-groupoid $(\mathcal G \ltimes E \rightrightarrows E; \mathcal G \rightrightarrows M)$ and the Lie groupoid cohomology with coefficients in $\tilde \nu$:
	\[
	H_{\mathrm{def}, \mathrm{lin}} (\mathcal G \ltimes E) \cong H (\mathcal G, \tilde \nu).
	\]
\end{proposition}

\subsection{Lie group actions on vector bundles}\label{sec:Lie_grp_vect}

Let $G$ be a Lie group with Lie algebra $\mathfrak g$. Assume that $G$ acts on a vector bundle $E \to M$ by vector bundle automorphisms. Then $G$ acts also on $M$ and $(G \ltimes E \rightrightarrows E; G \ltimes M \rightrightarrows M)$ is a trivial-core VB-groupoid. We want to discuss its linear deformation cohomology.

Of course the action of $G$ on $M$ induces an infinitesimal action of $\mathfrak g$ on $M$, and the Lie algebroid of $G \ltimes M$ is the action algebroid $\mathfrak g \ltimes M$. We recall from \cite{crainic:def2} that $G \ltimes M$ acts naturally on $\mathfrak g \ltimes M$, by extending the adjoint action of $G$ on $\mathfrak g$, and on $TM$ by differentiating the action of $G$ on $M$. Moreover, there is a short exact sequence of complexes:
\[
0 \longrightarrow C (G \ltimes M, TM) \longrightarrow C_{\mathrm{def}}(G \ltimes M)  \overset{p}{\longrightarrow} C (G \ltimes M, \mathfrak g \ltimes M)[1] \longrightarrow 0.
\]
The projection $p$ is defined as follows. Take a $k$-cochain $c: (G \ltimes M)^{(k+1)} \to T(G \ltimes M)$ in the deformation complex of $G \ltimes M$ and let $(h_0, \dots, h_k) \in (G \ltimes M)^{(k+1)}$, $h_0 = (g, x)$. Then
\[
c(h_0, \dots, h_k) \in T_{(g_0,x_0)}(G \ltimes M) \cong T_{g_0} G \times T_{x_0} M.
\]
Then we compose with the canonical isomorphism 
\[
TR_{g_0}^{-1} \times \mathrm{id}_{T_{x_0} M}: T_{g_0} G \times T_{x_0} M \overset{\cong}{\longrightarrow} \mathfrak g \times T_{x_0} M
\]
and finally with the projection $TM \to M$, to get an element in $\mathfrak g \ltimes M$. The kernel of $p$ is given by $TM$-valued cochains. Since the $TM$-component is the projection by the source map, it does not depend on the first component, and we conclude that the kernel is $C(G \ltimes M, TM)$.

We want to construct a similar sequence for $C_{\mathrm{def,lin}}(G \ltimes E)$ taking into account the linear nature of the action.  First of all, there is an obvious induced action of $G$ on $DE$ and the symbol map $\sigma: DE \to TM$ is $G$-equivariant. Hence there is a short exact sequence of cochain complexes:
\[
0 \longrightarrow C(G \ltimes M, \operatorname{End} E) \longrightarrow C(G \ltimes M, DE) \longrightarrow C(G \ltimes M, TM) \longrightarrow 0.
\]
In this case, the sequence (\ref{eq:coch_2}) reads
\[
0 \longrightarrow C(G \ltimes M, \operatorname{End} E) \longrightarrow C_{\mathrm{def,lin}}(G \ltimes E) \overset{\Pi}{\longrightarrow}  C_{\mathrm{def}} (G \ltimes M) \longrightarrow 0
\]
and, composing $\Pi$ with $p$, we get a cochain map
\begin{equation}\label{eq:Lie_vect_pr}
C_{\mathrm{def,lin}}(G \ltimes E) \longrightarrow C (G \ltimes M, \mathfrak g \ltimes M)[1] \longrightarrow 0.
\end{equation}
Applying the isomorphism (\ref{eq:nerve_act}), we get $(G \ltimes E)^{(k)} \cong G^k \times E$, and similarly $(G \ltimes M)^{(k)} \cong G^k \times M$. Then if $\tilde c \in C^k_{\mathrm{def,lin}}(G \ltimes E)$ is a linear cochain, the diagram (\ref{diag:c_tilde_2}) takes the following form:
\[
\begin{array}{c}
\xymatrix@C=0pt@R=7pt{
& G^{k+1} \times E \ar[dl] \ar[rr]^{\tilde c} \ar[dd]|!{[dl];[dr]}{\hole} & & TG \times TE \ar[dl] \ar[dd] \\
G^k \times E \ar[rr]^(.65){s_{\tilde c}} \ar[dd] & & TE \ar[dd] \\
& G^{k+1} \ar[dl] \ar[rr]^(.4)c |!{[ur];[dr]}{\hole} & & TG \times TM  \ar[dl] \\
G^k \times M \ar[rr]^(.65){s_c} & & TM  }
\end{array}.
\]
Moreover, $\tilde c$ is in the kernel of (\ref{eq:Lie_vect_pr}) if and only if its $TG$-component is $0$. In this case, it is clear that $\tilde c$ is determined by $s_{\tilde c}$, that is in turn equivalent to the map $\tilde c_2: G^k \times M \to DE$ defined by (\ref{eq:c_2}). Therefore, the kernel of (\ref{eq:Lie_vect_pr}) is $C(G \ltimes M, DE)$.

Summarizing there is an exact diagram of cochain complexes
\[
\begin{array}{c}
\xymatrix@C=15pt@R=15pt{
& 0 \ar[d] & 0 \ar[d] & & \\
0 \ar[r] & C (G \ltimes M, \operatorname{End} E) \ar[d] \ar@{=}[r] & C (G \ltimes M, \operatorname{End} E) \ar[r] \ar[d] & 0 \ar[d] \\
0 \ar[r] & C (G \ltimes M, DE) \ar[r] \ar[d] & C_{\mathrm{def}, \mathrm{lin}} (G \ltimes E) \ar[r] \ar[d] & C (G \ltimes M, \mathfrak g \ltimes M)[1] \ar[r] \ar@{=}[d] & 0 \\
0 \ar[r] & C (G \ltimes M, TM) \ar[r] \ar[d] & C_{\mathrm{def}}(G \ltimes M) \ar[r] \ar[d] & C (G \ltimes M, \mathfrak g \ltimes M)[1] \ar[r] \ar[d] & 0 \\
 & 0 & 0 & 0 &
}
\end{array}.
\]
This proves the following
\begin{proposition}
Let $G$ be a Lie group acting on a vector bundle $E \to M$ by vector bundle automorphisms. The linear deformation cohomology of the VB-groupoid $(G \ltimes E \rightrightarrows E, G \ltimes M \rightrightarrows M)$ fits in the exact diagram:
\[
\begin{array}{c}
\xymatrix@C=15pt@R=15pt{
& \vdots \ar[d] &  \vdots \ar[d] &  \vdots \ar[d] & \\
\cdots \ar[r] & H^{k} (G \ltimes M, \operatorname{End} E) \ar@{=}[r] \ar[d]& H^{k}(G \ltimes M, \operatorname{End} E) \ar[r] \ar[d]& 0 \ar[r] \ar[d] & \cdots \\
\cdots \ar[r] & H^{k} (G \ltimes M, DE) \ar[r] \ar[d]& H^k_{\mathrm{def}, \mathrm{lin}}(G \ltimes E) \ar[r] \ar[d]& H^{k+1} (G \ltimes M, \mathfrak g \ltimes M) \ar[r] \ar@{=}[d] & \cdots \\
\cdots \ar[r] & H^{k} (G \ltimes M, TM) \ar[r] \ar[d] & H^k_{\mathrm{def}}(G \ltimes M) \ar[r] \ar[d]& H^{k+1} (G \ltimes M, \mathfrak g \ltimes M) \ar[r] \ar[d] & \cdots \\
\cdots \ar[r] & H^{k+1} (G \ltimes M, \operatorname{End} E) \ar@{=}[r] \ar[d]& H^{k+1}(G \ltimes M, \operatorname{End} E) \ar[r] \ar[d]& 0 \ar[r] \ar[d] & \cdots \\
& \vdots  &  \vdots  &  \vdots  &
}
\end{array}.
\]
\end{proposition}

\backmatter
\cleardoublepage
\phantomsection 

\end{document}